\documentclass[leterpaper, 11pt]{article}
\usepackage[utf8]{inputenc}
\usepackage[T1]{fontenc}
\usepackage{lmodern}
\usepackage{graphicx,amsmath,amsfonts,latexsym,amsthm,amssymb,mathrsfs,mathtools}
\usepackage[abs]{overpic}
\usepackage[dvipsnames]{xcolor}
\usepackage[linktocpage=true]{hyperref}
\usepackage{caption}
\usepackage{subcaption}
\usepackage{dsfont}
\usepackage{stmaryrd}
\SetSymbolFont{stmry}{bold}{U}{stmry}{m}{n}
\usepackage{tikz-cd}
\usepackage{tikz}
\usepackage{tikz-3dplot}
\usepackage{hyperref}
\tikzset{
  symbol/.style={
    draw=none,
    every to/.append style={
      edge node={node [sloped, allow upside down, auto=false]{$#1$}}}
  }
}
\usepackage{comment}
\usepackage{enumitem}
\numberwithin{equation}{section}
\usepackage{titlesec}
\usepackage{microtype}
\usepackage{subfiles}
\usepackage{yhmath}
\usepackage{epigraph}
\usepackage{authblk}

\graphicspath{{./Figures/}}

\usepackage[margin=1in]{geometry}

\usepackage[backend=biber, maxnames=4, style=alphabetic, sorting=nyt]{biblatex}
\hypersetup{
  colorlinks,
  citecolor=teal,
  linkcolor=purple,
  urlcolor=teal}
  
\newtheorem{thmintro}{Theorem}

\newtheorem{corintro}{Corollary}

\newtheorem{propintro}{Proposition}

\newtheorem{conjintro}{Conjecture}

\newtheorem{thm}{Theorem}[section]
\newtheorem{techthm}[thm]{Technical Theorem}
\newtheorem{cor}[thm]{Corollary}
\newtheorem{lem}[thm]{Lemma}
\newtheorem{question}[thm]{Question}
\newtheorem{prop}[thm]{Proposition}
\newtheorem{conj}[thm]{Conjecture}
\newtheorem{defn}[thm]{Definition}
\newtheorem{ex}[thm]{Example}
\newtheorem{rem}[thm]{Remark}

\newtheorem{choice}{Choice}

\newtheorem*{mainthm}{Main Theorem}
\newtheorem*{thm*}{Theorem}
\newtheorem*{cor*}{Corollary}
\newtheorem*{question*}{Questions}
\newtheorem*{conj*}{Conjecture}

\newcommand{\N}{\mathbb{N}}
\newcommand{\Z}{\mathbb{Z}}
\newcommand{\R}{\mathbb{R}}

\newcommand{\Q}{\mathbb{Q}}
\newcommand{\F}{\mathbb{F}}

\title{Taut foliations and contact pairs in dimension three}

\author{Thomas Massoni\thanks{\href{mailto:tmassoni@stanford.edu}{tmassoni@stanford.edu}, ORCiD 0009-0000-9567-1330.}}
\affil{\small \textit{Department of Mathematics, Stanford University, \\
\textit{450 Jane Stanford Way, Building 380, Stanford, CA 94305, USA.}}}

\date{\today}

\begin{document}

\maketitle

\begin{abstract}

We present a new construction of codimension-one foliations from pairs of contact structures in dimension three. This constitutes a converse result to a celebrated theorem of Eliashberg and Thurston on approximations of foliations by contact structures. Under suitable hypotheses on the initial contact pairs, the foliations we construct are \emph{taut}, allowing us to characterize the existence of taut foliations entirely in terms of contact geometry. This viewpoint reveals some surprising \emph{flexibility} phenomena for taut foliations, and provides new insight into the $L$-space conjecture.

The first part of the proof builds upon the work of Colin and Firmo on positive contact pairs. The second part involves a wide generalization of a technical result of Burago and Ivanov on the construction of branching foliations tangent to continuous plane fields, and might be of independent interest.

\end{abstract}

\vspace{1em}
\begin{figure}[h]
    \centering
    \includegraphics[width=0.8\linewidth]{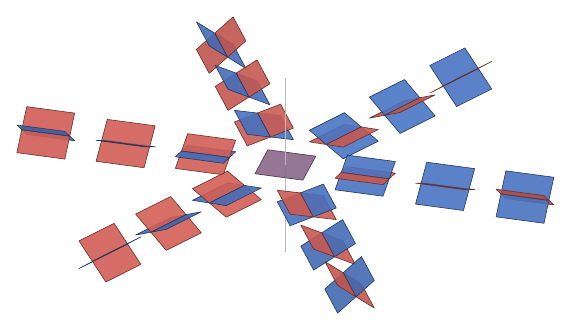}
    %\caption{Caption}
    %\label{fig:placeholder}
\end{figure}

%%%%%%%%%%%%%%%%%%%%%%%%%%%%%%%%%%%%%%%%%%%%%%%%%%%%
%%%%%%%%%%%%%%%%%%%%%%%%%%%%%%%%%%%%%%%%%%%%%%%%%%%%
\newpage
\tableofcontents

%%%%%%%%%%%%%%%%%%%%%%%%%%%%%%%%%%%%%%%%%%%%%%%%%%%%%%%%%%%%%%%%%%%%%%%%%%%%%%%
\newpage

\epigraph{Good things come in pairs.}{\textit{Chinese proverb}}
\section{Introduction}
%%%%%%%%%%%%%%%%%%%%%%%%%%%%%%%%%%%%%%%%%%%%%%%%%%%%%%%%%%%%%%%%%%%%%%%%%%%%%%%

Throughout this article, $M$ denotes a closed, connected, oriented three-dimensional manifold. The term foliation will always refer to a \emph{codimension one} foliation. All the foliations and contact structures under consideration will be \emph{(co)oriented}.

\medskip

The main result of this article is a construction of foliations from pairs of contact structures. This can be regarded as a converse result to a celebrated theorem of Eliashberg and Thurston~\cite{ET}.

\begin{mainthm}
Let $(\xi_-, \xi_+)$ be a positive contact pair on $M$, with a positively transverse vector field $\upsilon$. Assume either that one of $\xi_\pm$ is tight, or that they are everywhere transverse. Then there exists an aspherical $C^0$-foliation $\mathcal{F}$ on $M$ transverse to $\upsilon$.
\end{mainthm}

In particular, $T\mathcal{F}$ and $\xi_\pm$ are homotopic as oriented plane fields. If the vector field $\upsilon$ is volume preserving, then the foliation $\mathcal{F}$ is \emph{taut}. This allows us to completely characterize the existence of taut foliations on $3$-manifolds purely in terms of contact geometry. Therefore, this construction can be regarded as the first step of a program to address the \emph{$L$-space conjecture} from the point of view of contact topology.

%%%%%%%%%%%%%%%%%%%%%%%%%%%%%%%%%%%%%%%%%%%%%%%%%%%%%%%%%%%%%%%%%%%%%
        \subsection{Overview and motivations}
%%%%%%%%%%%%%%%%%%%%%%%%%%%%%%%%%%%%%%%%%%%%%%%%%%%%%%%%%%%%%%%%%%%%%

%%%
            \subsubsection{Foliations and contact structures}
%%%

We review some definitions and basic properties of (codimension-one) foliations and contact structures on $3$-manifolds.

\paragraph{Foliations.}

A $2$-dimensional foliation $\mathcal{F}$ on a $3$-dimensional manifold $M$ is, loosely speaking, a decomposition of $M$ by surfaces, or \emph{leaves}, which is locally modeled on the standard foliation by horizontal planes on $\R^3$. Alternatively, a smooth foliation may be described in terms of its tangent plane field $T \mathcal{F}$. Assuming that it is cooriented, it can be realized as the kernel of a smooth, nowhere vanishing $1$-form $\alpha$ satisfying
\begin{align}
    \alpha \wedge d\alpha = 0. \label{eq:frobalpha}
\end{align}
Conversely, any smooth nowhere vanishing $1$-form satisfying~\eqref{eq:frobalpha} defines a cooriented foliation $\mathcal{F}$ with $T \mathcal{F} = \ker \alpha$ in a unique way by the Frobenius integrability theorem.

\medskip

Foliations are among the most important structures on $3$-manifolds and their theory is extremely rich. We refer to the encyclopedic books~\cite{CC1,CC2} and~\cite{C07} for a particularly exhaustive account of this theory and for the standard definitions. We wish to think about foliations as \emph{fragile} objects which cannot be modified easily because of the integrability condition; a small perturbation of the $1$-form $\alpha$ typically does not satisfy~\eqref{eq:frobalpha} anymore.

\medskip

In this article, we will mostly consider foliations and plane fields which are not smooth but only continuous. Let us recall the important 

\begin{defn}
    A continuous plane field $\eta$ on $M$ is 
    \begin{itemize}
        \item \textbf{Locally integrable} if for every point $p \in M$, there exists a (germ of) $C^1$ surface tangent to $\eta$ passing through $p$,
        \item \textbf{Integrable} if there exists a topological foliation $\mathcal{F}$ on $M$ with $C^1$-immersed leaves tangent to $\eta$,
        \item \textbf{Uniquely integrable} if it is locally integrable, and at each point the germ of integral surface is unique.
    \end{itemize}
\end{defn}

We have the elementary implications
\begin{center}
    Uniquely integrable $\implies$ Integrable $\implies$ Locally integrable,
\end{center}
but the converse implications are not true in general. Even one dimension lower, there exist continuous (H\"{o}lder regular) vector fields on $\R^2$ which are \emph{not} tangent to foliations, or are tangent to many \emph{distinct} foliations, see~\cite{BF03}. However, for Lipschitz continuous plane fields (or vector fields), these three definitions are equivalent by the Picard--Lindel\"{o}f--Cauchy--Lipschitz theorem. Besides, a given integrable continuous plane field might be tangent to many distinct foliations.

We are interested in foliations whose leaves are tangent to a continuous plane field:

\begin{defn}
    A topological foliation $\mathcal{F}$ on $M$ is a \textbf{$C^0$-foliation} if its leaves are $C^1$ immersed surfaces tangent to a global continuous plane field on $M$.
\end{defn}

Essentially, this is the lowest regularity class for which $T \mathcal{F}$ is defined as a continuous plane field. However, $C^0$-foliations have \emph{no transverse regularity}. By a theorem of Calegari~\cite{C01}, any topological foliation on a $3$-manifold is topologically isotopic to a $C^0$-foliation with \emph{smooth} leaves. Therefore, we will often assume that $C^0$-foliations have smooth leaves. However, a $C^0$-foliation can typically \emph{not} be upgraded to a smooth---or even $C^1$---foliation in a natural way.

\medskip

Not every foliation on a $3$-manifold is interesting. Indeed, Thurston showed that for every homotopy class of plane fields on a closed $3$-manifold, there exists a smooth foliation whose tangent plane field is in that class, see~\cite{T76}. However, these foliations have \emph{Reeb components} by construction, making them too flexible to effectively measure the topology of the ambient manifold. Recall that a Reeb component of a foliation is a solid torus whose boundary is a leaf and such that the restriction of the foliation to the interior of the solid torus is a foliation by planes transverse to the core, see~\cite[Example 4.7]{C07}. A foliation without Reeb components is called \textbf{Reebless}. 

Another type of foliation was brought to prominence by the work of Thurston and Gabai: the class of \emph{taut foliations}.

\begin{defn} \label{def:taut}
    A (smooth) foliation $\mathcal{F}$ on $M$ is \textbf{(everywhere) taut} if for every $p \in M$, there exists an immersed loop in $M$ transverse to $\mathcal{F}$ and passing through $p$.
\end{defn}

There exist many equivalent definitions of tautness, see~\cite{C07} for a sample. This definition still makes sense for $C^0$-foliations. It can be adapted for topological foliations in different yet \emph{nonequivalent} ways, see~\cite{CKR19}. However, all these definitions are equivalent up to a $C^0$-small isotopy. While taut foliations are Reebless, the converse is not true in general but it holds on atoroidal manifolds. Importantly, the existence of a Reebless or taut foliation imposes strong restrictions on the topology and geometry of the ambient manifold $M$. For instance, Novikov's theorem implies that if $M$ carries a Reebless foliation $\mathcal{F}$, then the leaves of $\mathcal{F}$ are $\pi_1$-injective and the loops transverse to $\mathcal{F}$ are noncontractible, see~\cite{C07}. This implies that $M$ is irreducible and its universal cover is $\R^3$ (except if $M = S^1 \times S^2$).

\medskip

One could think of Reebless and taut foliations as \emph{brittle} objects, fragile and rigid at the same time. While it is known by the work of Gabai~\cite{G83} that every closed irreducible $3$-manifold with positive first Betti number carries a taut foliation, it is still a fundamental open problem to determine which \emph{rational homology spheres} carry a Reebless or a taut foliation. The \emph{$L$-space conjecture} described below provides some remarkable insight into that question.

\paragraph{Contact structures.}

A contact structure $\xi$ on $M$ is a smooth, \emph{maximally nonintegrable} plane field. It can be (locally) represented as the kernel of a smooth $1$-form $\alpha$ satisfying
\begin{align}
    \alpha \wedge d\alpha \neq 0 \label{eq:contalpha}
\end{align}
pointwise. Since $M$ is oriented and $\xi$ is assumed to be cooriented, we can find a globally defined contact form $\alpha$ that satisfies either $\alpha \wedge d\alpha > 0$ or $\alpha \wedge d\alpha < 0$ everywhere on $M$. This sign does not depend on the choice of $\alpha$. In the first case, $\xi$ is called \textbf{positive} and in the second case, it is called \textbf{negative}. We further assume that $\alpha$ agrees with the coorientation of $\xi$, in the sense that it evaluates positively on vectors transverse to $\xi$ with direction matching its coorientation.

\medskip

\noindent\fbox{%
\begin{minipage}{\dimexpr\linewidth-2\fboxsep-2\fboxrule\relax}
\textbf{We make the following orientation conventions:} the correspondence between the orientation and coorientation of a contact structure will be with respect to the orientation on $M$ for positive contact structures, and the \emph{opposite} orientation on $M$ for negative contact structures. Moreover, a vector field $Z$ is \emph{positively} transverse to a \emph{negative} contact structure $\xi_-$ if it is positively transverse to it as an oriented plane field, hence \emph{negatively} transverse to it as a cooriented plane field. In other words, a contact form for $\xi_-$ evaluates \emph{negatively} along $Z$.
\end{minipage}%
}

\medskip

By the Darboux theorem, a (positive) contact structure is locally modeled on the \emph{standard} contact structure on $\R^3$ defined by the $1$-form 
$$\alpha_\mathrm{std} \coloneqq r^2 d\theta + dz.$$
Here, $(r, \theta, z)$ denote the standard cylindrical coordinates on $\R^3$. Moreover, as opposed to foliations, contact structures are always smooth and quite \emph{elastic}: a small $C^1$ perturbation of $\alpha$ is still a contact form, and by Gray's stability theorem, the resulting contact structure is contact isotopic to the original one. We refer to the book~\cite{G09} for a thorough introduction to contact geometry and the standard definitions.

\medskip

As for foliations, not every contact structure on a $3$-manifold is interesting. Mirroring Thurston's existence theorem for foliations, Eliashberg showed that every homotopy class of oriented plane fields on $M$ contains a contact structure which is moreover \emph{overtwisted}. A contact structure $\xi$ is \textbf{overtwisted} if there exists an embedded disk in $M$ which is transverse to $\xi$, except along its boundary and its center where it is tangent to $\xi$. It is \textbf{tight} otherwise, and \textbf{universally tight} if its universal cover is tight. Overtwisted disks are to contact structures what Reeb components are to foliations, as they make contact structures flexible. On the other hand, tight contact structures are more rigid and more scarce.

\medskip

The choice of a contact form $\alpha$ induces a \textbf{Reeb vector field} $R$ on $M$ uniquely determined by 
\begin{align*}
    \alpha(R) =1, \qquad d\alpha(R, \cdot\, ) = 0.
\end{align*}
A contact structure is \textbf{hypertight} if it admits a contact form whose corresponding Reeb vector field has no contractible closed orbit. By a theorem of Hofer~\cite{H93}, hypertight contact structures are universally tight, but the converse does not hold in general.\footnote{The standard tight contact structure on $S^3$ is universally tight but not hypertight.}

\medskip

As opposed to Reebless or taut foliations, we would like to think of tight contact structures as \emph{ductile} objects, both supple and rigid. Moreover, the classification problem for tight contact structures is perhaps more reasonable than its counterpart for taut foliations.\footnote{It is not clear to us which equivalence relation on $C^0$-foliations to consider, whereas the answer is clear for contact structures: contact isotopy!} A \emph{coarse classification} was established by~\cite{CGH}, and a fine classification is known for a number of $3$-manifolds. Moreover, $3$-dimensional contact topology benefits from an arsenal of tools and invariants, from convex surfaces to pseudoholomorphic curves.

\medskip

Recently, Colin and Honda~\cite{CH20} showed a contact counterpart to Gabai's existence theorem for taut foliations: every closed irreducible $3$-manifold with positive first Betti number carries a hypertight contact structure. However, determining which rational homology spheres carry a (hyper)tight contact structure still remains an important open problem.

%%%
            \subsubsection{Eliashberg--Thurston approximation theorem}
%%%

At first glance, foliations and contact structures appear to be diametrically opposed notions---compare~\eqref{eq:frobalpha} with~\eqref{eq:contalpha}. However, the two theories share many key features, such as the dichotomy between \emph{flexibility} and \emph{rigidity}. A first breakthrough revealing a deep connection between foliations and contact structures was obtained by Eliashberg and Thurston:

\begin{thm*}[Eliashberg--Thurston~\cite{ET}]
    Let $\mathcal{F}$ be a $C^2$ aspherical foliation on $M$. Then $T\mathcal{F}$ can be $C^0$-approximated by positive and negative contact structures. Moreover, if $\mathcal{F}$ is taut, then its contact approximations are universally tight.
\end{thm*}

Here, $\mathcal{F}$ is \textbf{aspherical} if none of its leaves is homeomorphic to a sphere. By the Reeb Stability Theorem, $\mathcal{F}$ is aspherical if and only if it is not homeomorphic to the product foliation by spheres on $S^1 \times S^2$, see~\cite[Theorem 6.1.6]{CC1}.

Tautness can be characterized by the existence of a closed $2$-form $\omega$ that evaluates positively on leaves. If $\mathcal{F}$ is taut, Eliashberg and Thurston constructed from such a $2$-form a symplectic structure on $[-1,1] \times M$, which is a \emph{weak symplectic filling} of $(-M, \xi_-) \sqcup (M, \xi_+)$ where $\xi_\pm$ are negative and positive contact approximations of $\mathcal{F}$. They then appealed to a theorem of Eliashberg and Gromov to conclude that $\xi_-$ and $\xi_+$ are universally tight.

This result provides an efficient way to construct (tight) contact structures by approximating (taut) foliations. It has been extended and generalized in various directions; we now summarize some of them.

\paragraph{$C^0$-foliations.} The Eliashberg--Thurston theorem was generalized to $C^0$-foliations independently by Kazez--Roberts~\cite{KR17} and Bowden~\cite{B16a}.

\paragraph{Reebless foliations.} Colin~\cite{C02} showed that every $C^2$ aspherical Reebless foliation $\mathcal{F}$ can be approximated by universally tight contact structures. This has been strengthened by Bowden who showed that every $C^0$ small contact approximation of $\mathcal{F}$ is universally tight~\cite{B16b}, and that Colin's result holds for $C^0$-foliations~\cite{B16a}.\footnote{However, it is not known if every sufficiently small contact approximation of a Reebless $C^0$-foliation is universally tight.}

\paragraph{Foliations without invariant transverse measures.} More recently, Zung showed in~\cite{Z21} that a $C^2$ foliation $\mathcal{F}$ without invariant transverse measures can be approximated by contact structures whose Reeb vector fields are transverse to $\mathcal{F}$. Since $\mathcal{F}$ is taut, these contact approximations are hypertight. Moreover, Zung (implicitly) shows that $\mathcal{F}$ can be approximated by contact structures $\xi_-$ and $\xi_+$, negative and positive, respectively, supporting a \emph{Liouville pair}: there exist contact forms $\alpha_\pm$ for $\xi_\pm$ such that the $1$-form
\begin{align} \label{eq:liouv}
    \lambda \coloneqq e^{-s} \alpha_- + e^s \alpha_+
\end{align}
on $\R_s \times M$ is a \emph{Liouville form}, i.e., $d\lambda$ is a nondegenerate $2$-form. We call the pair $(\alpha_-, \alpha_+)$ a \textbf{Liouville pair} supported by $(\xi_-, \xi_+)$.

In $C^2$ regularity, the absence of invariant transverse measure is equivalent to the existence of an \emph{exact} $2$-form positive on the leaves. In view of Zung's theorem and in order to simplify the terminology, we propose the following definition:

\begin{defn} \label{def:hypertaut}
A cooriented $C^0$-foliation $\mathcal{F}$ is \textbf{hypertaut} if there exists an \emph{exact} $2$-form $\omega = d\beta$ such that $\omega_{\vert T \mathcal{F}}>0$.
\end{defn}

\paragraph{Anosov foliations.} An interesting class of hypertaut foliations arises as the weak-(un)stable foliations of \emph{Anosov flows}, see for instance the encyclopedic book~\cite{FH19} or the more accessible notes~\cite{B17}. In~\cite{H22a}, Hozoori showed a correspondence between Anosov flows on $3$-manifolds and a special class of Liouville pairs, that we call \textbf{Anosov Liouville pairs}. This was upgraded in~\cite{Mas22} to a \emph{homotopy equivalence} between the space of Anosov flows and the space of Anosov Liouville pairs. The latter structures can be thought of as very specific approximations of Anosov foliations that remember the hyperbolic nature of the flow.

\paragraph{Uniqueness of contact approximations.} In~\cite{V16}, Vogel proved that for most $C^2$ foliations without torus leaves, all the (positive) contact structures in a sufficiently small neighborhood of $T\mathcal{F}$ are contact isotopic. Therefore, the contact approximations of Eliashberg--Thurston are unique for these foliations. However, some foliations on $\mathbb{T}^3$ have infinitely many pairwise nonisotopic contact approximations.

\paragraph{Approximation versus deformation.} A \emph{contact deformation} of a $C^0$-foliation $\mathcal{F}$ is a $C^0$-continuous $1$-parameter family of plane fields $(\xi_t)$, $t \in [0,1]$, where $\xi_0 = T \mathcal{F}$, and for $t > 0$, $\xi_t$ is a contact structure. While $C^2$ foliations `with enough holonomy' can be deformed into contact structures (see~\cite[Proposition 2.9.2]{ET}), linear foliations without closed leaves on $\mathbb{T}^3$ do \emph{not} admit contact deformations (see~\cite[Corollary 9.11]{V16}). In contrast, Etnyre showed that any (cooriented) contact structure on a $3$-manifold is a smooth deformation of a foliation, see~\cite{E07}. However, the foliations constructed by Etnyre \emph{always} have Reeb components.\footnote{Etnyre starts with an \emph{open book decomposition} adapted to the contact structure, and inserts a suitable Reeb component at each component of the binding of the open book.}

%%%
            \subsubsection{A dictionary between foliations and contact pairs}
%%%

While $S^3$ carries a (unique) tight contact structure, namely the standard one, it does not carry a Reebless foliation. Therefore, there is no hope to construct a geometrically relevant foliation from a \emph{single} (tight) contact structure. Instead, one should consider a \emph{pair} of contact structures with opposite signs.

\begin{defn} \label{def:contpair}
A \textbf{contact pair} on $M$ is a pair $(\xi_-, \xi_+)$ where $\xi_-$ is a negative contact structure and $\xi_+$ is a positive contact structure. Both are cooriented. A \textbf{bicontact structure} is a contact pair $(\xi_-, \xi_+)$ such that $\xi_-$ and $\xi_+$ are everywhere transverse. A contact pair is \textbf{positive} if there exists a vector field $\upsilon$ positively transverse to $\xi_-$ and $\xi_+$.
\end{defn}

Alternatively, a contact pair $(\xi_-, \xi_+)$ is positive if and only if at every point $p \in M$ with $\xi_-(p) = \xi_+(p)$, these two planes have the same orientation (and opposite coorientations). A bicontact structure is automatically positive. Positive contact pairs were introduced by Colin and Firmo in~\cite{CF11}. See also~\cite{CH20} for a summary (in English!) of Colin--Firmo's paper.\footnote{Watch out: in~\cite{CH20}, the authors use the terminology \emph{contact pair} to mean \emph{positive contact pair}!}

The Eliashberg--Thurston theorem and its variations naturally produce positive contact pairs as approximations of foliations. We shall think of the positive and negative contact approximations, and the way they interact, as capturing some information about the original foliation. We will need both of them to be able to `reconstruct' the foliation, or rather produce a new one that is closely related to it.

\medskip

To each adjective characterizing the rigidity class of a foliation---\textit{Reebless, taut, hypertaut, Anosov}---we associate an adjective for positive contact pairs. Notice the implications:
\begin{center}
    Anosov $\implies$ Hypertaut $\implies$ Taut $\implies$ Reebless.
\end{center}

\begin{defn}
    Let $(\xi_-, \xi_+)$ be a contact pair.
\begin{itemize}
    \item It is \textbf{(universally) tight} if both $\xi_\pm$ are (universally) tight.
    \item It is \textbf{taut} if for every point $p \in M$, there exists a closed loop positively transverse to both $\xi_\pm$ and passing through $p$.\footnote{Those pairs were previously called \emph{strongly tight} in~\cite{CH20} and \emph{``fortement tendues''} in~\cite{CF11}, but we find that the name \emph{taut} is more appropriate as it parallels the foliation terminology in English.}
    \item It is \textbf{Liouville} if there exist contact forms $\alpha_\pm$ for $\xi_\pm$ such that $(\alpha_-, \alpha_+)$ is a Liouville pair, i.e., the $1$-form $\lambda$ on $\R \times M$ defined by~\eqref{eq:liouv} is a Liouville form.
    \item It is \textbf{Anosov Liouville} if there exist contact forms $\alpha_\pm$ for $\xi_\pm$ such that both $(\alpha_-, \alpha_+)$ and $(-\alpha_-, \alpha_+)$ are Liouville pairs.
\end{itemize}
\end{defn}
As we will see, the following implications hold:
\begin{center}
    Anosov Liouville $\implies$ Liouville $\implies$ Taut $\implies$ Universally tight positive.
\end{center}

Motivated by the previous results, we propose a dictionary between foliations and positive contact pairs in Table~\ref{table:dico}.
\begin{table}[h!]
\begin{center}
\begin{tabular}{|c | c|} 
 \hline 
 \bf{Foliations} & \bf{Positive contact} pairs \\ 
 \hline\hline
 Reebless & Universally tight \\ 
 \hline
 Taut & Taut \\
 \hline
Hypertaut & Liouville \\
 \hline
Anosov & Anosov Liouville \\
 \hline
\end{tabular}
\caption{A dictionary between foliations and positive contact pairs}
\label{table:dico}
\end{center}
\end{table}

An aspherical foliation satisfying a property in the left column is approximated by a positive contact pair with the corresponding property in the right column.\footnote{Strictly speaking, we do not know if Zung's result for hypertaut foliations holds in the $C^0$ setting. However, a hypertaut $C^0$-foliation is approximated by \emph{Liouville fillable} contact pairs, see Proposition~\ref{prop:fill}.} With the notable exception of the (universally) tight case, we will show that a contact pair satisfying a property in the right column `integrates' to a $C^0$-foliation with the corresponding property in the left column. In particular, a $3$-manifold $M \neq S^1 \times S^2$ carries a taut $C^0$-foliation \emph{if and only if} it carries a taut contact pair; see Corollary~\ref{corintro:tautchar} below.

%%%
            \subsubsection{The \texorpdfstring{$L$}{L}-space conjecture}
%%%

It is an open problem to determine which irreducible rational homology spheres carry a taut foliation. The famous $L$-space conjecture provides a hypothetical answer:

\begin{conj*}[$L$-space conjecture~\cite{BGW,J15}]
Let $M$ be an irreducible rational homology $3$-sphere. The following are equivalent:
\begin{enumerate}
    \item $M$ carries a cooriented taut foliation,
    \item $\pi_1(M)$ is left-orderable,
    \item $M$ is not an $L$-space.
\end{enumerate}
\end{conj*}

The second item can be thought of as an algebraic characterization, while the third one comes from gauge theory. Here, a rational homology sphere $M$ is an $L$-space if its `hat' Heegaard Floer homology is as small as possible, namely
$$\mathrm{rank} \widehat{HF}(M) = \vert H_1(M, \Z)\vert,$$
or alternatively if its reduced Heegaard Floer homology vanishes, $HF^\mathrm{red}(M) = 0$. Lens spaces are examples of $L$-spaces. We refer to~\cite{J15} for an overview of Heegaard Floer homology. In this article, we are mostly interested in the relation between items 1 and 3 of the $L$-space conjecture. We do not know how the (non) orderability of $\pi_1(M)$ could relate to contact pairs.

In~\cite{OS04}, Ozsváth and Szabó proved the implication $1 \implies 3$ by showing that the \emph{contact class} of a contact approximation of the foliation provided by the Eliashberg--Thurston theorem is nonzero in $HF^\mathrm{red}$. Interestingly, contact structures seem to have better interactions with gauge theory than foliations.\footnote{In~\cite{Z22}, Zhang associates to any \emph{smooth} hypertaut foliation $\mathcal{F}$ two nonzero invariants $c_\pm(\mathcal{F})$ in monopole Floer homology, without invoking the Eliashberg--Thurston theorem. However, this approach is still symplectic in nature as it relies on the construction of a suitable exact symplectic structure on $\R \times M$ depending on $\mathcal{F}$. This construction also requires a high degree of regularity for the foliation.}

\medskip

One key difficulty to reverse the Ozsváth--Szabó result and show that non-$L$-space irreducible rational homology spheres admit taut foliations is that very few general strategies are available to construct taut foliations. To our knowledge, all the known cases where the equivalence $1 \iff 3$ holds were obtained by rather ad hoc arguments, by considering a class of $3$-manifolds and independently determining which are $L$-spaces and which admit a taut foliation. See for instance~\cite{BC17, R17, RR17, HRRW20} for all graph manifolds, \cite{D20} for a large census of hyperbolic manifolds, and~\cite{K20, K23, S23a, S23b, S23c} for $3$-manifolds obtained by Dehn surgeries on suitable knots and links in $S^3$. Moreover, these approaches are quite combinatorial and perhaps a bit miraculous in nature. Instead, we propose a general program to construct Reebless and taut foliations that remains closer to the argument of Ozsváth--Szabó. Before that, let us state and discuss our main results.

%%%%%%%%%%%%%%%%%%%%%%%%%%%%%%%%%%%%%%%%%%%%%%%%%%%%%%%%%%%%%%%%%%%%%
        \subsection{Statement of results}
%%%%%%%%%%%%%%%%%%%%%%%%%%%%%%%%%%%%%%%%%%%%%%%%%%%%%%%%%%%%%%%%%%%%%

%%%
            \subsubsection{A new construction of foliations}
%%%

In~\cite{CF11}, the authors construct a \emph{locally integrable} plane field from a \emph{normal} positive contact pair. This technical condition can always be achieved after some \emph{nongeneric} deformation of the pair. In Section~\ref{sec:contpair}, we explain how to generalize this construction to \emph{any} positive contact pair, removing the normality condition. This will be important for some of our applications, in particular for Theorem~\ref{thmintro:liouv} below.

To a contact pair $(\xi_-, \xi_+)$, we will associate a \emph{canonical} vector field $X$ satisfying $X \in \xi_- \cap \xi_+$, and vanishing exactly on the set $\Delta \coloneqq \{ p \in M \ \vert \ \xi_-(p) = \xi_+(p)\}$. As a first step to obtain a foliation $\mathcal{F}$ from the contact pair, we construct a natural candidate for $T \mathcal{F}$:

\begin{thmintro} \label{thmintro:etau}
Let $(\xi_-, \xi_+)$ be \underline{any} positive contact pair on $M$ and let $(\phi_t)_t$ be the flow of the associated vector field $X$. 
    \begin{enumerate}
    \item The plane fields $\xi_\pm^t$ defined by
    \begin{align} \label{eq:xit}
        \xi_\pm^t(\phi_t(p)) \coloneqq d\phi_t [\xi_\pm(p)]
    \end{align}
    converge uniformly to a continuous plane field $\eta_u$ as $t$ goes to $+\infty$.
    \item For a generic pair $(\xi_-, \xi_+)$, $\eta_u$ is locally integrable on $M$.
    \end{enumerate}
Moreover, $\eta_u$ depends continuously on $(\xi_-, \xi_+)$ (for the smooth topology).
\end{thmintro}

Here, \emph{generic} means that there exists an unspecified comeagre set of positive contact pairs for which item $2$ holds. Actually, this genericity assumption is not absolutely necessary; one could show that $\eta_u$ is \emph{always} locally integrable. This follows from the next theorem and Remark~\ref{rem:etabranchfol} for specific positive contact pairs, but we expect this should hold more generally. We call $\eta_u$ the \textbf{(weak-)unstable plane field} of the contact pair $(\xi_-, \xi_+)$. It is only continuous, and it is not necessarily \emph{uniquely} integrable, even for bicontact structures; see for instance~\cite[Example 2.2.9]{ET}. However, we adapt the strategy of~\cite{BI08} to show:

\begin{thmintro} \label{thmintro:int}
Let $(\xi_-, \xi_+)$ be a positive contact pair on $M$. If either
\begin{itemize}
    \item $\xi_-$ or $\xi_+$ is tight, or
    \item $\xi_-$ and $\xi_+$ coincide along a generic link of saddle singularities (see Section~\ref{sec:locint}),
\end{itemize}
then the unstable plane field $\eta_u$ can be $C^0$-approximated by aspherical $C^0$-foliations.
\end{thmintro}

This immediately implies the Main Theorem from the beginning of the introduction. Let us briefly sketch the proof of this result. By Theorem~\ref{thmintro:etau}, we can \emph{canonically} associate to a (generic) positive contact pair $(\xi_-, \xi_+)$ a pair $(X, \eta_u)$, where $X$ is a smooth vector field and $\eta_u$ is a continuous, locally integrable plane field containing $X$ and invariant under the flow of $X$. We call such a pair a \textbf{polarized vector field}. After a generic perturbation of $(\xi_-, \xi_+)$, we can further simplify the dynamics of $X$ and precisely understand the behavior of $\eta_u$ near the singularities of $X$. We then adapt the methods of~\cite{BI08} to show that $\eta_u$ is tangent to a \emph{branching foliation}, namely, a collection of immersed, `maximal' surfaces in $M$ which cover $M$, and which may intersect but cannot topologically cross. Part~\ref{part:tech} of this article is entirely devoted to this very technical and subtle construction. We finally appeal to a result in~\cite{BI08} to separate the leaves of the branching foliation and obtain a genuine foliation, at the expense of a $C^0$-small perturbation of $\eta_u$. This solves Problem 5.3 from~\cite{CH20}.

\begin{rem} \label{rem:etabranchfol}
    In view of Lemma~\ref{lem:branchinglim} in the second part of this article, Theorem~\ref{thmintro:int} implies that the unstable plane field of \emph{any} positive contact pair satisfying the assumptions of the theorem is tangent to a \emph{branching foliation} in the sense of Burago--Ivanov~\cite{BI08} (see Definition~\ref{def:branchfol}). In particular, $\eta_u$ is locally integrable without any genericity assumption on the pair.
\end{rem}

In the special case where $(\xi_-, \xi_+)$ is a bicontact structure, $X$ is a \textbf{projectively (or conformally) Anosov flow} (see~\cite{ET,M95}) and $\eta_u$ coincides with the weak-unstable bundle of $X$. Our theorem implies that, without even modifying $(\xi_-, \xi_+)$, $\eta_u$ is tangent to a branching foliation which can further be approximated by $C^0$-foliations.

\medskip

With Theorem~\ref{thmintro:int} in hand, it is now relatively easy to show:

\begin{thmintro} \label{thmintro:taut}
    If $(\xi_-, \xi_+)$ is a taut contact pair, then $\eta_u$ can be $C^0$-approximated by taut aspherical $C^0$-foliations.
\end{thmintro}

Combined with the results of~\cite{CKR19} and the ($C^0$ version of the) Eliashberg--Thurston theorem, we obtain:

\stepcounter{thmintro}
\begin{corintro} \label{corintro:tautchar}
If $M \neq S^1 \times S^2$, then $M$ carries a taut (topological) foliation if and only if it carries a taut contact pair.
\end{corintro}

In Section~\ref{sec:liouv}, we show that a Liouville contact pair is automatically positive and supports a Liouville pair whose Reeb vector fields are both transverse to $\eta_u$. This implies:

\begin{thmintro} \label{thmintro:liouv}
If the contact pair $(\xi_-, \xi_+)$ supports a Liouville pair, then $\eta_u$ can be $C^0$-approximated by hypertaut $C^0$-foliations. Moreover, $\xi_\pm$ are hypertight. 
\end{thmintro}

Notice that supporting a Liouville pair is a $C^1$ open condition for contact pairs. However, unlike tautness, it is \emph{not} a $C^0$ open condition. Besides, it is quite surprising to us that Liouville contact pairs are \emph{automatically} hypertight; it is not true that a Liouville fillable contact structure is hypertight in general.

It turns out that the unstable plane field of a Liouville contact pair is quite rigid. For instance, it is $C^1$ for Anosov Liouville pairs (see~\cite{Mas22}) and uniquely integrable in certain cases:

\begin{thmintro}[Corollary~\ref{cor:transtaut}]\label{thmintro:unique}
Let $(\xi_-, \xi_+)$ be a Liouville contact pair on $M$. If either
\begin{itemize}
    \item $\xi_\pm$ are everywhere transverse, or 
    \item $\xi_\pm$ coincide along a generic link of source singularities (see Section~\ref{sec:locint}),
\end{itemize}
then $\eta_u$ is \underline{uniquely} integrable. Moreover, it is tangent to a hypertaut $C^0$-foliation.
\end{thmintro}

Here, the Liouville condition imposes strong geometric and dynamical restrictions on $\eta_u$ which force unique integrability. Perhaps more is true:

\begin{question} \label{quest:!int}
For a generic Liouville contact pair, is $\eta_u$ uniquely integrable?
\end{question}

\bigskip

Our techniques also provide a precise description of the \emph{skeleton} of Liouville pairs on $3$-manifolds: 

\stepcounter{thmintro}
\begin{propintro}[Proposition~\ref{prop:skel}]
    Let $(\alpha_-, \alpha_+)$ be a Liouville pair on $M$. The skeleton of the associated Liouville structure $\lambda$ on $V= \R_s \times M$ defined by~\eqref{eq:liouv} is of the form $$\mathfrak{skel}(V, \lambda) = \mathrm{graph}(\sigma) \coloneqq \big\{(\sigma(x), x) \ \vert \  x \in M\big\}$$ for some continuous function $\sigma : M \rightarrow \R$.
\end{propintro}

In particular, the skeleton is a continuous, separating hypersurface homeomorphic to $M$. Interestingly, it is not smooth in general.

%%%
        \subsubsection{Application: transverse surgeries on taut foliations}
%%%

Contact structures are more malleable objects than foliations. As an application of our machinery, we obtain that the existence of taut foliations is preserved under ``large slope'' surgeries along transverse links. For simplicity, we only discuss the case of knots, but our results immediately generalize to links. 

Let $K \subset M$ be a framed knot. If $s \in \Q$ is a rational number, we denote by $M_K(s)$ the manifold obtained from $M$ by performing a Dehn surgery of slope $s$ along $K$.

If $\mathcal{F}$ is a foliation on $M$ transverse to $K$, performing a nontrivial Dehn surgery along $K$ typically destroys $\mathcal{F}$, in the sense that $\mathcal{F}$ does not induce a foliation on $M_K(s)$. It rather induces a singular foliation, which is singular along the image of $K$ in $M_K(s)$. One could instead perform a \emph{turbulization}\footnote{Sometimes spelled \emph{turbularization}.} (see~\cite[Example 3.3.11]{CC1}) along $K$ by spinning the leaves of $\mathcal{F}$ along $K$ thereby creating a torus leaf around $K$, and inserting a Reeb component after surgery. Of course, this procedure \emph{never} yields a taut foliation.

However, it is much easier to perform a surgery on an approximating contact pair and obtain a taut contact pair on $M_K(s)$, provided that $\vert s\vert $ is large enough. Combined with Theorem~\ref{thmintro:taut}, this implies:

\begin{thmintro} \label{thmintro:transfol}
Let $\mathcal{F}$ be a taut aspherical $C^0$-foliation on $M$, and $K \subset M$ be a framed knot transverse to $\mathcal{F}$. There exists $s_0 = s_0(K) \geq 0$ such that for every rational number $s \in \Q$ satisfying $\vert s \vert \geq s_0$, $M_K(s)$ carries a taut $C^0$-foliation $\mathcal{F}'$. Moreover, the image of $K$ in $M_K(s)$ is transverse to $\mathcal{F}'$.
\end{thmintro}

See Section~\ref{sec:transurg} for a detailed proof. As a corollary, we obtain a generalization of the main result of~\cite{LR14}, with a somewhat easier and more natural proof. 

Let $K \subset S^3$ be a nontrivial knot. We denote by $\mathcal{S}_K \subset \Q$ the set of rational slopes $s \in \Q$ such that there exists a taut $C^0$-foliation on $S^3 \setminus K$ which intersects the boundary of a tubular neighborhood of $K$ transversally along a foliation by circles of slope $s$. By a celebrated theorem of Gabai~\cite{G87}, $\mathcal{S}_K$ contains $0$. By~\cite[Theorem 1.1]{LR14}, $\mathcal{S}_K$ contains a neighborhood of $0$. We actually get:

\stepcounter{thmintro}
\begin{corintro} \label{corintro:slope} $\mathcal{S}_K$ is an open subset of $\Q$.
\end{corintro}

We obtain this corollary as a rather brutal application of Theorem~\ref{thmintro:transfol} for which we have little control over the number $s_0$ in such generality. Instead, it would be interesting to construct a concrete contact pair adapted to $K$ and carefully analyze the slopes of the characteristic foliations of the contact structures on the boundary of a (potentially large) solid torus containing $K$. As a starting point, one could consider the case of \emph{fibered knots}. There is already an extensive literature on taut foliations in fibered knot complements, see~\cite{R01a, R01b, K20, K23}. Moreover, in view of the $L$-space conjecture and Heegaard Floer homology computations for Dehn surgeries on knots in $S^3$ (see~\cite{KMOS07, RR17}), one expects the following:

\begin{conj*}
If $K$ is a knot in $S^3$ with Seifert genus $g>0$, then either $\mathcal{S}_K = \Q$, $\mathcal{S}_K = (-\infty, 2g-1) \cap \Q$, or $\mathcal{S}_K = (-2g+1, +\infty) \cap \Q$.
\end{conj*}

However, to our knowledge, it is not known if $\mathcal{S}_K \subset \Q$ is an interval in general.

%%%%%%%%%%%%%%%%%%%%%%%%%%%%%%%%%%%%%%%%%%%%%%%%%%%%%%%%%%%%%%%%%%%%%
        \subsection{Discussion and further directions}
%%%%%%%%%%%%%%%%%%%%%%%%%%%%%%%%%%%%%%%%%%%%%%%%%%%%%%%%%%%%%%%%%%%%%

%%%
            \subsubsection{Towards a contact \texorpdfstring{$L$}{L}-space conjecture}
%%%

While we are able to construct taut foliations from taut contact pairs, our methods do not easily extend to tight contact pairs. However, we expect the following to be true:

\begin{conjintro}
If $M \neq S^1 \times S^2$, then $M$ carries a Reebless $C^0$-foliation if and only if it carries a tight positive contact pair.
\end{conjintro}

We give some evidence for this conjecture in Section~\ref{sec:tautight}. The main difficulty comes from the lack of regularity of the unstable plane field $\eta_u$, and the fact that Reeblessness is \emph{not} a natural property for plane fields (see the \emph{phantom Reeb components} from~\cite{CKR19}). Notice that this conjecture would immediately imply that on an \emph{atoroidal} closed $3$-manifold, the existence of a tight positive contact pair is equivalent to the existence of a taut foliation.

\medskip

A contact pair whose contact structures are homotopic as plane fields can always be deformed into a positive one when one of its contact structures is overtwisted, see Proposition~\ref{prop:pos} below. However, there is little hope to construct geometrically interesting foliations from such a pair. Upgrading tight contact pairs to positive ones seems much more difficult. At the very least, we should further assume that the contact structures are homotopic as oriented plane fields, but this might not be enough. Indeed, Lin recently remarked that the existence of a taut foliation on a rational homology sphere imposes more constraints on the (monopole) Floer homology of the manifold, see~\cite{L23}.\footnote{It is known that (the various flavors of) monopole Floer homology and Heegaard Floer homology are isomorphic, so we simply refer to these invariants as `Floer homology'.} In particular, the $\mathbb{F}[U]$-module $HM_*(M)$ has a \emph{direct $\F$-summand}, where $\F = \Z \slash 2\Z$. Moreover, the \emph{contact invariants} of the contact approximations of the foliation have to pair to $1$, for the natural perfect pairing
$$\langle \, \cdot \, , \, \cdot \, \rangle : HM_*(M) \otimes HM_*(-M) \longrightarrow \F$$
induced by the Poincar\'{e} duality isomorphism $HM^*(M) \cong HM_*(-M)$. This motivates the following:

\begin{defn}
A contact pair $(\xi_-, \xi_+)$ on $M$ is \textbf{algebraically tight} if the contact invariants $c(\xi_\pm) \in HM_*(\mp M)$ satisfy
$$\big \langle c(\xi_-), c(\xi_+)\big\rangle = 1.$$
\end{defn}

Notice that if $(\xi_-, \xi_+)$ is algebraically tight, then $\xi_-$ and $\xi_+$ are tight and homotopic as oriented plane fields. Indeed, the contact class of an overtwisted contact structure vanishes, and the grading of the contact class corresponds to the homotopy class of the contact structure (as an oriented plane field). We propose:

\begin{conjintro}
Let $(\xi_-, \xi_+)$ be an algebraically tight contact pair on $M$. Then $(\xi_-, \xi_+)$ is homotopic through contact pairs to a positive one.
\end{conjintro}

Here, we think of the contact classes as potential \emph{obstructions} to deforming a tight contact pair into a positive one. One possible approach to this conjecture would be to consider the set $\Delta_- \subset M$ where $\xi_-$ and $\xi_+$ coincide with \emph{opposite} orientations. It is generically an embedded link which is null-homologous if $\xi_-$ and $\xi_+$ are homotopic. Perhaps one could study the link Floer or sutured Floer homology of the pair $(M, \Delta_-)$ and use the algebraic tightness hypothesis to precisely understand the topology of $\Delta_-$. Additionally, there could be `elementary moves' that would inductively simplify $\Delta_-$. Some obstruction to performing these moves could be determined by Floer-theoretic invariants of the contact pair.

\medskip

It remains to understand how to construct suitable tight contact structures on rational homology spheres. This is of course a very delicate problem, but perhaps more tractable than constructing taut foliations.

\begin{question}[Contact realization] Assume that $M$ is an irreducible rational homology sphere.
\begin{itemize}
    \item Which nonzero elements in the Floer homology of $M$ can be realized as the contact class of a (necessarily tight) contact structure?
    \item If the Floer homology of $M$ has a direct $\F$-summand, does $M$ carry an algebraically tight contact pair?
\end{itemize}
\end{question}

Assuming the previous two conjectures, a positive answer to the second item of this question would imply:

\begin{conjintro}
If $M$ is an irreducible rational homology sphere such that $HM_*(M)$ has a direct $\F$-summand, then $M$ carries a Reebless $C^0$-foliation.
\end{conjintro}

Recently, Alfieri and Binns verified the existence of a direct $\mathbb{F}$-summand in Floer homology for a large class of non-$L$-space irreducible rational homology spheres~\cite{AB24}.

%%%
            \subsubsection{Variations and refinements}
%%%

We conclude this introduction with some potential generalizations of our main result.

\paragraph{Uniqueness of the foliation.}

While the polarized vector field $(X, \eta_u)$ associated to a positive contact pair $(\xi_-, \xi_+)$ is \emph{canonical}, the (branching) foliation tangent to $\eta_u$ that we construct is not. Indeed, it depends on a certain number of choices that are carefully described in Part~\ref{part:tech}. Understanding how the foliation depends on these choices seems a completely intractable problem. In particular, we have a very limited understanding of the topology of the leaves of the foliation. One could think of our approach as sacrificing some knowledge about the leaves in order to gain extra flexibility that ultimately allows us to integrate $\eta_u$. This should be contrasted with \emph{branched surface theory}, a well-established approach to construct foliations and laminations in dimension three. Besides, a given continuous plane field might be tangent to many distinct foliations which have drastically different holonomy properties. However, $\eta_u$ `remembers' many of the geometric properties of the contact pair.

\paragraph{Parametric version.} To simplify the construction of a (branching) foliation tangent to $\eta_u$, we impose some generic conditions on $(\xi_-, \xi_+)$ that make the dynamics of $X$ more tractable. We believe that the overall strategy of Part~\ref{part:tech} is quite robust and some of these conditions may be weakened. In particular, we ask:

\begin{question}
Let $(\xi^t_-, \xi^t_+)$, $t \in [0,1]$, be a smooth path of positive contact pairs on $M$, and assume that one of $\xi^0_\pm$ is tight. Does there exist a path of $C^0$-foliations $(\mathcal{F}_t)$, $t \in [0,1]$, such that $T\mathcal{F}_t$ is uniformly $C^0$-close to $\eta^t_u$?
\end{question}

It is not absolutely clear to us which topology to consider on the space of $C^0$-foliations, since those are not determined by their tangent plane field. To answer this question, one would have to precisely understand the type of singularities that the associated vector field $X$ can develop in a generic $1$-parameter family, and adapt the various choices made in the nonparametric construction to this setting.

\paragraph{Version with boundary.}

We now assume that $M$ is compact and has a nonempty boundary. Some of our constructions can be generalized to this setting. In particular, if $(\xi_-, \xi_+)$ is a positive contact pair such that the associated vector field $X$ is outward pointing along $\partial M$, then the construction of $\eta_u$ goes through and yields a locally integrable plane field. Moreover, the methods of Part~\ref{part:tech} should extend to this setting so that Theorem~\ref{thmintro:int} generalizes to manifolds with boundary. 

It would be interesting to first prescribe a $1$-dimensional (branching) foliation tangent to $\eta_u \cap \partial M$ along $\partial M$ and construct a (branching) foliation on $M$ which restricts to the chosen one along $\partial M$. A result along these lines will appear in~\cite{MassoniZung}, under stronger hypotheses, that in turn make it possible to control the boundary slope of the foliation.

%%%%%%%%%%%%%%%%%%%%%%%%%%%%%%%%%%%%%%%%%%%%%%%%%%%%%%%%%%%%%%%%%%%%%
        \subsection{Acknowledgments}
%%%%%%%%%%%%%%%%%%%%%%%%%%%%%%%%%%%%%%%%%%%%%%%%%%%%%%%%%%%%%%%%%%%%%

I am grateful to my PhD advisor John Pardon for his constant support and encouragement. This project greatly benefited from numerous discussions with Jonathan Bowden, Vincent Colin, Emmy Murphy, and Jonathan Zung. I would also like to thank Ian Agol, Danny Calegari, Yasha Eliashberg, John Etnyre, Peter Ozsv\'{a}th, Rachel Roberts, and B\"{u}lent Tosun for their interest, and Anshul Adve and Siddhi Krishna for stimulating conversations.

The figure on the front page was produced with the assistance of Claude Fable 5 (Anthropic), which was also used for proofreading and correcting typos. Beyond this, AI tools played no role in the mathematical content or proofs.

I am supported by a Stanford Science Fellowship.

%%%%%%%%%%%%%%%%%%%%%%%%%%%%%%%%%%%%%%%%%%%%%%%%%%%%%%%%%%%%%%%%%%%%%%%%%%%%%%%
\newpage
\part{From contact pairs to foliations}
%%%%%%%%%%%%%%%%%%%%%%%%%%%%%%%%%%%%%%%%%%%%%%%%%%%%%%%%%%%%%%%%%%%%%%%%%%%%%%%

This first part contains the proofs of the main theorems stated in the introduction, assuming a technical result whose proof will be deferred to Part~\ref{part:tech}.

%%%%%%%%%%%%%%%%%%%%%%%%%%%%%%%%%%%%%%%%%%%%%%%%%%%%%%%%%%%%%%%%%%%%%%%%%%%%%%%%%%%%%%%%%%%%%%%%%%%%%%%%%%%%%%%%%%%%%%%%%%%%%%%%%%%%%%%%%%%%%%%%%%%%%%%%%%%%%%
    \section{Contact pairs revisited} \label{sec:contpair}
%%%%%%%%%%%%%%%%%%%%%%%%%%%%%%%%%%%%%%%%%%%%%%%%%%%%%%%%%%%%%%%%%%%%%%%%%%%%%%%%%%%%%%%%%%%%%%%%%%%%%%%%%%%%%%%%%%%%%%%%%%%%%%%%%%%%%%%%%%%%%%%%%%%%%%%%%%%%%%

In this section, we extend the main construction of~\cite{CF11} and we prove Theorem~\ref{thmintro:etau}. It will be important to avoid using the normal form for contact pairs from~\cite{CF11}, as the deformation to a normal contact pair is \emph{not generic}.

%%%%%%%%%%%%%%%%%%%%%%%%%%%%%%%%%%%%%%%%%%%%%%%%%%%%%%%%%%%%%%%%%%%%%
            \subsection{The setup}
%%%%%%%%%%%%%%%%%%%%%%%%%%%%%%%%%%%%%%%%%%%%%%%%%%%%%%%%%%%%%%%%%%%%%

We first recall the following
        
\begin{defn}
A \textbf{contact pair} on $M$ is a pair of cooriented contact structures $(\xi_-, \xi_+)$, negative and positive, respectively. A \textbf{bicontact structure} is a contact pair $(\xi_-, \xi_+)$ such that $\xi_\pm$ are everywhere transverse.
\end{defn}

Let $(\xi_-, \xi_+)$ be a contact pair. We write
\begin{itemize}
\item $\Delta_+$ for the set of $p \in M$ such that $\xi_-(p)=\xi_+(p)$ as oriented planes,
\item $\Delta_-$ for the set of $p \in M$ such that $\xi_-(p)= - \xi_+(p)$ as oriented planes,
\end{itemize}
and we let $\Delta \coloneqq \Delta_- \sqcup \Delta_+$. The contact pair $(\xi_-, \xi_+)$ is \textbf{positive} (resp.~\textbf{negative}) if $\Delta_- = \varnothing$ (resp. $\Delta_+ = \varnothing$). Note that $\Delta$ has empty interior, since $\xi_\pm$ define opposite orientations on $M$.

\begin{rem}
A contact pair $(\xi_-, \xi_+)$ is positive if and only if there exists a (smooth) vector field $\upsilon$ which is positively transverse to both $\xi_-$ and $\xi_+$. In that case, $\xi_-$ and $\xi_+$ are homotopic as oriented plane fields, and their first Chern classes $c_1(\xi_-)$ and $c_1(\xi_+)$ agree. For a generic contact pair $(\xi_-, \xi_+)$, $\Delta$ is a smoothly embedded link. The Poincar\'{e} duals to the homology classes represented by $\Delta_\pm$ (with suitable orientations) are given by 
$$\mathrm{PD}[\Delta_\pm] = \frac{1}{2}\left(c_1(\xi_+) \pm c_1(\xi_-)\right) \in H^2(M; \Z).$$
In particular, those classes only depend on the homotopy classes of $\xi_\pm$. If $(\xi_-, \xi_+)$ can be deformed through contact pairs into a positive contact pair, then $\Delta_-$ is null-homologous in $M$.

Assume now that $\Delta_-$ is null-homologous in $M$. Then $\xi_-$ and $\xi_+$ are homotopic as plane fields if and only if the \emph{linking number} $\mathrm{lk}(\Delta_-, \Delta_+) \in \Z \slash d \Z $ of $\Delta_+$ with $\Delta_-$ vanishes. Here, $d \in \Z_{>0}$ is defined as $c_1(\xi_+) \cdot H_2(M; \Z) = d \Z$. This linking number is obtained as the algebraic intersection number of $\Delta_+$ (suitably oriented) with any oriented embedded surface $\Sigma \subset M$ with $\partial \Sigma = \Delta_-$; changing $\Sigma$ modifies the intersection number by a multiple of $d$. Let us choose an arbitrary Riemannian metric on $M$ as well as a trivialization $\tau$ of $TM$. Identifying plane fields with their positive unit normal vectors, $\xi_\pm$ become maps $\nu_\pm : M \rightarrow S^2$, which combine to a map $\nu_2 = (\nu_-, \nu_+) : M \rightarrow S^2 \times S^2$. We then have
$$\Delta_\pm = \nu_2^{-1}(\mathrm{Diag}_\pm),$$
where $\mathrm{Diag}_\pm$ denote the diagonal and antidiagonal in $S^2 \times S^2$, respectively. Since $c_1(\xi_+) = c_1(\xi_-)$, $\xi_-$ and $\xi_+$ are homotopic as oriented plane fields over the $2$-skeleton of $M$. The obstruction to extending the homotopy over $3$-cells is measured by a \emph{relative Hopf invariant} with value in $\Z$. The homotopy over the $2$-skeleton is not unique, and changing the homotopy modifies the relative Hopf invariant by a multiple of $d$. By an elementary version of the Pontryagin--Thom construction, this number equals the linking number $\mathrm{lk}(\Delta_-, \Delta_+)$ in $\Z \slash d \Z$.

To summarize, $\xi_-$ and $\xi_+$ are homotopic as oriented plane fields if and only if $[\Delta_-] = 0$ and $\mathrm{lk}(\Delta_-, \Delta_+)=0$.
\end{rem}

As discussed in the introduction, it is a central question to understand when a contact pair can be deformed to a positive one:

\begin{question} \label{quest:pair}
Let $(\xi_-, \xi_+)$ be a contact pair on $M$ such that $\xi_-$ and $\xi_+$ are homotopic as oriented plane fields. Can $(\xi_-, \xi_+)$ be deformed into a positive contact pair? Are there further Floer-theoretic obstructions?
\end{question}

The following proposition, suggested to us by Vincent Colin, gives an affirmative answer to this question when one of $\xi_\pm$ is \emph{overtwisted}:

\begin{prop} \label{prop:pos}
Let $\xi_+$ be a positive contact structure on $M$.
    \begin{enumerate}
    \item There exists an overtwisted negative contact structure $\xi_-$ on $M$ such that $(\xi_-, \xi_+)$ is a positive contact pair.
    \item Let $\xi_-$ be an overtwisted negative contact structure on $M$ which is homotopic to $\xi_+$ as an oriented plane field. There exists a path of negative contact structures $(\xi^\tau_-)_{0 \leq \tau \leq 1}$ such that $\xi^0_- = \xi_-$, and $(\xi^1_-, \xi_+)$ is a positive contact pair.
    \end{enumerate}
\end{prop}

\begin{proof}
By Etnyre~\cite{E07}, there exists a smooth foliation $\mathcal{F}$ on $M$ such that $\xi_+$ is a deformation of $T \mathcal{F} = \eta$, in the sense that there exists a (smooth) $1$-parameter family of smooth oriented plane fields $(\xi^\tau_+)_{0 \leq \tau \leq 1}$ such that $\xi^0_+ = \eta$, $\xi^1_+ = \xi_+$, and $\xi^\tau_+$ is a (positive) contact structure for $\tau > 0$. Moreover, the foliation $\mathcal{F}$ has Reeb components by construction, and is approximated by \emph{overtwisted} negative contact structures.\footnote{This essentially follows from an observation in~\cite[Section 4]{E07}: an approximation of $\mathcal{F}$ by a negative contact structure can equivalently be seen as an approximation by a positive contact structure after reversing the orientation of $M$; in that case, the leaves of $\mathcal{F}$ coming from the pages of the open book spiral towards the Reeb components in a \emph{counterclockwise} manner, which implies overtwistedness.} Therefore, there exist an overtwisted negative contact structure $\xi'_-$ and $0 < \tau \leq 1$ such that $(\xi'_-, \xi^{\tau}_+)$ is a positive contact pair. Applying Gray's stability yields 1. If $\xi_-$ is a negative contact structure as in $2$, then $\xi'_-$ and $\xi_-$ are two overtwisted negative contact structures which are homotopic as plane fields. Eliashberg's classification of overtwisted contact structures~\cite{E89} implies that they are contact isotopic, and $2$ follows easily.
\end{proof}

Of course, Question~\ref{quest:pair} is much more interesting (and difficult!) in the case where both $\xi_\pm$ tight.

\bigskip

We now investigate some basic quantities associated with contact pairs. Let $\mathrm{dvol}$ be a volume form on $M$ and $\alpha_\pm$ be contact forms for $\xi_\pm$ satisfying $$\alpha_+ \wedge d\alpha_+ = - \alpha_- \wedge d\alpha_- = \mathrm{dvol}.$$ There exists a unique (smooth) vector field $X$ on $M$ such that $$\alpha_\pm(X) = 0, \qquad \iota_X \mathrm{dvol} = \alpha_- \wedge \alpha_+.$$ Away from $\Delta$, $X$ is nonvanishing and spans the transverse intersection $\xi_- \pitchfork \xi_+$ with the correct orientation, and $X$ vanishes along $\Delta$. Note that $X$ does not depend on the choice of $\mathrm{dvol}$, so it only depends on the contact structures $\xi_\pm = \ker \alpha_\pm$. There exist smooth functions $f_0, g_0 : M \rightarrow \R$ defined by 
\begin{align*}
d(\alpha_- \wedge \alpha_+) &= f_0 \, \mathrm{dvol}, \\
\langle \alpha_-, \alpha_+ \rangle &= g_0 \, \mathrm{dvol},
\end{align*}
where
\begin{align*}
    \langle \alpha_-, \alpha_+ \rangle \coloneqq \alpha_- \wedge d\alpha_+ + \alpha_+ \wedge d\alpha_-.
\end{align*}
The function $g_0$ does not depend on the choice of $\mathrm{dvol}$, but the function $f_0$ does, since it is nothing more than the divergence of $X$ for $\mathrm{dvol}$: $$\mathcal{L}_X \mathrm{dvol} = f_0 \, \mathrm{dvol}.$$ If $\mathrm{dvol}'$ is another volume form on $M$, there exists a unique function $\rho : M \rightarrow \R$ such that $$\mathrm{dvol}' = e^\rho \, \mathrm{dvol}$$ and one easily checks that $$f'_0 = X \cdot \rho + f_0,$$ where $f'_0$ is the divergence of $X$ for $\mathrm{dvol}'$. However, $f_0$ does not depend on $\mathrm{dvol}$ along $\Delta$. Moreover:

\begin{lem} \label{lem:f0}
    Along $\Delta_+$, $f_0 \geq 2$ and along $\Delta_-$, $f_0 \leq -2$.
\end{lem}

\begin{proof}
Along $\Delta_+$, there exists $u_+ : \Delta_+ \rightarrow \R_{>0}$ such that $$\alpha_+(p) = -u_+(p) \alpha_-(p).$$
Hence, if $p \in \Delta_+$, $$f_0(p) = u_+(p) + \frac{1}{u_+(p)} \geq 2.$$
The case of $\Delta_-$ is similar.
\end{proof}

In particular, the divergence of $X$ (for any volume form or metric) along $\Delta_+$ (resp. $\Delta_-$) is positive (resp. negative). We also define smooth functions $g_\pm : M \rightarrow \R$ by
\begin{align*}
    \alpha_- \wedge d \alpha_+ &= g_+ \, \mathrm{dvol},\\
    \alpha_+ \wedge d \alpha_- &= g_- \, \mathrm{dvol},
\end{align*}
so that $g_0 = g_- + g_+$ and $f_0 = g_- - g_+$. One easily checks that
\begin{align}
    \mathcal{L}_X \alpha_- = g_- \, \alpha_- + \alpha_+, \label{eq:LX-}\\
    \mathcal{L}_X \alpha_+ = \alpha_- - g_+ \, \alpha_+. \label{eq:LX+}
\end{align}
Indeed, these formulae hold away from $\Delta$ since $\mathcal{L}_X \alpha_\pm$ vanish along $X$ and can be written as linear combinations of $\alpha_\pm$. The coefficients can be determined by wedging $\mathcal{L}_X \alpha_\pm$ with $\alpha_\pm$. Since $\Delta$ has empty interior and all of the quantities involved are smooth, these formulae hold along $\Delta$ as well.

\begin{rem} \label{rem:rank}
Let $p \in \Delta_+$. Since $X$ is tangent to $\xi_+$ (and $\xi_-$), the linearization $d_pX$ of $X$ at $p$ has its image contained in $\xi_+(p)$. Indeed:
$$\alpha_+ \circ d_pX = \mathcal{L}_X \alpha_+ (p) = d_p(\alpha_+(X)) + d \alpha_+\big(X(p),\cdot \,\big) = 0.$$
Let $L_p X : \xi_+(p) \rightarrow \xi_+(p)$ denote the restriction of $d_pX$ to $\xi_+(p)$. Then 
$$\mathrm{tr}(L_p X) = \mathrm{tr}(d_pX) = \mathrm{div}_p X = f_0(p) \geq 2$$
by Lemma~\ref{lem:f0}, hence $\mathrm{rank}(L_p X) \in \{1,2\}$. If $\mathrm{rank}(L_p X) = 2$, then the same holds for any nearby $p' \in \Delta_+$ and the constant rank theorem implies that $\Delta_+$ is a $1$-dimensional embedded submanifold of $M$ near $p$. Moreover, $\Delta_+$ is transverse to $\xi_+$ near $p$. If $\mathrm{rank}(L_p X)=1$, the situation is more complicated. Certainly, $\Delta_+$ is `at most $2$-dimensional' near $p$. It might be a $2$-dimensional embedded submanifold near $p$: the contact pair defined by the contact forms
\begin{align*}
    \alpha_\pm = \pm dz - y \, dx
\end{align*}
on $\R^3$ with standard coordinates $(x,y,z)$ has $\Delta_+ = \{y=0\}$. However, if $\mathrm{rank}\big(d_pX\big) = 2$, then $\Delta_+$ is $1$-dimensional near $p$ and \emph{tangent} to $\xi_+(p)$ at $p$.

Besides, for a \emph{generic} contact pair $(\xi_-, \xi_+)$, $\mathrm{rank}\big(d_pX\big) =2$ at every $p \in \Delta_+$, and $\mathrm{rank}(L_pX) = 2$ at all but finitely many points in $\Delta_+$. It follows that $\Delta_+$ is smoothly embedded link in $M$ that is transverse to $\xi_+$ away from finitely many isolated points. The same properties hold for $\Delta_-$.
\end{rem}

%%%%%%%%%%%%%%%%%%%%%%%%%%%%%%%%%%%%%%%%%%%%%%%%%%%%%%%%%%%%%%%%%%%%%
        \subsection{Invariant plane fields}
%%%%%%%%%%%%%%%%%%%%%%%%%%%%%%%%%%%%%%%%%%%%%%%%%%%%%%%%%%%%%%%%%%%%%
While $\Delta$ is a smoothly embedded link in $M$ for a generic contact pair $(\xi_-, \xi_+)$, we won't need to make such an assumption in the following proposition which generalizes~\cite[Proposition 2.4]{CF11}. We fix an auxiliary volume form $\mathrm{dvol}$ on $M$ which determines a preferred of contact form $\alpha_\pm$ for $\xi_\pm$ as in the previous section.

\begin{prop}\label{prop:distrib}
Let $(\xi_-, \xi_+)$ be \underline{any} contact pair on $M$. 

\begin{itemize}
    \item There exists a unique continuous plane field $\eta_u$ on $M \setminus \Delta_-$ satisfying
        \begin{enumerate}
            \item $X$ is tangent to $\eta_u$ and $\eta_u$ is invariant under the flow of $X$,
            \item Away from $\Delta$, $\eta_u\setminus \{0\}$ lies in the cone 
            $$C_+ \coloneqq \big\{ v \in T_pM \ \vert \  \alpha_\pm(v) > 0 \big\}\cup \big\{ v \in T_pM \ \vert \  \alpha_\pm(v) < 0 \big\}$$
            and along $\Delta_+$, $\eta_u$ coincides with $\xi_\pm$.
        \end{enumerate}
    Moreover, there exists a continuous $1$-form $\alpha_s$ on $M$ such that $\ker \alpha_s = \eta_u$ on $M \setminus \Delta_-$, $\alpha_s = 0$ on $\Delta_-$, $\alpha_s$ is differentiable along $X$ and $$\mathcal{L}_X \alpha_s = r_s \, \alpha_s,$$ where $r_s : M \rightarrow \R$ is a continuous function defined everywhere on $M$.

    \item Similarly, there exists a unique continuous plane field $\eta_s$ on $M \setminus \Delta_+$ satisfying
    \begin{enumerate}
            \item $X$ is tangent to $\eta_s$ and $\eta_s$ is invariant under the flow of $X$,
            \item Away from $\Delta$, $\eta_s\setminus \{0\}$ lies in the cone 
            $$C_- \coloneqq \big\{ v \in T_pM \ \vert \  \pm \alpha_\pm(v) > 0 \big\} \cup  \big\{ v \in T_pM \ \vert \  \pm \alpha_\pm(v) < 0 \big\}$$
            and along $\Delta_-$, $\eta_s$ coincides with $\xi_\pm$.
        \end{enumerate}
    Moreover, there exists a continuous $1$-form $\alpha_u$ on $M$ such that $\ker \alpha_u = \eta_s$ on $M \setminus \Delta_+$, $\alpha_u = 0$ on $\Delta_+$, $\alpha_u$ is differentiable along $X$ and $$\mathcal{L}_X \alpha_u = r_u \, \alpha_u,$$ where $r_u : M \rightarrow \R$ is a continuous function defined everywhere on $M$.
\end{itemize}

The $1$-forms $\alpha_u$ and $\alpha_s$ can be chosen so that the assignment 
$$(\xi_-, \xi_+) \longmapsto (\alpha_s, \alpha_u)$$ defines a continuous map from the space of (smooth) positive contact pairs on $M$ endowed with the $C^\infty$ topology to the space of pairs of continuous $1$-forms on $M$, endowed with the $C^0$ topology. We can further assume that
\begin{align*}
\alpha_s \wedge \alpha_u = \alpha_- \wedge \alpha_+, \qquad   r_s < r_u.
\end{align*}
\end{prop}

\begin{figure}[h]
    \centering
    \includegraphics[width=0.5\linewidth]{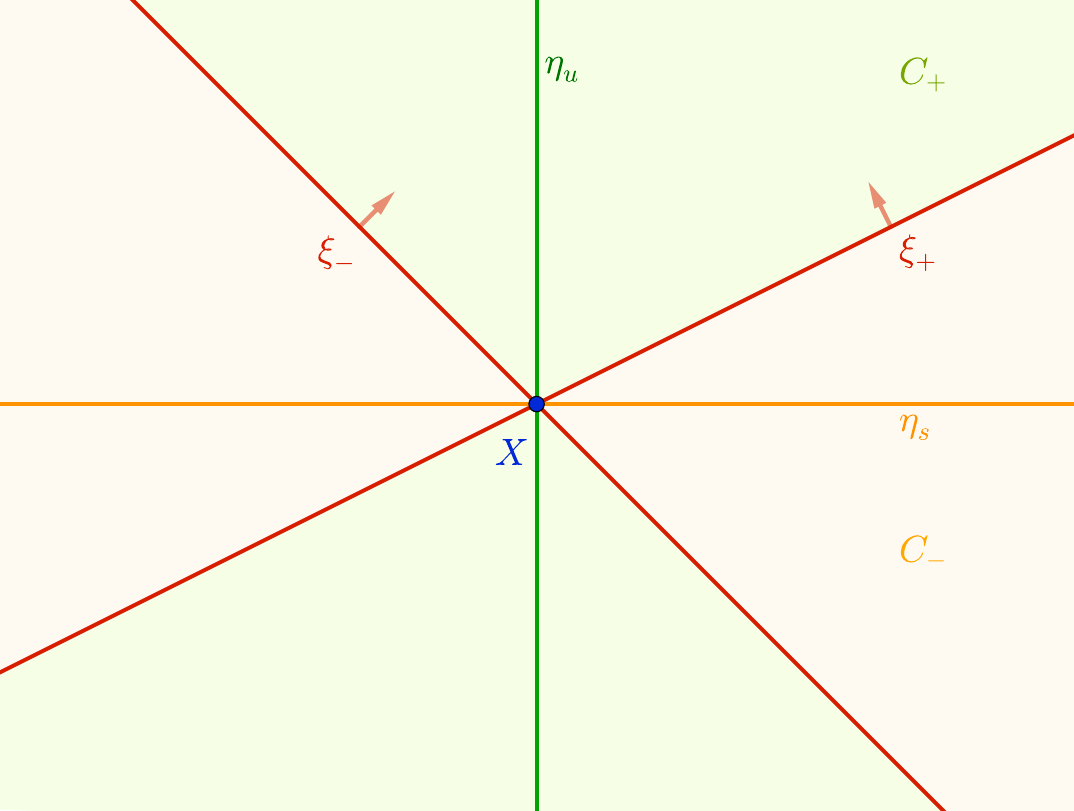}
    \caption{Away from $\Delta$: the contact structures $\xi_\pm$ with their respective coorientations, the cones $C_\pm$, and the plane fields $\eta_u$ and $\eta_s$. The vector field $X$ points toward the reader.}
    \label{fig:cones}
\end{figure}

We call $\eta_u$ (resp. $\eta_s$) the \textbf{(weak-)unstable} (resp. \textbf{(weak-)stable}) distribution of $X$. We call $r_u$ (resp. $r_s$) the \textbf{expansion rate} in the weak-unstable (resp. weak-stable) direction of $X$ (with respect to the specific choices of $\alpha_u$ and $\alpha_s$). The plane fields $\eta_u$ and $\eta_s$ are both locally (but not necessarily uniquely) integrable on $M \setminus \Delta$. See Figure~\ref{fig:cones} for a sketch of the relative positions of $\xi_\pm$, $\eta_u$ and $\eta_s$ away from $\Delta$.

\begin{rem}
    If $\Delta = \varnothing$, i.e., if $(\xi_-, \xi_+)$ is a bicontact structure, we recover the result of Eliashberg--Thurston~\cite{ET} and Mitsumatsu~\cite{M95} which states that the flow $\Phi$ of $X$ is conformally/projectively Anosov. In that case, as their names suggest, $\eta_u$ and $\eta_s$ are the weak-unstable and weak-stable distributions of $\Phi$, which are locally integrable but not necessarily uniquely integrable, see~\cite[Example 2.2.9]{ET}. In the more general setting where $\Delta \neq \varnothing$, we could call the flow $\Phi$ of $X$ a \textbf{singular projectively Anosov flow}.
\end{rem}

\begin{proof}[Proof of Proposition~\ref{prop:distrib}]
We look for $1$-forms $\alpha_u$ and $\alpha_s$ of the form
\begin{align}
\alpha_u &= e^{-\sigma_u} \alpha_- + e^{\sigma_u} \alpha_+ \label{eq:alphau}\\
\alpha_s &= e^{-\sigma_s} \alpha_- - e^{\sigma_s} \alpha_+ \label{eq:alphas}
\end{align}
where $\sigma_u, \sigma_s : M \rightarrow \R$ are continuous functions which are continuously differentiable along $X$. These $1$-forms have to satisfy
$$ \alpha_s \wedge \mathcal{L}_X \alpha_s = \alpha_u \wedge \mathcal{L}_X \alpha_u = 0,$$ which after some elementary computations is equivalent to 
\begin{align}
X \cdot \sigma_u &= \sinh(2\sigma_u) + \frac{1}{2} g_0, \label{eq:sigmau}\\
X \cdot \sigma_s &= - \sinh(2\sigma_s) + \frac{1}{2} g_0. \label{eq:sigmas}
\end{align}
For $p \in M$, we consider the ODE $$(E_p) : \quad \dot{y}(t) = \sinh(2y(t)) + \frac{1}{2} g_0\big(\phi^t_X(p)\big),$$
where $\phi^t_X$ denotes the flow of $X$. After embedding $M$ in some $\R^N$ by the Whitney embedding theorem, we obtain a system of ODEs of the form of the ones studied in Appendix~\ref{appsec:ODE}, with $$F(p,y,t) \coloneqq \sinh(2y) + \frac{1}{2} g_0(\phi^t_X(p)).$$
This function $F$ is easily seen to satisfy conditions (C1) and (C2) from Appendix~\ref{appsec:ODE}, so Lemma~\ref{lem:ODE} applies and we define $\sigma_u(p)$ as the unique initial value of $(E_p)$ such that the corresponding maximal solution is defined on $\R$ and is bounded. Moreover, no other solution to $(E_p)$ is defined for all positive times by Lemma~\ref{lem:ODE3}. We readily get that $t \mapsto \sigma_u(\phi^t_X(p))$ is differentiable and is the unique bounded solution to $(E_p)$ with initial value $\sigma_u(p)$, so $\sigma_u$ is differentiable along $X$ and solves~\eqref{eq:sigmau}. Note that on $\Delta$, $\sigma_u$ is given by $$\sigma_u(p) = - \frac{1}{2} \sinh^{-1} \left(\frac{g_0(p)}{2}\right).$$
Notice that this shows both the existence and uniqueness of $\eta_u$. A similar argument also shows that the assignment $(\xi_-, \xi_+) \mapsto \sigma_u \in C^0(M)$ is continuous. Indeed, it is enough to show that if $(\xi^n_-, \xi^n_+)$, $n \geq 0$, is a sequence of positive contact pairs on $M$ converging to $(\xi_-, \xi_+)$, then the corresponding sequence of functions $\sigma^n_u$, $n \geq 0$, converges to $\sigma_u$ uniformly. The corresponding sequence of functions $g^n_0$, $n \geq 0$, converges to $g_0$ and the corresponding sequence of vector fields $X_n$, $n \geq 0$, converges to $X$, both in the $C^\infty$ topology. Therefore, we obtain a sequence of functions $F_n : M \times \R \times \R \rightarrow \R$ defined by
$$F_n(p,y,t) \coloneqq \sinh(2y) + \frac{1}{2} g^n_0\big( \phi^t_{X_n}(p)\big)$$
which converges to $F$ in the $C^\infty_\mathrm{loc}$ topology. One easily checks that the hypotheses of Lemma~\ref{lem:ODE2} are satisfied, hence $\sigma^n_u$ converges to $\sigma_u$ uniformly. As a result, $\alpha_u$ depends continuously on $(\xi_-, \xi_+)$.

We proceed similarly for $\sigma_s$ by replacing $X$ with $-X$. 

Away from $\Delta_+$, $\alpha_u$ does not vanish and defines a continuous plane field $\eta_s$. Moreover, $\alpha_u$ vanishes along $\Delta_+$. Indeed, there exists a smooth function $u_+ : \Delta_+ \rightarrow \R_{>0}$ such that for $p \in \Delta_+$, $$\alpha_+(p) = - u_+(p) \, \alpha_-(p).$$
This implies that if $p \in \Delta_+$, $$g_0(p) = u_+(p) - \frac{1}{u_+(p)},$$
and since $$g_0(p) = e^{-2\sigma_u(p)} - e^{2\sigma_u(p)},$$
we obtain that $$e^{\sigma_u(p)} = \frac{1}{\sqrt{u_+(p)}},$$
hence $$\alpha_u(p) = \sqrt{u_+(p)} \, \alpha_-(p) - \frac{1}{\sqrt{u_+(p)}} u_+(p) \, \alpha_-(p) = 0.$$
Similarly, $\alpha_s$ does not vanish away from $\Delta_-$ and defines a continuous plane field $\eta_u$ on $M \setminus \Delta_-$, while it vanishes on $\Delta_-$. Analogous computations show that $\eta_u = \xi_\pm$ along $\Delta_+$ and $\eta_s = \xi_\pm$ along $\Delta_-$.

The expansion rates $r_u$ and $r_s$ can now be computed, using~\eqref{eq:LX-}, \eqref{eq:LX+} and~\eqref{eq:sigmau}, \eqref{eq:sigmas}:
\begin{align}
r_u &= \cosh(2\sigma_u) + \frac{1}{2}f_0, \label{eq:ru}\\
r_s &= -\cosh(2\sigma_s) + \frac{1}{2}f_0 \label{eq:rs}
\end{align}
Therefore, $r_u$ and $r_s$ are continuous on $M$, and $$r_u - r_s = \cosh(2\sigma_u) + \cosh(2 \sigma_s) \geq 2.$$
Finally, 
$$\alpha_s \wedge \alpha_u = 2 \cosh(\sigma_u - \sigma_s) \, \alpha_- \wedge \alpha_+,$$
but we can replace $\alpha_u$ and $\alpha_s$ with 
\begin{align*}
\alpha_u' \coloneqq \frac{1}{\sqrt{2 \cosh(\sigma_u - \sigma_s)}} \, \alpha_u, \\
\alpha_s' \coloneqq \frac{1}{\sqrt{2 \cosh(\sigma_u - \sigma_s)}} \, \alpha_s,
\end{align*}
so that $\alpha'_s \wedge \alpha'_u = \alpha_- \wedge \alpha_+$ and $r'_u - r'_s = r_u - r_s \geq 2$.
\end{proof}

\begin{rem}
    It can be shown that $r_u = 0$ along $\Delta_-$ and $r_s = 0$ along $\Delta_+$. This is independent of the choice of $\alpha_u$ and $\alpha_s$.
\end{rem}

%%%%%%%%%%%%%%%%%%%%%%%%%%%%%%%%%%%%%%%%%%%%%%%%%%%%%%%%%%%%%%%%%%%%%
            \subsection{Local integrability} \label{sec:locint}
%%%%%%%%%%%%%%%%%%%%%%%%%%%%%%%%%%%%%%%%%%%%%%%%%%%%%%%%%%%%%%%%%%%%%

To ensure local integrability of $\eta_u$ along $\Delta_+$ and of $\eta_s$ along $\Delta_-$, we will consider contact pairs satisfying some natural and \emph{generic} conditions.

\begin{defn} \label{def:reg}
The contact pair $(\xi_-, \xi_+)$ is \textbf{regular} if the following conditions hold.
\begin{itemize}
\item[(R1)] The singular set $\Delta \subset M$ is a smoothly embedded link,
\item[(R2)] There exist finitely many points $q_i$, $1 \leq i \leq n$, such that $\Delta$ has a quadratic contact with $\xi_\pm$ at $q_i$, and away from these points, $\Delta$ is transverse to $\xi_\pm$,
\item[(R3)] At every point of $\Delta$, the linearization of $X$ has rank $2$.
\end{itemize}
\end{defn}

It follows from standard transversality results that regular contact pairs form a comeagre set in the space of contact pairs. In particular, it is a dense subset by Baire's theorem. Writing $Q \coloneqq \{ q_i\}_{1 \leq i \leq n}$, it follows that points in $\Delta \setminus Q$ are of two types, depending on the sign of the determinant of the differential of $X$ in the normal direction:
\begin{itemize}
\item Either this determinant is positive, and the singularity is a \emph{source} (on $\Delta_+$) or a \emph{sink} (on $\Delta_-$), or
\item This determinant is negative, and the singularity is a \emph{saddle}.
\end{itemize}
Hence, we can write $$\Delta \setminus Q = \Delta_{\mathrm{so}} \sqcup \Delta_{\mathrm{si}} \sqcup \Delta_{\mathrm{sa}},$$ where $\Delta_{\mathrm{so}}$ (resp.
$\Delta_{\mathrm{si}}$, $\Delta_{\mathrm{sa}}$) is the disjoint union of open intervals in $\Delta$ corresponding to source (resp. sink, saddle) points.

Note that at every point $p \in \Delta_+$, the linearization of $X$ is tangent to $\eta_u(p) = \xi_+(p)$ by Remark~\ref{rem:rank}. Similarly, at every point $p \in \Delta_-$, the linearization of $X$ is tangent to $\eta_s(p)$.

The next proposition generalizes~\cite[Proposition 2.5]{CF11} by removing the normality assumption. We also fill some gaps and clarify some subtleties.

\begin{prop} \label{prop:locint}
Let $(\xi_-, \xi_+)$ be a regular contact pair. Then $\eta_u$ is locally integrable on $M \setminus  \Delta_- $, and $\eta_s$ is locally integrable on $M \setminus \Delta_+$. 
\end{prop}

\begin{proof} We show that $\eta_u$ is locally integrable at every point $p \in M \setminus \Delta_- $. We distinguish four cases.

\paragraph{Case 0:  $p \in M \setminus \Delta$.} Then $X(p) \neq 0$ and local integrability of $\eta_u$ at $p$ is standard, see~\cite[Proposition 2.5]{CF11}.

\paragraph{Case 1: $p \in \Delta_{\mathrm{so}}$.} By the Unstable Manifold Theorem, there exists a germ of smooth surface $S^u = S^u(p)$ passing through $p$ such that $X$ is tangent to $S^u$ and $T_p S^u = \eta_u(p) = \xi_\pm(p)$. Moreover, there exist smooth coordinates $(x,y,z)$ near $p$ in which the linearization of $X$ at $0$ is of one of the following forms:
\begin{itemize}
    \item \textit{(Node)} $X_1 = a x \partial_x + b y \partial_y$, where $a, b > 0$,
    \item \textit{(Degenerate node)} $X_1 = a \big(x \partial_x +  y \partial_y\big) + y\partial_x$, where $a> 0$,
    \item \textit{(Focus)} $X_1 = a \big(x \partial_x +  y \partial_y\big) + b \big( y \partial_x - x \partial_y \big)$, where $a> 0$, $b \neq 0$.
\end{itemize}
In any case, applying~\cite[Theorem (1.13)]{T74}, and noting that $X$ has a line of singularities passing through $0$, we obtain that $X$ is \emph{topologically} equivalent to the vector field $$X' = x \partial_x + y \partial_y.$$ In these last coordinates, $S^u \subset \{z=0\}$.

We claim that $\eta_u$ is tangent to $S^u$ in a neighborhood of $p$. Indeed, if $\eta_u(p')$ and $T_{p'}S^u$ disagree at a point $p'\in S^u$ near $p$, we consider a small disk $D$ passing through $p'$ and transverse to $X$. The intersection of $\eta_u$ with this disk yields a continuous vector field on $D$, for which we can find a small flow line $\gamma$ passing through $p'$. Because of the above (topological) normal form for $X$ near $p$, there exists a point $p'' \in \gamma$ which is contained in the unstable manifold $\widetilde{S}^u$ of a point $\widetilde{p} \in \Delta$ distinct from $p$. Under the backward flow of $X$, $p'$ converges to $p$ and $p''$ converges to $\widetilde{p}$. This and the topological normal form imply that the curve $\gamma_t \coloneqq \phi_X^{-t}(\gamma)$ obtained by flowing $\gamma$ along $X$ for time $-t$ uniformly converges to a segment on $\Delta_\mathrm{so}$ containing $p$ and $\widetilde{p}$ as $t \rightarrow +\infty$. Moreover, the length of $\gamma_t$ is uniformly bounded for $t \geq 0$. Indeed, normal hyperbolicity along $\Delta_\mathrm{so}$ implies that there exists a neighborhood $U$ of $p$ and a constant $C \geq 0$ such that for all $q \in U$ and $t \geq 0$, $\Vert d\phi_X^{-t}\Vert \leq C$. However, for $t$ large enough, such a $\gamma_t$ cannot be everywhere tangent to $\eta_u$ by the mean value theorem, which is a contradiction. As a result, $T S^u = \eta_u$ along $S^u$, and $S^u$ is a germ of integral surface of $\eta_u$ around $p$. This also shows that $\eta_u$ is \emph{uniquely} integrable at $p$, and the characteristic foliations of $\xi_\pm$ coincide on $S^u$ and are directed by $X$. The previous paragraph also implies the following property: if $S \subset M \setminus \Delta$ is a connected immersed surface tangent to $\eta_u$ and if the backward flow line of $X$ passing through some $p' \in S$ converges to $p \in \Delta_\mathrm{so}$, then $S$ is contained in the unstable manifold of $X$ at $p$, i.e., the saturation of $S^u(p)$ by the flow of $X$. This property will be crucial in Case 3 below.

\paragraph{Case 2: $p \in \Delta_{\mathrm{sa}}$.} The strategy of~\cite[Proof of Proposition 2.5]{CF11} can be adapted to this more general setting.\footnote{In~\cite{CF11}, the authors use that at a saddle point, $X$ is conjugated to its linearization; however, this conjugation is only \emph{continuous} in general!} By the Invariant Manifolds Theorem, there exist smooth coordinates $(x,y,z)$ near $p$ such that 
\begin{itemize}
\item The linearization $X_1$ of $X$ at $0$ is $X_1 = a x \partial_x - b y \partial_y$, where $a, b > 0$.\footnote{Since $X$ has positive divergence along $\Delta_+$, we even have $b < a$. Unlike in~\cite{CF11}, this won't be relevant in our proof.}  
\item $\Delta$ is contained in $\{x=y=0\}$,
\item The weak unstable manifold of $X$ is contained in $\{y=0\}$ and the strong unstable manifold of $X$ at $0$ is contained in $\{y=z=0\}$,
\item The weak stable manifold of $X$ is contained in $\{x=0\}$ and the strong stable manifold of $X$ at $0$ is contained in $\{x=z=0\}$.
\end{itemize}
It follows that in these coordinates, $\eta_u(0) = \mathrm{span} \{\partial_x, \partial_y \}$. We can also find small constants $\epsilon, \delta_x, \delta_y > 0$ such that the following properties hold. Both $\eta_u$ and $X$ are transverse to the disks
\begin{align*}
D_x^{\pm} &\coloneqq \left\{ x = \pm \epsilon, \ y^2 + z^2 \leq \delta_x^2 \right\} \\
D_y^{\pm} &\coloneqq \left\{ y = \pm \epsilon, \ x^2 + z^2 \leq \delta_y^2 \right\}
\end{align*}
and all the flow lines of $X$ starting on $D_y^{\pm}$ either intersect $D_x^{\pm}$ or converge to $\Delta$. By flowing $D_y^{\pm}$ along $X$, we obtain four smooth maps 
$$f^{i,j} : D_y^i \cap \{ j x > 0 \} \longrightarrow D_x^j,$$
where $i,j \in \{-, +\}$. Moreover, $df^{i,j}$ sends the intersection of $D_y^i$ with $\eta_u$ to the intersection of $D_x^j$ with $\eta_u$, and if we write $$p_x^{\pm} \coloneqq (\pm \epsilon, 0,0), \qquad p_y^{\pm} \coloneqq (0, \pm \epsilon, 0),$$
then
$$\lim_{p \rightarrow p_y^i} f^{i,j}(p) = p_x^j$$
for every $i,j \in \{-, +\}$. The latter can be seen by \emph{topologically} conjugating $X$ with its linearization at $0$, see~\cite[Theorem (1.13)]{T74}. See Figure~\ref{fig:hyperbolic} for a lower dimensional sketch.

\begin{figure}[h]
    \centering
    \includegraphics[width=0.65\linewidth]{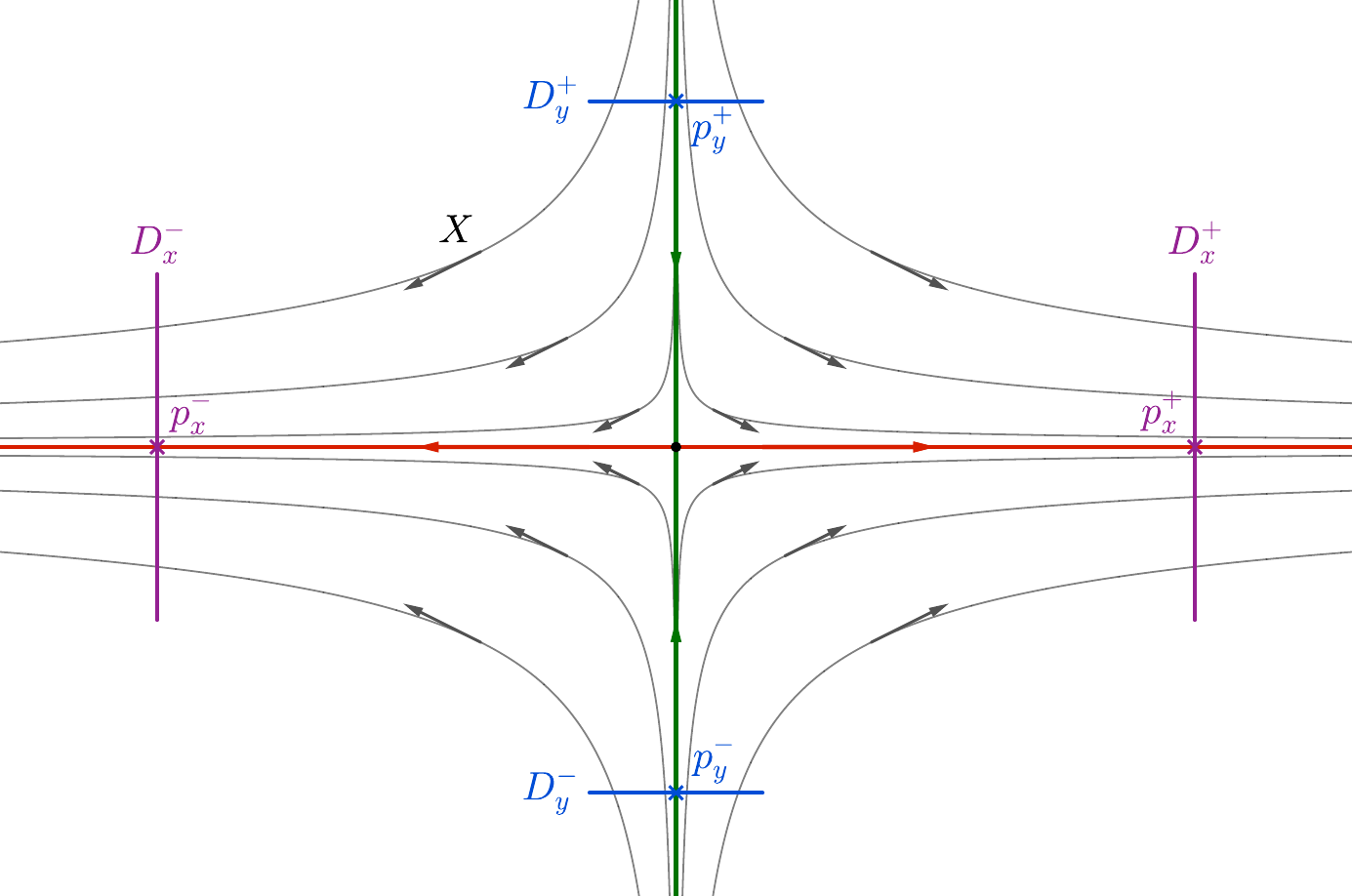}
    \caption{The $2$-dimensional analogue of Case 2.}
    \label{fig:hyperbolic}
\end{figure}

Since $\eta_u$ intersects $D_y^{\pm}$ along a continuous vector field, we can find small integral curves $\gamma_y^{\pm} : (-\delta, \delta) \rightarrow D_y^{\pm}$ tangent to $\eta_u$ and with $ \gamma_y^{\pm}(0) = p_y^{\pm}$, where $0 < \delta < \delta_y$. Flowing these curves along $X$, i.e., applying the maps $f^{i,j}$, we obtain four curves 
$$ \gamma_x^{i,j} : I^i \longrightarrow D_x^j,$$
where $I^- = (-\delta, 0)$ and $I^+ = (0, \delta)$, which are both tangent to $\eta_u$. Moreover, $$\lim_{t \rightarrow 0} \gamma_x^{i,j}(t) = p_x^j,$$ so the union of $\gamma_x^{i,j}$ and $\gamma_x^{-i,j}$ can be continuously extended to a curve $$\gamma_x^j : (-\delta, \delta) \longrightarrow D_x^j$$
with $\gamma_x^j(0) = p_x^j$. The key point is that since the $\gamma_x^{\pm}$'s are continuous, continuously differentiable away from $0$, and their derivatives extend continuously at $0$, these curves are in fact $C^1$.

The union of all of the flow lines of $X$ starting from $\gamma_y^{\pm}$ produces two open surfaces $S^\pm$ such that $$S \coloneqq S^- \cup \big\{\vert x \vert \leq \epsilon, \ y=z=0 \big\} \cup S^+$$ is a surface which can be written as the graph of a continuous function $f : U \subset \R^2 \rightarrow \R$ near $0$. Here, $U$ denotes an open neighborhood of $0$ in $\R^2$. Note that $S \cap D_x^{\pm} = \gamma_x^{\pm}$. By the previous paragraph, $f$ is differentiable away from $0$ and its graph is tangent to $\eta_u$. The differentiability of $f$ along $\big\{\vert x \vert \leq \epsilon, \ y=z=0 \big\} \setminus \{0\}$ follows from the differentiability of $\gamma_x^{\pm}$. Moreover, $f$ has (vanishing) partial derivatives at $0$, and its partial derivatives are continuous on $U$. As a result, $f$ is $C^1$ on $U$ and is everywhere tangent to $\eta_u$, so $S$ is a local integral surface for $\eta_u$ near $0$.

\paragraph{Case 3: $p \in Q$.} This is the most delicate case. For the sake of clarity, we divide it into three steps.

\begin{itemize}[leftmargin=*]
    \item \textit{Step 0 : setting things up.} We start with a \emph{topological} description of the flow lines of $X$ near $p$. Following~\cite[Proof of Lemme 2.2]{CF11}, there exist smooth coordinates $(x,y,z)$ around $p$ such that the linearization of $X$ at $0$ is of the form $$X_1 = a x \partial_x + z\partial_y,$$ where $a > 0$, and $\Delta$ is included in $\{ x=0, \ z+y^2 = 0\}$. We can further assume that the strong unstable manifold of $X$ at $0$ is included in $\{y = z=0\}$. By~\cite[Theorem (1.11)]{T74}, for every $k \geq 1$, there exists a germ $S^c$ of center manifold of class $C^k$ passing through $0$ which is tangent to the plane $\{x=0\}$ at $0$, and such that $X$ is tangent to $S^c$ along it. Let us take $k=10$ for safety. We denote by $Y = Y(y,z)$ the restriction of $X-X_1$ to $\{x=0\}$, which is of class $C^k$ and has vanishing $1$-jet at $0$. By~\cite[Theorem (1.13)]{T74}, $X$ is \emph{topologically} equivalent near $0$ to the vector field
    \begin{align*} %\label{eq:normal}
    X' = ax \partial_x + z \partial_y + Y(y,z).
    \end{align*}
    Notice that $S^c$ contains $\Delta$. Writing $Z(y,z) = z \partial_y + Y(y,z)$, $Z$ is a $C^k$ vector field in the $(y,z)$-plane whose linearization at $0$ is $z \partial_y$. Moreover, we can assume that in these new coordinates, $\Delta$ is the line $\{x=z=0\}$. This implies that for every $y$ sufficiently small, $Y(y,0) \equiv 0$, so we can write $Y(y,z) = zQ(y,z)$, where $Q$ is a $C^{k-1}$ vector field vanishing at $(0,0)$. We consider the vector field
    $$G(y,z) \coloneqq \partial_y + Q(y,z),$$
    so that
    $$Z = z\left( \partial_y + Q(y,z)\right).$$

    Here, $G$ is mostly horizontal. Moreover, the points in $\{y> 0, z=0\}$ correspond to singularities of $X$ in $\Delta_\mathrm{so}$, while the points in $\{y< 0, z=0\}$ correspond to singularities of $X$ in $\Delta_\mathrm{sa}$. Therefore, the flow lines of $G$ passing through a point $(y, 0)$ with $y<0$ are \emph{topologically} negatively transverse to the $y$-axis, while the flow lines of $G$ passing through $(y,0)$ with $y>0$ are topologically positively transverse to the $y$-axis. One easily deduces that the flow line of $G$ passing through $(0,0)$ only intersects the $y$-axis at $(0,0)$ and stays in the half-plane $\{z>0\}$ away from $(0,0)$. Moreover, every flow line of $G$ passing through $(y,0)$ for $y<0$ sufficiently close to $0$ intersects $\{y>0, z=0\}$. Those flow lines correspond to connections from $\Delta_\mathrm{so}$ to $\Delta_\mathrm{sa}$ in a neighborhood of $p$. The other stable branches of $\Delta_\mathrm{sa}$ near $p$ all escape the chosen coordinate neighborhood of $p$. 
    
    Since $G$ has a unique flow line passing through $0$, $Z$ has a unique flow line converging to $0$ in positive time, and a unique flow line converging to $0$ in negative time. In summary, we obtain that 
    \begin{itemize}
        \item There exists a unique flow line of $X$ converging to $0$, and
        \item The set of the points near $0$ which are converging to $0$ under the backward flow of $X$ is a \emph{topological} surface with boundary, whose boundary contains $0$ and is contained in $\{y=z=0\}$.
    \end{itemize}
    Unfortunately, we cannot conclude that this half-surface is tangent to $\eta_u$ since we don't know that it is $C^1$ yet.
    
    We now restrict our attention to the box 
    $$\mathsf{B}_\epsilon \coloneqq \left\{ \vert x \vert \leq \epsilon, \ \vert y \vert \leq \epsilon + \sqrt{\epsilon}, \  \vert z \vert \leq \epsilon \right\}$$ for $\epsilon > 0$ small enough. Still from~\cite[Proof of Lemme 2.2]{CF11}, for $\epsilon$ sufficiently small, $X$ is transverse to the vertical faces $$F_x^{\pm} \coloneqq \{ x = \pm \epsilon \} \cap \mathsf{B}_\epsilon, \qquad  F_y^{\pm} \coloneqq\{ y = \pm (\epsilon + \sqrt{\epsilon})\} \cap \mathsf{B}_\epsilon,$$ and enters $\mathsf{B}_\epsilon$ along the face $F_y^-$ and exits along the other vertical faces. Here, $X$ has source singularities along $\Delta \cap \{y>0\}$ and saddle singularities along $\Delta \cap \{y < 0\}$. Moreover, every backward flow line of $X$ starting on $\partial \mathsf{B}_\epsilon$ either converges to a point in $\Delta$, or exits $\partial \mathsf{B}_\epsilon$ through the face $F_y^-$ or the top face $F_z^+ \coloneqq \{z = \epsilon\} \cap \mathsf{B}_\epsilon$.
    
    We can find continuous simple curves $\gamma_u^\pm \subset F_x^\pm$ starting at $x^\pm \coloneqq ( \pm \epsilon, 0,0)$ and ending on the interior of $F_x^\pm \cap \{z = -\epsilon\}$, the bottom edge of $F_x^\pm$, such that under the flow of $-X$, every point in $\gamma_u^\pm \setminus \{x^\pm\}$ converges to a point in $\Delta \cap \{y> 0\}$. In other words, points in $\gamma_u^\pm \setminus \{x^\pm\}$ lie in the intersection of the unstable disks of the source points in $\Delta \cap \{y> 0\}$ with $F_x^\pm$. Similarly, we can find a continuous simple curve $\gamma_b \subset F_y^-$ intersecting the interior of the bottom edge of $F_y^-$ such that every point in $\gamma_b$ converges to a point in $\Delta \cap \{y < 0\}$ under the flow of $X$, except one of its endpoints $y_b$ which converges to $0$. In other words, every point of $\gamma_b$ lies in the stable manifold of a saddle point in $\Delta$, or in the stable branch of $X$ at $0$. Notice that $\gamma_b$ is contained in the intersection of $S^c$ with $F^-_y$.

    We define a smoothing $\Sigma$ of the union of the exit vertical faces $F_x^- \cup F_y^+ \cup F_x^+$ as follows. We first consider a $C^0$-small smoothing $\widetilde{c}$ of the piecewise linear planar curve
    $$c \coloneqq \big\{x =-\epsilon, \ \vert y \vert \leq \epsilon + \sqrt{\epsilon}\big\} \cup \big\{ \vert x \vert \leq \epsilon, \ y = \epsilon + \sqrt{\epsilon}\big\} \cup \big\{x = \epsilon, \ \vert y \vert \leq \epsilon + \sqrt{\epsilon}\big\}$$ such that $\widetilde{c}$ coincides with $c$ away from a small neighborhood of the corners of $c$ and is contained in the rectangle $\big\{ \vert x \vert \leq \epsilon, \ \vert y \vert \leq  \epsilon + \sqrt{\epsilon} \big\}$. We then set $$\Sigma \coloneqq \left\{ (x,y,z) \in \mathsf{B}_\epsilon \ \vert \  (x,y) \in \widetilde{c}, \ \vert z \vert \leq \epsilon \right\}.$$ An appropriate choice of $\widetilde{c}$ ensures that $X$ is positively transverse to $\Sigma$, and the curves $\gamma_u^\pm$ lie on $\Sigma$. If the curve $\widetilde{c}$ is parametrized by the variable $s \in [-1,1]$, we obtain coordinates $(s,z)$ on $\Sigma$ in which $\Sigma$ is simply the rectangle $[-1,1] \times [-\epsilon, \epsilon]$. We will assume that $x^+$ and $\gamma_u^+$ are contained in $[-1, 0) \times [-\epsilon, \epsilon]$, and $x^-$ and $\gamma_u^-$ are contained in $(0,1] \times [-\epsilon, \epsilon]$; see Figure~\ref{fig:sigma}.

\begin{figure}[t]
    \centering
    \includegraphics[width=1 \linewidth]{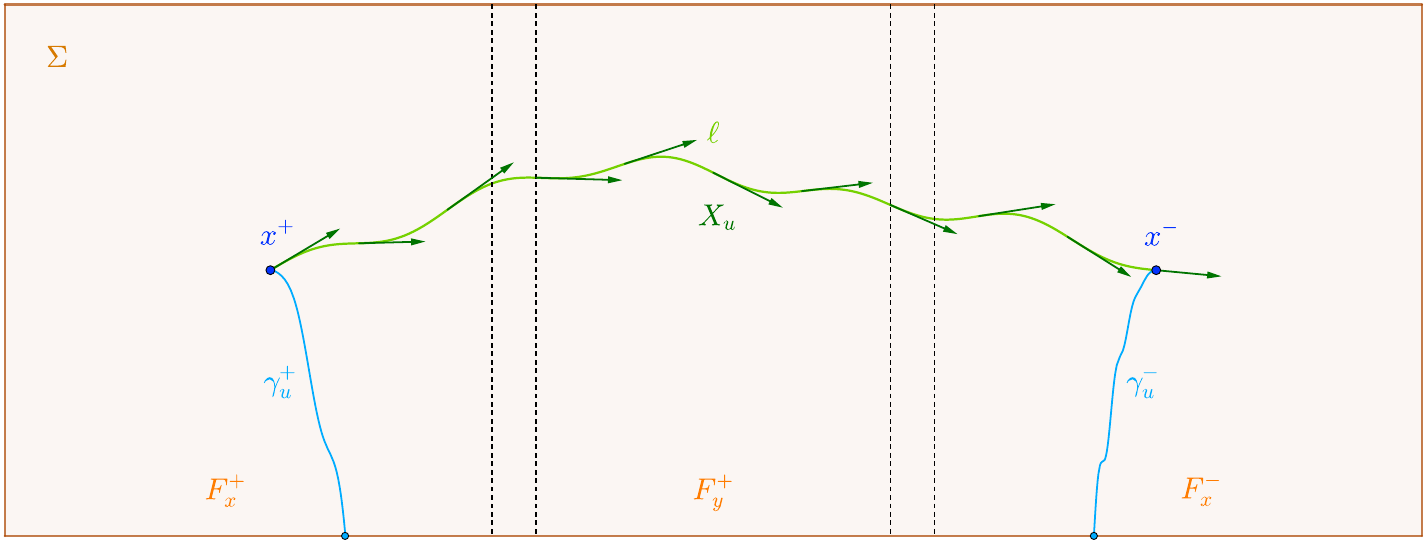}
    \caption{The surface $\Sigma$ and a pre-lasso $\ell$.}
    \label{fig:sigma}
\end{figure}

    In the next steps, we construct a germ of integral surface for $\eta_u$ passing through $0$ as follows. We construct two (germs of) integral surfaces with boundary $S^\pm \subset \mathsf{B}_\epsilon$, such that the unstable manifold of $X$ at $0$ is contained in the boundaries of $S^\pm$, and such that $S=S^- \cup S^+$ is the graph of a continuous function near $0$. Then, the same argument as in Case 2 will show that $S$ is $C^1$.

    \item \textit{Step 1: construction of $S^-$.} We adapt the reasoning of Case 2. Choosing $\epsilon$ sufficiently small, the center manifold $S^c$ intersects the face $F_y^-$ along a $C^1$ curve $\gamma_c$; this curve divides $F_y^-$ into two halves, and contains the curve $\gamma_b$. The specific form of the linearization of $X$ at $0$, and the topological normal form from Step 0 imply the following property: if $(x_n)_n$ is a sequence of points in $F_y^- \setminus \gamma_c$ converging to $y_b$ from the $x>0$ (resp. $x<0$) side at a sufficiently small asymptotic slope with respect to the plane $\{z=0\}$, then for $n$ large enough, the flow line of $X$ passing through $x_n$ exits $\mathsf{B}_\epsilon$ along $F_x^+$ (resp. $F_x^-$), and the intersection of this flow line with the corresponding face converges to $x^+$ (resp. $x^-$) when $n$ goes to $\infty$. Since $\eta_u$ is transverse to $S^c$ at $0$, it intersects $\gamma_c$ transversally at $y_b$ along the face $F_y^-$. We can further shrink $\epsilon$ to make the angle between $\eta_u(y_b)$ and the plane $\{z=0\}$ sufficiently small. As in Case 2, we then choose a flow line of $\eta_u \cap F_y^-$ passing through $y_b$. We flow it along $X$ to obtain a surface whose union with $\{y=z=0\}$, denoted by $S^-$, is a surface with boundary tangent to $\eta_u$.

    \item \textit{Step 2: construction of $S^+$.} The idea of the construction of $S^+$ is to find a $C^1$ curve $\ell$ on $\Sigma$ starting at $x^+$ and ending at $x^-$ so that $\ell$ is tangent to $\eta_u$, and every point on $\ell$ converges to $0$ under the backward flow of $X$. We will refer to such a curve $\ell$ as a \emph{lasso}. We first construct a \emph{pre-lasso}, i.e., a curve satisfying the same assumptions as a lasso except for the condition on the convergence to $0$; see Figure~\ref{fig:sigma}. We then construct a lasso from a suitable sequence of pre-lassos. The integral surface $S^+$ will be obtained by flowing a lasso along the backward flow of $X$.

    The intersection of $\eta_u$ with $\Sigma$ defines a continuous vector field $X_u$ on $\Sigma$ which is almost horizontal (its vertical coordinate is small compared to its horizontal coordinate). Coorienting $\eta_u$ so that $\partial_z$ is positively normal to $\eta_u$, this vector field essentially points ``from $x^+$ to $x^-$''. The key fact is the following: any integral curve of $X_u$ on $\Sigma$ passing through $x^+$ or $x^-$ is disjoint from $\gamma_u^+ \setminus \{x^+\}$ and $\gamma_u^- \setminus \{x^-\}$. Otherwise, a point on such a curve would converge to a point $p_s \in \Delta_\mathrm{so}$, and flowing this curve along the backward flow of $X$ would produce an integral surface for $\eta_u$ passing through $p_s$. By Case 1, we know that there exists a unique such integral surface near $p_s$, namely, the unstable manifold of $X$ at $p_s$. The topological normal form from Step 0 implies that our surface cannot contain $x_\pm$, a contradiction. We then argue that there exists a flow line of $X_u$ passing through $x^+$ which exits $\Sigma$ along 
    $$\partial_v^- \Sigma \coloneqq \partial \Sigma \cap \big\{ x = -\epsilon, \ y = -(\epsilon + \sqrt{\epsilon}) \big\}.$$ Indeed, if such a flow line $\gamma$ exits $\partial \Sigma$ along the bottom boundary
    $$\partial_b \Sigma \coloneqq \partial \Sigma \cap \{z = -\epsilon\},$$
    then it would intersect the unstable manifold of a source point, which is prohibited by the previous argument. If every such $\gamma$ exits $\Sigma$ along the interior of its top boundary
    $$\partial_t \Sigma \coloneqq \partial \Sigma \cap \{z = \epsilon\},$$
    we consider the set $\mathcal{S}$ of values of $s \in (-1,1)$ such that $(s, \epsilon)$ is an exit point of such a flow line $\gamma$, and we define $s_{\sup} \coloneqq \sup \mathcal{S}$. Picking a sequence $s_n \in \mathcal{S}$ converging to $s_{\sup}$, we consider a sequence of flow lines $\gamma_n$ for $X_u$ all starting at $x^+$ and exiting $\Sigma$ at $(s_n, \epsilon)$.  This sequence of curves is uniformly Lipschitz, since $X_u$ is continuous, so up to passing to a subsequence, we can assume that $(\gamma_n)_n$ converges uniformly to a curve $\gamma_{\sup}$ which is still a flow line of $X_u$ starting at $x^+$,\footnote{This is because a continuous function $\gamma : I \rightarrow \R$ defined on an interval containing $0$ is a $C^1$ solution to the ODE $\dot{y}(t) = F(y(t))$ if and only if it satisfies $y(t) = y(0) +\int_0^t F(y(s)) \, ds$ for every $t\in I$.} and which exits $\Sigma$ at $(s_{\sup}, \epsilon)$. By assumption, we have $s_{\sup} \leq 1$. There exists a small $\delta > 0$ such that the point $x'\coloneqq(s_{\sup}, \epsilon - \delta)$ escapes $\mathsf{B}_\epsilon$ when flown backward along $X$, and escapes $\Sigma$ along its top or lateral boundary when flown along $X_u$. We now consider a flow line $\gamma'$ of $-X_u$ starting at $x'$. It cannot intersect $\gamma_{\sup}$, since otherwise, we could merge a portion of $\gamma_{\sup}$ with a portion of $\gamma'$ and obtain a flow line of $X_u$ starting at $x^+$ and escaping $\Sigma$ along a point on the top boundary further right than $(s_{\sup}, \epsilon)$, or on the lateral boundary of $\Sigma$, and both of these are excluded by hypothesis. Therefore, $\gamma'$ has to pass below $x^+$ and intersect $\gamma_u^+$. This is also excluded, since some point on $\gamma'$ converges to $0$ under the backward flow of $X$ by the topological normal form of $X$ near $p$, which contradicts a previous argument.

    To summarize, we now know that there exists a flow line of $X_u$ passing through $x^+$ and which exits $\Sigma$ on its lateral boundary $\partial_v^- \Sigma$. Such a flow line is not allowed to intersect $\gamma_u^-$, so it necessarily passes above $x^-$. We consider the set $\mathcal{Z}$ of $z \in (-\epsilon, \epsilon)$ such that there exists a flow line of $X_u$ starting at $x^+$, passing through $(s^-, z)$, and escaping $\Sigma$ along $\partial_v^- \Sigma$. Here, $x^- = (s^-, z^-)$. Similarly as before, we define $z_{\inf} \coloneqq \inf \mathcal{Z}$, and we can obtain a flow line $\gamma_{\inf}$ starting at $x^+$ and passing through $(s^-, z_{\inf})$. We know that $z_{\inf} \geq z^-$. If $z_{\inf} > z^-$, arguing as before, we could obtain a flow line $\gamma''$ of $-X_u$ starting at $x^-$ and intersecting $\gamma_u^+$, which is impossible. Therefore, $z_{\inf} = z^-$ and $\gamma_{\inf}$ passes through $x^-$. In other terms, we have shown the existence of a pre-lasso. 

    We now explain how to construct a lasso, i.e., a pre-lasso such that all of its points converge to $0$ under the backward flow of $X$. Let $\mathcal{L}$ be the set of (unparametrized) pre-lassos $\ell \subset \Sigma$, i.e., the set of $C^1$ curves on $\Sigma$ tangent to $X_u$ starting at $x^+$ and ending at $x^-$. We can assume without loss of generality that $s^+ = -1/2$ and $s^- = 1/2$. For $\ell \in \mathcal{L}$, we consider its intersection point $(0, z_0)$ with the vertical line $\{s=0\} \cap \Sigma$. Arguing as before, by taking the infimum of such values of $z_0$, we obtain a pre-lasso $\ell_0$ passing through a point $(0, z_0)$ which converges to $0$ under the backward flow of $X$. Arguing inductively, we can construct a sequence $(\ell_n)_n$ of pre-lassos such that for every $n$ and every dyadic number $\mathfrak{d}$ of height at most $n$ in $[-1/2, 1/2]$, the intersection point of $\ell_n$ with $\{s=\mathfrak{d}\} \cap \Sigma$ converges to $0$ under the backward flow of $X$, and $\ell_m$ intersects  $\{s=\mathfrak{d}\} \cap \Sigma$ at this same point when $m \geq n$. Passing to a subsequence, these pre-lassos converge to a pre-lasso $\ell_\infty$ such that a dense subset of its points converge to $0$ under the backward flow of $X$. Since the points on $\ell$ not converging to $0$ form an open subset of $\ell$, we obtain that $\ell_\infty$ is a lasso.

    As explained before, we obtain a surface $S^+$ by flowing a lasso $\ell$ under the backward flow of $X$. The boundary of $S^+$ is then $\ell \cup \left( \{y=z=0\} \cap \mathsf{B}_\epsilon \right)$, and $S^+$ is everywhere tangent to $\eta_u$.
\end{itemize}

Finally, the union $S$ of $S^-$ and $S^+$ inside of $\mathsf{B}_\epsilon$ is a continuous surface which can be written as the graph of a continuous function $f$. Arguing as in Case 2, by looking at the intersection of this surface with $F_x^\pm$ near $x^\pm$, we conclude that $f$ is in fact $C^1$, so $S$ is an integral surface for $\eta_u$ passing through $0$.
\end{proof}

\begin{rem}
Handling Case $3$ is particularly tricky because we are not aware of a sufficiently nice $C^1$ normal form for $X$ near a quadratic point $p \in Q$.
\end{rem}

\begin{rem}
It is worth noting that in~\cite{CF11}, the authors assume that $(\xi_-, \xi_+)$ is in normal form to show the existence of $(\eta_u, \eta_s)$, as they need to show that the limits of $\xi_\pm$ along the flow of $X$ exist. With the present method, we show the existence without any assumption on $(\xi_-, \xi_+)$. Then, even though \cite{CF11} also uses the normal form to show local integrability along $\Delta$, what is really needed for the proof is a precise understanding of the qualitative behavior of $X$ near $\Delta$. This can only be achieved under some genericity assumption on $(\xi_-, \xi_+)$, the point being that a plane field invariant under a possibly singular but sufficiently nice vector field tangent to it is locally integrable, even at the singularities.
\end{rem}

The following lemma provides a more refined description of the integral surfaces of $\eta_u$ passing through a quadratic point $p \in Q$:

\begin{lem} \label{lem:quad}
Let $p \in Q$ be a quadratic contact point. There exists a unique flow line $\gamma_s = \gamma_s(p)$ of $X$ converging to $p$. Moreover, there exists a unique (germ of) surface with boundary $S^u=S^u(p)$ whose boundary contains $p$ and is contained in the unstable manifold of $X$ at $p$, and such that $X$ is tangent to $S^u$. The plane field $\eta_u$ is tangent to $S^u$ and is uniquely integrable along $S^u$ in the following sense. If $(S, q)$ is a germ of integral surface of $\eta_u$ at $q \in M$, then
\begin{itemize}
    \item If $q$ lies in the interior of $S^u$, then $S \subset S^u$.
    \item If $q$ lies in the boundary of $S^u$, then the unstable branches of $X$ at $p$ divide $S$ into two surfaces with boundary, one of which is entirely contained in $S^u$.
\end{itemize}
    
Therefore, every (germ of) integral surface $S$ of $\eta_u$ at $p$ is the union of two surfaces with boundary $S^+$ and $S^-$, where
    \begin{itemize}
    \item The boundary of $S^\pm$ contains $p$ and is contained in the unstable manifold of $X$ at $p$,
    \item $S^+$ is included in the saturation of $S^u$ by the flow of $X$,
    \item $S^-$ contains $\gamma_s$.
    \end{itemize}
Hence, as a germ at $p$, $S$ is entirely determined by its intersection with a small disk transverse to $\gamma_s$. 
\end{lem}

We call $\gamma_s(p)$ the \textbf{stable branch} of $X$ at $p$, and $S^u(p)$ the \textbf{unstable half-disk} of $X$ at $p$.

\begin{proof}
In the proof of Proposition~\ref{prop:locint}, we saw that there is a unique flow line of $X$ converging to $p$, and the flow lines converging to $p$ in backward times form a topological half-disk containing $p$ and the unstable branches of $X$ at $p$. We also constructed a $C^1$ half-disk near $p$ tangent to $\eta_u$ satisfying these properties, so these two half-disks coincide. In particular, this unstable half-disk is $C^1$. Moreover, the lasso $\ell$ constructed in the proof of Proposition~\ref{prop:locint} is unique, and $y_b \in F^-_y$ is the unique point in $F^-_y$ converging to $0$ under the flow of $X$. See Figure~\ref{fig:sigma2} for a detailed picture of $\Sigma$.

\begin{figure}[h]
    \centering
    \includegraphics[width=1 \linewidth]{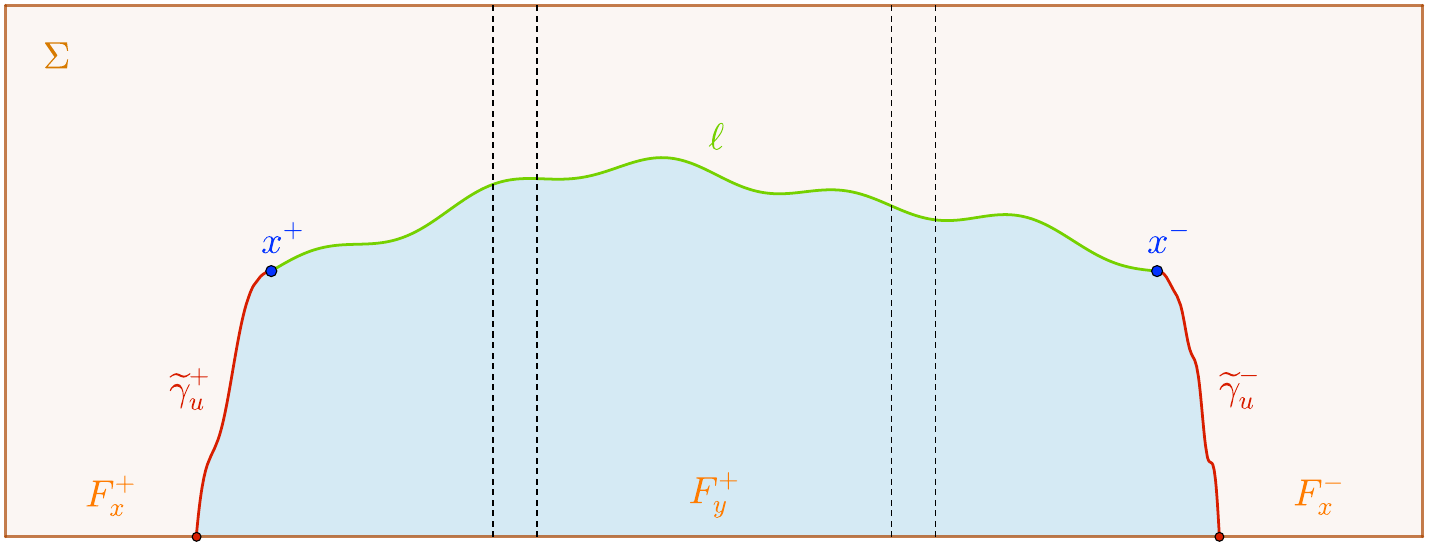}
    \caption{The surface $\Sigma$ and its unique lasso $\ell$. The curves $\widetilde{\gamma}^\pm_u$ are the intersection of $\Sigma$ with the unstable manifolds of the saddle points on $\Delta$. The blue region is the intersection of $\Sigma$ with the unstable manifolds of the source points on $\Delta$. All of the points in the orange region leave the box $\mathsf{B}_\epsilon$ when flowed backward along $X$.}
    \label{fig:sigma2}
\end{figure}

We now show that $\eta_u$ is uniquely integrable along $S^u$ in the sense of the statement of the lemma. Otherwise, there would exist a point $p_0 \in \ell \subset \Sigma$ and a small integral curve $\gamma_0$ of $X_u$ containing $p_0$ which is not contained in $\ell$. After deleting a portion of $\gamma_0$ and without loss of generality, we can assume that $\gamma_0 \cap \ell = \{p_0\}$, $p_0$ is the left boundary point of $\gamma_0$, and $p_0 \neq x^-$. Note that $\gamma_0$ necessarily lies above $\ell$, as points below $\ell$ lie in the unstable manifolds of source singularities in $\Delta$. Moreover, if $\gamma_0$ is sufficiently short, then all of the points in $\gamma_0 \setminus \{p_0\}$ are flown to the back face $F_y^-$ of the box $\mathsf{B}_\epsilon$ under the backward flow of $X$. The intersection points of $F_y^-$ with these flow lines form a continuous curve $\gamma_0^-$ which is tangent to $\eta_u \cap F_y^-$, and one of its ends converges to $y_b$. It follows that $\gamma_0^- \cup \{y_b\}$ can be described near $y_b$ as the graph of a $C^1$ curve. However, by Step 1 in Case 3 of the proof of Proposition~\ref{prop:locint}, the (forward) flow lines of $X$ starting at points in $\gamma_0^-$ near $\{y_b\}$ would intersect $\partial \mathsf{B}_\epsilon$ away from $\ell$, which is a contradiction.

The description of an integral surface $S$ near $p$ then follows from the proof of Proposition~\ref{prop:locint}.
\end{proof}

\begin{rem}
In the proof of Proposition~\ref{prop:locint}, we showed that there is a unique germ of integral surface for $\eta_u$ at a source point $p \in \Delta_\mathrm{so}$, and a germ of integral surface at a saddle point $p \in \Delta_\mathrm{sa}$ is entirely determined by its intersection with small disks transverse to the unstable (or stable) branches of $X$ at $p$.
\end{rem}

%%%%%%%%%%%%%%%%%%%%%%%%%%%%%%%%%%%%%%%%%%%%%%%%%%%%%%%%%%%%%%%%%%%%%
        \subsection{Proof of Theorem~\ref{thmintro:etau}}
%%%%%%%%%%%%%%%%%%%%%%%%%%%%%%%%%%%%%%%%%%%%%%%%%%%%%%%%%%%%%%%%%%%%%

First of all, the plane fields $\xi^t_\pm$ defined by~\eqref{eq:xit} converge pointwise to (not necessarily continuous) plane fields $\eta^{\infty}_\pm$ when $t$ goes to $+\infty$. This can be obtained from the standard argument explained in~\cite[Section 2.4]{CF11}. We choose an arbitrary Riemannian metric $g$ on $M$ and we denote by $\nu$ the orthogonal plane field to $X$ defined away from $\Delta$. We can then measure the angle $\theta^t_\pm(p)$ between $\xi_\pm(p) \cap \nu(p)$ and $\xi^t_\pm(p) \cap \nu(p)$, which satisfy $$-\pi < \theta^t_- \leq 0 \leq  \theta^t_+ < \pi.$$
The contact conditions for $\xi_\pm$ imply $$\frac{d}{dt} \theta^t_+(p) > 0, \qquad \frac{d}{dt} \theta^t_-(p) < 0,$$ for $p \in M \setminus \Delta$. Moreover, if $p \in \Delta$, then for every $t\in \R$, $\xi^t_+(p) = \xi^t_-(p) = \xi_+(p) = \xi_-(p)$. This shows that the limits $\eta^\infty_\pm$ exist. If $\gamma$ is a (possibly trivial) flow line of $X$, then by definition, the restrictions of $\eta^\infty_\pm$ to $\gamma$ are invariant under the flow of $X$. In particular, they are differentiable along $X$. Moreover, $\eta^\infty_\pm$ are ``sandwiched'' between $\xi_-$ and $\xi_+$ as in Proposition~\ref{prop:distrib}. We do not know that $\eta^\infty_\pm$ is continuous yet, but the proof of the latter readily implies that $\eta^\infty_\pm$ coincide with $\eta_u$ along $\gamma$,\footnote{This is because only one solution to $(E_p)$ is defined for all positive times.} hence these plane fields coincide on $M$. Now that we know that $\eta^\infty_\pm = \eta_u$ is continuous, we can argue as in~\cite[Lemme 3.1]{CF11} and apply Dini's theorem to show that the convergence of $\xi^t_\pm$ to $\eta_u$ is uniform. The first item of Theorem~\ref{thmintro:etau} is proved, and the second item immediately follows from Proposition~\ref{prop:locint}.  \qed

%%%%%%%%%%%%%%%%%%%%%%%%%%%%%%%%%%%%%%%%%%%%%%%%%%%%%%%%%%%%%%%%%%%%%%%%%%%%%%%%%%%%%%%%%%%%%%%%%%%%%%%%%%%%%%%%%%%%%%%%%%%%%%%%%%%%%%%%%%%%%%%%%%%%%%%%%%%%%%
\section{Constructing foliations from contact pairs} \label{sec:c0fol}
%%%%%%%%%%%%%%%%%%%%%%%%%%%%%%%%%%%%%%%%%%%%%%%%%%%%%%%%%%%%%%%%%%%%%%%%%%%%%%%%%%%%%%%%%%%%%%%%%%%%%%%%%%%%%%%%%%%%%%%%%%%%%%%%%%%%%%%%%%%%%%%%%%%%%%%%%%%%%%

In this section, we prove Theorem~\ref{thmintro:int} from the introduction, up to a (very) technical result (Theorem~\ref{thm:tech}) whose proof occupies Part~\ref{part:tech}.

%%%%%%%%%%%%%%%%%%%%%%%%%%%%%%%%%%%%%%%%%%%%%%%%%%%%%%%%%%%%%%%%%%%%%
        \subsection{Polarized vector fields}
%%%%%%%%%%%%%%%%%%%%%%%%%%%%%%%%%%%%%%%%%%%%%%%%%%%%%%%%%%%%%%%%%%%%%

\begin{defn}
Let $X$ be a smooth, possibly singular vector field on $M$. A \textbf{polarization} of $X$ is a continuous, cooriented plane field $\eta$ which contains $X$ and is invariant under the flow of $X$. We call the pair $(X, \eta)$ a \textbf{polarized vector field}.\footnote{This terminology is motivated by the following observation: away from the singular set of $X$, $\eta$ defines an oriented line subbundle of the plane bundle $TM \slash \langle X \rangle$ which is invariant under the flow of $X$. After choosing a (locally defined) transverse plane field to $X$, this line bundle can be identified with an oriented line field transverse to $X$.}
\end{defn}

In the previous sections, we saw that a positive contact pair $(\xi_-, \xi_+)$ gives rise to a polarized vector field $(X, \eta_u)$ which depends continuously on $(\xi_-, \xi_+)$.

By analogy with Definition~\ref{def:reg}, we say that a polarized vector field $(X, \eta)$ on $M$ is \textbf{regular} if the following conditions are satisfied:
\begin{enumerate}
    \item The singular set $\Delta$ of $X$ is a smoothly embedded link,
    \item The linearization of $X$ has rank $2$ and positive trace along $\Delta$, 
    \item $\Delta$ is transverse to $\eta$ away from a finite set of points $Q \subset \Delta$ where it has ``quadratic tangencies'' with $\eta$, in the following sense. At each $p \in Q$, there exist smooth coordinates $(x,y,z)$ in which the linearization of $X$ at $0$ is of the form
    $$X_1 = a x \partial x + z \partial_y,$$
    and $\eta(0) = \mathrm{span} \{\partial_x, \partial_y\}$. Moreover, $\Delta = \{ x=0, z = \pm y^2\}$ near $p$ in these coordinates.
    \item $(X, \eta)$ satisfies the conclusions of Lemma~\ref{lem:quad} at these quadratic points.
\end{enumerate}

The proof of Proposition~\ref{prop:locint} implies that $\eta$ is locally integrable. Similarly, we have a decomposition 
$$\Delta = \Delta_\mathrm{so} \sqcup Q \sqcup \Delta_\mathrm{sa}$$ in terms of source, quadratic, and saddle singularities.

The analysis of the different types of singularities of $X$ along $\Delta$ in the proofs of Proposition~\ref{prop:locint} and Lemma~\ref{lem:quad} readily implies the following lemma which will be useful in Part~\ref{part:tech}.

\begin{lem} \label{lem:nbd}
There exists a tubular neighborhood $\mathcal{N}$ of $\Delta$, made of a disjoint union of tubular neighborhoods of the connected components of $\Delta$, and satisfying the following property. If a nontrivial flow line $\gamma$ of $X$ is entirely contained in $\mathcal{N}$, then $\gamma$ connects a source singularity to a saddle singularity, both in the same component of $\Delta$.
\end{lem}

To simplify some of our proofs, especially in Part~\ref{part:tech}, we will make further generic assumptions on the contact pairs under consideration and their associated polarized vector fields. Recall that a \textbf{connection} of a vector field $X$ is a nontrivial flow line $\gamma$ converging to singularities of $X$ in positive and negative times.

\begin{defn} \label{def:admcon}
Let $(X, \eta)$ be a regular polarized vector field. A connection $\gamma$ of $X$ between singularities $\lim_{\pm \infty} \gamma =: p_\pm \in \Delta$ is \textbf{admissible} if it is of one of the four following types:
\begin{enumerate}
    \item[(A1)] $p_-$ is a source singularity and $p_+$ is a saddle singularity,
    \item[(A2)] $p_-$ is a source singularity and $p_+$ is a quadratic singularity,
    \item[(A3)] $p_-$ is a quadratic singularity, $p_+$ is a saddle singularity, and $\gamma$ is not an unstable branch of $p_-$,
    \item[(A4)] Both $p_\pm$ are saddle singularities, $p_- \neq p_+$, and $\gamma$ is the only saddle-saddle connection with endpoint $p_-$ or $p_+$.
\end{enumerate}
\end{defn}

\begin{figure}[h]
    \renewcommand{\thesubfigure}{A\arabic{subfigure}}
    \setcounter{subfigure}{0}
    \centering
    \begin{subfigure}{0.5\textwidth}
        \centering
        \includegraphics{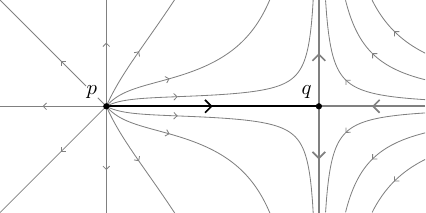}
        \caption{Source--saddle connection.}
        \label{fig:sosa}
    \end{subfigure}%
    \begin{subfigure}{0.5\textwidth}
        \centering
        \includegraphics{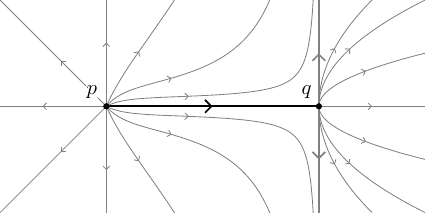}
        \caption{Source--quadratic connection.}
        \label{fig:soqu}
    \end{subfigure}
    \par\bigskip
    \begin{subfigure}{0.5\textwidth}
        \centering
        \includegraphics{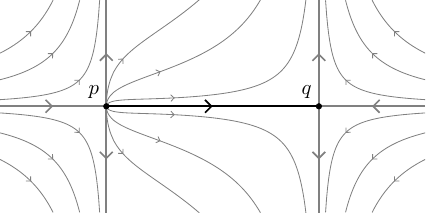}
        \caption{Quadratic--saddle connection.}
        \label{fig:qusa}
    \end{subfigure}%
    \begin{subfigure}{0.5\textwidth}
        \centering
        \includegraphics{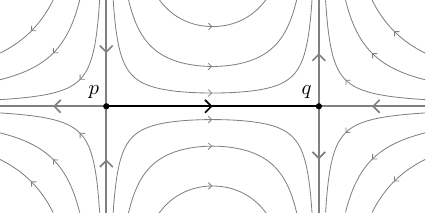}
        \caption{Saddle--saddle connection.}
        \label{fig:sasa}
    \end{subfigure}
    \caption{The four types of admissible connections.}
    \label{fig:admconn}
\end{figure}

The four types of admissible connections are depicted in Figure~\ref{fig:admconn}. We will also require a strengthening of case (A4) above, which will only become relevant in Section~\ref{sec:envorder} of Part~\ref{part:tech}.

\begin{defn} \label{def:brokensad}
    A \textbf{broken triple saddle connection} is a finite sequence $(\gamma_0, \dots, \gamma_{n+1})$, $n \geq 0$, of pairwise distinct unoriented connections of $X$ such that 
    \begin{itemize}
        \item $\gamma_0$ and $\gamma_{n+1}$ are saddle-saddle connections,
        \item For every $0 \leq i \leq n$, $\gamma_i$ and $\gamma_{i+1}$ share an endpoint, 
        \item For every $1 \leq i \leq n$, $\gamma_i$ belongs to the unstable manifold of a source or quadratic singularity of $X$.
    \end{itemize}
\end{defn}

\begin{defn} \label{def:adm}
A regular polarized vector field $(X, \eta)$ is \textbf{admissible} if $X$ only has admissible connections, and no broken triple saddle connections. A regular positive contact pair $(\xi_-, \xi_+)$ is \textbf{admissible} if its associated polarized vector field $(X, \eta_u)$ is admissible.
\end{defn}

The next lemma follows from standard general position arguments applied to the stable and unstable manifolds of the singularities of $X$:

\begin{lem} \label{lem:pertu}
Any regular positive contact pair can be (smoothly) approximated by an admissible one with the same singular set $\Delta$. Therefore, admissible positive contact pairs form a comeagre subspace of the space of positive contact pairs.
\end{lem}

\medskip

We will need further assumptions on the contact pairs under consideration to prevent the existence of \emph{disks of tangency} for their associated polarized vector fields.

\begin{defn} \label{def:disk}
A \textbf{disk of tangency} for a polarized vector field $(X, \eta)$ is an immersion $f: D \rightarrow M$ of the open $2$-disk in $M$ such that
    \begin{itemize}
    \item $f$ is tangent to $\eta$,
    \item $f$ extends to a continuous map $\bar{f} : \overline{D} \rightarrow M$ whose restriction to $\partial \overline{D} \cong S^1$ maps to a closed orbit of $X$.
    \end{itemize}
\end{defn}

The presence of such disks could obstruct some of the main constructions in Part~\ref{part:tech} below. However, if one of the contact structures in the pair is tight, or if $\Delta_\mathrm{so} = \varnothing$ (see the hypothesis of Theorem~\ref{thmintro:int}), then these disks cannot exist:

\begin{prop} \label{prop:disk}
Let $(\xi_-, \xi_+)$ be an admissible positive contact pair. If its associated polarized vector field $(X, \eta_u)$ has a disk of tangency, then $\Delta_{\mathrm{so}} \neq \varnothing$ and both $\xi_-$ and $\xi_+$ are overtwisted.
\end{prop}

\begin{proof} Let $f : D \rightarrow M$ be a disk of tangency that extends to $\bar{f} : \overline{D} \rightarrow M$, and let $\gamma$ be the closed orbit of $X$ such that $\bar{f}\big(\partial \overline{D}\big) \subset \gamma$.

Notice that $\bar{f}$ is an immersion as well: intersect $\gamma$ with a small transversal surface $\Sigma$ and consider the intersection of $f(D)$ with $\Sigma$. It is a collection of flow lines of the continuous vector field $X_u$ obtained by intersecting $\eta_u$ with $\Sigma$. These flow lines converge to the point $p_0 = \gamma \cap \Sigma$, so they extend in a $C^1$ way at that point. This implies that $\bar{f}$ is a $C^1$ immersion tangent to $\eta_u$.

Let $\widetilde{X}$ denote the pullback of $X$ to $\overline{D}$. It is a continuous, uniquely integrable vector field tangent to $\partial \overline{D}$ with isolated singularities in $D$. These singularities are of three types: source (index $+1$), birth-death (index $0$), or saddle (index $-1$). By the Poincar\'{e}--Hopf theorem, there must be a source singularity in $D$, i.e., $f(D) \cap \Delta_\mathrm{so} \neq \varnothing$.

Let us temporarily assume that $\bar{f}$ is an embedding and write $\mathcal{D} \coloneqq \bar{f}\big(\overline{D}\big)$. Since $L \coloneqq \partial \mathcal{D}$ is Legendrian for both $\xi_-$ and $\xi_+$, and since $X$ is tangent to $L$ and to $\eta_u$, the Thurston--Bennequin numbers $\mathrm{tb}_\pm(L)$ of $L$ for $\xi_\pm$ are both zero. Since a Legendrian knot $K$ bounding an embedded disk in a tight contact $3$-manifold has $\mathrm{tb}(K) \leq - \chi(\overline{D}) = -1$ by~\cite{E92}, this implies that $\xi_-$ and $\xi_+$ are both overtwisted.

We now explain how to reduce to the case where $\bar{f}$ is an embedding. We will repeatedly use the Poincar\'{e}--Bendixson theorem (see~\cite[Theorem 7.16]{T12}) and the following fact: if a domain in $D$ is bounded by a closed orbit of $\widetilde{X}$, then it contains a singularity of $\widetilde{X}$, by the Poincar\'{e}--Hopf theorem (or even the hairy ball theorem by doubling the disk).

Let us consider the oriented graph $G = (V,E)$ on $D$ whose vertices are the singularities of $\widetilde{X}$ and whose oriented edges are the flow lines of $\widetilde{X}$ connecting two singularities; see Figure~\ref{fig:disk}. 

\begin{figure}[h]
    \centering
    \includegraphics[width=0.5 \linewidth]{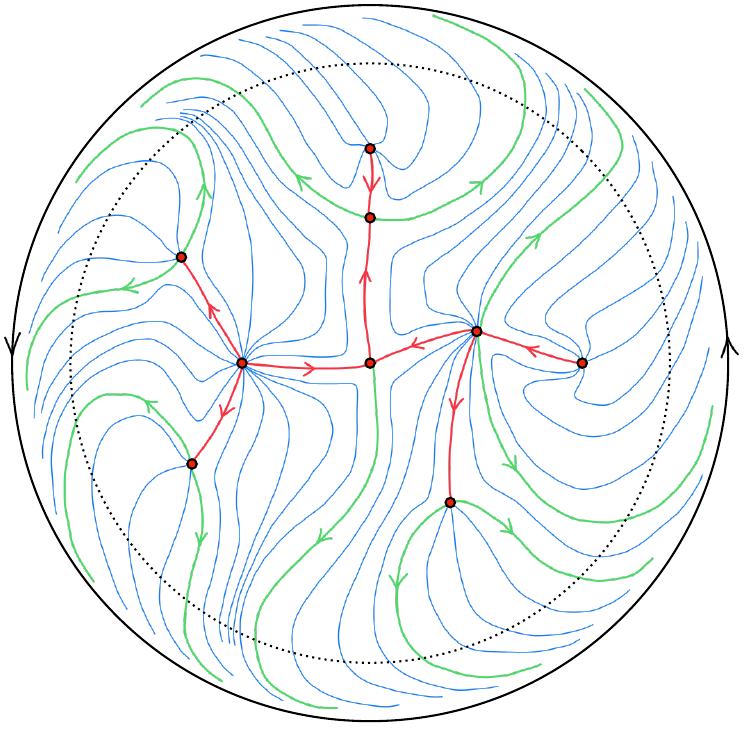}
    \caption{A tree on a leaf. In red: the tree $G$. In green: the unstable branches limiting to the boundary orbit. In blue: flow lines of $\widetilde{X}$. The outer circle is a closed orbit of $\widetilde{X}$, the dotted circle is transverse to $\widetilde{X}$. Notice the nonadmissible connection between two quadratic singularities...}
    \label{fig:disk}
\end{figure}

Notice that these flow lines are stable branches of the saddle and birth-death singularities, since there are no sinks, and there are finitely many of them. This implies that $G$ is a finite graph. We now argue that \emph{after possibly restricting to a smaller disk of tangency}, $G$ is a tree. We first show by contradiction that $G$ has no cycle, after possible restriction. If $G$ has a cycle, we choose an innermost one and we denote it by $c$. We also consider $\overline{D}_c \subset D$ the closed domain bounded by $c$ in $D$. Let $\delta$ be an edge in $c$ connecting the singularity $p$ to the singularity $q$. We distinguish three cases:
\begin{enumerate}[leftmargin=*]
    \item \textit{$p$ and $q$ are saddles.} By the admissibility condition, we have $p \neq q$ and the two unstable branches of $q$ cannot converge to singularities in positive time. One of these two branches is contained inside of $\overline{D}_c$. We denote it by $u$ and we distinguish two subcases:
    \begin{itemize}
        \item If $u$ limits to a closed orbit, the latter bounds a domain which is a smaller disk of tangency. We restart the argument with that disk, which strictly reduces the graph $G$.
        \item Otherwise, the Poincar\'{e}--Bendixson theorem implies that $u$ limits to a connected, \emph{strict} subgraph $G'$ of $G$ included in $\overline{D}_c$ with at least one edge. Moreover, by elementary properties of $\omega$-limi sets, the graph $G'$ has the following property: if $\delta'$ is an edge in $G'$ limiting to a singularity $z$, which is a saddle or a birth-death, then at least one of the two unstable branches of $z$ is also an edge in $G'$. This readily implies that $G'$ contains a cycle, since there are finitely many singularities. Recall that $c$ is an innermost cycle, so $G'$ necessarily contains $c$, and in particular it contains $\delta$; however, it is impossible for an unstable branch of $q$ to accumulate on $\delta$ as it would also accumulate on itself, a contradiction.
    \end{itemize}
    \item \textit{$p$ is a source or birth-death and $q$ is a saddle.} None of the two unstable branches of $q$ can be an edge of the cycle $c$, since the admissibility condition would imply that it converges to another saddle singularity, and this situation is prohibited by the first case. Therefore, the stable branch of $q$ different than $\delta$ is an edge in $c$. If the unstable branch of $q$ contained in $\overline{D}_c$ does not limit to a singularity, the argument of the first case applies and leads to a contradiction. This unstable branch could limit to a saddle $z$ in $\overline{D}_c$, but the unstable branches of $z$ cannot limit to singularities, and we similarly obtain a contradiction.
    \item \textit{$p$ is a source and $q$ is a birth-death.} The unstable branches of $q$ cannot converge to singularities by admissibility. Therefore, there is an unstable flow line in the unstable half-disk of $q$, different than the unstable branches, which is an edge in $c$. Once again, the behavior of the unstable branch of $q$ contained in $\overline{D}_c$ leads to a contradiction. 
\end{enumerate}

To summarize, we have shown that after possibly restricting to a smaller disk of tangency, the graph $G$ is a forest, namely, a disjoint union of trees. Notice that since there is no cycle of singularities, every flow line of $\widetilde{X}$ either converges to a singularity or to a closed orbit. This allows us to restrict to the case where $D$ does not contain any nontrivial closed orbit of $\widetilde{X}$. We can argue by induction on the number of connected components of $G$. Indeed, if $G_0$ denotes a connected component of $G$, then every leaf of $G_0$ admits an unstable flow line $u$ (unstable branch or a flow line contained in the unstable manifold) which is not contained in $G$. This flow line must limit to a closed orbit $\gamma$ of $\widetilde{X}$. If the disk $D_\gamma$ bounded by $\gamma$ is disjoint from $G_0$, then we restrict $f$ to $D_\gamma$, reducing the number of connected components of the graph of singularities of $\widetilde{X}$. Otherwise, $G_0$ is contained in $D_\gamma$, and every other nontrivial closed orbit of $\widetilde{X}$ contained in $D_\gamma$ bounds a disk which is \emph{disjoint} from $G_0$; otherwise this disk would contain $u$, hence contain $D_\gamma$. If such an orbit exists, we restrict $f$ to the disk it bounds, reducing the number of connected components. Otherwise, we restrict $f$ to $D_\gamma$, and we are done.

We now show that $G$ is connected. Notice that the argument in the previous paragraph readily implies that $\widetilde{X}$ has contracting holonomy along $\partial D$. The Poincar\'{e}--Bendixson theorem also implies that all of the stable branches of singularities are edges in $G$. Let $G_1, \dots, G_n$, $n \geq 1$, denote the connected components of $G$. For every $1 \leq i \leq n$, there is an open neighborhood $\mathcal{N}_i$ of $G_i$ in $D$ with smooth boundary such that $\widetilde{X}$ is outward transverse to $\partial \mathcal{N}_i$. We can further assume that these neighborhood are pairwise disjoint. To construct such neighborhoods, we consider the union of small bands around each edge in $G$ and small disks around the vertices. The conditions on $G$ imply that we can modify this neighborhood so that $\widetilde{X}$ becomes outward transverse to its boundary. We also consider a slightly smaller disk $D' \subset D$ such that $\widetilde{X}$ is outward transverse to $\partial D'$. If $n \geq 2$, we would obtain a nonsingular vector field on $D' \setminus \bigcup_{1 \leq i \leq n} \mathcal{N}_i$ which is impossible for degree reasons. Therefore, $G$ is connected and $D'$ retracts onto $G$ under the backward flow of $\widetilde{X}$.

We now show that $f$ is injective. It is enough to show that the restriction of $f$ to the vertices of $G$ is injective. Indeed, if two distinct points $p, q \in D$ are such that $f(p) = f(q)$, then our previous arguments imply that they lie on two distinct flow lines of $\widetilde{X}$ which converge to two singularities in negative time. These singularities are necessarily distinct since $f$ is an immersion. Therefore, we obtain two distinct vertices of $G$ with the same image under $f$. We define an equivalence relation $\sim$ on the vertices of $G$ by $$ p \sim q \iff f(p) = f(q).$$ Since $f$ is an immersion, this relation extends to the edges of $G$. We denote the quotient graph by $H \coloneqq G\slash\!\!\sim$ and the associated quotient map by $\tilde{f} : G \rightarrow H$. The graph $H$ can also be described as the graph in $M$ whose vertices are the singularities $f(p)$, $p \in V(G)$, and the edges are the flow lines of $X$ connecting them. It is then easy to see that $\tilde{f}$ is a covering map in the graph theoretical sense. Since $H$ is covered by the \emph{finite} tree $G$, $H$ itself a tree and $\tilde{f}$ is an isomorphism. This implies that $f_{\vert G}$ is injective, thus $f$ is injective.

We are left to show that the restriction of $\bar{f}$ to $\partial \overline{D}$ is injective. If not, then this map winds around the closed orbit $\gamma$ at least twice. With the above notations, the intersection of the transversal $\Sigma$ with $\bar{f}(\overline{D})$ is a finite collection of flow lines of $X_u$ which only intersect at $p_0$. If $\gamma^0_u$ is such a flow line, we consider its image $\gamma^1_u$ under the return map of $X$ on $\Sigma$. By our assumptions, $\gamma^0_u$ and $\gamma^1_u$ are distinct and only intersect at $p_0$. After choosing coordinates $(x,y)$ on $\Sigma$ in which $X_u$ is almost flat, we can represent $\gamma^0_u$ and $\gamma^1_u$ as graphs of $C^1$ functions defined for $x \leq 0$ and passing through $(0,0) \cong p_0$. Without loss of generality, we can assume that $\gamma^0_u \leq \gamma^1_u$. We now consider the iterate images of $\gamma^0_u$ under the return map and obtain a sequence of flow lines $\gamma^i_u$, $i \geq 0$, satisfying $\gamma^i_u \leq \gamma^{i+1}_u$ for every $i \geq 0$. Here, we are using that $\eta_u$ is cooriented and invariant under the flow of $X$. The $\gamma^i_u$'s have to be pairwise disjoint, which is impossible since $\bar{f}_{\vert \partial \overline{D}}$ winds around $\gamma$ finitely many times.
\end{proof}

\begin{rem}
Disks of tangency should be compared to the \emph{(twisted) disks of contact} from~\cite{GO89}, and the \emph{sink disks} from~\cite{L02}, which play a crucial role in the theory of branched surfaces and laminations. The presence of such disks obstructs the construction of a lamination from a branched surface, and it is often a delicate combinatorial problem to modify a branched surface and remove those disks. While it can be particularly hard to certify that a given branched surface is sink disk free, or if a lamination extends to a foliation in its complement, Proposition~\ref{prop:disk} provides a very natural and tractable way to preclude the existence of disks of tangency.
\end{rem}

%%%%%%%%%%%%%%%%%%%%%%%%%%%%%%%%%%%%%%%%%%%%%%%%%%%%%%%%%%%%%%%%%%%%%
        \subsection{Proof of Theorem~\ref{thmintro:int}}
%%%%%%%%%%%%%%%%%%%%%%%%%%%%%%%%%%%%%%%%%%%%%%%%%%%%%%%%%%%%%%%%%%%%%

Let $(\xi_-, \xi_+)$ be a positive contact pair satisfying the hypotheses of Theorem~\ref{thmintro:int}. After a generic $C^\infty$-small perturbation away from $\Delta$, which induces a small $C^0$ perturbation of $\eta_u$, we can assume that this contact pair is admissible by Lemma~\ref{lem:pertu}. By Proposition~\ref{prop:disk}, the polarized flow $(X, \eta_u)$ has no disks of tangency. We now appeal to the 

\begin{techthm} \label{thm:tech}
Let $(X, \eta)$ be an admissible polarized vector field on $M$ without disks of tangency. For every $\epsilon > 0$, there exists a $C^0$-foliation $\mathcal{F}_\epsilon$ such that $d_{C^0}(\eta, T \mathcal{F}_\epsilon) < \epsilon$. Moreover, there exists a continuous map $h_\epsilon : M \rightarrow M$, $\epsilon$-close to the identity, such that for every leaf $L$ of $\mathcal{F}_\epsilon$, the restriction ${h_{\epsilon}}_{\vert L} : L \rightarrow M$ is a $C^1$ immersion tangent to $\eta$.
\end{techthm}

The proof of this highly nontrivial result is deferred to Part~\ref{part:tech} and will occupy the whole second half of this paper. 

Let $\mathcal{F}$ be a $C^0$-foliation approximating $\eta_u$ together with a continuous map $h : M \rightarrow M$ as above. By~\cite[Theorem 4.1]{KR19} (see also~\cite[Remark 3.3]{B16a}), we can further assume that the leaves of $\mathcal{F}$ are smoothly immersed. We show by contradiction that $\mathcal{F}$ has no spherical leaves. Otherwise, we readily obtain an embedded sphere $S$ in $M$ whose tangent bundle is homotopic to ${\eta_u}_{\vert S}$, and since $\xi_\pm$ are homotopic to $\eta_u$, $$\vert \langle e(\xi_\pm), [S] \rangle \vert  = \vert \langle e(\eta_u), [S] \rangle \vert = \chi(S) = 2.$$ However, by Bennequin--Eliashberg's inequality (see~\cite{E92}), $$\langle e(\xi), [S] \rangle = 0$$ for a tight contact structure $\xi$, hence both $\xi_\pm$ are overtwisted. By applying the map $h$, we also obtain a $C^1$-immersed sphere $f : S^2 \rightarrow M$ tangent to $\eta_u$, and pulling back $X$ along $f$, we get a (continuous, uniquely integrable) vector field $\widetilde{X}$ on $S^2$ with isolated singularities which are sources (index $+1$), saddles (index $-1$), or birth-death singularities (index $0$). The Poincar\'{e}--Hopf theorem readily implies that $\widetilde{X}$ has at least one source singularity, so $\Delta_{\mathrm{so}} \neq \varnothing$. This contradicts either hypothesis of Theorem~\ref{thmintro:int}, and finishes the proof. \qed

\begin{rem}
    Technically speaking, the hypothesis that $(X, \eta)$ has no disk of tangency is not absolutely necessary for Theorem~\ref{thm:tech} to hold. It could be weakened, e.g., by requiring that every disk of tangency is contained in an open subset of $M$ on which $\eta$ is uniquely integrable. However, we do not know if the latter condition is achieved for the polarized vector field of a generic positive contact pair.
\end{rem}

%%%%%%%%%%%%%%%%%%%%%%%%%%%%%%%%%%%%%%%%%%%%%%%%%%%%%%%%%%%%%%%%%%%%%%%%%%%%%%%%%%%%%%%%%%%%%%%%%%%%%%%%%%%%%%%%%%%%%%%%%%%%%%%%%%%%%%%%%%%%%%%%%%%%%%%%%%%%%%
    \section{Taut contact pairs}
%%%%%%%%%%%%%%%%%%%%%%%%%%%%%%%%%%%%%%%%%%%%%%%%%%%%%%%%%%%%%%%%%%%%%%%%%%%%%%%%%%%%%%%%%%%%%%%%%%%%%%%%%%%%%%%%%%%%%%%%%%%%%%%%%%%%%%%%%%%%%%%%%%%%%%%%%%%%%%

%%%%%%%%%%%%%%%%%%%%%%%%%%%%%%%%%%%%%%%%%%%%%%%%%%%%%%%%%%%%%%%%%%%%%
        \subsection{Tautness and tightness} \label{sec:tautight}
%%%%%%%%%%%%%%%%%%%%%%%%%%%%%%%%%%%%%%%%%%%%%%%%%%%%%%%%%%%%%%%%%%%%%

The following definition already appears in~\cite{CF11} and~\cite{CH20}.

\begin{defn}
A contact pair $(\xi_-, \xi_+)$ on $M$ is \textbf{taut} if for all $p, q \in M$, there exists a smooth arc from $p$ to $q$ positively transverse to both $\xi_\pm$.
\end{defn}

A taut contact pair is necessarily positive. It will be convenient to reformulate this definition in various ways that are similar to some of the many characterizations of taut foliations. In Appendix~\ref{appsec:taut}, we generalize some of these definitions to (continuous) plane fields, which might also be of independent interest.

\begin{prop} \label{prop:strongtight}
Let $(\xi_-, \xi_+)$ be a positive contact pair on $M$. The following are equivalent:
\begin{enumerate}
    \item $(\xi_-, \xi_+)$ is taut,
    \item For every $p \in M$, there exists a smooth loop positively transverse to $\xi_\pm$ passing through $p$,
    \item There exists a smooth closed $2$-form $\omega$ such that $\omega_{\vert \xi_\pm} > 0,$
    \item There exists a smooth volume preserving vector field positively transverse to $\xi_\pm$,
    \item $\eta_u$ is a taut plane field (see Definition~\ref{def:tautdist}).
\end{enumerate}
In particular, $\xi_\pm$ are weakly semi-fillable, hence universally tight.
\end{prop}

\begin{proof}
The implications 
$$\mathit{1} \implies \mathit{2} \implies \mathit{3} \implies \mathit{4} \implies \mathit{5}$$
are standard. The last implication follows from Proposition~\ref{prop:tautdistrib}. If $\eta_u$ is taut, since the distributions $\xi^t_\pm$ defined by~\eqref{eq:xit} converge uniformly to $\eta_u$ as $t$ goes to $+\infty$, it is not hard to show that $\big(\xi^T_-, \xi^T_+\big)$ is taut for $T$ sufficiently large; see the proof of~\cite[Proposition 3.7]{CF11}. Hence, $(\xi_-, \xi_+)$ is taut as well. Weak semi-fillability follows from the argument in~\cite[Corollary 3.2.8]{ET}.
\end{proof}

The third or fourth characterization immediately implies that tautness is an open condition in the $C^0$ topology. Moreover, if $\mathcal{F}$ is a taut $C^0$-foliation on $M$ and if $\xi_\pm$ are positive and negative contact structures on $M$ which are sufficiently $C^0$-close to $T \mathcal{F}$, then $(\xi_-, \xi_+)$ is taut. Conversely, we have:

\begin{thm}[Theorem~\ref{thmintro:taut}]
Let $(\xi_-, \xi_+)$ be a taut contact pair. Then $\eta_u$ can be $C^0$-approximated by taut $C^0$-foliations.
\end{thm}

\begin{proof}
The tautness of $(\xi_-, \xi_+)$ implies that $\eta_u$ is a taut plane field by Proposition~\ref{prop:strongtight}. Moreover, $\xi_\pm$ are both tight so Theorem~\ref{thmintro:int} applies. Since (everywhere) tautness is an open condition for continuous plane fields, $C^0$-foliations sufficiently close to $\eta_u$ are also taut.
\end{proof}

\begin{rem}
In~\cite[Proposition 3.6]{CF11}, it is shown that taut contact pairs can be deformed into positive pairs without quadratic singularities, i.e., for which $\Delta = \Delta_+$ is transverse to $\xi_\pm$ \emph{everywhere}. It would be interesting to better understand the structure of $\eta_s$, the stable plane field of the pair, near $\Delta$ where it `vanishes'.
\end{rem}

The tautness assumption on the pair $(\xi_-, \xi_+)$ is quite restrictive and hard to work with in practice. A natural weakening would be to only assume that $(\xi_-, \xi_+)$ is a positive contact pair and that both $\xi_\pm$ are tight (in particular, Theorem~\ref{thmintro:int} applies). In that case, Colin--Firmo~\cite[Th\'{e}or\`{e}me 3.3]{CF11} showed that \emph{if $\eta_u$ is uniquely integrable}, then it is tangent to a foliation without Reeb components with nullhomologous core, and moreover $M$ carries a Reebless foliation. The unique integrability assumption is crucial for their argument and is unfortunately impractical. As a plausible generalization of this result, we propose:

\begin{conj} \label{conj:tight} If $(\xi_-, \xi_+)$ is a \emph{tight} positive contact pair, then $\eta_u$ can be $C^0$-approximated by a $C^0$-foliation without rationally nullhomologous Reeb components, and $M$ carries a Reebless $C^0$-foliation.
\end{conj}

This conjecture would have the following important consequences:
\begin{itemize}
    \item A closed oriented $3$-manifold different from $S^1 \times S^2$ carries a Reebless foliation if and only if it carries a positive pair of tight contact structures.
    \item An atoroidal closed oriented $3$-manifold different from $S^1 \times S^2$ carries a taut foliation if and only if it carries a positive pair of tight contact structures.
    \item On an atoroidal rational homology sphere, the unstable plane field of a positive pair of tight contact structures can be $C^0$-approximated by taut $C^0$-foliations.
    \item An atoroidal rational homology sphere carries a pair of transverse taut foliations if and only if it carries a tight bicontact structure.
\end{itemize}

As noted in~\cite{CF11}, under the hypotheses of Conjecture~\ref{conj:tight}, Bennequin's inequality for transverse curves implies that \emph{no closed curve transverse to $\eta_u$ can bound an embedded disk in $M$.} Unfortunately, this property is \emph{not} open in the $C^0$ topology! If it were, it would be possible to approximate $\eta_u$ by a foliation $\mathcal{F}$ with trivial meridional holonomy along its Reeb components. Following the strategy of the proof of~\cite[Th\'{e}or\`{e}me 3.3]{CF11}, it would then be possible to remove the Reeb components of $\mathcal{F}$ and get a Reebless foliation $\mathcal{F}'$. Of course, the homotopy class of $T \mathcal{F}'$ might be different than the one of $\eta_u$. It seems more credible that some foliations approximating $\eta_u$ might have Reeb components with nontrivial meridional holonomy, and some might not (compare this to the \emph{phantom Reeb components} from~\cite{CKR19}). In Appendix~\ref{sec:ODEper}, we prove a two dimensional result showing that closed transverse curves are the only obstructions to finding foliations by circles, see Proposition~\ref{prop:ODEper}. This serves as evidence for the above conjecture and might be of independent interest.

\begin{figure}[h]
    \centering
    \includegraphics[width=0.5\linewidth]{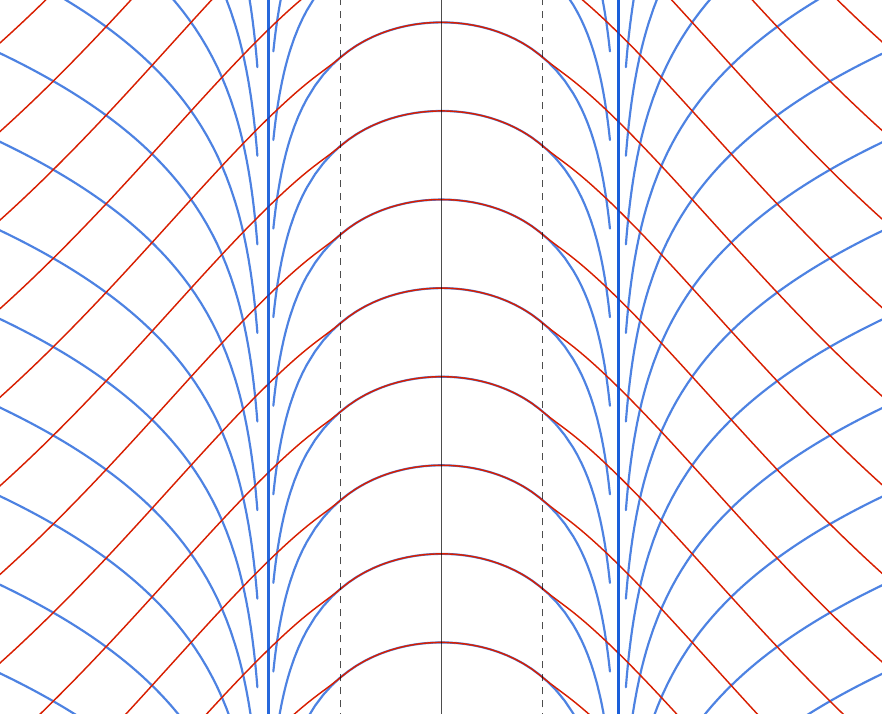}
    \caption{A $\theta$-slice of the foliations $\mathcal{F}$ (in blue) and $\mathcal{G}$ (in red).}
    \label{fig:reeb}
\end{figure}

Interestingly, Etnyre~\cite{E07} shows that a Reeb component with \emph{nontrivial} meridional holonomy can be approximated by a contact pair $(\xi_-, \xi_+)$ where one of $\xi_\pm$ is universally tight while the other one is overtwisted. Moreover, a Reeb component with \emph{trivial} meridional holonomy can be approximated by a pair of \emph{universally tight} contact structures:

\begin{ex}
On $\R^2 \times S^1$ with coordinates $(r, \theta, z)$, where $(r, \theta)$ denote standard polar coordinates on $\R^2$, we define
\begin{align}
    \alpha &\coloneqq (1-r^2) dz + r dr,\\
    \beta &\coloneqq \frac{1}{2} r^2 d\theta.
\end{align}
These are both smooth $1$-forms. The kernel of $\alpha$ is tangent to a smooth foliation $\mathcal{F}$ on $\R^2 \times S^1$ for which $D^2 \times S^1$ is a \emph{Reeb component}. Moreover, 
\begin{align*}
\langle \alpha, \beta \rangle &\coloneqq \alpha \wedge d \beta + \beta \wedge d\alpha \\
&= r dr\wedge d\theta \wedge dz \\
&> 0,
\end{align*}
and
\begin{align*}
\alpha \wedge d\alpha = \beta \wedge d\beta = 0,
\end{align*}
so $\mathcal{F}$ can be \emph{linearly deformed} into two contact structures with opposite signs $$\xi^\tau_\pm \coloneqq \ker \left( \pm \alpha + \tau \beta \right)$$
for any $\tau > 0$, and $(\xi^\tau_-, \xi^\tau_+)$ is a positive contact pair. Both $\xi^\tau_\pm$ are \emph{universally tight} since their lifts to $\R^2 \times \R$ embed into the standard positive and negative tight contact structures on $\R^3$, respectively. This can be seen by considering a $(z, \theta)$-invariant foliation by planes $\mathcal{G}$ on $\R^3$ which coincides with (the lift of) $\mathcal{F}$ near $\{r = 0\}$ and is transverse to $\mathcal{F}$ elsewhere (see Figure~\ref{fig:reeb}). We might further arrange that this foliation is tangent to $\mathrm{span}\{\partial_r, \partial_\theta\}$ outside a sufficiently large compact set. The characteristic foliation of $\xi^\tau_\pm$ along a leaf of $\mathcal{G}$ has a unique singularity at $0$, is radial near $\{r=0\}$, and is diffeomorphic to the standard radial foliation on $\R^2$ (it cannot have a closed characteristic since $\xi^\tau_\pm$ is transverse to $\partial_\theta$). Therefore, we may apply a $(z, \theta)$-equivariant diffeomorphism of $\R^3$ which maps $\mathcal{G}$ to the standard foliation of $\R^3$ by horizontal planes, to obtain a $(z, \theta)$-invariant positive contact structure which is tangent to $\partial_r$, and is only tangent to a horizontal plane at the origin. It easily follows that $\widetilde{\xi}^\tau_+$ embeds into the standard tight contact structure $\ker \big( dz + r^2 d\theta\big)$ on $\R^3$. A similar argument applies to $\xi^\tau_-$.

This shows that on an open $3$-manifold, some foliations with Reeb components \emph{can} be linearly deformed into a positive pair of universally tight contact structures.
\end{ex}

%%%%%%%%%%%%%%%%%%%%%%%%%%%%%%%%%%%%%%%%%%%%%%%%%%%%%%%%%%%%%%%%%%%%%
        \subsection{Transverse surgeries} \label{sec:transurg}
%%%%%%%%%%%%%%%%%%%%%%%%%%%%%%%%%%%%%%%%%%%%%%%%%%%%%%%%%%%%%%%%%%%%%

We now show that taut contact pairs are preserved under transverse surgeries of sufficiently large slopes. Theorem~\ref{thmintro:transfol} will follow easily.

\begin{prop} \label{prop:tautsurgery}
Let $(\xi_-, \xi_+)$ be a taut contact pair on $M$ and $K \subset M$ be a knot positively transverse to $\xi_\pm$ and equipped with some framing. There exists $s_0 > 0$ such that for any rational number $s$ satisfying $\vert s \vert  \geq s_0$, $M_K(s)$ carries a taut contact pair $(\xi^s_-, \xi^s_+)$. Moreover, the surgery dual knot of $K$ in $M_K(s)$ is positively transverse to $\xi^s_\pm$.
\end{prop}

\begin{proof}
After possibly perturbing $K$ to make it disjoint from $\Delta$, there exists a small tubular neighborhood $\mathcal{N}_K \cong S^1 \times D^2$ of $K \cong S^1 \times \{0\}$ such that
\begin{itemize}
    \item The vector field $\partial_t$ is positively transverse to $\xi_\pm$ in $\mathcal{N}_K$, and the framing of $K$ in the standard coordinates on $\mathcal{N}_K$ becomes trivial,
    \item For every $t \in S^1$, $X$ is tangent to $\{t\} \times D^2$ and $$\xi_\pm \cap T\big(\{t\} \times D^2\big) = \langle X \rangle.$$
\end{itemize}
We can apply~\cite[Proposition 3.5]{CF11} to find a new positive pair $(\xi'_-, \xi'_+)$ with associated vector field $X'$ and singular set $\Delta'$ such that
\begin{itemize}
    \item $\xi'_\pm$ is isotopic to $\xi_\pm$ rel $M \setminus \mathcal{N}_K$,
    \item The vector field $\partial_t$ is positively transverse to $\xi'_\pm$ in $\mathcal{N}_K$,
    \item For every $t \in S^1$, the characteristic foliation of $\xi'_\pm$ on $\{t\} \times D^2$ is spanned by $X'$ and has two singularities, a source and a saddle, in canceling position.
\end{itemize}
It follows that $\Delta'$ coincides with $\Delta$ outside of $\mathcal{N}_K$,  and $\Delta' \cap \mathcal{N}_K$ has two parallel components, made of source and saddle singularities, respectively. We can further arrange that the new source component coincides with $K$ and that there exists a small $\epsilon > 0$ such that in $S^1 \times D^2_\epsilon \subset S^1 \times D^2 \cong \mathcal{N}_K$, $\xi'_\pm$ is the kernel of the $1$-form $$\alpha^0_\pm \coloneqq \pm dt + r^2 d\theta,$$ where $(r,\theta)$ denotes the polar coordinates on $D^2$. Notice that $(\xi'_-, \xi'_+)$ is still taut.

We now modify $\xi^0_\pm \coloneqq \ker \alpha^0_\pm$ in $S^1 \times D^2_\epsilon$ relative to $T_\epsilon \coloneqq S^1 \times \partial D_\epsilon$ and perform an \emph{admissible contact surgery} on the new contact structures simultaneously. To this end, we consider a sufficiently small radius $0 < r_0 < \epsilon$ such that $r_0^2 \in \mathbb{Q}$, and we write $r_0^2 = p/q$ for some coprime positive integers $p, q \in \Z_{>0}$. Along $T_{r_0} \coloneqq S^1 \times \partial D_{r_0}$, the characteristic foliation of $\xi^0_+$ has slope $s \coloneqq -1/r_0^2=-q/p$. We modify $\xi^0_-$ inside of $S^1 \times D_\epsilon$ so that its characteristic foliation also has slope $s$ on $T_{r_0}$. We choose a smooth function $f : [r_0, \epsilon] \rightarrow \R$ satisfying
\begin{itemize}
\item $f' > 0$,
\item $f(r) = r^2$ near $\epsilon$,
\item $f(r) = r^2 - 2 r_0^2$ near $r_0$.
\end{itemize}
See Figure~\ref{fig:f}. We then define $$\widetilde \alpha^0_- \coloneqq -dt + f(r) d\theta$$ inside of $S^1 \times \left([r_0, \epsilon] \times S^1 \right)$, which coincides with $\alpha^0_-$ near $T_\epsilon$. It is a negative contact form whose kernel ``twists more'' than the one of $\alpha^0_-$. Moreover, its characteristic foliation on $T_{r_0}$ coincides with the one of $\alpha^0_+$. It is then easy to perform a \emph{contact cut} (see~\cite{L01,BE13}) along $T_{r_0}$ for both $\alpha^0_+$ and $\widetilde \alpha^0_-$. This is possible since their moment maps for the $S^1$-action on a neighborhood of $T_{r_0}$ given by a translation of slope $s$ in the $(\theta, t)$ direction coincide. After this operation, we obtain two contact structures $\xi^s_\pm$ on the Dehn filling of $S^1 \times \left( [r_0, \epsilon] \times S^1\right)$ of slope $s$ along $T_{r_0}$, which coincide with $\xi'_\pm$ near $T_\epsilon$. After extending them by $\xi'_\pm$ outside of $S^1 \times D_\epsilon$, we obtain two contact structures, still denoted by $\xi^s_\pm$, on $M_K(s)$. If $r_0$ is sufficiently small, i.e., if $\vert s \vert$ is sufficiently large, then $(\xi^s_-, \xi^s_+)$ is taut in $S^1 \times D_\epsilon$. Since $(\xi'_-, \xi'_+)$ is taut in $M$, it is also taut in $M \setminus \left( S^1 \times D_\epsilon \right)$ and $(\xi^s_-, \xi^s_+)$ is taut in $M_K(s)$. To perform a surgery of large positive slope, we can apply a similar argument to $\xi'_-$, or equivalently, reverse the orientation of $M$ and consider the positive contact pair $(\xi'_+, \xi'_-)$.

\begin{figure}[h]
    \centering
    \includegraphics[width=0.5\linewidth]{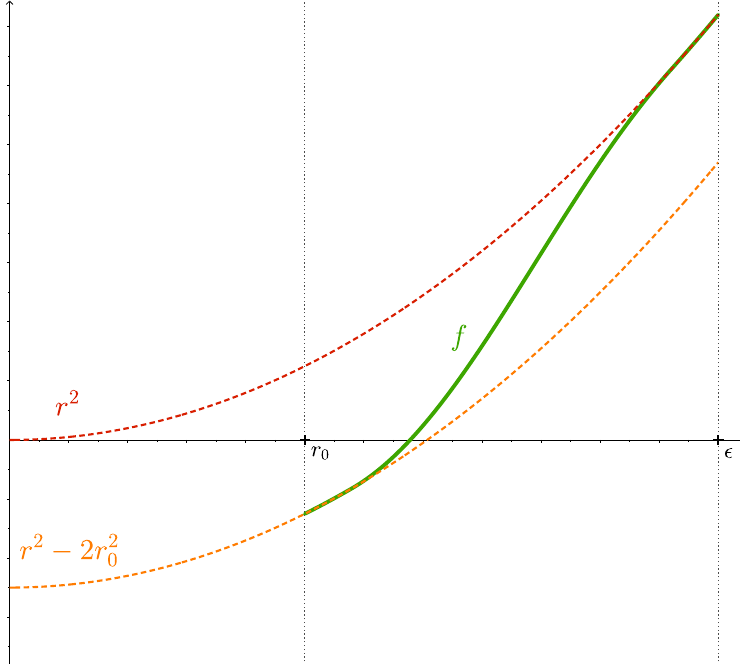}
    \caption{The function $f$.}
    \label{fig:f}
\end{figure}

\end{proof}

\begin{proof}[Proof of Theorem~\ref{thmintro:transfol}] 
Let $\mathcal{F}$ be a taut $C^0$-foliation on $M$ different from the standard foliation on $S^1 \times S^2$, and let $K$ be a framed knot transverse to $\mathcal{F}$. By~\cite[Theorem 1.2]{KR17}, $T \mathcal{F}$ can be approximated by a taut contact pair $(\xi_-, \xi_+)$ such that $K$ is positively transverse to $\xi_\pm$. By Proposition~\ref{prop:tautsurgery}, for any $s \in \mathbb{Q}$ such that $\vert s \vert$ is sufficiently large, $M_K(s)$ carries a taut contact pair $(\xi'_-, \xi'_+)$ such that $K'$, the image of $K$ after surgery, is positively transverse to $\xi'_\pm$. This implies that $K'$ is transverse to $\eta'_u$, the unstable plane field of $(\xi'_-, \xi'_+)$. Applying Theorem~\ref{thmintro:taut}, we readily obtain a taut $C^0$-foliation $\mathcal{F}'$ in $M_K(s)$ such that $K'$ is transverse to $\mathcal{F}'$.
\end{proof}

\begin{proof}[Proof of Corollary~\ref{corintro:slope}]
Note that $s \in \mathcal{S}_K$ if and only if $S^3_K(s)$ admits a taut $C^0$-foliation transverse to $K'_s$, the image of $K$ after surgery. If $s = p/q$, for coprime $p$ and $q$ (or $p=0$ and $q=1$ if $s=0$), the choice of a framing for $K'_s$ is equivalent to the choice of integers $u, v \in \mathbb{Z}$ such that $u p - v q = 1$, which is also equivalent to the choice of a matrix in $\mathrm{SL}(2,\mathbb{Z})$ representing a mapping class of the torus $T^2 = S^1 \times S^1$ sending the meridian curve (the second component) to a curve of slope $s$. This matrix can be written as
$$A = \begin{pmatrix}
u & q \\
v & p 
\end{pmatrix},$$
and sends a curve of slope $s' = p'/q'$ to a curve of slope $$r = \frac{v q' + pp'}{u q' + qp'} = s \cdot \frac{s' + v/p }{s'+ u/q}$$ if $s \neq 0$, and of slope
$$r = - \frac{1}{u+s'}$$
otherwise (in which case $v=-1$). In both cases, $r$ converges to $s$ as $s'$ goes to $\pm \infty$, and every slope $r$ sufficiently close to $s$ can be written as such for some $s' \in \mathbb{Q}$ with $\vert s' \vert$ sufficiently large. Since an $s$-slope surgery on a nontrivial knot in $S^3$ never produces $S^1 \times S^2$ by Gabai's \emph{Property R} theorem~\cite{G87} (and some elementary homological considerations for $s \neq 0$), the proof now easily follows from Theorem~\ref{thmintro:transfol}.
\end{proof}

%%%%%%%%%%%%%%%%%%%%%%%%%%%%%%%%%%%%%%%%%%%%%%%%%%%%%%%%%%%%%%%%%%%%%%%%%%%%%%%%%%%%%%%%%%%%%%%%%%%%%%%%%%%%%%%%%%%%%%%%%%%%%%%%%%%%%%%%%%%%%%%%%%%%%%%%%%%%%%
    \section{Liouville pairs}
%%%%%%%%%%%%%%%%%%%%%%%%%%%%%%%%%%%%%%%%%%%%%%%%%%%%%%%%%%%%%%%%%%%%%%%%%%%%%%%%%%%%%%%%%%%%%%%%%%%%%%%%%%%%%%%%%%%%%%%%%%%%%%%%%%%%%%%%%%%%%%%%%%%%%%%%%%%%%%

%%%%%%%%%%%%%%%%%%%%%%%%%%%%%%%%%%%%%%%%%%%%%%%%%%%%%%%%%%%%%%%%%%%%%
        \subsection{Elementary properties}
%%%%%%%%%%%%%%%%%%%%%%%%%%%%%%%%%%%%%%%%%%%%%%%%%%%%%%%%%%%%%%%%%%%%%

We recall that a Liouville pair is a pair of contact forms $(\alpha_-, \alpha_+)$ on $M$, negative a positive, such that the $1$-form $$\lambda = e^{-s} \alpha_- + e^s \alpha_+$$ is a positive Liouville form on $\R_s \times M$. In other terms, $\omega \coloneqq d\lambda$ is a nondegenerate $2$-form on $\R \times M$, and $\omega \wedge \omega$ is a positive volume form.

If $(\alpha_-, \alpha_+)$ is a Liouville pair and $\sigma : M \rightarrow \R$ is a smooth function, then $(e^{-\sigma} \alpha_-, e^\sigma \alpha_+)$ is also a Liouville pair. We say that these two Liouville pairs are \textbf{equivalent}. Up to equivalence, we can always assume that $(\alpha_-, \alpha_+)$ is \textbf{balanced}, i.e.,
$$\alpha_+ \wedge d\alpha_+ = - \alpha_- \wedge d\alpha_-.$$
In that case, we denote by $$\mathrm{dvol} \coloneqq \alpha_+ \wedge d\alpha_+$$ the associated volume form. Writing $$d(\alpha_- \wedge \alpha_+) = f_0 \, \mathrm{dvol}$$ as before, where $f_0 : M \rightarrow \R$ is a smooth function, elementary calculations show that the Liouville condition is equivalent to
\begin{align} \label{eq:condf0}
f_0 > -2.
\end{align}

The following lemma was already observed in~\cite{Mas22}.

\begin{lem}
Let $(\alpha_-, \alpha_+)$ be a Liouville pair on $M$. Then the underlying contact pair $(\xi_-, \xi_+) = (\ker \alpha_-, \ker \alpha_+)$ is positive, and $r_u > 0$ (see Proposition~\ref{prop:distrib} for the notation).
\end{lem}

\begin{proof}
The fact that $(\xi_-, \xi_+)$ is positive immediately follows from Lemma~\ref{lem:f0} and inequality~\eqref{eq:condf0}. For $\alpha_u$ defined by~\eqref{eq:alphau}, the expression~\eqref{eq:ru} together with inequality~\eqref{eq:condf0} imply $r_u > 0$.
\end{proof}

If $(\xi_-, \xi_+)$ is a contact pair and $(\alpha_-, \alpha_+)$ is a Liouville pair such that $(\ker \alpha_-, \ker \alpha_+) = (\xi_-, \xi_+)$, we say that $(\xi_-, \xi_+)$ \textbf{supports} $(\alpha_-, \alpha_+)$. The following lemma shows that such a Liouville pair (if it exists) is unique up to homotopy.

\begin{lem} Let $(\xi_-, \xi_+)$ be a contact pair. The space of Liouville pairs $(\alpha_-, \alpha_+)$ supported by $(\xi_-, \xi_+)$ is either empty or contractible.
\end{lem}

\begin{proof}
Let $(\alpha_-, \alpha_+)$ be a Liouville pair for $(\xi_-, \xi_+)$. We can assume that it is balanced. It defines a volume form $\mathrm{dvol}$ and a vector field $X$. It is enough to show that the space of \emph{balanced} Liouville pairs for $(\xi_-, \xi_+)$ is contractible. Any such Liouville pair can be written as
$$\alpha'_\pm = e^\sigma \alpha_\pm$$
for some smooth function $\sigma : M \rightarrow \R$. Moreover, the Liouville condition is equivalent to
\begin{align} \label{eq:condliouv2}
    2 X \cdot \sigma + f_0 > -2.
\end{align}
As a result, the space of balanced Liouville pairs for $(\xi_-, \xi_+)$ is diffeomorphic to the space of functions $\sigma$ satisfying~\eqref{eq:condliouv2}. The latter is obviously convex, hence contractible.
\end{proof}

%%%%%%%%%%%%%%%%%%%%%%%%%%%%%%%%%%%%%%%%%%%%%%%%%%%%%%%%%%%%%%%%%%%%%
        \subsection{From hypertaut foliations to Liouville pairs} \label{sec:liouv}
%%%%%%%%%%%%%%%%%%%%%%%%%%%%%%%%%%%%%%%%%%%%%%%%%%%%%%%%%%%%%%%%%%%%%

By work of Sullivan~\cite{S76}, a cooriented $C^1$-foliation admits no nontrivial invariant transverse measure if and only if it is hypertaut in the sense of Definition~\ref{def:hypertaut}, i.e., it admits an \emph{exact} two-form positive on its leaves. For $C^0$-foliations, Sullivan's argument still shows that hypertaut foliations have no invariant transverse measures, but the converse fails in general. Indeed, as we mentioned earlier, the presence of phantom Reeb components can obstruct even tautness. However, as we will show in work in preparation, the converse holds up to an arbitrarily small isotopy.

\begin{lem}[\cite{MasHypertaut}] \label{lem:hypertaut}
Let $\mathcal{F}$ be a cooriented $C^0$-foliation on $M$, and assume that it has no nontrivial invariant transverse measure. Then it is isotopic, via an arbitrarily small isotopy, to a hypertaut $C^0$-foliation.
\end{lem}

Equivalently, $\mathcal{F}$ is hypertaut if for some/any Riemannian metric $g$ on $M$, there exists a smooth nonsingular \emph{zero-flux} vector field $\upsilon$ transverse to $\mathcal{F}$. Here, the \emph{flux} of $\upsilon$ through a closed oriented surface $\Sigma \subset M$ is defined as
$$\mathrm{Flux}(\upsilon) [\Sigma] \coloneqq \int_\Sigma \langle \upsilon, \nu_\Sigma \rangle \, \omega_\Sigma,$$
where $\nu_\Sigma$ is the positive normal vector to $\Sigma$ and $\omega_\Sigma$ is the area form on $\Sigma$ induced by $g$.

The next proposition is a version of~\cite[Corollary 3.2.5]{ET}.

\begin{prop} \label{prop:fill}
Let $\mathcal{F}$ be a hypertaut $C^0$-foliation on $M$. Then $\mathcal{F}$ can be $C^0$-approximated by contact structures $\xi_-$ and $\xi_+$, negative and positive, respectively, such that there exists a Liouville structure on $[-1,1] \times M$ which is an exact filling of $\big(-\!M, \xi_- \big) \sqcup \big(M, \xi_+ \big)$.
\end{prop}

\begin{proof}
Let $\beta$ be a smooth $1$-form on $M$ such that $d \beta$ dominates $\mathcal{F}$. Since $\mathcal{F}$ has no closed (spherical) leaf, it can be approximated by negative and positive contact structures $\xi_-$ and $\xi_+$ by~\cite{KR17}. We can further assume that $d\beta_{\vert \xi_\pm} > 0$. Let $\alpha$ be a continuous $1$-form defining $T\mathcal{F}$, and $\widetilde{\alpha}$ be a smoothing of $\alpha$ satisfying $\widetilde{\alpha} \wedge d\beta > 0$. For $\epsilon > 0$, we consider the $1$-form
$$\lambda \coloneqq \epsilon t \widetilde{\alpha} + \beta$$
on $[-1,1]_t \times M$. If $\epsilon$ is sufficiently small, then $\omega \coloneqq d \lambda$ is symplectic and dominates $\xi_\pm$ on $\{\pm1\} \times M$. We now appeal to~\cite[Proposition 3.1]{E91} (see also~\cite[Proposition 4.1]{E04}) to modify $\lambda$ near $\{\pm1\} \times M$ and obtain the desired Liouville filling.
\end{proof}

It is not clear to us how to recover a hypertaut foliation from a Liouville fillable positive contact pair. It is also unclear whether the Liouville structure obtained from a hypertaut $C^0$-foliation by the previous proposition is unique up to homotopy (the contact approximations of $\mathcal{F}$ might not be unique). However, for $C^2$ hypertaut foliations, it is possible to construct a more specific type of Liouville structure, namely, one induced by a Liouville pair. This construction relies on the following proposition which is the main result of Jonathan Zung's article~\cite{Z21}:

\begin{prop}[\cite{Z21}] \label{prop:Zung}
If $\mathcal{F}$ is a $C^2$ hypertaut foliation, then there exist $1$-forms $\alpha$ and $\beta$ such that 
\begin{align*}
\ker \alpha = T \mathcal{F}, \qquad \alpha \wedge d\beta > 0, \qquad \beta \wedge d\alpha \geq 0.
\end{align*}
\end{prop}
The first inequality simply means that $d\beta$ is a dominating $2$-form for $\mathcal{F}$.

\begin{cor} \label{cor1}
If $\mathcal{F}$ is a hypertaut $C^2$ foliation, then there exists a Liouville pair $(\alpha_-, \alpha_+)$ on $M$ such that the contact structures $\xi_\pm = \ker \alpha_\pm$ are $C^0$-close to $T\mathcal{F}$.
\end{cor}

\begin{proof}
For $\epsilon > 0$, we define
\begin{align} \label{eq:liouvpair}
    \alpha_\pm \coloneqq \epsilon \beta \pm \alpha,
\end{align}
where $\alpha$ and $\beta$ are as in Proposition~\ref{prop:Zung}. Following~\cite{ET}, we write $$\langle \alpha, \beta \rangle \coloneqq \alpha \wedge d \beta + \beta \wedge d \alpha.$$ One easily computes 
\begin{align*}
    \alpha_+ \wedge d\alpha_+ &= \epsilon \langle \alpha, \beta \rangle  + O(\epsilon^2), \\
    \alpha_- \wedge d\alpha_- &= - \epsilon \langle \alpha, \beta \rangle + O(\epsilon^2), \\
    d(\alpha_- \wedge \alpha_+) &= - 2 \epsilon \, d(\alpha \wedge \beta).
\end{align*}
For $\epsilon$ small enough, $\alpha_\pm$ are contact forms with opposite orientations. Moreover, $(\alpha_-, \alpha_+)$ is a Liouville pair if and only if $$ - d(\alpha_- \wedge \alpha_+) < 2 \sqrt{- (\alpha_+ \wedge d\alpha_+) \cdot (\alpha_- \wedge d\alpha_-)},$$ which is equivalent to $$d(\alpha \wedge \beta) < \langle \alpha, \beta \rangle + O(\epsilon).$$ This condition is satisfied for $\epsilon$ small enough, since $$d(\alpha \wedge \beta) = \beta \wedge d \alpha - \alpha \wedge d\beta < \beta \wedge d\alpha + \alpha \wedge d\beta = \langle \alpha, \beta \rangle.$$
The $1$-forms $\alpha_\pm$ might only be $C^1$, but they can easily be smoothed to yield a smooth Liouville pair.
\end{proof}

%%%%%%%%%%%%%%%%%%%%%%%%%%%%%%%%%%%%%%%%%%%%%%%%%%%%%%%%%%%%%%%%%%%%%
        \subsection{Proof of Theorem~\ref{thmintro:liouv}}
%%%%%%%%%%%%%%%%%%%%%%%%%%%%%%%%%%%%%%%%%%%%%%%%%%%%%%%%%%%%%%%%%%%%%

Let $(\xi_-, \xi_+)$ be a contact pair on $M$ supporting a Liouville pair $(\alpha_-, \alpha_+)$. By Proposition~\ref{prop:Reeb} below, we can assume that the Reeb vector field of $\alpha_+$ is positively transverse to $\eta_u$, i.e., ${d\alpha_+}_{\vert \eta_u} > 0$. Since $\xi_\pm$ are (strongly) semi-fillable, they are both tight. By Theorem~\ref{thmintro:int}, $\eta_u$ can be $C^0$-approximated by a $C^0$-foliation $\mathcal{F}$, and we can further arrange that ${d\alpha_+}_{\vert T \mathcal{F}} > 0$. Therefore, $\mathcal{F}$ is hypertaut. In particular, it is taut and the above Reeb vector fields have no contractible closed orbits, hence $\xi_\pm$ are hypertight.

\begin{prop} \label{prop:Reeb}
Let $(\alpha_-, \alpha_+)$ be a Liouville pair on $M$, with underlying positive contact pair $(\xi_-, \xi_+)$ and unstable plane field $\eta_u$. Then $(\alpha_-, \alpha_+)$ is equivalent to a Liouville pair $(\alpha'_-, \alpha'_+)$ whose  Reeb vector fields $R'_\pm$ are both transverse (with opposite signs) to $\eta_u$.
\end{prop}

\begin{proof}
As usual, we assume that $(\alpha_-, \alpha_+)$ is balanced, and we use the notations of the proof of Proposition~\ref{prop:distrib}. By~\eqref{eq:alphau} and~\eqref{eq:alphas}, we can write
\begin{align}
\alpha_- &= \frac{1}{2\cosh(\sigma_u - \sigma_s)} \left( e^{\sigma_s} \alpha_u + e^{\sigma_u} \alpha_s \right), \\
\alpha_+ &= \frac{1}{2\cosh(\sigma_u - \sigma_s)} \left( e^{-\sigma_s} \alpha_u - e^{-\sigma_u} \alpha_s \right).
\end{align}
Let $\sigma : M \rightarrow \R$ be some smooth function and $(\alpha'_-, \alpha'_+) \coloneqq (e^{-\sigma} \alpha_-, e^\sigma \alpha_+)$. We compute $\alpha_s\big(R'_\pm\big)$ in terms of $\sigma$. From the identities
\begin{align}
\alpha'_\pm\big( R'_\pm \big) &= 1, \\
\mathcal{L}_X \alpha'_\pm \big( R'_\pm \big) &=0,
\end{align}
we get (after some slightly tedious computations)
\begin{align}
\alpha_s \big( R'_+ \big) &= - \overset{> 0}{\overbrace{\frac{2 e^{\sigma_u - \sigma} \cosh(\sigma_u - \sigma_s)}{X \cdot (\sigma_u - \sigma_s) + r_u - r_s}}} \, \left( X \cdot (\sigma - \sigma_s) - X \cdot \ln \cosh (\sigma_u - \sigma_s) + r_u \right), \\
\alpha_s \big( R'_- \big) &= - \underset{> 0}{\underbrace{\frac{2 e^{\sigma - \sigma_u} \cosh(\sigma_u - \sigma_s)}{- X \cdot (\sigma_u - \sigma_s) + r_u - r_s}}} \, \left( X \cdot (\sigma - \sigma_s) + X \cdot \ln \cosh (\sigma_u - \sigma_s) - r_u \right).
\end{align}
From~\eqref{eq:sigmau}, \eqref{eq:sigmas}, \eqref{eq:ru}, and~\eqref{eq:rs}, one easily computes
\begin{align*}
    \pm X \cdot \big(\sigma_u - \sigma_s\big) + r_u - r_s = e^{\pm 2\sigma_u} + e^{\pm 2 \sigma_s} > 0.
\end{align*}
Furthermore, we claim that $$- X \cdot\ln \cosh (\sigma_u - \sigma_s) +r_u > 0.$$
Indeed, by~\eqref{eq:sigmau} and~\eqref{eq:sigmas},
\begin{align*}
X \cdot\ln \cosh (\sigma_u - \sigma_s) &= \tanh(\sigma_u - \sigma_s) \, X \cdot (\sigma_u - \sigma_s) \\
&= \tanh(\sigma_u - \sigma_s) \big( \sinh(2\sigma_u) + \sinh(2\sigma_s) \big) \\
&= 2\tanh(\sigma_u - \sigma_s) \big(\sinh(\sigma_u + \sigma_s) \cosh(\sigma_u - \sigma_s) \big)\\
&= 2 \sinh(\sigma_u-\sigma_s) \sinh(\sigma_u + \sigma_s)\\
&= \cosh(2\sigma_u) - \cosh(2\sigma_s).
\end{align*}
Here, we used the classical formulae
\begin{align*}
\sinh(x) + \sinh(y) &= 2\sinh\left( \frac{x+y}{2} \right) \cosh\left( \frac{x-y}{2} \right)\\
2\sinh(x) \sinh(y) &=  \cosh(x+y) - \cosh(x-y).
\end{align*}
Therefore, using~\eqref{eq:ru},
\begin{align*}
- X \cdot\ln \cosh (\sigma_u - \sigma_s) +r_u &= \cosh(2\sigma_s) - \cosh(2\sigma_u) + \cosh(2\sigma_u) + f_0/2\\
&= \cosh(2\sigma_s) + f_0/2\\
&>0,
\end{align*}
since $f_0 > -2$. Let us fix some $\epsilon > 0$ such that $$- X \cdot\ln \cosh (\sigma_u - \sigma_s) + r_u  > \epsilon.$$
By Lemma~\ref{lem:smooth}, we can choose $\sigma$ so that $$\vert X \cdot (\sigma - \sigma_s)\vert_{C^0} < \epsilon.$$ This ensures 
\begin{align*}
\alpha_s \big( R'_+ \big) &< 0, \\
\alpha_s \big( R'_- \big) &> 0,
\end{align*}
so $R'_\pm$ are both transverse to $\eta_u = \ker \alpha_s$ as desired. Note that we do not need $\sigma$ to be $C^0$-close to $\sigma_s$.
\end{proof}

%%%%%%%%%%%%%%%%%%%%%%%%%%%%%%%%%%%%%%%%%%%%%%%%%%%%%%%%%%%%%%%%%%%%%
        \subsection{The skeleton} \label{sec:skel}
%%%%%%%%%%%%%%%%%%%%%%%%%%%%%%%%%%%%%%%%%%%%%%%%%%%%%%%%%%%%%%%%%%%%%

For a general Liouville manifold $(V, \lambda)$ with Liouville vector field $Z$ and Liouville flow $(\phi^t_Z)$, the \textbf{skeleton} is defined as
$$\mathfrak{skel}(V, \lambda) \coloneqq \bigcup_{\substack{K \subset V \\\mathrm{compact}}} \bigcap_{t > 0} \phi^{-t}_Z(K).$$

In this section, we determine the geometric structure of the skeleton of a Liouville pair, defined as the skeleton of the associated Liouville structure. This result applies in particular to \emph{Anosov Liouville structures}, defined in~\cite{Mas22}.

\begin{prop} \label{prop:skel}
    Let $(\alpha_-, \alpha_+)$ be a Liouville pair on $M$. The skeleton of the associated Liouville structure $$\lambda = e^{-s} \alpha_- + e^s \alpha_+$$ on $V=\R_s \times M$ is of the form $$\mathfrak{skel}(V, \lambda) = \mathrm{graph}(\sigma) = \big\{(\sigma(p), p) \ \vert \  p \in M\big\}$$ for some continuous function $\sigma : M \rightarrow \R$. Moreover, $\sigma$ is continuously differentiable along the vector field $X$ associated with the contact pair $(\xi_-, \xi_+)=(\ker \alpha_-, \ker \alpha_+)$, and the restriction of the Liouville flow to the skeleton is orbit equivalent to the flow of $X$.\footnote{Two flows are orbit equivalent if there exists a homeomorphism sending the oriented flow lines of the first to oriented flow lines of the second, possibly with a different parametrization.}
\end{prop}

\begin{proof}
    First of all, we can assume that $(\alpha_-, \alpha_+)$ is balanced, after applying a graphical diffeomorphism $\varphi : (s,p) \mapsto (s+ \sigma_0(p), p)$ for a suitable smooth function $\sigma_0 : M \rightarrow \R$. Then, with the previous notations, the Liouville vector field for $\lambda$ is given by
    \begin{align*}
        Z(s,p) &= \frac{1}{\cosh(2s) + f_0(p)/2} \left(  \big( \sinh(2s) + g_0(p)/2\big) \, \partial_s + X(p) \right).
    \end{align*}
    Writing 
    \begin{align*}
        F(s,p) &= \sinh(2s) + g_0(p)/2, \\
        \widetilde{Z} &=  F \partial_s + X,
    \end{align*}
    the vector fields $Z$ and $\widetilde{Z}$ share the same skeleton. The flow line of $\widetilde{Z}$ starting at $(s_0, p_0)$ is given by
    \begin{align*}
    \phi^t_{\widetilde{Z}}(s_0, p_0) &= \left(\sigma(t), \phi^t_X(p_0) \right),
    \end{align*}
    where $\sigma$ satisfies the ODE
    \begin{equation*}
    E_{(s_0, p_0)} :  \left\{
    \begin{aligned}
    \dot{\sigma}(t) &= F\big(\sigma(t), \phi^t_X(p_0)\big), \\
    \sigma(0) &= s_0.
    \end{aligned}
  \right.
\end{equation*}
    By definition, $(s_0, p_0) \in \mathfrak{skel} = \mathfrak{skel}(V, \lambda)$ if and only if the corresponding maximal solution to $E_{(s_0, p_0)}$ is bounded in positive time. Lemma~\ref{lem:ODE} implies that for every $p_0 \in M$, there exists a unique $\sigma(p_0) \in \R$ such that $(\sigma(p_0), p_0) \in \mathfrak{skel}$; with the notation of the proof of Proposition~\ref{prop:distrib}, $\sigma = \sigma_u$, which is continuous and continuously differentiable along $X$. Therefore, under the balancing condition, $\mathfrak{skel} = \mathrm{graph}(\sigma_u)$. Finally, the restriction of the Liouville flow to the skeleton is clearly orbit equivalent to the flow of $X$ because of the form of $\widetilde{Z}$.
\end{proof}

\begin{rem}
We see from the proof that the regularity of the skeleton, i.e., the regularity of the function $\sigma$, is the same as the regularity of the function $\sigma_u$. If $X$ is Anosov, then it is as regular as the weak-stable bundle of $X$, which is always $C^1$ but not necessarily $C^2$.
\end{rem}

Proposition~\ref{prop:skel} shows in particular that the skeleton of the Liouville structure associated to a Liouville pair is a codimension one topological submanifold homeomorphic to $M$. It would be interesting to investigate the converse:

\begin{question}
    Let $\lambda_0$ be a Liouville structure on $V = \R \times M$. If $\mathfrak{skel}(V, \lambda_0) \subset V$ is a codimension one compact embedded topological submanifold, does there exist a Liouville pair $(\alpha_-, \alpha_+)$ on $M$ with associated Liouville structure $\lambda$ on $V$ such that $(V, \lambda_0)$ and $(V, \lambda)$ are exact symplectomorphic? 
\end{question}

%%%%%%%%%%%%%%%%%%%%%%%%%%%%%%%%%%%%%%%%%%%%%%%%%%%%%%%%%%%%%%%%%%%%%
        \subsection{Semi-Anosov flows and unique integrability} \label{sec:semiano}
%%%%%%%%%%%%%%%%%%%%%%%%%%%%%%%%%%%%%%%%%%%%%%%%%%%%%%%%%%%%%%%%%%%%%

In this section, we define and study a certain type of $3$-dimensional flows whose behavior interpolates between projectively Anosov and Anosov flows.

\begin{defn} Let $\Phi = (\phi_t)_t$ be a smooth flow on $M$ generated by a vector field $X$, and $\Lambda \subset M$ be a $\Phi$-invariant compact subset disjoint from the singular set of $X$. A \textbf{semi-Anosov splitting} for $\Phi$ on $\Lambda$ is a continuous, $\Phi$-invariant splitting $$T_\Lambda M = E^{ws} \oplus E^u$$ along $\Lambda$, where $E^u$ is $1$-dimensional, $E^{ws}$ is $2$-dimensional, $X \in E^{ws}$, and for some/any metric $g$ on $M$, there exist constants $C > 0$ and $a, b > 0$ such that
\begin{itemize}
\item For all $v_u \in E^u$ and $t \geq 0$, 
$$\Vert d\phi_t(v_u) \Vert \geq C e^{at} \Vert v_u \Vert,$$
\item For all $v_s \in E^{ws}\setminus\{0\}$, $v_u \in E^u \setminus \{0\}$, and $t \geq 0$, 
$$\frac{\Vert d\phi_t(v_s) \Vert}{\Vert v_s \Vert} \leq C e^{-bt} \frac{\Vert d\phi_t(v_u) \Vert}{\Vert v_u \Vert}.$$
\end{itemize}
We call $E^u$ the \textbf{(strong-)unstable bundle} of $\Phi$, and $E^{ws}$ the \textbf{weak-stable bundle} of $\Phi$. We also call $E^{wu} \coloneqq \langle X \rangle \oplus E^u$ the \textbf{weak-unstable bundle} of $\Phi$. We call $\Phi$ \textbf{semi-Anosov} if it is nonsingular and admits a semi-Anosov splitting on $M$.
\end{defn}

By definition, semi-Anosov flows are projectively Anosov, as the second bullet point in the definition is exactly the existence of a dominated splitting. Moreover, closed orbits or semi-Anosov flows are not contracting. However, a semi-Anosov flow might fail to be Anosov.

\begin{question} Is there a $3$-manifold $M$ carrying a semi-Anosov flow, but no Anosov flow?
\end{question}

The next proposition is very similar to the Anosov case, and the proof is essentially the same.

\begin{prop} \label{prop:semiint}
If $\Phi$ admits a semi-Anosov splitting on $\Lambda$, then $E^u$ is uniquely integrable and its leaves are smooth.
\end{prop}

\begin{proof} It follows from a version of the Hadamard--Perron theorem~\cite[Theorem B.5.2]{FH19} which implies a version of the Unstable Manifold theorem~\cite[Theorem 6.1.1]{FH19}.
\end{proof}

We now establish a connection between Liouville pairs and semi-Anosov flows.

\begin{prop} \label{prop:semitrans}
Let $(\alpha_-, \alpha_+)$ be a Liouville pair with associated vector field $X$, and let $\Lambda \subset M$ be a compact subset disjoint from $\Delta$ and invariant by $X$. Then the flow of $X$ admits a semi-Anosov splitting along $\Lambda$.
\end{prop}

\begin{proof}
Since $\Delta \cap \Lambda = \varnothing$, $\eta_s$ and $\eta_u$ are defined on $\Lambda$, and there exist continuous, nonvanishing $1$-forms $\alpha_s$, $\alpha_u$ such that $\ker \alpha_{u/s} = \eta_{s/u}$ and $$\mathcal{L}_X \alpha_{s/u} = r_{s/u} \alpha_{s/u},$$ where $r_s < r_u$ and $r_u > 0$. We let $E^{ws} \coloneqq \eta_s$ and we want to find a suitable $1$-dimensional, $\Phi$-invariant subbundle $E^u \subset \eta_u$. On $N_X \coloneqq TM \slash  \langle X \rangle$ restricted to $\Lambda$, we have a dominated splitting $$N_X = \overline{\eta}_s \oplus \overline{\eta}_u,$$ where $\overline{\eta}_{s/u} \coloneqq \eta_{s/u} \slash \langle X \rangle$. Dualizing $(\alpha_s, \alpha_u)$, we get continuous sections $\overline{e}_{s/u}$ of $\overline{\eta}_{s/u}$ and we define a (continuous) metric $\overline{g}$ on $N_X$ for which $(\overline{e}_s, \overline{e}_u)$ is orthonormal. Moreover, since $r_u > 0$, there exist constants $C > 0$ and $a > 0$ such that for all $v \in \overline{\eta}_u$ and $t \in \R$, $$\Vert d\phi_t(v) \Vert \geq C e^{at} \Vert v \Vert.$$
The proof of~\cite[Proposition 1.1]{D87} readily implies that there exists a $\Phi$-invariant $1$-dimensional subbundle $E^u \subset \eta_u$ defined along $\Lambda$ satisfying the desired conditions.
\end{proof}

We conclude with the following result, which is essentially the only instance in which we can ensure that $\eta$ is a uniquely integrable plane field:

\begin{cor} \label{cor:transtaut}
Let $(\xi_-, \xi_+)$ be a regular contact pair on $M$ supporting a Liouville pair $(\alpha_-, \alpha_+)$. If $\Delta$ has no saddle singularities, then $\eta_u$ is \underline{uniquely integrable}, and its integral foliation is a hypertaut $C^0$-foliation.
\end{cor}

\begin{proof}
If $\Delta = \varnothing$, combining Proposition~\ref{prop:semiint}, Proposition~\ref{prop:semitrans} and Proposition~\ref{prop:Reeb}, $\eta_u$ is uniquely integrable, its integral foliation is a $C^0$-foliation, and there exists a smooth contact form $\alpha_+'$ whose Reeb vector field is transverse to $\mathcal{F}$. Equivalently, $d\alpha_+'$ is positive on the leaves of $\mathcal{F}$, which implies hypertautness.

If $\Delta \neq \varnothing$ and $\Delta = \Delta_\mathrm{so}$, we consider a small tubular neighborhood $U$ of $\Delta$ such that $X$ is transverse to $\partial U$ and outward pointing. We then define $$\mathcal{V} \coloneqq \bigcup_{t > 0} \phi_t(U), \qquad \Lambda \coloneqq M \setminus \mathcal{V},$$ so that $\Lambda$ is compact and invariant under the flow of $X$. By the proof of Proposition~\ref{prop:locint}, $\eta_u$ is uniquely integrable on $\mathcal{V}$, since it is uniquely integrable near $\Delta$. Moreover, the reasoning of the previous paragraph shows that $\eta_u$ is uniquely integrable on $\Lambda$ as well.
\end{proof}

Unfortunately the hypotheses of this corollary are quite restrictive. We ask:

\begin{question}
Can the condition $\Delta_\mathrm{sa} = \varnothing$ in Corollary~\ref{cor:transtaut} be removed? To show this, it would be sufficient to find a uniquely integrable vector field $X_u$ defined away from $\Delta$ such that $\eta_u = \mathrm{span}\{ X, X_u\}$.
\end{question}

 \label{part:I}

%%%%%%%%%%%%%%%%%%%%%%%%%%%%%%%%%%%%%%%%%%%%%%%%%%%%%%%%%%%%%%%%%%%%%%%%%%%%%%%
\newpage
\part{Branching foliations} \label{part:tech}
%%%%%%%%%%%%%%%%%%%%%%%%%%%%%%%%%%%%%%%%%%%%%%%%%%%%%%%%%%%%%%%%%%%%%%%%%%%%%%%

In this part, $(X, \eta)$ denotes a polarized vector field on $M$ which is admissible (see Definition~\ref{def:adm}) and without disks of tangency (see Definition~\ref{def:disk}). We also choose some arbitrary Riemannian metric $\boldsymbol{g}$ on $M$.

The goal of this part is to prove the Technical Theorem~\ref{thm:tech}. It is completely independent of the first part and does not rely on any notion from contact or symplectic topology.

%%%%%%%%%%%%%%%%%%%%%%%%%%%%%%%%%%%%%%%%%%%%%%%%%%%%%%%%%%%%%%%%%%%%%
        \section{Preliminaries}
%%%%%%%%%%%%%%%%%%%%%%%%%%%%%%%%%%%%%%%%%%%%%%%%%%%%%%%%%%%%%%%%%%%%%

Because of the lack of smoothness and unique integrability of $\eta$, we will not be able to construct a foliation tangent to $\eta$ directly. However, we will construct a weaker structure tangent to $\eta$ called a \emph{branching foliation}.

Branching foliations are defined in~\cite[Definition 4.1]{BI08}. Roughly speaking, those are collections of `maximal' surfaces which are allowed to touch but not to cross, and whose union is $M$. We recall the precise definition for the convenience of the reader.

\begin{defn} \label{def:branchfol}
A \textbf{leaf collection} tangent to $\eta$ is a set of maps  $\mathscr{F} = \{ f_i : \Sigma_i \rightarrow M\}_{i \in I}$ indexed by a set $I$ such that for every $i \in I$, 
    \begin{itemize}
        \item $\Sigma_i$ is a smooth surface (i.e., a $2$-dimensional manifold without boundary),
        \item $f_i : \Sigma_i \rightarrow M$ is a $C^1$ immersion tangent to $\eta$.
    \end{itemize}
A leaf collection $\mathscr{F}$ is a \textbf{branching foliation} if
    \begin{itemize}
        \item For every $i \in I$, the length metric for the intrinsic $C^0$ Riemannian metric $\boldsymbol{g}_i \coloneqq f_i^*\boldsymbol{g}$ on $\Sigma_i$ is complete,
        \item For all $i,j \in I$, $f_i$ and $f_j$ have no topological crossings (see~\cite[Definition 4.1]{BI08}),
        \item The leaves of $\mathscr{F}$ cover $M$:  $$\bigcup_{i \in I} f_i(\Sigma_i) = M.$$
    \end{itemize}
\end{defn}

Notice that this definition does not depend on the choice of the metric $\boldsymbol{g}$ since $M$ is compact. Crucially, branching foliations can be $C^0$ perturbed into genuine foliations by separating the leaves without introducing crossings. More precisely, the combination of \cite[Lemma 7.1]{BI08} and~\cite[Theorem 7.2]{BI08} implies the following

\begin{thm} \label{thm:separ}
If $\eta$ is tangent to a branching foliation, then it can be $C^0$-approximated by a $C^0$-foliation $\mathcal{F}$. Moreover, there exists a continuous map $h: M \rightarrow M$ such that for every leaf $L \in \mathcal{F}$, the restriction $h_{\vert L} : L \rightarrow M$ is a $C^1$ immersion tangent to $\eta$.
\end{thm}

We can further assume that $h$ is close to the identity, hence of degree $1$ and surjective, and that the leaf collection 
$$\mathscr{F} \coloneqq \big\{ h_{\vert L} : L \rightarrow M \ \big\vert \  L \in \mathcal{F} \big\}$$
is a branching foliation tangent to $\eta$. However, we will not need these additional properties.

One can think of branching foliations as the closure of the space of $C^0$-foliations, in the following sense. By the previous theorem, any (oriented) plane field tangent to a branching foliation is an accumulation point of the space of plane fields tangent to $C^0$-foliations for the $C^0$ topology. On the other hand, we have:

\begin{lem} \label{lem:branchinglim}
    Let $(\eta_n)_{n \geq 0}$ be a sequence of continuous plane fields tangent to branching foliations, which converges to a plane field $\eta$ in the $C^0$-topology. Then $\eta$ is tangent to a branching foliation.
\end{lem}

\begin{proof}
    We only sketch the main ideas since the strategy is quite standard.

    First, we consider a point $p \in M$ as well as a sequence of complete immersed surfaces $f_n : \Sigma_n \rightarrow M$, $n \geq 0$, such that $f_n$ is a leaf of a branching foliation tangent to $\eta_n$ passing through $p$. Then, by Arzelà--Ascoli and after passing to a suitable subsequence, $(f_n)_n$ converges to a complete immersion $f : \Sigma \rightarrow M$ tangent to $\eta$.

    Now, we consider a dense subset $\mathcal{D} \subset M$, and applying the previous strategy (and passing to further subsequences using a diagonal argument), we can construct a leaf collection $\mathscr{F}$ made of complete, non topologically crossing surfaces tangent to $\eta$, which moreover covers $\mathcal{D}$. We can then argue as in~\cite[Lemma 7.1]{BI08} and enlarge $\mathscr{F}$ by taking all the limits of its leaves to obtain a branching foliation tangent to $\eta$.
\end{proof}

It remains to show that $\eta$ is tangent to a branching foliation. To that end, we will essentially follow the strategy of~\cite{BI08} and we explain how to adapt their proof to our setup. Recall that the structure of interest is a pair $(X, \eta)$ where
\begin{itemize}
\item $X$ is a smooth vector field which plays the role of $E^s$ from~\cite{BI08},
\item $\eta$ is a continuous plane field containing $X$ and invariant under its flow, which plays the role of $E^{cs}$ in~\cite{BI08}.
\end{itemize}
However, there are two major differences between our vector field $X$ and the vector field $E^s$ from~\cite{BI08}: unlike $E^s$, $X$ may have closed orbits and singularities. A careful inspection of the proofs in~\cite{BI08} reveals that the presence of closed orbits (in the absence of singularities) is not a serious issue. However, the presence of singularities dramatically complicates the arguments and requires more general notions of \emph{surfaces}, \emph{patches}, \emph{prefoliations}, etc. Furthermore, the hypothesis that $(X, \eta)$ has no disks of tangency will be essential, as these disks could obstruct the construction of prefoliations, and more specifically the construction of \emph{upper enveloping surfaces}. The (forward) structure of these surfaces is also more involved in our singular setting, and we will need to be particularly careful about their behavior near $\Delta$.

%%%%%%%%%%%%%%%%%%%%%%%%%%%%%%%%%%%%%%%%%%%%%%%%%%%%%%%%%%%%%%%%%%%%%
        \section{General construction}
%%%%%%%%%%%%%%%%%%%%%%%%%%%%%%%%%%%%%%%%%%%%%%%%%%%%%%%%%%%%%%%%%%%%%

Because of the presence of saddle and quadratic singularities, we will need to consider more general surfaces in $M$ tangent to $\eta$ than the ones in~\cite{BI08}. In particular, we will consider surfaces with corners with suitable local forms near $\Delta$. We will call these surfaces \emph{tiles}. We then generalize the main technical concept from~\cite{BI08}, namely, the one of a \emph{prefoliation}. These structures are intermediate steps in the construction of a genuine branching foliation, in the sense that their leaves are possibly nonmaximal tiles. Prefoliations, and their weaker version defined below, also come with the structure of a (pre)order at every point of $M$, allowing us to distinguish which tiles are above and below others, even though they geometrically coincide. The main technical---and most difficult!---step is to extend a (weak) prefoliation by suitably enlarging its tiles and adding new ones, while also extending the order structure. See Proposition~\ref{prop:prefoil} below, whose proof will occupy Sections~\ref{sec:forenv},~\ref{sec:extorder}, and~\ref{sec:proof}. As in~\cite{BI08}, this extension is particularly nontrivial because of the global behavior of the vector field $X$ (see~\cite[Section 6.1]{BI08}, in particular Figures 1 and 2).

%%%
            \subsection{Tiles}
%%%

We now fix a neighborhood $\mathcal{N}$ of $\Delta$ as in Lemma~\ref{lem:nbd}. We say that smooth local coordinates $(x,y,z)$ around a point $p \in M$ are \textbf{adapted} to $\eta$ if $\partial_z$ is positively transverse to $\eta$ and $\eta(0) = \mathrm{span}\{\partial_x, \partial_y\}$ in these coordinates. Every point $p \in M$ admits such local coordinates. If moreover $p \in \Delta_\mathrm{sa}$, then the projections of the stable and unstable branches of $X$ at $p$ on the $(x,y)$-plane determine four local quadrants at $0$ in $\R^2$. If $p \in Q$, then the projections of the unstable branches of $X$ at $p$ on the $(x,y)$-plane determine two local half-planes at $0$ in $\R^2$, one of which corresponds to the unstable half-disk of $p$. The projection of the stable branch at $p$ divides the other half-plane into two quadrants.

\begin{defn}[See~\cite{BI08}, Definition 4.1] A \textbf{surface} in $M$ is a $C^1$ immersion $f : U \rightarrow M$, where $U$ is a connected surface with boundary and \underline{convex and concave corners}.
\end{defn}

The rest of~\cite[Definition 4.1]{BI08} is still valid for this more general notion of surface.

\begin{defn}[See~\cite{BI08}, Definition 4.2]
An \textbf{$\eta$-surface} is a surface $f: U \rightarrow M$ tangent to $\eta$, invariant under the flow of $X$, and such that $U$ is homeomorphic to a simply connected subset of $\R^2$.
\end{defn}

The invariance condition means that the flow of the pullback $\widetilde{X} \coloneqq f^*X$ of $X$ to $U$ is complete, or that the surface ``consists of whole trajectories of $X$''. Here, $\widetilde{X}$ is a continuous and uniquely integrable vector field on $U$. We call a nonsingular flow line of $\widetilde{X}$ a \textbf{$X$-line} on $U$. Since $(X, \eta)$ has no disks of tangency, $\widetilde{X}$ has no nontrivial closed orbit on $U$, and $X$-lines are homeomorphic to $\R$.

\begin{defn} \label{def:locgraph}
Let $f: U \rightarrow M$ be an $\eta$-surface, and $D_0$ be a connected open subset of $U$. If $p \in \overline{f(D_0)}$, we say that $f$ is \textbf{locally graphical near $p$ over $D_0$} if there exists a small open ball $B$ centered at $p$ with smooth coordinates $(x,y,z)$ adapted to $\eta$, a connected domain $D \subset \R^2$, an embedding $\varphi : D \rightarrow D_0 \cap f^{-1}(B)$ which is a diffeomorphism onto the connected component containing $p$, and a map $\widetilde{f} :  D \rightarrow \R$ such that $f \circ \varphi : D \rightarrow B \subset \R^3$ is the graph of $\widetilde{f}$ over $D \subset \R^2$ in the following sense:
$$\forall (x, y) \in D, \quad f \circ \varphi (x,y) = \big(x,y,\widetilde{f}(x,y)\big) \in \R^3.$$

Furthermore, we say that $f$ has a \textbf{puncture} at $p$ in $D_0$ if $D$ is the complement of $0$ in a neighborhood of $0$, and $f$ has a \textbf{slit} at $p$ in $D_0$ if $D$ is the complement of a $C^1$ embedded curve starting at $0$ in a neighborhood of $0$, after possibly shrinking $B$ in each case.
\end{defn}

Notice that if $D_0$ is a sufficiently small disk centered at $a \in U$, then $f$ is automatically locally graphical near $f(a)$ over $D_0$, since $f$ is an immersion tangent to $\eta$. The local behavior of $f$ near $a$ can be described as follows. 

\begin{enumerate}
    \item If $f(a) \in \Delta$, then either $a$ is in the interior of $U$, hence $0$ is in the interior of $D$, or $a$ is a boundary point or a corner of $U$, and the same holds for $0$ in $D$. In the latter case, we further distinguish three cases.
        \begin{enumerate}
            \item $f(a) \in \Delta_\mathrm{so}$. The intersection of $D$ with a small open neighborhood of $0$ is a sector bounded by two embedded $C^1$ curves meeting at $0$, corresponding to two flow lines of $X$ emanating from $a$.
            \item $f(a) \in \Delta_\mathrm{sa}$. Recall that the projections of the stable and unstable branches at $f(a)$ on the $(x,y)$-plane in the above coordinates determine four (local) quadrants at $0$. If $a$ is a boundary point of $U$, the intersection of $D$ with a small open neighborhood $V$ of $0$ is exactly the intersection with $V$ of the union of two consecutive quadrants. If $a$ is a convex corner of $U$, then $D \cap V$ is exactly the intersection with $V$ of one of the quadrants. If $a$ is a concave corner of $U$, then $D \cap V$ is exactly the intersection with $V$ of three consecutive quadrants.
            \item $f(a) \in Q$. Recall that the projections of the unstable branches at $f(a)$ on the $(x,y)$-plane determine two (local) half-planes at $0$, and the projection of the stable branch divides one of these half-planes into two quadrants. The other half-plane corresponds to the projection of the unstable half-disk of $f(a)$. The description of the different possibilities for $D \cap V$ when $a$ is a boundary point or a corner of $U$ is left as an exercise to the reader.
        \end{enumerate}
    \item If $f(a) \notin \Delta$, meaning that $X(f(a)) \neq 0$, we can further assume that $X$ coincides with $\partial_x$ in the above local coordinates. By the invariance of $\eta$ under $X$, $\widetilde{f} : D \rightarrow \R$ is a function of $y$ only. Note that $a$ cannot be a corner of $U$. If $a$ is in the interior of $U$, then $D$ can be chosen of the form $(-\epsilon, \epsilon)^2$, and if $a$ is a boundary point of $U$, $D$ can be chosen of the form $(-\epsilon, \epsilon) \times [0, \epsilon)$ or  $(-\epsilon, \epsilon) \times (-\epsilon, 0]$, where $\epsilon > 0$ is sufficiently small. We can further distinguish two cases:
        \begin{enumerate}
            \item The orbit of $X$ passing through $f(a)$ is not periodic. In that case, $f$ maps the flow line of $\widetilde{X}$ in $U$ passing through $a$ to the flow line of $X$ in $M$ passing through $f(a)$ bijectively.
            \item The orbit of $X$ passing through $f(a)$ is periodic. In that case, the flow line of $\widetilde{X}$ in $U$ passing through $a$ is not closed, since $U$ is simply connected and $(X, \eta)$ has no disks of tangency. The restriction of $f$ to this flow line covers the closed orbit of $X$ in $M$ passing through $f(a)$ infinitely many times, in positive and negative times.
        \end{enumerate}
\end{enumerate}
Notice that $f$ sends the corners of $U$ to $\Delta$, and $f^{-1}(\Delta)$ is a discrete subset of $U$. Moreover, $f$ sends the boundary components of $U \setminus f^{-1}(\Delta)$ to entire flow lines of $X$.

\medskip

As in~\cite{BI08}, we can consider the \textbf{completion} of an $\eta$-surface $f: U \rightarrow M$. The completion $\overline{U}$ is the completion of $U$ for the intrinsic metric on $U$, and $f$ uniquely extends to a continuous map $\bar{f} : \overline{U} \rightarrow M$. However, the completion of an $\eta$-surface in our singular setting is not necessarily an $\eta$-surface, since $\overline{U}$ might fail to be a smooth surface and $\bar{f}$ might fail to be an immersion. This occurs when $f$ is not locally graphical near a singularity $p \in \Delta$ or has a slit at $p$. Because of this issue, we need to be particularly careful about the way the $\eta$-surfaces we consider intersect $\mathcal{N}$ and behave near $\Delta$.

If $p \in \Delta$ and $D_0$ is a connected component of $f^{-1}(\mathcal{N})$ such that $p \in \overline{f(D_0)}$, we call $D_0$ a \textbf{singular region} of $f$ for $p$.

\begin{defn} \label{def:tile}
Let $f: U \rightarrow M$ be an $\eta$-surface.
\begin{itemize}
    \item $f$ is \textbf{well-cornered} if for every $p \in \Delta$ and every singular region $D_0 \subset U$ for $p$, $f$ is locally graphical near $p$ over $D_0$ and has no puncture nor slit at $p$ in $D_0$.
    \item $f$ is \textbf{unstable-complete} if for every $a \in U$,
        \begin{enumerate}
            \item If $f(a)$ lies in the unstable disk of a source singularity $p \in \Delta_\mathrm{so}$, then there exists an open subset $V \subset U$ containing the point $a$ such that the restriction $f_{\vert V} : V \rightarrow M$ is a diffeomorphism onto this unstable disk.
            \item If $f(a)$ lies in the interior of the unstable half-disk of a quadratic singularity $q \in Q$, then there exists an open subset $V \subset U$ containing the point $a$ such that the restriction $f_{\vert V} : V \rightarrow M$ is a diffeomorphism onto the interior of this unstable half-disk.
        \end{enumerate}
    \item $f$ is an \textbf{$\eta$-tile}, or $\textbf{tile}$ for simplicity, if it is well-cornered and unstable-complete.
\end{itemize}
\end{defn}

More concretely, a tile $f$ has the following behavior near a singularity $p \in \Delta$:
\begin{itemize}
    \item If $p \in \Delta_\mathrm{so}$, then the domain $D$ from Definition~\ref{def:locgraph} contains a neighborhood of $0$.
    \item If $p \in \Delta_\mathrm{sa}$ and if $0$ is not in the interior of $D$, then after possibly shrinking $B$, the interior of $D$ is either the interior of a quadrant, or the interior of the union of two or three adjacent quadrants.
    \item If $p \in Q$ and if $0$ is not in the interior of $D$, then after possibly shrinking $B$, the interior of $D$ is the interior of a quadrant, the interior of the union of the two quadrants, the interior of the (projection of the) unstable half-plane, or the interior of the union of this half-plane with one quadrant.
\end{itemize}
We only have described the possibilities for the \emph{interior} of $D$ (intersected with a neighborhood of $0$). The precise description of $D$ is left as an exercise to the reader.

\begin{lem}[See~\cite{BI08}, Lemma 4.3] The completion $\bar{f} : \overline{U} \rightarrow M$ of a tile $f : U \rightarrow M$ is a tile. Moreover, $\overline{U} \setminus U$ is contained in $\partial \overline{U}$.
\end{lem}

\begin{proof}
If $a \in \overline{U} \setminus U$ is such that $\bar{f}(a) \notin \Delta$, then the proof of~\cite[Lemma 4.3]{BI08} shows that $a$ is in the set-theoretic boundary of $\overline{U}$ and near $a$, $\overline{U}$ is obtained from $U$ by adding an $X$-line. Otherwise, if $\bar{f}(a) \in \Delta$, then the effect of the completion near $a$ corresponds to taking the closure of the domain $D$ from Definition~\ref{def:locgraph}. Then $\overline{U}$ differs from $U$ near $a$ by adding the point $a$ together with boundary components which are $X$-lines. It easily follows that $\overline{U}$ is a surface with boundary and corners, and that $\bar{f}$ is a tile. The details of the proof are left to the reader.
\end{proof}

We can now extend the definitions in~\cite[Section 4]{BI08} to tiles as follows. An \textbf{edge} of a tile $f : U \rightarrow M$ is a connected component of $\partial \overline{U} \setminus \bar{f}^{-1}(\Delta)$. A point in $\bar{f}^{-1}(\Delta) \cap \partial \overline{U}$ is a \textbf{boundary singularity}. It is either a corner of $\overline{U}$, or a point between two edges which are $X$-lines with opposite orientations. The latter will be called a \textbf{fake corner}. Notice $\bar{f}$ sends fake corners to points in $Q \cup \Delta_\mathrm{sa}$. An edge or boundary singularity (corner or fake corner) is \textbf{proper} if it is contained in $U$, and \textbf{improper} otherwise. Two edges are \textbf{adjacent} if they meet at a corner or at a fake corner. As in~\cite{BI08}, an edge has a forward or backward coorientation. This notion can be extended to (genuine) corners, but not to fake corners: two edges adjacent to a fake corner have opposite coorientations. We warn the reader that we will use the terminology \emph{forward/backward corner} in a different way below.

The definition~\cite[Definition 4.4]{BI08} also extends in a straightforward way: a \textbf{marked tile} is a pair $(f, a)$ where $f : U \rightarrow M$ is a tile and $a \in U$. It is a \textbf{forward (resp.~backward)} marked tile if the marking $a$ belongs to a forward (resp.~backward) edge or corner, and \textbf{passing} marked tile if the marking $a$ belongs to the interior of $U$. Unlike in~\cite{BI08}, there is no natural notion of ``forward or backward half-tile'' at a point $a \in U$. For instance, $a$ could belong to an $X$-line connecting two singularities in the interior of $U$, and this $X$-line would not divide $U$ into two half-surfaces.

%%%
            \subsection{Weak prefoliations}
%%%

As in~\cite[Definition 4.5]{BI08}, if $\mathscr{A}$ is a collection of tiles and $p \in M$, we set 
$$\mathscr{A}_p \coloneqq \{ (f,a) \ \vert \, f \in \mathscr{A}, \ f(a) = p\}$$
and
$$\mathscr{A}_\star \coloneqq \bigsqcup_{p \in M} \mathscr{A}_p.$$
If $f \in \mathscr{A}$, we will typically denote a corresponding element in $\mathscr{A}_\star$ by $f_\star$. 

The definition of the \emph{geometric order} in~\cite[Definition 4.6]{BI08} extends in a straightforward way to a collection $\mathscr{A}$ of tiles. The notion of \emph{patch} from Definition 4.7 of~\cite{BI08} is harder to adapt and is not crucial for our argument.

We now introduce a natural equivalence relation on marked tiles.

\begin{defn}
Two marked tiles $f_\star = (f : U \rightarrow M, a)$ and $g_\star = (g: V \rightarrow M, b)$ are \textbf{equivalent} if there exists a pointed diffeomorphism $\varphi : (U,a) \rightarrow (V,b)$ such that $g \circ \varphi = f$. In that case, we write $f_\star \approx g_\star$. 
\end{defn}

It is straightforward to check that $\approx$ is an equivalence relation on marked tiles. Notice that if $(f, a) \approx (g,b)$, then $f(a) = g(b)$, and the diffeomorphism $\varphi$ as above is unique, since $f$ is an immersion and $U$ is connected. 

\begin{rem}
It might happen that a single tile $f : U \rightarrow M$ admits two distinct markings $a, b \in U$ such that $(f, a) \approx (f,b)$. In that case, it is easy to check that the corresponding diffeomorphism $\varphi$ has no fixed points and acts properly discontinuously on $U$. Therefore, the quotient $\Sigma \coloneqq U \slash \varphi$ has a natural manifold structure and $f$ induces an immersion $\Sigma \rightarrow M$ which is tangent to $\eta$. In other words, $f$ is a nontrivial cover of a surface tangent to $\eta$, but this surface is not a tile since $\Sigma$ is not simply connected. Furthermore, if $U$ only has one proper edge, then $\Sigma$ has a closed boundary component which is mapped to a closed orbit of $X$ in $M$.
\end{rem}

\begin{defn}[See~\cite{BI08}, Definition 4.8] \label{def:weakpref}
A \textbf{weak prefoliation} $(\mathscr{A}, \precsim)$ is a collection $\mathscr{A}$ of tiles together with a (partial) preorder $\precsim$ on $\mathscr{A}_\star$ satisfying the following axioms.\footnote{A preorder on a set is a transitive and reflexive binary relation. It is an order if and only if it is antisymmetric. Dropping the antisymmetry requirement makes weak prefoliations more flexible than the prefoliations from~\cite{BI08}. See also the next footnote.}
\begin{enumerate}
    \item[0.] Two marked tiles $f_\star \in \mathscr{A}_p$ and $g_\star \in \mathscr{A}_q$ are comparable by $\precsim$ if and only if $p=q$. Moreover, if $f_\star \precsim g_\star$ and $g_\star \precsim f_\star$ then $f_\star \approx g_\star$.
    \item[1.] The preorder $\precsim$ is a refinement of the geometric order in the following sense: if $f_\star \precsim g_\star$, then $g_\star$ is geometrically locally above $f_\star$.
    \item[2.] The preorder $\precsim$ is coherent along the intersection of tiles in the following sense: if $\gamma : [0,1] \rightarrow M$ is a continuous path, and if $\gamma_1 : [0,1] \rightarrow U_1$ and $\gamma_2 : [0,1] \rightarrow U_2$ are lifts of $\gamma$ to tiles $f_1 : U_1 \rightarrow M$ and $f_2 : U_2 \rightarrow M$, respectively, then $(f_1, \gamma_1(0)) \precsim (f_2, \gamma_2(0))$ if and only if $(f_1, \gamma_1(1)) \precsim (f_2, \gamma_2(1))$.
\end{enumerate}
\end{defn}

In Axiom \textit{0}, we \emph{do not} assume that $f_\star \precsim g_\star$ if $f_\star \approx g_\star$. In other words, two equivalent tiles might be strictly comparable by $\precsim$, while two nonequivalent tiles are necessarily strictly comparable by $\precsim$. We write $f_\star \simeq g_\star$ if $f_\star \precsim g_\star$ and $g_\star \precsim f_\star$, and $f_\star \prec g_\star$ if $f_\star \precsim g_\star$ and $f_\star \not\simeq g_\star$. If $f_\star \simeq g_\star$, we say that $f_\star$ and $g_\star$ are \textbf{order-equivalent}.\footnote{Obviously, $\precsim$ induces a genuine order on the quotient $\mathscr{A}_\star \slash \simeq$. Alternatively, we could circumvent the lack of antisymmetry of $\precsim$ by considering a more general definition of tiles, allowing their domains to be arbitrary surfaces different than the closed disk and the sphere.}

Axioms \textit{0} and \textit{1} readily imply that if $g_\star \in \mathscr{A}_p$ is strictly locally above $f_\star \in \mathscr{A}_p$, then $f_\star \prec g_\star$. Combined with Axiom \textit{2}, this implies that there are no topological crossings between tiles in $\mathscr{A}$; see the discussion after~\cite[Definition 4.8]{BI08}.

Axiom \textit{2} can be modified by further requiring that $\gamma$ is $C^1$, tangent to $\eta$, and either tangent or transverse to $X$. A slight adaptation of the discussion after~\cite[Definition 4.8]{BI08} shows that this seemingly weaker version of Axiom \textit{2} is in fact equivalent to it. Indeed, if $\gamma : [0,1] \rightarrow M$ is only continuous, then $\gamma_1$ can be approximated by a piecewise $C^1$ curve $\widetilde{\gamma}_1$ in $U$ with the same endpoints as $\gamma_1$, which intersects the same $X$-lines and singularities as $\gamma_1$, and whose $C^1$ segments are either tangent or transverse to $\widetilde{X}$.\footnote{To see that, one can cover $U$ with open disks which are either flow boxes for $\widetilde{X}$ or standard neighborhoods of the isolated singularities of $\widetilde{X}$.} Therefore, $\widetilde{\gamma} \coloneqq f \circ \widetilde{\gamma}_1$ admits a lift $\widetilde{\gamma}_2 : [0,1] \rightarrow U_2$ to $f_2$ with the same endpoints as $\gamma_2$.

\medskip

A weak prefoliation is \textbf{forward (resp.~backward)} if its tiles only have edges with forward (resp.~backward) coorientations, and it is \textbf{complete} if its tiles are all complete surfaces. It \textbf{covers} $M$ if for every $p \in M$, $\mathscr{A}_p \neq \varnothing$.

\begin{ex} \label{ex:prefoil}
A branching foliation $\mathscr{F}$ has a natural (weak) prefoliation structure, constructed in the proof of~\cite[Theorem 7.2]{BI08}. We give a slightly different construction. First of all, the leaves in $\mathscr{F}$ are obviously tiles. Now, if $f_\star = (f : U \rightarrow M, a)$ and $g_\star = (g : V \rightarrow M, b)$ are two marked leaves in $\mathscr{F}_p$, $p \in M$, we define the relation $\precsim^{\mathscr{F}}_p$ as follows: $f_\star \precsim^\mathscr{F}_p g_\star$ if for every (continuous, or $C^1$) path $\gamma : [0,1] \rightarrow U$ starting at $\gamma(0) = a$ which lifts to a path $\widetilde{\gamma}$ to $g_\star$, $(g, \widetilde{\gamma}(1))$ is locally above $(f, \gamma(1))$. Since $f_\star$ and $g_\star$ have no topological crossings, they are comparable by $\precsim_p$. Moreover, it is not hard to show that $f_\star \precsim^\mathscr{F}_p g_\star $ and $g_\star \precsim^\mathscr{F}_p f_\star$ if and only if $f_\star \approx g_\star$, so axiom $0$ of Definition~\ref{def:weakpref} is satisfied. Axioms $1$ and $2$ immediately follow from the definition of $\precsim^\mathscr{F}$.
\end{ex}

\begin{defn}[See~\cite{BI08}, Definition 4.9] \label{def:ext}
Let $\mathscr{A}$ be a weak prefoliation. An \textbf{extension} of $\mathscr{A}$ is the data of
\begin{itemize}
    \item A weak prefoliation $\mathscr{B}$,
    \item A map $i : \mathscr{A} \rightarrow \mathscr{B}$, called an \textbf{extension map},
    \item For every $f \in \mathscr{A}$ with $f: U \rightarrow M$ and $i(f) : U' \rightarrow M$, an embedding $\varphi_f : U \hookrightarrow U'$ satisfying $i(f) \circ \varphi_f = f$,
\end{itemize}
such that $i$ preserves the preorders. More precisely, the map $i$ induces a map $i_\star : \mathscr{A}_\star \rightarrow \mathscr{B}_\star$ defined by $$i_\star(f_\star) \coloneqq \big(i(f), \varphi_f(a)\big),$$
where $f_\star = (f,a) \in \mathscr{A}$. We require that the following holds:
$$\forall f_\star, g_\star \in \mathscr{A}_\star, \quad f_\star \precsim^\mathscr{A} g_\star \iff i_\star(f_\star) \precsim^\mathscr{B} i_\star(g_\star).$$
\end{defn}

Unfortunately, Proposition 4.11, Proposition 4.12 and Proposition 4.13 from~\cite{BI08} do not individually extend in a straightforward way to our setup. Instead, we will combine modified versions of these statements and show:

\begin{prop} \label{prop:prefoil}
Let $\mathscr{A}$ be a complete and backward (resp.~forward) weak prefoliation. There exists an extension $\mathscr{A} \rightarrow \mathscr{B}$ where $\mathscr{B}$ is complete, forward (resp.~backward), and covers $M$.
\end{prop}

The proof of this proposition occupies the next four sections. We now show how it implies the Technical Theorem~\ref{thm:tech}.

\begin{proof}[Proof of the Technical Theorem~\ref{thm:tech}]
As in~\cite{BI08}, there is $r_0>0$ such that every ball of radius $r_0$ centered at a point in $M \setminus \mathcal{N}$ admits smooth coordinates $(x,y,z)$ in which $X = \partial_x$, and the angle between $\eta$ and $\mathrm{span} \{ \partial_x, \partial_y\}$ is at most $10^{-10}$. We can rescale the metric $\boldsymbol{g}$ so that $r_0 \geq 100$.

Applying Proposition~\ref{prop:prefoil} inductively, we construct a sequence of complete weak prefoliations $(\mathscr{A}^n)_{n \geq 0}$ such that 
\begin{itemize}
    \item $\mathscr{A}^0 = \varnothing$,
    \item $\mathscr{A}^1$ is forward and covers $M$,
    \item For every $n \geq 0$, $\mathscr{A}^{n+1}$ is an extension of $\mathscr{A}^n$, with an extension map $i^n : \mathscr{A}^n \rightarrow \mathscr{A}^{n+1}$,
    \item For every $k \geq 0$, $\mathscr{A}^{2k}$ is backward and $\mathscr{A}^{2k+1}$ is forward.
\end{itemize}
For every $f \in \mathscr{A}^1$, we consider the sequence of tiles $(f^n)_{n \geq 1}$ inductively defined by
$$f^1 \coloneqq f, \qquad \forall n \geq 1, \ f^{n+1} = i^n(f^n).$$
If $f^n : U_n \rightarrow M$, then by definition there exists a sequence of embeddings 
$$U_1 \overset{\varphi^1}{\longrightarrow} U_2 \overset{\varphi^2}{\longrightarrow} \dots \overset{\varphi^{n-1}}{\longrightarrow}  U_n \overset{\varphi^n}{\longrightarrow} U_{n+1} \overset{\varphi^{n+1}}{\longrightarrow} \dots$$
such that for every $n \geq 1$, $f^{n+1} \circ \varphi^n = f^n$. Taking the colimit of the above diagram yields a surface $U_\infty$, a sequence of embeddings $(\psi^n : U_n \rightarrow U_\infty)_{n\geq 1}$, and an immersion $f^\infty : U_\infty \rightarrow M$ such that for every $n \geq 1$, $f^\infty \circ \psi^n = f^n$ and $\psi^{n+1} \circ \varphi^n = \psi^n$. In particular, $f^\infty$ is tangent to $\eta$. If $a_\infty \in U_\infty$, then for $n \geq 1$ sufficiently large, there exists $a_n \in U_n$ such that $\psi^n(a_n) = a_\infty$. We say that the marked tile $(f^n, a_n) \in \mathscr{A}^n_\star$ is an \emph{ancestor} of $(f^\infty, a_\infty)$. For $1 \leq m \leq n$, we write 
\begin{align*}
    \varphi^m_n &\coloneqq \varphi^{n-1} \circ \dots \circ \varphi^m, \\
    U^m_n &\coloneqq \varphi^m_n(U_m) \subset U_n, \\ 
    U^m_\infty &\coloneqq \psi^m(U_m) \subset U_\infty, \\
    a^m_n &\coloneqq \varphi^m_n(a_m) \in U_n.
\end{align*}

We claim that the surface $f^\infty$ is open and complete. Indeed, let $a_\infty \in U_\infty$ and $(f^n, a_n)$ be an ancestor of $(f^\infty, a_\infty)$. If $f^n(a_n) \notin \Delta$, then $a_{n+1} = \varphi^n(a_n)$ is in the interior of $U_{n+1}$, since either $a_n$ is in the interior of $U_n$, or it is in an edge of $U_n$ and $a_{n+1}$ is not in an edge of $U_{n+1}$ because these surfaces have edges with opposite coorientations. If $f^n(a_n) \in \Delta$ is a source, then it is already in the interior of $U_n$, and if it is a saddle or a quadratic singularity, then $a_{n+3} = a^n_{n+3}$ is in the interior of $U_{n+3}$ by a similar argument as above. In both cases, $a_\infty = \psi^{n+1}(a_{n+1}) = \psi^{n+3}(a_{n+3})$ is in the interior of $U_\infty$. To show completeness, it is enough to show that every closed ball for the intrinsic metric $f^* \boldsymbol{g}$ on $U_\infty$ is compact, by the Hopf--Rinow--Cohn-Vossen Theorem~\cite[Theorem 2.5.28]{BBI01}. The key observation is that if $a_\infty \in U_\infty$ and $(f^n, a_n)$ is an ancestor of $(f^\infty, a_\infty)$, then for every $k \geq 1$, the point $a_{n+3k} = a^n_{n+3k} \in U_{n+3k}$ is at distance at least $k$ from the boundary of $U_{n+3k}$. For $k=1$, we can argue as follows. If $\gamma : [0,1] \rightarrow \overline{U}_{n+3}$ is a piecewise $C^1$ (or rectifiable) curve from $a_{n+3}$ to an edge of $U_{n+3}$, and if $t_0 \in [0,1]$ is the last time such that $\gamma(t_0) \in \partial U^n_{n+3}$, then there exists $t_1 \in [t_0,1]$ such that $f^{n+3}(\gamma(t_1)) \notin \mathcal{N}$.\footnote{For this to be true, we crucially need to apply the forward or backward extensions thrice, because of the possible presence of connections between two saddles.} The endpoints of the restriction of $\gamma$ to $\overline{U_{n+3}} \setminus U^n_{n+3}$ belong to (the closures of) $X$-lines with the \emph{same} coorientations, and the normalization of $r_0$ implies that the length of either $\gamma_{\vert [t_0, t_1]}$ or $\gamma_{\vert [t_1, 1]}$ is at least $1$, so the length of $\gamma$ is at least $1$. The general case $k \geq 1$ can be obtained by induction. Therefore, if $B = B(a_\infty, r) \subset U_\infty$ is a closed intrinsic ball of radius $r > 0$ centered at $a_\infty$, then $B \subset U^N_\infty$ for a sufficiently large $N \geq 1$, and $B_N = \big(\psi^N\big)^{-1}(B) \subset U_N$ coincides with a closed intrinsic ball of radius $r$ in $U_N$. By~\cite[Theorem 2.5.28]{BBI01}, $B_N$ is compact since $U_N$ is complete, so $B = \psi^N(B_N)$ is compact as well.

We now define
$$\mathscr{A}^\infty \coloneqq \big\{f^\infty \ \vert \ f \in \mathscr{A}^1 \big\}.$$
For every $p \in M$, $\mathscr{A}^\infty_p \neq \varnothing$ since $\mathscr{A}^1$ covers $M$. If $\precsim^n$ denotes the preorder on $\mathscr{A}^n_\star$, we define a preorder $\precsim^\infty$ on $\mathscr{A}^\infty_\star$ as follows: if $f^\infty_\star, g^\infty_\star \in \mathscr{A}^\infty_\star$, we set $f^\infty_\star \precsim^\infty g^\infty_\star$ if and only if there exists $n\geq 1$ and ancestors $f^n_\star,  g^n_\star \in \mathscr{A}^n_\star$ of $f^\infty_\star$ and $g^\infty_\star$, respectively, satisfying $f^n_\star \precsim^n g^n_\star$. It is straightforward to check that $\precsim^\infty$ is well-defined and that $(\mathscr{A}^\infty, \precsim^\infty)$ is a weak prefoliation. In particular, two surfaces $f^\infty, g^\infty \in \mathscr{A}^\infty$ have no topological crossings. Therefore, we have shown that $\mathscr{A}^\infty$ is a branching foliation tangent to $\eta$, and Theorem~\ref{thm:separ} finishes the proof.
\end{proof}

\begin{rem}
    Let us define a \textbf{partial branching foliation} in the same way as a branching foliation, without the condition that the leaves cover $M$. Then, Example~\ref{ex:prefoil} readily shows that a partial branching foliation induces a complete and backward/forward weak prefoliation. Proposition~\ref{prop:prefoil} and the previous proof immediately imply the following: for any partial branching foliation $\mathscr{F}^\circ$, there exists a branching foliation $\mathscr{F}$ with $\mathscr{F}^\circ \subset \mathscr{F}$.
\end{rem}
%%%%%%%%%%%%%%%%%%%%%%%%%%%%%%%%%%%%%%%%%%%%%%%%%%%%%%%%%%%%%%%%%%%%%
        \section{Upper enveloping tiles} \label{sec:forenv}
%%%%%%%%%%%%%%%%%%%%%%%%%%%%%%%%%%%%%%%%%%%%%%%%%%%%%%%%%%%%%%%%%%%%%

In this section, we extend the notions of \emph{upper enveloping surface} and \emph{forward envelope} from~\cite[Section 6.1]{BI08} to our setup. We will assume that the reader is familiar with the notion of a one-dimensional prefoliation on the plane from~\cite[Section 5]{BI08}. Our proof of the existence of forward envelopes essentially follows the strategy of~\cite{BI08}, but we find it convenient to phrase it in a slightly more abstract way.

%%%
        \subsection{Preliminaries} \label{sec:prelim}
%%%

Let $(\mathscr{A}, \precsim)$ be a (possibly empty) weak prefoliation which is backward and complete. Our goal is to extend the tiles in $\mathscr{A}$ to tiles with only forward edges. The first step is to define new neighboring tiles, the \emph{forward envelopes}, which play the role of ``upper neighbors'' of the existing tiles (or are new tiles if $\mathscr{A}$ was empty). In the present section, we extend the definition of forward envelopes to our setup. In the next section, we extend the preorder to these tiles and to their completions. In the section after, we glue the tiles in $\mathscr{A}$ to their newly created neighbors and we further extend the preorder one last time. Before getting to the heart of the matter, let us make some remarks and introduce some terminology.

\bigskip

The flow of $X$ induces a natural $\R$-action on $\mathscr{A}_\star$ by ``sliding the marked point'' along an $X$-line. More precisely, for every $f_\star = (f:U \rightarrow M, a) \in \mathscr{A}_\star$ and $t \in \R$, we define the \textbf{$t$-translation} of $f_\star$ by
$$f^t_\star := \big(f, \widetilde{\phi}^t(a)\big),$$
where $\widetilde{\phi}^t$ denotes the flow of $\widetilde{X} = f^* X$. For every $p \in M$, we have
$$f_\star \in \mathscr{A}_p \iff f^t_\star \in \mathscr{A}_{\phi^t(p)},$$
where $\phi^t$ denotes the flow of $X$. Therefore, this operation extends to a bijection
$$\Phi^t : \mathscr{A}_p \rightarrow \mathscr{A}_{\phi^t(p)}$$
which obviously defines an $\R$-action on $\mathscr{A}_\star$.

\bigskip 

If $p \in M \setminus \Delta$ and $D \subset M$ is an open $2$-disk transverse to $X$ centered at $p$, then $\mathscr{A}$ induces a \emph{one-dimensional prefoliation} (in the sense of~\cite[Section 5.2]{BI08}) on $D$ for the (continuous) oriented line field obtained as the trace of $\eta$ on $D$. More precisely, if $X_\eta$ denotes a unit vector field directing this line field, then the intersections of tiles in $\mathscr{A}$ with $D$ form a collection $\mathcal{A}$ of flow lines of $X_\eta$, and for every $q \in D$, the preorder $\precsim$ induces a total order on the set of flow lines of $X_\eta$ passing through $q$ as in~\cite[Section 5.2]{BI08} (see also the beginning of~\cite[Section 5.1]{BI08}). We call the data of $\mathcal{A}$ equipped with this order the \textbf{trace} of $\mathscr{A}$ on $D$.

\bigskip

We say that a singularity $p\in \Delta_\mathrm{sa} \cup Q$ is \textbf{surrounded by unstable manifolds} if either
\begin{itemize}
    \item Every stable branch at $p$ intersects the unstable manifold of a source or quadratic singularity, or
    \item $p \in \Delta_\mathrm{sa}$, a stable or unstable branch of $p$ is connected to $q \in \Delta_\mathrm{sa}$, and the three stable branches of $p$ and $q$ other than the connection intersect unstable manifolds of source or quadratic singularities.
\end{itemize}

We define the \textbf{unstable locus} $M^\mathrm{un} \subset M$ as the union of interiors of unstable manifolds of source and quadratic singularities, together with the singularities in $\Delta_\mathrm{sa} \cup Q$ surrounded by unstable manifolds and their unstable branches.

We define the \textbf{regular locus} of $M$ as $M^\mathrm{reg} \coloneqq M \setminus \Delta$. We further partition the singularity locus $\Delta$ into two sets $$\Delta^\mathrm{good} \coloneqq \Delta \cap M^\mathrm{un}, \qquad \Delta^\mathrm{bad} \coloneqq \Delta \setminus \Delta^\mathrm{good},$$ and we define the \textbf{good locus} by $M^\mathrm{good} \coloneqq M \setminus \Delta^\mathrm{bad} = M^\mathrm{reg} \cup \Delta^\mathrm{good}$. Roughly speaking, this  is the set of points at which it will make sense to define upper enveloping tiles.

Finally, we define the \textbf{pure locus} of $M$ as $M^\mathrm{pure} \coloneqq M^\mathrm{reg} \setminus M^\mathrm{un}$. It is foliated by the flow lines of $X$ disjoint from the unstable locus. We say that such a flow line is \textbf{pure}. Notice that we have $$M^\mathrm{good} = M^\mathrm{un} \sqcup M^\mathrm{pure}.$$

Let $p \in \Delta_\mathrm{sa} \cup Q$. If $p \in \Delta_\mathrm{sa}$ and no stable branch of $X$ at $p$ is in $M^\mathrm{un}$, we say that $p$ is a \textbf{pure} saddle singularity. If exactly one stable branch of $X$ at $p$ is contained in $M^\mathrm{pure}$, we say that $p$ is a \textbf{quadratic-like} singularity. We denote by $\Delta^\mathrm{pure}_\mathrm{sa}$ the set of pure saddle singularities of $X$, and by $Q^*$ the set of quadratic-like singularities of $X$, which also includes quadratic singularities.

\bigskip

We now describe the behavior of $\eta$ on $M^\mathrm{un}$.

\begin{lem} \label{lem:uniqint}
The plane field $\eta$ is uniquely integrable on $M^\mathrm{un}$.
\end{lem}

\begin{proof}
We need to show that $\eta$ is uniquely integrable at the singularities surrounded by unstable manifolds and along their unstable branches. This will immediately follow from the two slightly stronger claims:
\begin{enumerate}[leftmargin=*]
\item \textit{If there exists a connection $\gamma$ from $p \in \Delta_{\mathrm{so}} \cup Q$ to $q \in \Delta_{\mathrm{sa}} \cup Q$, then $\eta$ is uniquely integrable at $q$ and along the unstable branches at $q$ \underline{on the side of $\gamma$}.}
\item \textit{If moreover there exists a connection $\gamma'$ from $q$ to $q' \in \Delta_\mathrm{sa}$, then $\eta$ is uniquely integrable at $q'$ and along the unstable branch at $q'$ \underline{on the quadrant at $q'$ determined by $\gamma$ and $\gamma'$}.}
\end{enumerate}
The meaning of the underlined parts will be clarified in the proof. Let $S^u$ denote the interior of the unstable manifold of $X$ at $p$.

For the first claim, let $S$ be a small open $C^1$ disk tangent to $\eta$ and passing through a point $q_0$ lying on one of the unstable branches $\gamma^u$ at $q$, and let $S_\gamma$ be the closure of the connected component of $S \setminus \gamma_u$ on the side of $\gamma$. This side can be characterized as follows: if $D$ is a small smooth disk in $M$ transverse to $X$ centered at a point $p_0 \in \gamma$ close to $q$, then the points in $S_\gamma$ sufficiently close to $\gamma^u$ all get sent to $D$ or to $q$ under the negative flow of $X$, and the corresponding flow lines stay close to $\gamma \cup \{q\} \cup \gamma^u$. Let $\widehat{S}_\gamma$ be the union of the saturation of $S_\gamma$ by the flow of $X$ with $\gamma$, so that $\widehat{S}_\gamma$ is a $C^1$ surface with boundary and corner tangent to $\eta$. Since $\widehat{S}_\gamma$ intersects $S^u$ along $\gamma$, and $\eta$ is uniquely integrable there, the interior of $\widehat{S}_\gamma$ is entirely contained in $S^u$. Therefore, $S_\gamma$ is contained in $S^u \cup \gamma^u$. If $q_0 = q$, then $S_\gamma$ is defined in a similar way, and its interior contains points in $\gamma$, which belong to the unstable manifold of $p$. This implies that $S_\gamma$ is contained in the union of $S^u$ with $q$ and the two unstable branches at $q$. The desired unique integrability property of $\eta$ follows.

For the second claim, $\gamma^u$ will denote the unstable branch at $q'$ lying on the side determined by $\gamma$ and $\gamma'$. Similarly, let $S$ be a small open $C^1$ disk tangent to $\eta$ and passing through a point $q'_0 \in \gamma^u$, and let $S_{\gamma, \gamma'}$ be the closure of the connected component of $S \setminus \gamma^u$ in the side determined by $\gamma$ and $\gamma'$. This side can be characterized as follows: if $D$ is as before, then the points in $S_{\gamma, \gamma'}$ sufficiently close to $\gamma^u$ all get sent to $D$ or to $q'$ under the negative flow of $X$, and these flow lines stay close to $\gamma \cup \{q\} \cup \gamma' \cup \{q'\} \cup \gamma^u$. Let $\widehat{S}_{\gamma, \gamma'}$ be the union of the saturation of $S_{\gamma, \gamma'}$ by the flow of $X$ with $\gamma'$, so that $\widehat{S}_{\gamma, \gamma'}$ is a $C^1$ surface with boundary and corners tangent to $\eta$. Since $\widehat{S}_{\gamma, \gamma'}$ intersects $\gamma'$ and lies on the side of $\gamma'$ determined by $\gamma$, the previous point implies that $\widehat{S}_{\gamma, \gamma'}$ is entirely contained in $S^u \cup \gamma' \cup \gamma^u$, hence $S_{\gamma, \gamma'}$ is contained in $S^u \cup \gamma^u$. Finally, if $q'_0 =q'$, we denote by $S_{\gamma, \gamma'}$ the closure of the connected component of $S \setminus \big(\gamma' \cup \gamma^u\big)$ lying in the quadrant at $q'$ determined by $\gamma$ and $\gamma'$. Arguing as before, we deduce that $S_{\gamma, \gamma'}$ is entirely contained in $S^u \cup \gamma' \cup \{q'\} \cup \gamma^u$. The desired unique integrability property of $\eta$ follows.
\end{proof}

Therefore, $\eta$ induces a foliation $\mathcal{F}^\mathrm{un}$ on $M^\mathrm{un}$. Notice that $M^\mathrm{un}$ is invariant under the flow of $X$, and contains no nontrivial closed orbits of $X$ (independently of the hypothesis that $(X, \eta)$ has no disk of tangency).

\begin{lem} \label{lem:unstabtile}
Let $f : L \rightarrow M$ be a (complete) leaf of $\eta$ in $M^\mathrm{un}$. Its universal cover $\widetilde{f} : \widetilde{L} \rightarrow M$ is a tile. We call it an \textbf{unstable tile}.
\end{lem}

\begin{proof}
By construction, $f$ and $\widetilde{f}$ are made of  unstable manifolds of source singularities and interiors of unstable manifolds of quadratic singularities, glued together along unstable connections of saddle or quadratic singularities surrounded by them. One easily deduces, using the Poincar\'{e}--Bendixson and the Poincar\'{e}--Hopf theorems, that $L$ is not a sphere, hence $\widetilde{L}$ is homeomorphic to a simply connected subset of $\R^2$. Moreover, $\widetilde{f}$ is unstable-complete by definition. We are left to show that $\widetilde{f}$ is well-cornered at every singularity in $\Delta_\mathrm{sa} \cup Q$. It is clearly the case at a singularity surrounded by unstable manifolds, by construction. Let $p \in \Delta_\mathrm{sa}$ be surrounded by unstable manifolds, and $\gamma$ be a connection from $p$ to $q \in \Delta_\mathrm{sa}$ which is \emph{not} surrounded by unstable manifolds. Then, the unstable branches at $q$ both belong to $M^\mathrm{pure}$ and every singular region $D_0 \subset \widetilde{L}$ for $q$ contains an open half-disk $D'_0$ mapping to a half-disk $D$ in $M$ containing a portion of $\gamma$, and whose `straight' part of the boundary is contained in the union of $q$ and its unstable branches. This shows that $\widetilde{f}$ is well-cornered at $q$.

Since $X$ only has admissible connections, it is now enough to show the following claim: if $g : U \rightarrow M$ is a $C^1$ parametrization of the (interior of) the unstable manifold of singularity $p \in \Delta_\mathrm{so} \cup Q$, then $g$ is well-cornered. We will only treat the case of a source singularity, the case of a quadratic singularity being similar. Let $q \in \Delta_\mathrm{sa} \cup Q$, and let $D_0 \subset U$ be a singular region of $g$ for $q$. Since all of the $X$-lines in $g$ converge to $g^{-1}(p)$ in negative time, $D_0$ cannot intersect $X$-lines corresponding to the unstable branches at $q$. Therefore, if $D_0$ contains an $X$-line corresponding to a connection from $p$ to $q$, then $D_0$ contains a small open half-disk $D'_0$ similar as before. Otherwise, $D_0$ contains a small open quadrant $C_0$ mapped to an open quadrant $C$ in $M$ whose `corner' boundary is sent to the union of $q$ with two consecutive unstable and stable branches. In both cases, $g$ is well-cornered at $q$ over $D_0$.
\end{proof}

%%%
        \subsection{Forward curves}
%%%

To construct and study forward envelopes, it will be crucial to understand the `forward' structure of the tiles under consideration with respect to the pullback of the vector field $X$. The forward structure of surfaces in~\cite{BI08} is much easier to understand due to the absence of singularities of $X$. In the context of this article, we will need to carefully subdivide tiles into various regions, and extend the definition of \emph{forward curves}. A crucial step is to understand how two forward curves with the same endpoints are related to each other, which is achieved in Lemma~\ref{lem:forcurv2} below. We also collect some technical results that will be useful later.

\begin{defn} \label{def:loci}
Let $f : U \rightarrow M$ be a tile. 
\begin{itemize}
    \item The \textbf{regular locus} of $f$ is the open subset $U^\mathrm{reg} \coloneqq  f^{-1}(M^\mathrm{reg}) \subset U$. It is foliated by $X$-lines.
    \item The \textbf{good locus} of $f$ is the open subset $U^\mathrm{good} \coloneqq f^{-1}(M^\mathrm{good}) \subset U$.
    \item The 
    \textbf{unstable locus} of $f$ is the open subset $U^\mathrm{un} \coloneqq f^{-1}(M^\mathrm{un}) \subset U$. A connected component of $U^\mathrm{un}$ is called an \textbf{unstable component} of $f$.
    \item The \textbf{pure locus} of $f$ is the subset $U^\mathrm{pure} \coloneqq f^{-1}(M^\mathrm{pure}) \subset U$. An $X$-line in $U^\mathrm{pure}$ is called \textbf{pure}.
\end{itemize}
\end{defn}

We now fix a marked tile $f_\star = (f : U \rightarrow M, a)$ for the rest of this section.

\begin{defn}
A \textbf{forward (resp.~backward) curve} in $U$ is a $C^1$ immersed curve $\gamma : [0,1] \rightarrow U$ such that for every $t \in [0,1]$, either
    \begin{enumerate}
        \item $\gamma(t) \in U^\mathrm{un}$, or
        \item $\gamma(t) \in U^\mathrm{reg}$ and $\dot{\gamma}(t)$ is positively (resp.~negatively) transverse to the $X$-line passing through $\gamma(t)$.\footnote{Positive transversality means that $\big(X(\gamma(t)), \dot{\gamma}(t)\big)$ forms a positively oriented basis of $T_{\gamma(t)} U$.}
    \end{enumerate}

Furthermore, we say that $\gamma$ is \textbf{strictly forward (resp.~backward)} if it satisfies condition 2 above at every $t \in [0,1]$.
\end{defn}

We emphasize that we do not impose any constraints on the part of a forward curve which is contained in $U^\mathrm{un}$. Outside of that part, the curve has to be transverse to the $X$-lines it meets. If $\gamma : [0,1] \rightarrow U$ is a curve, we denote by $\check{\gamma} : [0,1] \rightarrow U$ the corresponding reversed curve defined by $\check{\gamma}(t) \coloneqq \gamma(1-t)$, $t \in [0,1]$. Obviously, $\gamma$ is forward if and only if $\check{\gamma}$ is backward.

The following lemma will be particularly useful in Section~\ref{sec:extorder} below and serves as a warm-up.

\begin{lem} \label{lem:forcurv}
Let $\gamma : [0,1] \rightarrow U$ be a forward curve. 
\begin{enumerate}
    \item If $\ell \subset U^\mathrm{pure}$ is a pure $X$-line, then $\gamma$ intersects $\ell$ at most once.
    \item If $\gamma(0) = \gamma(1)$, then $\gamma$ is entirely contained in an unstable component of $f$.
\end{enumerate}
\end{lem}

\begin{proof}
We show 1 by contradiction. Without loss of generality, let us assume that $\gamma(0) \in \ell$ and $\gamma(1) \in \ell$. We can further modify $\gamma$ so that $\gamma(0) = \gamma(1)$ and $\gamma'(0) = \gamma'(1)$, by flowing $\gamma$ along $\widetilde{X}$ near its starting point, since $\gamma$ is positively transverse to $\ell$ at $t =0$. After a small perturbation, we can further assume that $\gamma$ has finitely many (and transverse) self-intersections. The latter can be resolved (whether they occur in $U^\mathrm{un}$ or $U^\mathrm{pure}$) in order to obtain a finite collection of \emph{embedded} forward curves. We can therefore assume that $\gamma$ is embedded, and bounds a domain $D \subset U$ homeomorphic to a closed disk. We now analyze the possible behaviors of $\ell$ in $D$. Observe that $\gamma$ cannot intersect $\ell$ except at $x=\gamma(0)=\gamma(1)$, so one of the two connected components of $\ell \setminus \{x\}$ is entirely contained in $D$. Assume first the `negative' component is contained in $D$. Then by Poincar\'{e}--Bendixson and the absence of closed orbits of $\widetilde{X}$, $\ell$ necessarily converges in negative times to a singularity $x_0 \in U$ contained in the interior of $D$. Moreover, since $\ell$ is pure, then $x_0$ is not a source and at least one stable branch of $\widetilde{X}$ at $x_0$ lies in $U^\mathrm{pure}$. The latter cannot intersect $\gamma$ since $\gamma$ is forward, so it must converge in negative times to a saddle singularity $x_1$. Once again, one of the stable branches at $x_1$ must lie in $U^\mathrm{pure}$, and cannot intersect $\gamma$ nor converge to a singularity in negative times, a contradiction. The same argument applies when the `positive' component of $\ell \setminus \{x\}$ is contained in $D$ by considering unstable branches instead of stable branches.

The second item follows immediately, as a forward curve which is not entirely contained in $U^\mathrm{un}$ must intersect a pure $X$-line.
\end{proof}

Using similar arguments, one can easily show that the unstable components of a tile are simply connected. The proof is left as an exercise for the reader. Combined with the previous lemma, this implies that a forward curve is homotopic relative to $U^\mathrm{pure}$ (hence homotopic through forward curves) to an \emph{embedded} one. 

\begin{figure}[h]
    \centering
            \includegraphics{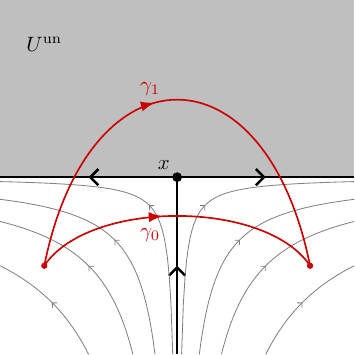}
    \caption{A (full) $T$-move.}
    \label{fig:Tmove}
\end{figure}

\begin{defn}
Let $\gamma_0 : [0,1] \rightarrow U$ and $\gamma_1 : [0,1] \rightarrow U$ be two forward curves with the same endpoints, namely, $\gamma_0(0) = \gamma_1(0)$ and $\gamma_0(1) = \gamma_1(1)$.
\begin{itemize}
    \item We say that $\gamma_0$ and $\gamma_1$ are \textbf{forward-homotopic} if they are homotopic relative to their endpoints through forward curves.
    \item We say that $\gamma_1$ is obtained from $\gamma_0$ by a \textbf{full $T$-move} if these two curves only intersect at their endpoints and bound a compact domain $D$ in $U$ containing exactly one singularity of $\widetilde{X}$. We further require that this singularity is quadratic-like, and $\gamma_0$ intersects its pure stable branch. See Figure~\ref{fig:Tmove}.
    \item We say that $\gamma_1$ is obtained from $\gamma_0$ by a \textbf{(local) $T$-move} if there exists a segment $I \subset [0,1]$ such that $\gamma_0$ and $\gamma_1$ coincide on $[0,1] \setminus I$, and ${\gamma_1}_{\vert I}$ is obtained from ${\gamma_0}_{\vert I}$ by a full $T$-move.
\end{itemize}
\end{defn}

Notice that two forward-homotopic curves intersect the same pure $X$-lines and the same unstable components in $U$. Two forward curves that differ by a $T$-move across a singularity $x$ intersect the same pure $X$-lines except for the pure stable branch and unstable branches at $x$, and they intersect the same unstable components in $U$ except for the unstable component at $f(x) \in Q^*$.

Two forward curves with the same endpoints are not necessarily forward-homotopic, but they become forward-homotopic after performing finitely many local $T$-moves:

\begin{lem} \label{lem:forcurv2}
Let $\gamma_0$ and $\gamma_1$ be two forward curves in $U$ with the same endpoints. There exists forward curves $\widetilde{\gamma}_0$ and $\widetilde{\gamma}_1$ such that
    \begin{itemize}
        \item For $i \in \{0,1\}$, $\gamma_i$ and $\widetilde{\gamma}_i$ are forward-homotopic,
        \item $\widetilde{\gamma}_0$ and $\widetilde{\gamma}_1$ differ by a finite sequence of local $T$-moves.
    \end{itemize}
\end{lem}

\begin{proof}
After applying forward-homotopies, we can assume that $\gamma_0$ and $\gamma_1$ are embedded and intersect at finitely many points. We can further reduce to the case where $\gamma_0$ and $\gamma_1$ only intersect at their endpoints. Let $D$ be the compact domain in $U$ bounded by $\gamma_0$ and $\gamma_1$. It contains finitely many singularities of $\widetilde{X}$. We will show the desired statement by induction on the number of these singularities.

If $D$ contains no singularity of $\widetilde{X}$, then every $X$-line in the interior of $D \cap U^\mathrm{pure}$ must intersect both $\gamma_0$ and $\gamma_1$ by Poincar\'{e}--Bendixson. If there is no such $X$-line, then (the interiors of) $\gamma_0$ and $\gamma_1$ are contained in the same unstable component and the result is immediate. Otherwise, since $D$ is compact, we can find an open neighborhood $U_0$ of $D \cap U^\mathrm{pure}$ foliated by $X$-lines intersecting $\gamma_0$ and $\gamma_1$ once each, and such that $D \setminus U_0$ intersects finitely many unstable components of $U$. We then use the flow of $\widetilde{X}$ to forward-homotope $\gamma_0$ to a forward curve $\widetilde{\gamma}_0$ which coincides with $\gamma_1$ on $U_0$, and we further homotope $\widetilde{\gamma}_0$ on each of the finitely many components of $D \setminus U_0$ to match $\gamma_1$.

Let us now assume that $x \in D$ is a singularity of $\widetilde{X}$. We analyze the three different possible cases.
\begin{itemize}[leftmargin=*]
    \item \textbf{Case 1:} $f(x) \in \Delta_{\mathrm{so}}$. Then at least one $X$-line in the unstable disk of $x$ intersects $\gamma_0$ or $\gamma_1$, say $\gamma_0$. We can forward-homotope $\gamma_0$ to $\widetilde{\gamma}_0$ so that $\widetilde{\gamma}_0$ and $\gamma_1$ bound a compact domain $D' \subset D$ which does not contain $x$, hence reducing the number of singularities of $\widetilde{X}$ between the two forward curves.  
    \item \textbf{Case 2:} $f(x) \in Q$. If $x \in U^\mathrm{un}$, we argue as in Case 1. Otherwise, the stable branch at $x$ lies in $U^\mathrm{pure}$ and must intersect $\gamma_0$ or $\gamma_1$, say $\gamma_0$. We then perform a local $T$-move on $\gamma_0$ along this stable branch to reduce the number of singularities of $\widetilde{X}$ in $D$.
    \item \textbf{Case 3:} $f(x) \in \Delta_\mathrm{sa}$. If $x \in U^\mathrm{un}$, then we can argue as before. Otherwise, at least one stable branch at $x$ lies in $U^\mathrm{pure}$ and intersects $\gamma_0$ or $\gamma_1$, say $\gamma_0$. If the other stable branch at $x$ lies in $U^\mathrm{un}$, then $f(x) \in Q^*$ and we argue as in Case 2. Otherwise, it lies in $U^\mathrm{pure}$ and either intersects $\gamma_0$, or converge in negative times to another saddle singularity $x'$ in $D$. In the latter case, at least one stable branch at $x'$ is contained in $U^\mathrm{pure}$, and it must intersect $\gamma_0$. In both cases, one of the two unstable branches at $x$ is `trapped' in a domain bounded by $\gamma_0$ and by stable and unstable branches at $x$ and at $x'$ if it exists. If this unstable branch converges to another saddle singularity, the unstable branches of the latter are also trapped and cannot converge to singularities. However, this is impossible by Poincar\'{e}--Bendixson and the absence of closed orbits of $\widetilde{X}$.
\end{itemize}
In all three cases, we constructed a forward-homotopy of $\gamma_0$ or $\gamma_1$, or performed some local $T$-move, reducing the number of singularities between $\gamma_0$ and $\gamma_1$. Note that for embedded forward curves, $T$-moves and forward-homotopies commute, so these operations can be performed in the order prescribed by the statement.
\end{proof}

\begin{figure}[h]
    \centering
            \includegraphics{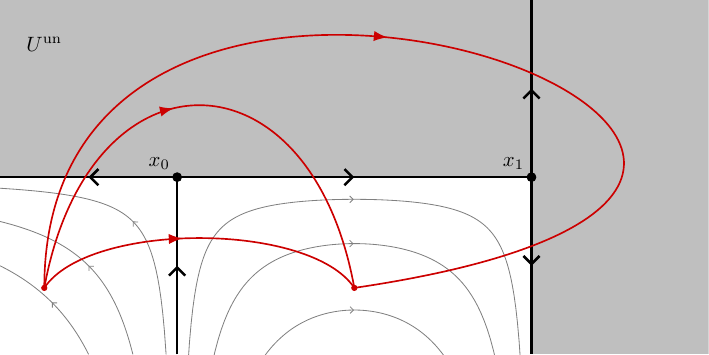}
    \caption{A double $T$-move.}
    \label{fig:2Tmove}
\end{figure}

Let us now consider the following situation. Assume that there exist two saddle singularities $x_0$ and $x_1$ in $U$ such that 
\begin{itemize}
    \item There exists a connection $\gamma$ from $x_0$ to $x_1$,
    \item One stable branch at $x_0$ is contained in $U^\mathrm{pure}$ while the other is contained in $U^\mathrm{un}$,
    \item The stable branch at $x_1$ different than $\gamma$ is contained in $U^\mathrm{un}$.
\end{itemize}
Then one can perform a sequence of two consecutive $T$-moves, one across $x_0$ followed by another one across $x_1$, as depicted in Figure~\ref{fig:2Tmove}. We call this operation a \textbf{double $T$-move}. By our assumptions on $X$, there cannot be longer sequences of $T$-moves.

\medskip

Let now $g_\star = (g : V \rightarrow M, b)$ be another marked tile satisfying $g(b)=f(a)$. The following lemma will be useful in the next sections.

\begin{lem} \label{lem:forlifts}
Let $\gamma_0, \gamma_1 : [0,1] \rightarrow U$ be two forward curves such that $\gamma_0(0) = \gamma_1(0) = a$ and $\gamma_0(1) = \gamma_1(1)$. If both $\gamma_0$ and $\gamma_1$ admit lifts $\widetilde{\gamma}_0$ and $\widetilde{\gamma}_1$ to $g_\star$, then $\widetilde{\gamma}_0(1) = \widetilde{\gamma}_1(1)$.
\end{lem}

\begin{proof}
If a forward curve $\gamma$ is obtained from $\gamma_0$ by a $T$-move across some singularity $x \in U$, we say that it is a \emph{negative} $T$-move if $\gamma$ intersects the pure stable branch at $x$. In that case, $\gamma$ lifts to $g_\star$ to a curve $\widetilde{\gamma}$ with $\widetilde{\gamma}(1) = \widetilde{\gamma}_0(1)$. Here, we used that $g$ has no slits.

The result is clear when $\gamma_0$ and $\gamma_1$ are forward-homotopic. In general, it is a consequence of the following immediate refinement of Lemma~\ref{lem:forcurv2}: there exists a forward curve $\gamma$ in $U$ such that for every $i \in \{0,1\}$, $\gamma$ is obtained from $\gamma_i$ by performing a forward-homotopy followed by finitely many \emph{negative} $T$-moves. Therefore, $\gamma$ lifts to $g_\star$ to a curve $\widetilde{\gamma}$ satisfying $\widetilde{\gamma}_0(1)=\widetilde{\gamma}(1)=\widetilde{\gamma}_1(1)$.
\end{proof}

\medskip

In the next sections, we will need a strengthening of the unstable-complete condition from Definition~\ref{def:tile}:

\begin{defn} \label{def:strongun}
A tile $f: U \rightarrow M$ is \textbf{strongly unstable-complete} if for every connected component $V$ of $U^\mathrm{un} \subset U$, there exists an unstable tile $g : W \rightarrow M$ and a diffeomorphism $\varphi : V \rightarrow W$ such that $g \circ \varphi = f$.
\end{defn}

Notice that a complete and forward tile is automatically strongly unstable-complete. 

%%%
        \subsection{Forward envelopes}
%%%

We consider a (possibly empty) weak prefoliation $(\mathscr{A}, \precsim^\mathscr{A})$ which is backward and complete. If $p \in M$, a \textbf{section} of $(\mathscr{A}_p, \precsim^\mathscr{A}_p)$ is a subset $\Sigma_0 \subset \mathscr{A}_p$ which is \emph{downward closed}: for all $f_\star \in \Sigma_0$ and $g_\star \in \mathscr{A}_p$, if $g_\star \precsim^\mathscr{A}_p f_\star$ then $g_\star \in \Sigma_0$. Typical examples of sections are:
\begin{itemize}
    \item $\Sigma_0 = \varnothing,$
    \item $\Sigma_0 = \big\{ g_\star \in \mathscr{A}_p \ \vert \ g_\star \precsim^\mathscr{A}_p f_\star\big\}$, where $f_\star \in \mathscr{A}_p$.
\end{itemize}

We now extend the definition of an \emph{upper enveloping surface} from~\cite{BI08}:

\begin{defn}[See~\cite{BI08}, Definition 6.1] \label{def:upenv}
An \textbf{upper enveloping tile} is a pair $(f, \Sigma)$, where $f:U \rightarrow M$ is a tile and $\Sigma : U \rightarrow \mathcal{P}\big( \mathscr{A}_\star \big)$ is a map which assigns a section $\Sigma(a) \subset \mathscr{A}_{f(a)}$ to every $a \in U$, such that the following conditions hold.
\begin{enumerate}
\setcounter{enumi}{-1}
    \item Either $f$ has exactly one proper edge $e_0$, and this edge is a pure $X$-line with forward coorientation, or $f$ has no proper edge but has a distinguished unstable component $U_0$. In the first case, we call $e_0$ the \textbf{initial edge} of $f$ and we say that $f$ is of \textbf{pure type}, and in the second case, we call $U_0$ the \textbf{initial unstable component} of $f$ and we say that $f$ is of \textbf{unstable type}. In both cases, we refer to $e_0$ or $U_0$ as the \textbf{initial locus} of $f$, and we denote it by $I_0$.
    \item $f$ is strongly unstable-complete, in the sense of Definition~\ref{def:strongun}.
    \item Every point in $U^\mathrm{good}$ can be reached from the initial edge or unstable component by a forward curve.
    \item Let $a \in U^\mathrm{reg}$ and $D \subset M$ be a smooth open $2$-disk transverse to $X$ and containing $f(a)$. Let $\gamma : [0,1] \rightarrow U$ be a strictly forward curve starting at $a$ and contained in $f^{-1}(D)$, and let $\mathcal{A}$ be the trace of $\mathscr{A}$ on $D$. We require that the curve $f \circ \gamma$ on $D$ is contained in the one-dimensional upper envelope of $(f(a), \Sigma(a) \cap D)$ with respect to $\mathcal{A}$ (as defined in~\cite[Section 5.1]{BI08}). Moreover, we require that for every $t \in [0,1]$, $\Sigma(\gamma(t))$ is the shadow of $(f(a), \Sigma(a) \cap D)$ at $f(\gamma(t))$ (as defined in~\cite[Section 5.1]{BI08}).
    \item $\Sigma$ is coherent along $X$-lines, in the following sense: for every $a \in U$ and $t \in \R$, $\Sigma\big(\widetilde{\phi}^t(a)\big) = \Phi^t\big( \Sigma(a) \big)$.
    \item $\Sigma$ is coherent at singularities, in the following sense. Let $a \in f^{-1}(\Delta)$ and $g_\star \in \mathscr{A}_{f(a)}$, where $g_\star = (g : V \rightarrow M, b)$, and assume that there exist
        \begin{itemize} 
        \item an $X$-line $\gamma$ in $U$ converging to the point $a$ in positive (resp.~negative) times, and
        \item two points $a' \in \gamma$ and $b' \in V$ such that $f(a')=g(b')$ and the $X$-line in $V$ passing through $b'$ converges to $b$ in positive (resp.~negative) times.
        \end{itemize}
    Then the following holds: $$g_\star \in \Sigma(a) \iff (g, b') \in \Sigma(a').$$
\end{enumerate}
\end{defn}

These conditions easily imply that $\Sigma$ is coherent along unstable components in the following sense. If $a, a' \in U$ belong to the same unstable component of $f$, then for every $g_\star = (g : V \rightarrow M, b) \in \mathscr{A}_{f(a)}$, $b \in V^\mathrm{un}$ and there is a unique $b' \in V^\mathrm{un}$ in the same unstable component as $b$ such that $g(b') = f(a')$. In that case, we have: $$(g, b) \in \Sigma(a) \iff (g,b') \in \Sigma(a').$$

Moreover, if $a_0 \in U$ belongs to the initial locus $I_0$ of $f$, then the map $\Sigma$ is entirely determined by $\Sigma(a_0) \subset \mathscr{A}_{f(a_0)}$. If $p = f(a_0)$ and $\Sigma(a_0) = \Sigma_0$, we say that the upper enveloping tile $F = (f, \Sigma)$ \textbf{starts at $(p, \Sigma_0)$}. Notice that $f$ is of pure type if and only if $p \in M^\mathrm{pure}$, and $f$ is of unstable type if and only if $p \in M^\mathrm{un}$. It will be convenient to include the choice of an initial \emph{basepoint} $a_0$ and consider \emph{pointed} upper enveloping tiles $F_\circ = (f, \Sigma, a_0)$. We will typically denote upper enveloping tiles by capital letters, such as $F$ for $(f, \Sigma)$, $G$ for $(g, \Sigma')$, etc. and triples $(f,\Sigma, a_0)$ as above by $F_\circ$.\footnote{We choose the subscript $\circ$ and not $\star$ to emphasize that the basepoint lies in the initial locus $I_0$ of the tile, and plays a different role than a marked point.} Following~\cite{BI08}, we call the map $\Sigma$ the \textbf{shadow} of $F$.

\medskip

For the rest of this section, we fix a point $p \in M^\mathrm{good}$ and a section $\Sigma_0 \subset \mathscr{A}_p$. We denote by $\mathscr{E}_p(\Sigma_0)$ the collection of pointed upper enveloping tiles starting at $(p, \Sigma_0)$.

The following lemma is a generalization of~\cite[Lemma 6.2]{BI08}. Although the main idea is the same, the proof is much more complicated in our setup.

\begin{lem} \label{lem:exists}
There exists an upper enveloping tile starting at $(p, \Sigma_0)$.
\end{lem}

\begin{proof}
If $p \in M^\mathrm{un}$, we simply consider the unstable tile $f$ passing through $p$ with an appropriate basepoint, and we extend $\Sigma_0$ in the only possible way (recall that $\eta$ is uniquely integrable along unstable tiles).

If $p \in M^\mathrm{pure}$, we construct an upper envelope starting at $p$ in 3 steps.
\begin{itemize}[leftmargin=*]
    \item \textit{Step 1: constructing a forward strip.} By definition, there is a pure flow line of $X$ passing through $p$. Let $D \subset M$ be a small open smooth disk transverse to $X$ and passing through $p$, and $\mathcal{A}$ be the trace of $\mathscr{A}$ on $D$. It is a two-dimensional prefoliation in the sense of~\cite[Section 5.1]{BI08}. Let $f_0 : [0,\epsilon) \rightarrow D$ be the (forward) upper envelope for $(p, \Sigma_0 \cap D)$ with respect to $\mathcal{A}$, with shadow $\overline{\Sigma}_0$, defined up to some $\epsilon > 0$. As opposed to~\cite[Lemma 6.2]{BI08}, we do not impose any lower bound on $\epsilon$. Let now $f: [0,\epsilon) \times \R \rightarrow M$ be the extension of $f_0$ under the flow of $X$, defined by $f(x, t) \coloneqq \varphi^t_X(f_0(x))$. The one-dimensional shadow $\overline{\Sigma}_0$ extends to a shadow $\Sigma$ for $f$. Notice that $f$ is not quite a tile yet, since it might fail to be strongly unstable-complete and well-cornered. To remedy this, we will further extend $f$.
    \item \textit{Step 2: unstable-completion.} Since the unstable locus of $f$ is open, there exists countably many open intervals $I_n \subset [0, \epsilon)$, $n \in \N$, such that this unstable locus is exactly the disjoint union of the subsets $I_n \times \R$. For every $n \in \N$, there exists an unstable tile $g_n : V_n \rightarrow M$ and a $C^1$ map $\varphi_n : I_n \times \R \rightarrow V_n$ satisfying $g_n \circ \varphi_n = f_{\vert I_n \times \R}$. If $I_n$ is of the form $(a,b)$ with $b< \epsilon$, then $\varphi_n$ is injective: the curve $t \in (a,b) \mapsto \varphi_n(t,0)$ in $V_n$ does not intersect the same $X$-line twice, as it reaches two distinct edges of $V_n$ on both sides. Therefore, $\varphi_n$ is an embedding. If $b=\epsilon$, $\varphi_n$ might fail to be injective, e.g., if the  transverse disk $D$ `spirals' around a source singularity. However, one can easily shrink $\epsilon$ (or modify $D$) so that $\varphi_n$ becomes an embedding. Recall that the completion of an unstable tile does not have closed boundary components, thanks to our hypothesis that $(X, \eta)$ has no disk of tangency.
        
    We denote by $\widetilde{f} : \widetilde{U} \rightarrow M$ the extension of $f$ obtained by gluing the $g_n$'s to $f$ along the $\varphi_n$'s. There is a natural embedding $\varphi : [0, \epsilon) \times \R \hookrightarrow \widetilde{U}$ satisfying $\widetilde{f} \circ \varphi = f$. By construction, $\widetilde{U}$ is homeomorphic to a simply connected subset of $\R^2$, and $\widetilde{f}$ is strongly unstable-complete. The shadow $\Sigma$ extends in a unique way to a shadow $\widetilde{\Sigma}$ for $\widetilde{f}$, and the conditions 0--5 in Definition~\ref{def:upenv} are satisfied. However, $\widetilde{f}$ might not be a tile yet, and we need to extend it a bit further in order to make it well-cornered.

    \item \textit{Step 3: filling the slits.} Let $q \in \Delta_\mathrm{sa} \cup Q$, not surrounded by unstable manifolds, and $D_0 \subset \widetilde{U}$ be a singular region of $\widetilde{f}$ for $q$. Notice that if $\widetilde{f}$ fails to be graphical near $q$ over $D_0$, or has a puncture or a slit at $q$ in $D_0$, then one of the two cases happen: 
        \begin{itemize}[leftmargin=*]
            \item \textbf{Case 1:} There exist two $X$-lines intersecting $D_0$ and mapping to opposite stable branches at $q$.
            \item \textbf{Case 2:} There exist two $X$-lines intersecting $D_0$ and mapping to opposite unstable branches at $q$.
        \end{itemize}

    We claim that Case 1 never happens. Otherwise, we could find a third $X$-line $\ell_u$ intersecting $D_0$ and mapping to an unstable branch at $q$. Since $q$ is not surrounded by unstable manifolds, $\ell_u$ and at least one of the stable $X$-lines lie in $\widetilde{U}^\mathrm{pure}$. We denote the latter one by $\ell_s$. By definition of $\widetilde{f}$, there exists an embedded and strictly forward curve $\gamma$ in $\widetilde{U}$ intersecting both $\ell_s$ and $\ell_u$ \emph{positively}. Moreover, $\gamma$ belongs to $\varphi\big([0, \epsilon) \times \R\big)$. Let us assume without loss of generality (up to reverting $X$ and considering backward curves instead) that $\gamma$ starts on $\ell_u$ and ends on $\ell_s$. We then `slide' $\gamma$ near its endpoint along $\ell_s$ towards $q$ and then along $\ell_u$ to obtain a \emph{closed} strictly forward curve intersecting $\ell_u$; see Figure~\ref{fig:slide}. By construction, this curve belongs to $\varphi\big([0, \epsilon) \times \R\big)$ which is impossible.

\begin{figure}[h]
    \centering
            \includegraphics{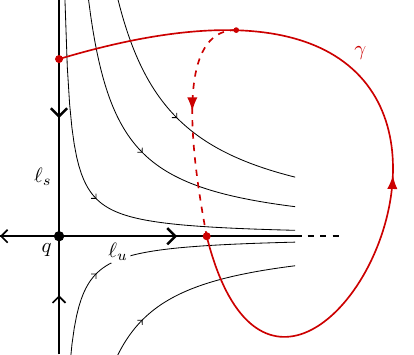}
    \caption{Sliding $\gamma$ near its endpoint in Case 1.}
    \label{fig:slide}
\end{figure}

    We now focus on Case 2. There are two subcases depending on whether $q \in \Delta_\mathrm{sa}$ or $q \in Q$. If $q \in \Delta_\mathrm{sa}$, the previous case implies that $D_0$ intersects a unique $X$-line mapping to a stable branch at $q$, and this $X$-line lies in $\widetilde{U}^\mathrm{un}$. Therefore, if $\widetilde{f}$ is not well-cornered at $q$ over $D_0$, then it has a slit corresponding to the other stable branch at $q$. Similarly, if $q \in Q$ and if $\widetilde{f}$ is not well-cornered at $q$ over $D_0$, then arguing as before shows that no $X$-line intersecting $D_0$ is mapped to the stable branch at $q$, and $\widetilde{f}$ has a slit corresponding to the stable branch at $q$. Note that it might happen that one of the unstable branches at $q$ is the initial edge of $\widetilde{f}$, i.e., the $X$-line corresponding to $f_{\vert \{0\} \times \R}$. In that case, it is easy to see that $\widetilde{f}$ is well-cornered at $q$ over $D_0$.

    To summarize, in both subcases, if $\widetilde{f}$ is not well-cornered at $q$ over $D_0$, then it has a slit corresponding to a stable branch of $q$, and the other stable branch lies in $M^\mathrm{un}$; see Figure~\ref{fig:slit}. Moreover, this stable branch might possibly be connected to another singularity $q_-$ in negative times, in which case $q_- \in \Delta_\mathrm{sa}$.

    We now explain how to fill the slits to make $\widetilde{f}$ well-cornered, and how to extend $\widetilde{\Sigma}$. With the previous notations, let $V \subset \widetilde{U}$ be the saturation of $D_0$ by $X$-lines. It is an open subset homeomorphic to an open disk. Clearly, there exists an embedding $\widehat{\varphi}_{q, D_0} = \widehat{\varphi} : V \hookrightarrow \widehat{V}$ where $\widehat{V}$ is homeomorphic to an open disk and $\widehat{\varphi}(V) \subset \widehat{V}$ is the complement of a proper embedding of $[0, \infty)$ into $\widehat{V}$, and such that $\widetilde{f}_{\vert V}$ extends on $\widehat{V}$ to a map $\widehat{f}_{q, D_0}$ satisfying:
        \begin{itemize}
        \item $\widehat{f}_{q, D_0} \circ \widehat{\varphi} = \widetilde{f}_{\vert V}$,
        \item $\widehat{f}_{q, D_0}$ is an injective immersion of an open disk into $M$ tangent to $\eta$ whose image contains $q$. 
        \end{itemize}

\begin{figure}[h]
    \centering
            \includegraphics{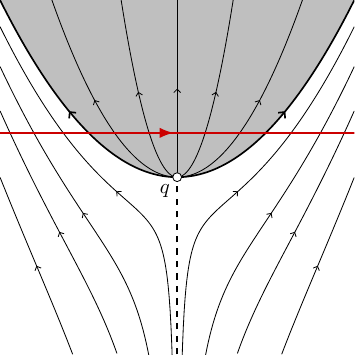}
    \caption{A slit at $q$. The red line corresponds to the forward curve induced by $f_0$.}
    \label{fig:slit}
\end{figure}

    We can now glue $\widetilde{f}$ and $\widehat{f}_{q, D_0}$ along $V$ to remove the slit at $q$ over $D_0$. Performing this operation at every singularity with a slit, we obtain an immersion $\widehat{f} : \widehat{U} \rightarrow M$ tangent to $\eta$, where $\widehat{U}$ is homeomorphic to a simply connected domain of $\R^2$. We further have to extend $\widetilde{\Sigma}$ to $\widehat{U}$ so that the properties of Definition~\ref{def:upenv} are satisfied. We fix $q$ and $D_0$ as before, so that $\widetilde{f}$ has a simple slit at $q$ over $D_0$, and we will view $D_0 \subset V \subset \widehat{V}$ as subsets of $\widehat{U}$. Let $x_0 \in \widehat{V}$ lying on the $X$-line mapped to the stable branch $\ell_s$ at $q$, and let $p_0 \coloneqq \widehat{f}(x_0)$. Let $D \subset M$ be a small embedded disk in $M$ passing through $p_0$ and transverse to $X$. Let $\gamma_0:[-1,1] \rightarrow \widehat{V}$ be a strictly forward curve with $\gamma_0(0) = x_0$ and such that $\widehat{f} \circ \gamma_0$ parametrizes the intersection of $D$ with $\widehat{f}\big(\overline{D}_0\big)$, where $\overline{D}_0$ denotes the closure of $D_0$ in $\widehat{U}$. We have to show that $\widehat{f} \circ \gamma_0$ is contained in the one-dimensional upper envelope of $\big(\widehat{f} \circ \gamma_0(-1), \widetilde{\Sigma}(\gamma_0(-1)) \cap D \big)$ in $D$, and that $\widetilde{\Sigma}(\gamma_0(t)) \cap D$ coincides with its shadow at $\widehat{f} \circ \gamma_0(t)$ for every $t \in [-1, 1] \setminus \{0\}$. In particular, this provides the desired extension of $\widetilde{\Sigma}$ at $x_0$, by considering the shadow of the former upper envelope at that point. To show the claim, we first observe that since $\mathscr{A}$ is backward and complete by hypothesis, a tile in $\mathscr{A}$ passing through $q$ is of one of the following forms near $q$: a quadrant with backward boundary and convex corner at $q$, three quadrants (if $q \in \Delta_\mathrm{sa}$) or the union of a quadrant and a half plane (if $q \in Q$) with backward boundary and concave corner at $q$, or it contains $q$ in the interior of its domain; see Figure~\ref{fig:backcor}. This induces a partition of $\mathscr{A}_q$ into three subsets $\mathscr{A}^\mathrm{cvx}_q$, $\mathscr{A}^\mathrm{ccv}_q$, and $\mathscr{A}^0_q$, depending on whether a marked tile has a convex corner, concave corner, or no corner at $q$. Moreover, let $p_+$ be a point of $\mathcal{N}$ on the right unstable branch at $q$, and let $x_+$ be the preimage in $\widehat{V}$ of $p_+$ by $\widehat{f}$. Then `sliding the marked point along stable/unstable branches at $q$' induces natural injective maps $\varsigma_0 : \mathscr{A}_{p_0} \rightarrow \mathscr{A}_q$ and $\varsigma_+:\mathscr{A}_{p_+} \rightarrow \mathscr{A}_q$, which preserve the preorder $\precsim$. Moreover, because of the form of the tiles in $\mathscr{A}_{q}$, $\varsigma_0$ is bijective, and $\varsigma_+$ is a bijection onto $\mathscr{A}^0_q \sqcup \mathscr{A}^\mathrm{ccv}_q$. Let $\gamma_+ : [0,1] \rightarrow \widehat{V}$ be an embedded strictly forward curve in $\widehat{V}$ from $x_+$ to $\gamma_0(1)$. If $D_+$ is an embedded disk in $M$ transverse to $X$ and containing $\widehat{f} \circ \gamma_+$, then $\widehat{f} \circ \gamma_+$ is contained in the one-dimensional upper envelope for $\big(p_+, \widetilde{\Sigma}(x_+) \cap D_+\big)$, and $\widetilde{\Sigma}(\gamma_+(t))\cap D_+$ coincides with its shadow at $\widehat{f} \circ \gamma_+(t)$ for every $t \in [0,1]$. We define a section $\widehat{\Sigma}_+ \subset \mathcal{A}_{p_0} = \mathscr{A}_{p_0} \cap D$ as the smallest section of $\mathcal{A}_{p_0}$ (with respect to the inclusion) containing $\widetilde{\Sigma}_+ \coloneqq \varsigma^{-1}_0 \circ \varsigma_+\big( \widetilde{\Sigma}(x_+)\big) \cap D$.\footnote{Concretely, if $(A, \leq)$ is a preordered set and $\Sigma \subset A$, then the smallest section of $A$ containing $\Sigma$ is $$\widehat{\Sigma} = \bigcup_{s \in \Sigma} \{ a \in A \ \vert \ a \leq s\}.$$} Notice that $\widehat{\Sigma}_+ \setminus \widetilde{\Sigma}_+$ only contains curves with a right boundary at $p_0$, corresponding to tiles which have a convex corner at $q$, i.e., whose images under $\varsigma_0$ lie in $\mathscr{A}^\mathrm{cvx}_q$. Therefore, the invariance of the various structures of interest under the flow of $X$ implies the curve ${\gamma_0}_{\vert [0,1]}$ is contained in the one-dimensional upper envelope of $\big(p_0, \widehat{\Sigma}_+ \big)$ in $D$ and its shadow at $\widehat{f} \circ \gamma_0(t)$ for $t \in (0,1]$ is exactly $\widetilde{\Sigma}(\gamma_0(t)) \cap D$. Now, if $\widehat{\Sigma}_- \subset \mathcal{A}_{p_0}$ denotes the shadow at $p_0$ of the one-dimensional upper envelope of $\big(\widehat{f} \circ \gamma_0(-1), \widetilde{\Sigma}(\gamma_0(-1)) \cap D \big)$ in $D$, then $\widehat{\Sigma}_+ \subset \widehat{\Sigma}_-$, and $\widehat{\Sigma}_- \setminus \widehat{\Sigma}_+$ only contains curves with a right boundary at $p_0$. This implies that ${\gamma_0}_{\vert (0,1]}$ is contained in the latter upper envelope, and its shadow at $\widehat{f} \circ \gamma_0(t)$ for $t \in (0,1]$ is $\widetilde{\Sigma}(\gamma_0(t)) \cap D$, as desired.

\begin{figure}[h]
            \centering
            \begin{subfigure}{0.33\textwidth}
                    \includegraphics{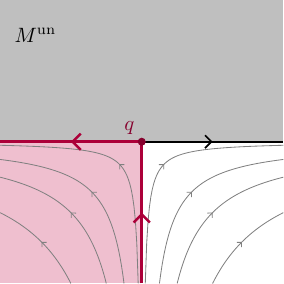}
            \centering
            \caption{A tile in $\mathscr{A}^\mathrm{cvx}_q$.}
            \end{subfigure}%
            \begin{subfigure}{0.33\textwidth}
                    \includegraphics{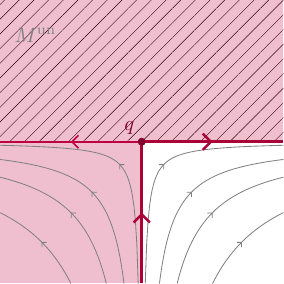}
            \centering
            \caption{A tile in $\mathscr{A}^\mathrm{ccv}_q$.}
            \end{subfigure}%
            \begin{subfigure}{0.33\textwidth}
                    \includegraphics{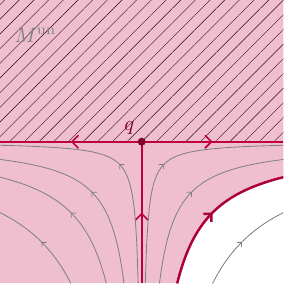}
            \centering
            \caption{A tile in $\mathscr{A}^0_q$.}
            \end{subfigure}
    \caption{The three types of tiles in $\mathscr{A}_q$.}
    \label{fig:backcor}
\end{figure}

    Moreover, if $x \in \widehat{V}$ is the preimage of $q$ under $\widehat{f}$, there is a unique way to extend $\widehat{\Sigma}$ at $x$ so that condition 5 of Definition~\ref{def:upenv} is satisfied.

    Unfortunately, $\widehat{f}$ might still fail to be well-cornered, as new slits can appear in the extension process. This situation arises exactly when $\widetilde{f}$ has a slit corresponding to a stable branch of a singularity $q$ which is connected to another singularity $q_- \in \Delta_\mathrm{sa}$ in negative times, one of the stable branches at $q_-$ lies in $M^\mathrm{un}$, and the unstable branch at $q_-$ not connected to $q$ is already contained in $\widetilde{f}$. In that case, filling the slit at $q$ creates a slit at $q_-$ corresponding to the other stable branch at $q_-$. Our admissibility conditions on $X$ guarantee that the latter stable branch is not connected to another singularity of $X$. We can then apply the previous argument at these `secondary slits' to further extend $\big(\widehat{f}, \widehat{\Sigma}\big)$ to a pair $\big(\wideparen{f}, \wideparen{\Sigma}\big)$, where $\wideparen{f}$ is now a \emph{tile}. By construction, all the conditions of Definition~\ref{def:upenv} are satisfied, hence $\big(\wideparen{f}, \wideparen{\Sigma}\big)$ is an upper enveloping tile starting at $(p, \Sigma_0)$. \qedhere
\end{itemize}
\end{proof}

We now define a natural preorder, the \emph{extension preorder}, on $\mathscr{E}_p(\Sigma_0)$.

\begin{defn}
Let $F_\circ =(f : U \rightarrow M, \Sigma, a_0)$ be a pointed upper enveloping tile. An \textbf{extension} of $F$ is a pointed upper enveloping tile $G_\circ =(g : V \rightarrow M, \Sigma', b_0)$ together with an embedding $\varphi : U \rightarrow V$ satisfying
\begin{itemize}
    \item $g \circ \varphi = f$,
    \item $\Sigma' \circ \varphi = \Sigma$,
    \item $\varphi(a_0) = b_0$.
\end{itemize}
If $\varphi$ is a diffeomorphism, we say that $F_\circ$ and $G_\circ$ are \textbf{equivalent}. Otherwise, we say that $G_\circ$ is a \textbf{strict extension} of $F_\circ$. We call $\varphi$ the \textbf{extension morphism}, and we write $\varphi : F_\circ \hookrightarrow G_\circ$, or simply $F_\circ \hookrightarrow G_\circ$.

If $F$ and $G$ are upper enveloping tiles without basepoints, we say that \textbf{$G$ is an extension of $F$} (resp.~$F$ and $G$ are \textbf{equivalent}) if there exist basepoints $a_0$ and $b_0$ for $F$ and $G$, respectively, such that $(G, b_0)$ is an extension of $(F, a_0)$ (resp.~$(F, a_0)$ and $(G, b_0)$ are equivalent).
\end{defn}

Observe that the extension morphism $\varphi$ is \emph{unique}: if $\varphi$ and $\psi$ are two extension morphisms from $F_\circ$ to $G_\circ$, then the set $\{ a \in U \ \vert \ \varphi(a)=\psi(a)\}$ is nonempty and closed, and it is open since $g \circ \varphi = g \circ \psi$ and $g$ is an immersion. This implies that if $\varphi : F_\circ \hookrightarrow G_\circ$ and $\psi : G_\circ \hookrightarrow F_\circ$, then $\phi$ and $\psi$ are diffeomorphism which are inverse of each other. In particular, $F_\circ$ and $G_\circ$ are equivalent. This shows that the preorder $\hookrightarrow$ on $\mathscr{E}_p(\Sigma_0)$ is a partial order up to equivalence.

We will be interested in the maximal elements for this preorder. By the following lemma, those are exactly the upper enveloping tiles without backward edges:

\begin{lem} \label{lem:ext}
An upper enveloping tile admits a strict extension if and only if it has an edge with backward coorientation.
\end{lem}

\begin{proof}
Let $F = (f: U \rightarrow M, \Sigma)$ be an upper enveloping tile with a choice of initial basepoint $a_0 \in I_0$. 

Let us first assume that $F_\circ = (F, a_0)$ admits a strict extension $\varphi : F_\circ \hookrightarrow G_\circ$, where $G_\circ = (g:V \rightarrow M, \Sigma', b_0)$. We identify $U$ with the subset $\varphi(U) \subset V$. Since the extension is strict, the latter inclusion is strict. Let $x \in V \setminus U$, and $\gamma$ be a forward curve from $a_0$ to $x$ in $V$. Let $C_0 \subset [0,1]$ be the connected component of $\gamma^{-1}(U)$ containing $0$, and $t_0 \coloneqq \sup C_0$. Since $U$ is open in $V$, $x_0 \coloneqq \gamma(t_0) \notin U$. Moreover, $x_0 \notin V^\mathrm{un}$, since $f$ and $g$ are strongly unstable-complete, so $x_0 \in V^\mathrm{pure}$. If $\ell$ denotes the pure $X$-line passing through $x_0$, then $\ell$ corresponds to a backward (improper) edge of $f$.

We now assume that $F$ has an improper backward edge corresponding to a flow line $\gamma$ of $X$ in $M^\mathrm{pure}$. Let $p \in \gamma$. By condition 3 of Definition~\ref{def:upenv}, $\Sigma$ induces a section $\Sigma_p \subset \mathscr{A}_p$ in a natural way, by extending the one-dimensional envelope and its shadow at $p$. By Lemma~\ref{lem:exists}, we consider an upper enveloping tile $H_\circ = (h : W \rightarrow M, \Sigma'', c_0)$ starting at $(p, \Sigma_p)$. We can then extend $F$ across $\gamma$ by attaching $H$ to obtain an immersion $g : V \rightarrow M$ tangent to $\eta$ together with a map $\Sigma': V \rightarrow \mathcal{P}(\mathscr{A}_\star)$ extending $\Sigma$ and $\Sigma''$. If $\gamma$ is not connected to a singularity of $X$ in positive or negative times, then one easily checks that $(g, \Sigma')$ is an upper enveloping tile, and that it is a strict extension of $F$. In general, $g$ might have \emph{slits}. This happens precisely when $\gamma$ is part of an improper \emph{concave} corner of $f$, and/or $\gamma$ is connected in negative times to a quadratic-like singularity $q \in Q^*$. The methods from Step 3 in the proof of Lemma~\ref{lem:exists} can be adapted to fill these slits (potentially in two steps for `double slits') to obtain the desired extension of $F$. The (tedious yet straightforward) details are left to the reader.
\end{proof}

The previous lemma motivates the following

\begin{defn}[See~\cite{BI08}, Definition 6.10]
A \textbf{forward envelope} is an upper enveloping tile such that all of its edges have forward coorientation.
\end{defn}

The construction of a forward envelope in our setup will be more delicate than in~\cite{BI08}. A crucial step is to show that the preorder $\hookrightarrow$ on $\mathscr{E}_p(\Sigma_0)$ is \emph{directed}.

\begin{lem} \label{lem:comext}
Two pointed upper enveloping tiles $F_\circ$ and $G_\circ$ starting at $(p, \Sigma_0)$ admit a common extension. More precisely, there exist a pointed upper enveloping tile $H_\circ$ starting at $(p, \Sigma_0)$ and two extension morphisms $\varphi : F_\circ \hookrightarrow H_\circ$ and $\psi : G_\circ \hookrightarrow H_\circ$.
\end{lem}

\begin{proof}
Let us write $F_\circ = (f : U \rightarrow M, \Sigma, a_0)$ and $G_\circ = (g : V \rightarrow M, \Sigma', b_0)$. We can first extend $F_\circ$ and $G_\circ$ so that both $f$ and $g$ are \emph{$T$-complete}, in the following sense: if $U$ contains an $X$-line mapped to a stable branch of a quadratic-like singularity, then $U$ contains the corresponding unstable tile as well. This can be achieved by applying the strategy of Lemma~\ref{lem:ext} at the improper backward edges of $f$ and $g$ corresponding to such stable branches, and applying it a second time in the situations corresponding to double $T$-moves.

We denote by $\mathcal{U}$ the collection of open subset $U' \subset U$ containing $a_0$ such that
\begin{itemize}
    \item $U'$ is \emph{forward}, namely, for every point in $x \in U' \cap U^\mathrm{good}$, there exists a forward curve from $a_0$ to $x$ contained in $U'$,
    \item There exists an embedding $\iota : U' \hookrightarrow V$ such that $\iota(a_0) = b_0$ and $g \circ \iota = f_{\vert U'}$.
\end{itemize}
The embedding $\iota$ is then uniquely determined as follows: if $\gamma : [0,1] \rightarrow U$ is forward curve in $U'$ from $a_0$ to $x$, then it lifts to $g_\star$ to a forward curve $\widetilde{\gamma} = \iota \circ \gamma$ in $V$ from $b_0$ to $\iota(x)$. Since such a lift is unique, any other embedding $\iota'$ as in the second item would satisfy $\iota'(x) = \iota(x)$. 

We argue that $\mathcal{U}$ is not empty. It is obvious if $p \in M^\mathrm{un}$, since $f$ and $g$ are strongly unstable-complete. If $p \in M^\mathrm{pure}$, it follows from condition 3 in Definition~\ref{def:upenv}, since $1$-dimensional upper envelopes are unique (see also~\cite[Corollary 6.4]{BI08}). Moreover, Lemma~\ref{lem:forlifts} readily implies that $\mathcal{U}$ is stable under taking unions. We then define
$$U_g \coloneqq \bigcup_{U' \in \mathcal{U}} U',$$
with a corresponding embedding $\iota : U_g \rightarrow M$. Note that $U_g$ is saturated by the flow of $\widetilde{X}$ and strongly unstable-complete. We make the following claim, which will be proved later.
\textit{ \begin{flushleft}
\noindent \textbf{Claim.} Let $\gamma$ be a forward curve in $U$ from $x_0$ to $x_1$. If $x_1 \in U_g$, then $\gamma \subset U_g$, hence $x_0 \in U_g$.
\end{flushleft}
}

Specializing at $x_0 = a_0$, this claim together with the conditions in Definition~\ref{def:upenv} readily imply that $f$ and $g$ have the same shadow on $U_g$, in the sense that $\Sigma' \circ \iota = \Sigma_{\vert U_g}$.

We construct a common extension of $F_\circ$ and $G_\circ$ by gluing $f$ and $g$ along $U_g$ via $\iota$. More precisely, we define $$W \coloneqq \big( U \sqcup V \big) \slash \sim_\iota,$$
where $x \sim_\iota y$ if and only if $x \in U_g$, $y \in V$, and $y= \iota(x)$. We have natural injective maps $\varphi : U \hookrightarrow M$ and $\psi : V \hookrightarrow M$, together with a unique map $h : W \rightarrow M$ satisfying $h \circ \varphi = f$ and $h \circ \psi = g$. We extend $\Sigma$ and $\Sigma'$ to a map $\Sigma'' : W \rightarrow \mathcal{P}\big(\mathscr{A}_\star\big)$ satisfying $\Sigma'' \circ \varphi = \Sigma$ and $\Sigma'' \circ \psi = \Sigma'$.

So far, $W$ is only defined as a \emph{topological space}. We will show that it has a natural structure of $2$-manifold, possibly with boundary, such that $\varphi$ and $\psi$ are smooth, and $h$ is a $C^1$ immersion tangent to $\eta$. Note that $W$ is simply connected by the Seifert--Van Kampen theorem.

We analyze how $f$ and $g$ overlap near a point $x \in \partial U_g \subset U$. We first assume that $\underline{x \in U^\mathrm{pure}}$, and we denote by $\ell$ the $X$-line passing through $x$. Let us consider a small neighborhood $O_q$ of $q\coloneqq f(x)$ in $M$ such that the connected component of $f^{-1}(O_q)$ containing $x$ is diffeomorphic to an open disk. Denoting this connected component by $U_x$, we can further assume that $U_x \subset U^\mathrm{reg}$. We will show that $U_x \cap U_g$ is the left open half-disk in $U_x$ bounded by $\ell$, and that the connected component of $g^{-1}(O_q)$ intersecting $V_x \coloneqq \iota(U_x \cap U_g)$ is precisely $V_x$. In particular, $g$ has an improper edge bounding the side of $V_x$ corresponding to the flow line of $X$ passing through $q$. Informally, $g$ is entirely contained in $f$ near $x$, and we can use $f$ to define a chart for $W$ near $\varphi(x)$.

To prove this, we consider two points $x_0, x_1 \in U_x$, where $x_0$ is to the left of $\ell$ and $x_1$ is to the right of $\ell$, as well as a (strictly) forward curve $\gamma$ from $x_0$ to $x_1$ contained in $U_x$.  Since $x \in \partial U_g$, there must exist points along $\gamma$ arbitrarily close to $x$ belonging to $U_g$, and the claim above implies that $x_0 \in U_g$. Moreover, if $x_1 \in U_g$, then still by the claim, $\gamma$ would be entirely contained in $U_g$, contradicting that $x \in \partial U_g \subset U \setminus U_g$. Therefore, we obtained that $x_0 \in U_g$ and $x_1 \notin U_g$. This implies the desired description for $U_x \cap U_g$, and the description of $g$ easily follows.

We now assume that $\underline{x \in f^{-1}(\Delta)}$. Then $f(x) \in \Delta_\mathrm{sa} \cup Q$. With the above notations, the previous argument can be adapted to show that $U_x \cap U_g$ is an open quadrant or the interior of the union of two or three consecutive quadrants near $x$, and $g$ is entirely contained in $f$ near $x$. The case of two opposite open forward quadrants can be ruled out by considering two forward curves in $U_g$ starting at $a_0$ and ending in these two quadrants. These curves can be extended past consecutive stable and unstable branches at $x$ in order to obtain two forward curves from $a_0$ to a common $x' \in U$. None of the stable branches at $x$ are in $U^\mathrm{un}$, since otherwise $U_g$ would contain two consecutive quadrants at $x$. Moreover, one of these stable branches is not connected to any other singularity, and we assume that it is the stable branch involved in the previous extension. Since our two forward curves are related by forward-homotopies and $T$-moves, and that the chosen stable branch cannot be involved in a $T$-move, it must intersect both forward curves. A contradiction quickly follows by analyzing which portions of the forward curves are contained in $U_g$. The other possible cases for the shape of $U_g$ at $x$ are easier to investigate, and the details are left to the reader.

We can define a collection of open subsets $\mathcal{V}$ of $V$ in the same manner as for $U$, and the corresponding maximal element $V_f$ coincides with $\iota(U_g)$, with associated embedding $\iota^{-1} : \iota(U_g) \hookrightarrow V$. By symmetry, a similar description of the overlap between $g$ and $f$ holds near $\partial V_f \subset V$. Therefore, there exists a (unique) smooth structure on $W$ making $\varphi$ and $\psi$ smooth and for which $h$ is $C^1$.

From the above arguments, it is easy to check that $h$ is an upper enveloping tile starting at $(p, \Sigma_0)$. If $c_0 \coloneqq \varphi(a_0) = \psi(b_0)$, then $H_\circ \coloneqq (h : W \rightarrow M, \Sigma'', c_0)$ is a common extension of $F_\circ$ and $G_\circ$, via $\varphi$ and $\psi$.

We finish the proof by showing the claim. We first make the following observation: if $\gamma_0$ and $\gamma_1$ are two forward curves in $U$ with the same endpoints, then $\gamma_0$ is contained in $U_g$ if and only if $\gamma_1$ is. Indeed, by Lemma~\ref{lem:forcurv2}, $\gamma_0$ and $\gamma_1$ differ by a finite sequence of $T$-moves and a forward-isotopy. Since the latter preserves the set of pure $X$-lines and unstable components intersected by a forward curve, we can reduce to the case where $\gamma_0$ and $\gamma_1$ differ by a full $T$-move. Then, the desired equivalence easily follows from the conditions in Definition~\ref{def:upenv} and the fact that $f$ and $g$ are $T$-complete and without slits.

Now, with the notations from the claim, we assume that $x_1 \in U_g$ and we denote by $\gamma_0$ a forward curve in $U$ from $a_0$ to $x_0$. The concatenation $\widetilde{\gamma}_1 \coloneqq \gamma_0 \ast \gamma$ is then a forward curve from $a_0$ to $x_1$. Let also $\gamma_1$ be a forward curve in $U_g$ from $a_0$ to $x_1$. Since $\gamma_1$ and $\widetilde{\gamma}_1$ are two forward curves from $a_0$ to $x_1$ and $\gamma_1 \subset U_g$, the previous paragraph implies that $\widetilde{\gamma}_1 \subset U_g$, hence $\gamma \subset U_g$ as desired.
\end{proof}

We can now prove the main result of this section:

\begin{prop}[See ~\cite{BI08}, Lemma 6.11] \label{prop:forenv}
There exists a forward envelope starting at $(p, \Sigma_0)$. It is unique up to equivalence.
\end{prop}

\begin{proof}
By the previous observations and lemmas, it suffices to show that $\mathscr{E}_p(\Sigma_0)$ admits a maximal element with respect to the extension preorder. Indeed, a maximal element is a forward envelope starting at $(p, \Sigma_0)$, and it is unique up to equivalence since the preorder is directed. In particular, it is a greatest element.

We consider the disjoint union 
$$\underline{U} \coloneqq \bigsqcup_{\substack{F_\circ \in \mathscr{E}_p(\Sigma_0),  \\ f : U \rightarrow M}} U,$$
which ranges over all the pointed upper enveloping tiles $F_\circ = (f : U \rightarrow M, \Sigma, a_0)$ starting at $(p, \Sigma_0)$. We denote the component corresponding to $F_\circ$ by $U_{F_\circ}$. To avoid set-theoretic issues, we can restrict to the ones with $U \subset \R^2$. By definition, every tile is equivalent to one whose domain is a subset of $\R^2$.

There is a natural map $\underline{f} : \underline{U} \rightarrow M$ whose restriction on the component corresponding to $F_\circ = (f : U \rightarrow M,\Sigma, a_0)$ is $f$. We define an equivalence relation $\equiv$ on $\underline{U}$ as follows: if $a \in U_{F_\circ}$ and $b \in U_{G_\circ}$, then $a \equiv b$ if and only if there exists a common extension $H_\circ$ of $F_\circ$ and $G_\circ$ with extension morphisms $\varphi : F_\circ \hookrightarrow H_\circ$ and $\psi : G_\circ \hookrightarrow H_\circ$ such that $\varphi(a) = \psi(b)$. By the existence of common extensions and the uniqueness of extension morphisms, it is easy to check that $\equiv$ is a well-defined equivalence relation on $\underline{U}$. Moreover, if $a \equiv b$ then $\underline{f}(a) = \underline{f}(b)$. Therefore, $\underline{f}$ induces a map 
$$f^\infty : U^\infty \rightarrow M, \qquad U^\infty \coloneqq \underline{U} \slash \equiv.$$
Here, we endow $U^\infty$ with the quotient topology.

We now show that $U^\infty$ is a $2$-dimensional manifold, possibly with boundary. By definition, $U^\infty$ admits a smooth atlas for which $f^\infty$ is a $C^1$ immersion tangent to $\eta$, and $U^\infty$ is saturated by the flow of $X$. Note that $U^\infty$ has a distinguished base point $a^\infty_0$ corresponding to the class of the initial base point of any element in $\mathscr{E}_p(\Sigma_0)$, and every point in $U^\infty$ can be connected to $a^\infty_0$ by a continuous path. Moreover, for every $F_\circ = (f : U \rightarrow M, \Sigma, a_0) \in \mathscr{E}_p(\Sigma_0)$, there is a natural smooth embedding $q_{F_\circ} : U \hookrightarrow U^\infty$ satisfying $f^\infty \circ q_{F_\circ} = f$. Lemma~\ref{lem:comext} easily implies that $U^\infty$ is Hausdorff and satisfies the following property : for every compact $K \subset U^\infty$, there exists $F_\circ \in \mathscr{E}_p(\Sigma_0)$ and a compact $K' \subset U$ such that $K = q_{F_\circ}(K')$. One deduces that $U^\infty$ is simply connected and aspherical, since every $U_{F_\circ}$ is. Observe that if $p \in M^\mathrm{un}$ then $U^\infty$ has no boundary. Otherwise, $U^\infty$ has a unique boundary component and it contains $a^\infty_0$. Furthermore, $U^\infty$ admits a (continuous) Riemannian metric (by either pulling back the Riemannian metric $\boldsymbol{g}$ on $M$ along $f^\infty$, or patching together the intrinsic Riemannian metrics on the $U_{F_\circ}$'s). Since it is Hausdorff, it is metrizable and since it is connected, it is second-countable; see~\cite[Appendix A, Theorem 1]{S99}.\footnote{We are grateful to Ryan Unger for pointing out this reference to us.} Therefore, we have shown that $U^\infty$ is a smooth surface, possibly with boundary, which embeds into $\R^2$.

The $\eta$-surface $f^\infty : U^\infty \rightarrow M$ is obviously strongly unstable-complete. One easily checks that it is well-cornered, using that every path $\gamma : [0,1]  \rightarrow U^\infty$ admits a lift $\widetilde{\gamma} : [0,1] \rightarrow U_{F_\circ}$ for some $F_\circ \in \mathscr{E}_p(\Sigma_0)$, i.e., it satisfies $q_{F_\circ} \circ \widetilde{\gamma} = \gamma$. Therefore, $f^\infty$ is a tile.

The collection of shadows of elements $F_\circ \in \mathscr{E}_p(\Sigma_0)$ naturally induces a shadow $\Sigma^\infty$ for $f^\infty$. It is uniquely characterized by the following property : for every $F_\circ \in \mathscr{E}_p(\Sigma_0)$, $\Sigma^\infty \circ q_{F_\circ} = \Sigma$, where $\Sigma$ denotes the shadow of $F_\circ$. It is straightforward to check that $(f^\infty, \Sigma^\infty)$ is an upper enveloping tile, using that any path in $U^\infty$ can be lifted to a path in some $U_{F_\circ}$. Therefore, $F^\infty_\circ \coloneqq (f^\infty, \Sigma^\infty, a^\infty_0)$ is an element of $\mathscr{E}_p(\Sigma_0)$. By construction, $F^\infty_\circ$ is an extension of every $F_\circ \in \mathscr{E}_p(\Sigma_0)$, which means that it is a greatest element for the extension order.
\end{proof}

%%%%%%%%%%%%%%%%%%%%%%%%%%%%%%%%%%%%%%%%%%%%%%%%%%%%%%%%%%%%%%%%%%%%%
        \section{Extending the order} \label{sec:extorder}
%%%%%%%%%%%%%%%%%%%%%%%%%%%%%%%%%%%%%%%%%%%%%%%%%%%%%%%%%%%%%%%%%%%%%

Let $\mathscr{U}$ denote the collection of forward envelopes starting at the pairs $(p, \Sigma_0)$ where 
\begin{itemize}
    \item $p \in M^\mathrm{good}$, and \textbf{$p$ does not belong to the pure stable branch of a quadratic-like singularity,}\footnote{This condition is not strictly necessary, but it will greatly simplify the construction of the preorder on $\mathscr{U}$ by drastically reducing the number of cases to consider in several proofs.}
    \item $\Sigma_0 \subset \mathscr{A}_p$ is any section.
\end{itemize}
To avoid set-theoretic issues, we can restrict our attention to the forward envelopes whose domains are subsets of $\R^2$, which is always the case up to equivalence.

We will extend the results of~\cite[Section 6.2]{BI08} to our setup. More precisely, we will define a preorder $\precsim$ on $\mathscr{U}_\star$ making $(\mathscr{U}, \precsim)$ a weak prefoliation. This preorder has to satisfy the following checklist:
\begin{itemize}
    \item It is compatible with the geometric order and coherent along intersections,
    \item It extends to the set of \emph{completed} forward envelopes $\overline{\mathscr{U}}$,
    \item Every marked tile in $\mathscr{A}_\star$ with marking on an edge has an \emph{upper neighbor} (see~\cite[Definition 4.10]{BI08} and Section~\ref{sec:upneigh} below) in $\overline{\mathscr{U}}$. 
\end{itemize}
Unfortunately, these requirements are not enough to entirely specify $\precsim$ and make it a total relation. We will adapt the methods of~\cite{BI08} and define $\precsim$ by successive `layers', depending on how the forward envelopes intersect. This will crucially rely on the \emph{forward structure} with respect to $X$ on forward envelopes, which essentially allows us to reduce the construction to a $1$-dimensional problem. However, this forward structure is more complicated in our setup because of the presence of singularities.

%%%
        \subsection{Forward order on pure \texorpdfstring{$X$}{X}-lines}
%%%

Let $f : U \rightarrow M$ be a tile.

\begin{defn}
We denote by $\mathcal{L} = \mathcal{L}_f$ the set of \underline{pure} $X$-lines of $f$ (see Definition~\ref{def:loci}). We define a (strict) preorder $\lhd$ on $\mathcal{L}$ as follows: if $\ell_0, \ell_1 \in \mathcal{L}$, then $\ell_0 \lhd \ell_1$ if and only if there exists a forward curve from $\ell_0$ to $\ell_1$. We denote by $\unlhd$ the corresponding (nonstrict) preorder.
\end{defn}

\begin{figure}[h]
    \centering
        \includegraphics{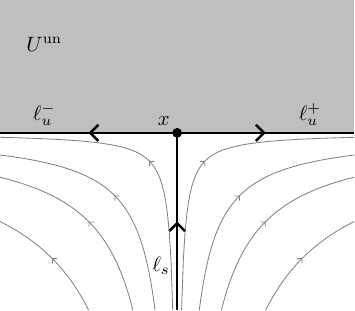}
    \caption{A special triple.}
    \label{fig:spetrip}
\end{figure}

Lemma~\ref{lem:forcurv} immediately implies that $\unlhd$ is a \emph{partial order}. We call it the \textbf{forward order} on $\mathcal{L}$. Recall that two forward-homotopic curves intersect the same $X$-lines in $\mathcal{L}$. If they differ by a $T$-move across a singularity $x \in U$, the sets of $X$-lines in $\mathcal{L}$ that they intersect differ as follows. One intersects the stable branch $\ell_s \in \mathcal{L}$ of $x$ but not its unstable branches $\ell_{u}^-, \ell_{u}^+ \in \mathcal{L}$, and the other intersects $\ell_{u}^-$ and $\ell_{u}^+$ but not $\ell_s$. Otherwise, they intersect the same $X$-lines in $\mathcal{L} \setminus \{ \ell_s, \ell_u^-, \ell_u^+\}$. Here, we chose our notations so that $\ell_u^- \lhd \ell_u^+$, see Figure~\ref{fig:spetrip}. We call $\big(\ell_s, \ell_u^-, \ell_u^+\big)$ a \textbf{special triple} in $\mathcal{L}$.

The next lemma is an easy consequence of Lemma~\ref{lem:forcurv}. The proof is left to the reader.
\begin{lem} \label{lem:propT}
We have the following properties:
\begin{enumerate}
    \item $\ell_s$ and $\ell_u^\pm$ are not comparable with respect to $\unlhd$.
    \item If $\ell \in \mathcal{L}$, then 
\begin{align*}
\ell_s \lhd \ell &\iff \ell^+_u \lhd \ell,\\
\ell \lhd \ell_s &\iff \ell \lhd \ell^-_u, \\
\ell^-_u \unlhd \ell \unlhd \ell^+_u &\iff \ell \in \{ \ell^-_u, \ell^+_u\}.
\end{align*}
\end{enumerate}
\end{lem}

Notice that if $f$ is the underlying tile of a forward envelope with a (necessarily unique and forward) proper edge, then this edge corresponds to \emph{the smallest} element of $(\mathcal{L}, \unlhd)$. If $f$ has no edge and has an initial unstable component $U_0$, then the minimal elements of $(\mathcal{L}, \unlhd)$ are exactly the $X$-lines in $\partial U_0 \subset U$ with backward coorientation with respect to $U_0$.

%%%
            \subsection{Backward domains of overlap}
%%%

In this section, we adapt the definitions from~\cite[Section 6.2.1]{BI08}. Our goal is to precisely describe how marked forward envelopes overlap in the backward region determined by the marked point. This will be crucial in the definition of the preorder on forward envelopes in the next section. As in~\cite[Section 6.2]{BI08}, this preorder will depend on how two marked forward envelopes overlap in this backward region.

%%%
                \subsubsection{Backward domain of an envelope}
%%%

We fix a point $p \in M^\mathrm{good}$ together with a marked forward envelope $f_\star = (f : U \rightarrow M, a) \in \mathscr{U}_p$, which implicitly comes with the additional data of a shadow $\Sigma$ and an initial locus $I_0$.

We denote by $\mathcal{I}_{f_\star} \subset \mathcal{L}_f$ the set of all $X$-lines in $\mathcal{L}_f$ intersected by the \emph{backward} curves starting at $a$, or alternatively, by the forward curves from $I_0$ to $a$. If $f$ is of pure type and $a$ lies on the initial edge $e_0$ of $f$, we set $\mathcal{I}_{f_\star} \coloneqq \{e_0\}$ instead. This set is well-defined by the second item of Definition~\ref{def:upenv}. It is empty exactly when $f$ is of unstable type and $a \in U_0$. 

Because of the presence of unstable components and special triples in $\mathcal{I}_{f_\star}$, this set is not necessarily an interval for the preorder $\unlhd$, unlike in~\cite{BI08}. To remedy this, we define an equivalence relation $\sim_T$ on $\mathcal{I}_{f_\star}$. We first define a homogeneous relation $\sim^\mathrm{pre}_T$ as follows: if $\ell_0, \ell_1 \in \mathcal{I}_{f_\star}$, then $\ell_0 \sim^\mathrm{pre}_T \ell_1$ if and only if $\ell_0 = \ell_1$, or $\ell_0$ and $\ell_1$ belong to the same special triple. We then define $\sim_T$ as the transitive closure of $\sim^\mathrm{pre}_T$. Notice that an $X$-line $\ell$ belongs to at most two distinct special triples, in which case $\ell$ coincides with the positive or negative unstable branch of one and with the stable branch of the other. The \textbf{reduction} of $\mathcal{I}_{f_\star}$ is defined by 
$$\mathcal{I}^\mathrm{red}_{f_\star} \coloneqq \mathcal{I}_{f_\star} \slash \sim_T.$$
By Lemma~\ref{lem:propT}, the partial order $\unlhd$ naturally induces a partial order $\unlhd^\mathrm{red}$ on $\mathcal{I}^\mathrm{red}_{f_\star}$, and we have:

\begin{lem} \label{lem:embed}
The poset $\big(\mathcal{I}^\mathrm{red}_{f_\star}, \unlhd^\mathrm{red}\big)$ embeds into $([0,1], \leq)$ as a closed subset.
\end{lem}

\begin{proof}
Let $\gamma : [0,1] \rightarrow U$ be a forward curve from $I_0$ to $a$, and let $C \coloneqq \gamma^{-1}(U^\mathrm{pure})$. Then $C$ is a closed subset of $[0,1]$ in bijection with the set of pure $X$-lines intersected by $\gamma$. Let $\lambda: C \hookrightarrow \mathcal{I}_{f_\star}$ be the increasing map sending $t \in C$ to the $X$-line passing through $\gamma(t)$. It induces a nondecreasing map $\lambda^\mathrm{red} :  C \rightarrow \mathcal{I}^\mathrm{red}_{f_\star}$ which is surjective by Lemma~\ref{lem:forcurv2}. We define an equivalence relation $\sim$ on $C$ as follows: if $t, t' \in C$, then $t \sim t'$ if and only if $\lambda^\mathrm{red}(t) = \lambda^\mathrm{red}(t')$. One easily checks that if $t<t'$ and $t \sim t'$, then $\{t,t'\} \subset \partial C$ and $(t,t') \cap C$ is either empty or a single point, depending on whether $\lambda(t)$ and $\lambda(t')$ are negative and positive branches of a special triple or of a double special triple, respectively. The order $\leq$ induces an order $\leq^\mathrm{red}$ on the quotient $C^\mathrm{red} \coloneqq C \slash \sim$. Since $\lambda^\mathrm{red}$ induces an isomorphism between $\big(C^\mathrm{red}, {\leq}^\mathrm{red}\big)$ and $\big(\mathcal{I}^\mathrm{red}_{f_\star}, \unlhd^\mathrm{red}\big)$, it is enough to show that $\big(C^\mathrm{red}, {\leq}^\mathrm{red} \big)$ embeds into $([0,1], \leq)$ as a closed subset. 

Let $\mathcal{U}$ be the collection of connected components of $[0,1] \setminus C$ of the form $(t,t')$ with $t \sim t'$. We consider the closed set $$\widehat{C} \coloneqq [0,1] \setminus \bigcup_{U \in \mathcal{U}} U \subset [0,1].$$
Notice that $C \subset \widehat{C}$ is closed, and $\sim$ extends to $\widehat{C}$ so that ${\leq}^\mathrm{red}$ extends to a total order on the quotient $\widehat{C}^\mathrm{red} \coloneqq \widehat{C}\slash \sim$, and the induced map $C^\mathrm{red} \rightarrow \widehat{C}^\mathrm{red}$ is (strictly) order preserving. Moreover, the image of the latter map is closed for the order topology. We can loosely think of quotienting $\widehat{C}$ by $\sim$ as `collapsing the gaps' between the connected components of $\widehat{C}$. Recall the following characterization due to Cantor (see~\cite[Theorem 4.3]{J03} for the unbounded case): every totally ordered set which is bounded, dense, separable, and complete is either a singleton or isomorphic to $([0,1], \leq)$. It is straightforward to check that $\big(\widehat{C}^\mathrm{red}, {\leq}^\mathrm{red} \big)$ satisfies these properties, hence embeds into $([0,1], \leq)$ as a closed subset. The composition $C^\mathrm{red} \rightarrow \widehat{C}^\mathrm{red} \rightarrow [0,1]$ yields the desired embedding.
\end{proof}

A priori, the embedding of Lemma~\ref{lem:embed} depends on the choice of forward curve. However, one can easily show that a different forward curve induces the same embedding up to post-composition with an increasing homeomorphism of $[0,1]$. Alternatively, one can appeal to the following elementary fact: if a poset $(I, \preceq)$ embeds into $([0,1], \leq)$ as a closed subset, then this embedding is continuous for the order topology, and it is unique up to post-composition with an increasing homeomorphism of $[0,1]$.\footnote{If $\iota: I \rightarrow [0,1]$ is such an embedding with closed image $C$, then $\iota$ is open so $\iota^{-1} : C \rightarrow I$ is continuous. Since $C$ is compact and the order topology on $I$ is Hausdorff, $\iota^{-1}$ is a homeomorphism, hence $\iota$ is continuous. If $\iota' : I \rightarrow [0,1]$ is another embedding with closed image $C'$, we simply extend $\iota' \circ \iota^{-1} : C \rightarrow [0,1]$ to an increasing homeomorphism $\phi : [0,1] \rightarrow [0,1]$ such that $\phi$ is affine on the connected components of $[0,1] \setminus C$, and $\iota' = \phi \circ \iota$.}

We now fix an embedding as in the previous lemma, and we define $I_{f_\star} \subset [0,1]$ to be the image of $\mathcal{I}^\mathrm{red}_{f_\star}$ under this embedding. We denote by $\mathcal{D}_{f_\star}$ the domain in $U$ defined as the union of all of the $X$-lines in $\mathcal{I}_{f_\star}$, together with all of the unstable components of $f$ crossed by the forward curves from $I_0$ to $a$. We call $\mathcal{D}_{f_\star}$ the \textbf{backward domain} of $f_\star$.\footnote{Beware: this domain has a more complicated topology than the domain 
$D_X$ from~\cite[Section 6.2.1]{BI08}. Because of the way we defined it, $\mathcal{D}_{f_\star}$ might not be simply connected as it can have punctures. This will not be relevant for our arguments, and these missing points can be added for psychological comfort.}

%%%
                \subsubsection{Overlap of two envelopes}
%%%

Let us consider a second marked forward envelope $g_\star = (g : V \rightarrow M, b) \in \mathscr{U}_p$ with corresponding poset $\mathcal{I}_{g_\star} \subset \mathcal{L}_g$ and backward domain $\mathcal{D}_{g_\star} \subset V$. If $\gamma : [0,1] \rightarrow U$ is a \emph{backward} curve starting from $a$, a \textbf{lift} of $\gamma$ to $g_\star$ is a curve $\widetilde{\gamma} : [0,1] \rightarrow V$ starting from $b$ satisfying $f \circ \gamma = g \circ \widetilde{\gamma}$. If it exists, such a lift is unique and is a backward curve in $V$. We define $\mathcal{I}_{f_\star \cap g_\star} \subset \mathcal{I}_{f_\star}$ to be the set of all $X$-lines in $\mathcal{L}_{f_\star}$ intersected by some backward curve starting at $a$ which lifts to $g_\star$. If $g$ is of pure type and $b$ belongs to the initial edge $e'_0$ of $g$, we set $\mathcal{I} \coloneqq \{e'_0\}$ instead. We define the \textbf{backward domain of overlap} $\mathcal{D}_{f_\star \cap g_\star} \subset \mathcal{D}_{f_\star}$ to be the union of all of the $X$-lines in $\mathcal{I}_{f_\star \cap g_\star}$ together with all of the unstable components of $f$ crossed by the aforementioned backward curves lifting to $g_\star$. Notice that $\mathcal{I}_{f_\star \cap g_\star}$ is an \emph{upward section} of $\mathcal{I}_{f_\star}$, in the following sense: for all $\ell_0 \in \mathcal{I}_{f_\star \cap g_\star}$ and $\ell_1 \in \mathcal{I}_{f_\star}$, if $\ell_0 \unlhd \ell_1$ then $\ell_1 \in \mathcal{I}_{f_\star \cap g_\star}$. The same is true for the reduction
$$\mathcal{I}^\mathrm{red}_{f_\star \cap g_\star} \coloneqq \mathcal{I}_{f_\star \cap g_\star} \slash \sim_T, \qquad  \mathcal{I}^\mathrm{red}_{f_\star \cap g_\star} \subset \mathcal{I}^\mathrm{red}_{f_\star}.$$
Let $I_{f_\star \cap g_\star}\subset I_{f_\star}$ denote the image of $\mathcal{I}^\mathrm{red}_{f_\star \cap g_\star}$ under the chosen embedding from Lemma~\ref{lem:embed}. Notice that $I_{f_\star \cap g_\star}$ can be empty, for instance if $g$ is of unstable type and $b$ belongs to the initial unstable component of $g$. We refer to the situation $I_{f_\star \cap g_\star} = \varnothing$ as case (I0). Otherwise, if $I_{f_\star \cap g_\star} \neq \varnothing$, we define $$m \coloneqq \inf I_{f_\star \cap g_\star} \in I_{f_\star},$$ and we distinguish several cases: 
\begin{enumerate}
    \item[(I1)] $m \in I_{f_\star \cap g_\star}$, and $I_{f_\star \cap g_\star} = [m, 1] \cap I_{f_\star}$. Then $\mathcal{I}_{f_\star \cap g_\star}$ is of one of the following three types:
        \begin{enumerate}
            \item[(i)] $\mathcal{I}_{f_\star \cap g_\star} = \big\{ \ell \in \mathcal{I}_{f_\star} \ \vert \ \ell_m \unlhd \ell \big\}$, where $\ell_m \in \mathcal{I}_{f_\star}$. Here, $\ell_m$ is the least element of $\mathcal{I}_{f_\star \cap g_\star}$. If $\ell_m$ is the negative unstable branch of a special triple in $\mathcal{I}_{f_\star}$, then it corresponds to the initial edge of $g_\star$.\footnote{Recall that we exclude the case where $\ell_m$ is the stable branch of a quadratic-like singularity and corresponds to the initial edge of $g_\star$.}
            \item[(ii)] $\mathcal{I}_{f_\star \cap g_\star} = \{\ell_s\} \cup \big\{ \ell \in \mathcal{I}_{f_\star} \ \vert \ \ell_u^- \unlhd \ell \big\}$, where $(\ell_s, \ell_u^-, \ell_u^+)$ is a special triple in $\mathcal{I}_{f_\star}$. Here, $\ell_s$ and $\ell^-_u$ are the two minimal elements of $\mathcal{I}_{f_\star \cap g_\star}$.
            \item[(iii)] $\mathcal{I}_{f_\star \cap g_\star} = \{\ell_s, \ell'_s\} \cup \big\{ \ell \in \mathcal{I}_{f_\star} \ \vert \ \ell_u'^- \unlhd \ell \big\}$, where $(\ell_s, \ell_u^-, \ell_u^+)$ and $(\ell'_s, \ell'^-_u, \ell'^+_u)$ are a special triples in $\mathcal{I}_{f_\star}$ with $\ell'_s = \ell^-_u$. Here, $\ell_s$, $\ell'_s$ and $\ell'^-_u$ are the three minimal elements of $\mathcal{I}_{f_\star \cap g_\star}$.
        \end{enumerate}
    In these three cases, the minimal elements of $\mathcal{I}_{f_\star}$ might be forward edges of an unstable component of $U$ contained in $\mathcal{D}_{f_\star \cap g_\star}$, except in case (i) if $\ell_m$ corresponds to the initial edge of $g_\star$.
    \item[(I2)] $m \notin I_{f_\star \cap g_\star}$, and $I_{f_\star \cap g_\star} = (m, 1] \cap I_{f_\star}$. In that case, there exists a (possibly non unique) $\ell_m \in \mathcal{I}_{f_\star}$ such that $\mathcal{I}_{f_\star \cap g_\star} = \big\{ \ell \in \mathcal{I}_{f_\star} \ \vert \ \ell_m \lhd \ell \big\}$. We can further distinguish three cases:
    \begin{enumerate}
            \item[(i)] $\ell_m$ does not belong to a special triple in $\mathcal{I}_{f_\star}$, and $\ell_m$ is the infimum of $\mathcal{I}_{f_\star \cap g_\star}$ in $\mathcal{I}_{f_\star}$,
            \item[(ii)] There exists a special triple $(\ell_s, \ell^-_u, \ell^+_u)$ such that $\ell^+_u$ is not the stable branch of another special triple, and $\ell_m \in \{\ell_s, \ell^+_u\}$. In that case, $\mathcal{I}_{f_\star \cap g_\star}$ has exactly two maximal lower bounds in $\mathcal{I}_{f_\star}$, namely $\ell_s$ and $\ell^+_u$.
            \item[(iii)] There exists a special triple $(\ell_s, \ell^-_u, \ell^+_u)$ such that $\ell^+_u$ is the stable branch of another special triple $(\ell'_s, \ell'^-_u, \ell'^+_u)$, and $\ell_m \in \{ \ell_s, \ell'_s, \ell^+_u\}$. In that case, $\mathcal{I}_{f_\star \cap g_\star}$ has exactly three maximal lower bounds in $\mathcal{I}_{f_\star}$, namely $\ell_s$, $\ell'_s$, and $\ell'^+_u$.
        \end{enumerate}
\end{enumerate}

Notice that the case $\mathcal{I}_{f_\star \cap g_\star} = \{\ell_s\} \cup \big\{ \ell \in \mathcal{I}_{f_\star} \ \vert \ \ell_u^+ \unlhd \ell \big\}$, where $(\ell_s, \ell_u^-, \ell_u^+)$ is a special triple with $\ell_s, \ell_u^+ \in \mathcal{I}_{f_\star}$, cannot occur. Indeed, if both $\ell_s$ and $\ell_u^+$ are in $\mathcal{I}_{f_\star \cap g_\star}$, then the unstable locus in $U$ corresponding to this special triple is in $\mathcal{D}_{f_\star \cap g_\star}$, since forward envelopes have at most one proper edge,\footnote{We also implicitly use the unique integrability properties of $\eta$ established in the proof of Lemma~\ref{lem:uniqint}.} and it follows that $\ell_u^-$ is also in $\mathcal{I}_{f_\star \cap g_\star}$. Other possible cases can be ruled out by similar arguments. This discussion implies:

\begin{lem} \label{lem:casesI}
Exactly one of the following holds:
\begin{enumerate}
    \item[(I0)] $\mathcal{I}_{f_\star \cap g_\star} = \varnothing$.
    \item[(I1)] $\mathcal{I}_{f_\star \cap g_\star}$ has one, two, or three minimal elements.
    \item[(I2)] $\mathcal{I}_{f_\star \cap g_\star}$ has no minimal element but has one, two, or three maximal lower bounds in $\mathcal{I}_{f_\star}$.
\end{enumerate}
Moreover, a unique minimal element of $\mathcal{I}_{f_\star \cap g_\star}$ is a least element, and a unique maximal lower bound of $\mathcal{I}_{f_\star \cap g_\star}$ in $\mathcal{I}_{f_\star}$ is an infimum.
\end{lem}

\begin{figure}[h]
    \centering
        \begin{subfigure}{0.2\textwidth}
                \includegraphics{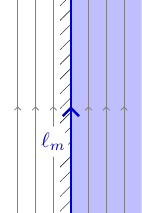}
        \centering
        \caption{Case (I1) (i).}
        \end{subfigure}%
        \begin{subfigure}{0.3\textwidth}
                \includegraphics{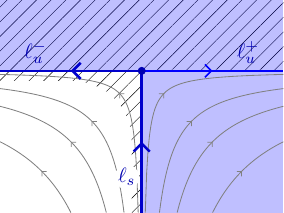}
        \centering
        \caption{Case (I1) (ii).}
        \end{subfigure}%
        \begin{subfigure}{0.5\textwidth}
                \includegraphics{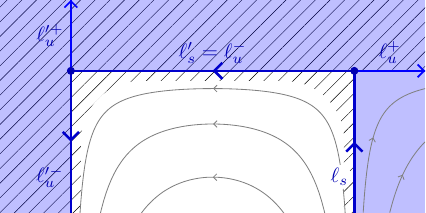}
        \centering
        \caption{Case (I1) (iii).}
        \end{subfigure}
        \par\bigskip
        \begin{subfigure}{0.2\textwidth}
                \includegraphics{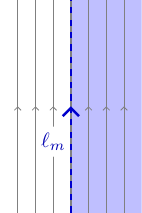}
        \centering
        \caption{Case (I2) (i).}
        \end{subfigure}%
        \begin{subfigure}{0.3\textwidth}
                \includegraphics{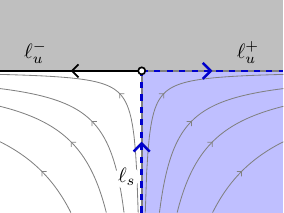}
        \centering
        \caption{Case (I2) (ii).}
        \end{subfigure}%
        \begin{subfigure}{0.5\textwidth}
                \includegraphics{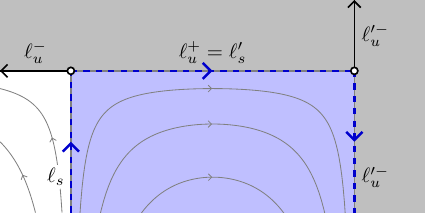}
        \centering
        \caption{Case (I2) (iii).}
        \end{subfigure}
    \caption{Double overlaps. In each picture, the whole rectangle corresponds to $\mathcal{D}_{f_\star}$, and the blue region corresponds to $\mathcal{D}_{f_\star \cap g_\star}$. The blue hatched regions are possible unstable regions contained in $\mathcal{D}_{f_\star \cap g_\star}$. The gray hatching to the left of thick blue lines are `shadows' of $g_\star$. Thick blue lines belong to $\mathcal{D}_{f_\star}$, while dashed blue lines do not. The gray regions are unstable regions in $\mathcal{D}_{f_\star}$.}
    \label{fig:overlap}
\end{figure}

Figure~\ref{fig:overlap} summarizes all the possibilities for these cases, as well as the corresponding domains $\mathcal{D}_{f_\star \cap g_\star}$.

\medskip

There is a natural (strictly) order preserving map $\alpha_{f_\star \cap g_\star} : \mathcal{I}_{f_\star \cap g_\star} \rightarrow \mathcal{I}_{g_\star}$ defined as follows. If $p \in M^\mathrm{pure}$ and $\ell$ is the $X$-line in $U$ passing through $a$, we define $\alpha_{f_\star \cap g_\star}(\ell)$ as the $X$-line in $V$ passing through $b$. Otherwise, if $\ell \in \mathcal{I}_{f_\star \cap g_\star}$ and $\gamma : [0,1] \rightarrow U$ is a forward curve from $\ell$ to $a$ which admits a lift $\widetilde{\gamma} : [0,1] \rightarrow V$ to $g$ with $\widetilde{\gamma}(1)=b$, then $\alpha_{f_\star \cap g_\star}(\ell)$ is the $X$-line in $V$ passing through $\widetilde{\gamma}(0)$. By Lemma~\ref{lem:forlifts}, this definition does not depend on the choice of the forward curve $\gamma$. By construction, $\alpha_{f_\star \cap g_\star}$ is strictly order preserving. 

There is also a natural map $\varphi_{f_\star \cap g_\star} : \mathcal{D}_{f_\star \cap g_\star} \rightarrow \mathcal{D}_{g_\star}$ satisfying 
\begin{align} \label{eq:phifg}
g \circ \varphi_{f_\star \cap g_\star} = f_{\vert \mathcal{D}_{f_\star \cap g_\star}},
\end{align}
defined as follows. If $p \in M^\mathrm{pure}$ and $x \in \mathcal{D}_{f_\star \cap g_\star}$ belongs to the $X$-line passing through $a$, then there exists a unique $t \in \R$ such that $x = \phi_{f^*X}^t(a)$ and $\varphi_{f_\star \cap g_\star}(x) \coloneqq \phi_{{g^*}X}^t(b)$. Otherwise, if $x \in \mathcal{D}_{f_\star \cap g_\star}$ and $\gamma : [0,1] \rightarrow U$ is a forward curve from $x$ to $a$ and $\widetilde{\gamma} : [0,1] \rightarrow V$ is a lift of $\gamma$ to $g$ ending at $b$, then we define $\varphi_{f_\star \cap g_\star}(x) \coloneqq \widetilde{\gamma}(0)$. Lemma~\ref{lem:forlifts} readily implies that $\varphi_{f_\star \cap g_\star}$ is well-defined and independent of the choice of $\gamma$. It satisfies the equality~\eqref{eq:phifg} by construction.

Exchanging $g_\star$ and $f_\star$, we similarly obtain a poset $\mathcal{I}_{g_\star \cap f_\star} \subset \mathcal{I}_{g_\star}$ and a domain $\mathcal{D}_{g_\star \cap f_\star} \subset \mathcal{D}_{g_\star}$, as well as maps $\alpha_{g_\star \cap f_\star} : \mathcal{I}_{g_\star \cap f_\star} \rightarrow \mathcal{I}_{g_\star}$ and $\varphi_{g_\star \cap f_\star} : \mathcal{D}_{g_\star \cap f_\star} \rightarrow \mathcal{D}_{g_\star}$. By definition, the image of $\varphi_{f_\star \cap g_\star}$ is contained in $\mathcal{D}_{g_\star \cap f_\star}$, and $\varphi_{g_\star \cap f_\star} \circ \varphi_{f_\star \cap g_\star} = \mathrm{id}_{\mathcal{D}_{f_\star \cap g_\star}}$. Similarly, $\varphi_{f_\star \cap g_\star} \circ \varphi_{g_\star \cap f_\star} = \mathrm{id}_{\mathcal{D}_{g_\star \cap f_\star}}$, hence $\varphi_{f_\star \cap g_\star}$ is a diffeomorphism between $\mathcal{D}_{f_\star \cap g_\star}$ and $\mathcal{D}_{g_\star \cap f_\star}$. Likewise, $\alpha_{f_\star \cap g_\star}$ induces an order preserving bijection between $\mathcal{I}_{f_\star \cap g_\star}$ and $\mathcal{I}_{g_\star \cap f_\star}$ whose inverse is given by $\alpha_{g_\star \cap f_\star}$.

Notice that if $\mathcal{I}_{f_\star \cap g_\star} = \varnothing$, then $p \in M^\mathrm{un}$ and $U_c = \mathcal{D}_{f_\star \cap g_\star} \subset U$ is an unstable component of $f$. Similarly, if $\mathcal{I}_{f_\star \cap g_\star}$ has a least element $\ell_m$ which is not a forward boundary component of $\mathcal{D}_{f_\star \cap g_\star}$ in $\mathcal{D}_{f_\star}$, then there exists an unstable component $U_c \subset \mathcal{D}_{f_\star \cap g_\star}$ such that $\ell_m$ is a backward boundary component of $U_c$. In both cases, we say that $U_c$ is the \textbf{initial unstable component of overlap} of $f_\star$ and $g_\star$.

%%%
                \subsubsection{Initial corners} \label{sec:initcor}
%%%

Let $q \in \Delta_\mathrm{sa} \cup Q$. Recall that $q$ is a \textbf{pure} saddle singularity if no stable branch of $X$ at $q$ is in $M^\mathrm{un}$, and $q$ is a \textbf{quadratic-like} singularity if exactly one stable branch of $X$ at $q$ is contained in $M^\mathrm{pure}$. The set of pure saddle singularities of $X$ is denoted by $\Delta^\mathrm{pure}_\mathrm{sa}$, and the set of quadratic-like singularities of $X$ is denoted by $Q^*$.

If $q \in \Delta^\mathrm{pure}_\mathrm{sa} \cup Q^*$, a \textbf{corner at $q$} is a pair $c = (\gamma_s, \gamma_u)$, where $\gamma_s$ and $\gamma_u$ are a stable and unstable branches of $X$ at $q$, respectively, both contained in $M^\mathrm{pure}$. We denote by $\mathfrak{C}_q$ the set of corners at $q$, and $\mathfrak{C}$ the set of all corners. Note that $\mathfrak{C}_q$ has four elements if $q \in \Delta^\mathrm{pure}_\mathrm{sa}$, and two elements if $q \in Q^*$.

A corner $c = (\gamma_s, \gamma_u) \in \mathfrak{C}$ is \textbf{forward} (resp.~\textbf{backward}) if $\gamma_u$ lies to the right (resp.~to the left) of $\gamma_s$. To make sense of this orientation, we consider the (unit) vector $v_s, v_u \in \eta(q)$ defined as
\begin{align*}
v_s \coloneqq \lim_{t \rightarrow +\infty} \frac{\dot{\gamma}_s(t)}{\vert \dot{\gamma}_s(t)\vert}, & & v_u \coloneqq \lim_{t \rightarrow -\infty} \frac{\dot{\gamma}_u(t)}{\vert \dot{\gamma}_u(t)\vert}.
\end{align*}
If $(v_s, v_u)$ is positively (resp.~negatively) oriented, with respect to the orientation on $\eta(q)$, then $\gamma_u$ is to the left (resp.~to the right) of $\gamma_s$. We denote by $\mathfrak{C}^+_q \subset \mathfrak{C}_q$ the set of forward corners at $q$, by $\mathfrak{C}^+ \subset \mathfrak{C}$ the set of all forward corners. Note that $\mathfrak{C}^+_q$ has two elements if $q \in \Delta^\mathrm{pure}_\mathrm{sa}$, and one element if $q \in Q^*$. Moreover, there is a canonical `projection' map
$$\mathfrak{p} : \mathfrak{C}^+ \rightarrow \Delta^\mathrm{pure}_\mathrm{sa} \cup Q^*.$$

\begin{figure}[h]
    \centering
        \begin{subfigure}[c]{0.5\textwidth}
            \raisebox{-\height}{\includegraphics{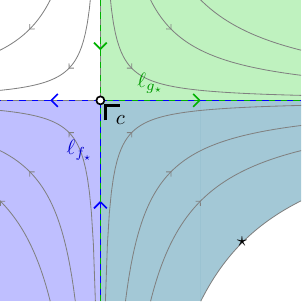}
                }
        \centering
        \caption{Case (C1) (i).}
        \end{subfigure}%
        \begin{subfigure}[c]{0.5\textwidth}
            \raisebox{-\height}{\includegraphics{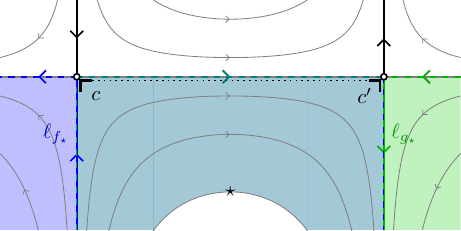}
            }%
            \vspace{2ex}
            \raisebox{-\height}{\includegraphics{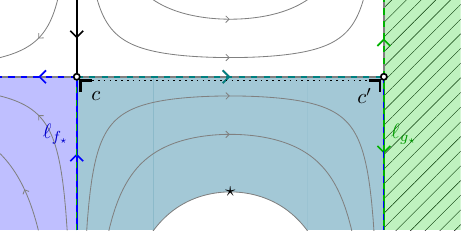}   
            }
        \centering
        \caption{Case (C1) (ii). The right stable branch at $c'$ might be in $M^\mathrm{pure}$ (top) or in $M^\mathrm{un}$ (bottom).}
        \label{fig:corC1iiun}
        \end{subfigure}
        \par\bigskip
        \begin{subfigure}{0.5\textwidth}
                    \includegraphics{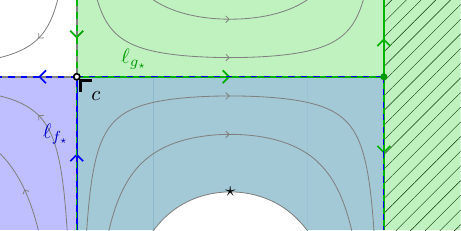}
        \centering
        \caption{Case (C2)(i).}
        \label{fig:corC2i}
        \end{subfigure}%
        \begin{subfigure}{0.5\textwidth}
                    \includegraphics{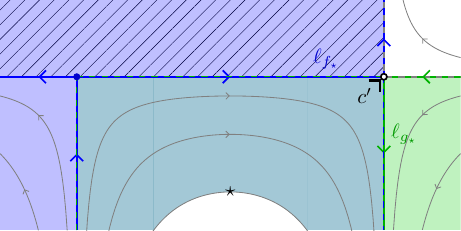}
        \centering
        \caption{Case (C2)(ii).}
        \end{subfigure}
    \caption{Possible cases for initial (double) corners of overlap. The regions in blue correspond to $\mathcal{D}_{f_\star}$, the regions in green correspond to $\mathcal{D}_{g_\star}$, and their region of overlap is in teal. The thinly hatched regions are in $M^\mathrm{un}$.}
    \label{fig:corners}
\end{figure}

We now refine case (I2) above in the description of the backward overlap of $f_\star$ and $g_\star$. We assume that $\mathcal{I}_{f_\star \cap g_\star}$ (hence $\mathcal{I}_{g_\star \cap f_\star}$) satisfies case (I2) of Lemma~\ref{lem:casesI}. We say that a corner $c \in \mathfrak{C}^+$ is the \textbf{initial (forward) corner of overlap} of $f_\star$ and $g_\star$ if, after possibly switching $f_\star$ and $g_\star$, there exists a maximal lower bounds $\ell_s \in \mathcal{I}_{f_\star}$ of $\mathcal{I}_{f_\star \cap g_\star}$ and a maximal lower bound $\ell_u \in \mathcal{I}_{g_\star}$ of $\mathcal{I}_{g_\star \cap f_\star}$ such that $c = \big(f \circ \ell_s, g \circ \ell_u\big)$. A more subtle situation can occur. If $c=(\gamma_s, \gamma_u) \in \mathfrak{C}^+_q$ and $c'=(\gamma'_s, \gamma'_u) \in \mathfrak{C}^+_{q'}$ are two forward corners such that $\gamma_u = \gamma'_s \eqqcolon \gamma_0$, in which case $q$ and $q'$ are saddle singularities connected by $\gamma_u$, we say that $(c,c')$ is the \textbf{initial (forward) double corner of overlap} of $f_\star$ and $g_\star$ if, after possibly switching $f_\star$ and $g_\star$, $\mathcal{I}_{f_\star \cap g_\star} \subset \mathcal{I}_{f_\star}$ has an infimum $\ell_s$, $\mathcal{I}_{g_\star \cap f_\star} \subset \mathcal{I}_{g_\star}$ has an infimum $\ell'_u$, and $\gamma_s = f \circ \ell_s$ and $\gamma'_u = g \circ \ell'_u$. Then, after possibly switching $f_\star$ and $g_\star$, one of the following cases holds:

\begin{enumerate}
    \item[(C1)] Both $\mathcal{I}_{f_\star \cap g_\star} \subset \mathcal{I}_{f_\star}$ and $\mathcal{I}_{g_\star \cap f_\star} \subset \mathcal{I}_{g_\star}$ satisfy case (I2) (i) above, and have respective infima $\ell_{f_\star}$ and $\ell_{g_\star}$. Then, either
        \begin{enumerate}
            \item[(i)] $f_\star$ and $g_\star$ have an initial corner of overlap $c = \big(f \circ \ell_{f_\star}, g \circ \ell_{g_\star}\big)$, or
            \item[(ii)] $f_\star$ and $g_\star$ have an initial double corner of overlap $(c,c')$, where $c = \big(f \circ \ell_{f_\star}, \gamma\big)$ and $c'= \big(\gamma, g \circ \ell_{g_\star}\big)$ for some saddle-saddle connection $\gamma$ of $X$.
        \end{enumerate}
    \item[(C2)] Exactly one of $\mathcal{I}_{f_\star \cap g_\star}$ or $\mathcal{I}_{g_\star \cap f_\star}$ satisfies case (I2) (i) while the other satisfies case (I2) (ii). More precisely, 
        \begin{enumerate}
        \item[(i)] $\mathcal{I}_{f_\star \cap g_\star}$ has an infimum $\ell_s$ in $\mathcal{I}_{f_\star}$, while $\mathcal{I}_{g_\star \cap f_\star}$ has exactly two maximal lower bounds $\ell'_s$ and $\ell'^+_u$ in $\mathcal{I}_{g_\star}$, such that $c= \big(f \circ \ell_s, g \circ \ell'_s \big)$ and $c' = \big(g \circ \ell'_s, g \circ \ell'^+_u\big)$ form a double corner. In that case, $c$ is the initial corner of overlap of $f_\star$ and $g_\star$.
        \item[(ii)] $\mathcal{I}_{f_\star \cap g_\star}$ has exactly two maximal lower bounds $\ell_s$ and $\ell^+_u$ in $\mathcal{I}_{f_\star}$, while $\mathcal{I}_{g_\star \cap f_\star}$ has an infimum $\ell'^+_u$ in $\mathcal{I}_{g_\star}$, such that $c= \big(f \circ \ell_s, f \circ \ell^+_u \big)$ and $c' = \big(f \circ \ell^+_u, g\circ\ell'^+_u\big)$ form a double corner. In that case, $c'$ is the initial corner of overlap of $f_\star$ and $g_\star$.
        \end{enumerate}
\end{enumerate}
See Figure~\ref{fig:corners} for illustrations of all these cases.

\medskip

We now describe how the shadows of $f_\star$ and $g_\star$ interact at an initial corner. This will be relevant in the definition of the extended preorder on $\mathscr{U}$. We begin with some definition.

If $c = (\gamma_s, \gamma_u) \in \mathfrak{C}^+_q$ is a forward corner at $q$, then for any points $q_s \in \gamma_s$ and $q_u \in \gamma_u$, there exist well-defined order preserving map $\varsigma_s : \mathscr{A}_{q_s} \hookrightarrow \mathscr{A}_q$ and $\varsigma_u : \mathscr{A}_{q_u} \hookrightarrow \mathscr{A}_q$ obtained by sliding the marked point along the stable and unstable branch at $c$, respectively, to the corner point $q$. Here, we use that $\mathscr{A}$ is complete. Similar maps already appeared in the proof of Lemma~\ref{lem:exists}. If $c$ is the initial corner of overlap of $f_\star$ and $g_\star$, and $\ell_s \in \mathcal{I}_{f_\star}$ and $\ell_u \in \mathcal{I}_{g_\star}$ are the corresponding `entrance $X$-lines' of $f_\star$ and $g_\star$ along $c$, respectively, we choose arbitrary points $x_s \in \ell_s$ and $x_u \in \ell_u$ so that $q_s \coloneqq f(x_s) \in \gamma_s$ and $q_u \coloneqq g(x_u) \in \gamma_u$. We then set
\begin{align*}
    \Sigma(c) &\coloneqq \varsigma_s\big(\Sigma(x_s) \big) \subset \mathscr{A}_q, \\
    \Sigma'(c) &\coloneqq \varsigma_u\big(\Sigma'(x_u) \big) \subset \mathscr{A}_q.
\end{align*}
Note that these sets are independent of the points $x_s$ and $x_u$ by axiom 4 of Definition~\ref{def:upenv}. We then define $\widehat{\Sigma}(c) \subset \mathscr{A}_q$ and $\widehat{\Sigma}'(c)\subset \mathscr{A}_q$ as the smallest sections of $\mathscr{A}_q$ (with respect to the inclusion) containing $\Sigma(c)$ and $\Sigma'(c)$, respectively. Since those are two sections, and the preorder on $\mathscr{A}_q$ is total, one is included in the other. If $\widehat{\Sigma}(c) = \widehat{\Sigma}'(c)$, we say that $f_\star$ and $g_\star$ are \textbf{geometrically coherent at $c$}. Otherwise, we say that they are \textbf{geometrically noncoherent}.

We extend this definition to a double corner $(c,c')$. In the above notations, we consider points $q_s \in \gamma_s$, $q'_u \in \gamma'_u$, and $q_0 \in \gamma_0$. Sliding the marked points as before induces order preserving maps organized as follows:
$$\begin{tikzcd}[column sep=scriptsize, row sep=scriptsize]
	{\mathscr{A}_{q}} & {\mathscr{A}_{q_0}} & {\mathscr{A}_{q'}} \\
	{\mathscr{A}_{q_s}} && {\mathscr{A}_{q'_u}}
	\arrow["{\varsigma'_s}", hook, from=1-2, to=1-3]
	\arrow["{\varsigma_u}"', hook', from=1-2, to=1-1]
	\arrow["{\varsigma_s}", hook, from=2-1, to=1-1]
	\arrow["{\varsigma'_u}"', hook, from=2-3, to=1-3]
\end{tikzcd}$$
We then consider $x_s \in \ell_s$, $x'_u \in \ell'_u$, and we define
\begin{align*}
    \Sigma(c) &\coloneqq \varsigma_s\big(\Sigma(x_s) \big) \subset \mathscr{A}_q, \\
    \Sigma'(c') &\coloneqq \varsigma'_u\big(\Sigma'(x'_u) \big) \subset \mathscr{A}_{q'}.
\end{align*}
As before, these sets do not depend on the choices of $x_s$ and $x'_u$. Denoting $\widehat{\Sigma}(c) \subset \mathscr{A}_q$ and $\widehat{\Sigma}'(c') \subset \mathscr{A}_{q'}$ the smallest sections containing $\Sigma(c)$ and $\Sigma'(c')$, respectively, we define
\begin{align*}
    \Sigma_0(c) &\coloneqq \varsigma^{-1}_u\big( \widehat{\Sigma}(c) \big) \subset \mathscr{A}_{q_0}, \\
    \Sigma'_0(c') &\coloneqq {(\varsigma'_s)}^{-1}\big( \widehat{\Sigma}'(c') \big) \subset \mathscr{A}_{q_0},
\end{align*}
and we denote by $\widehat{\Sigma}_0(c)$ and $\widehat{\Sigma}'_0(c')$ the smallest sections of $\mathscr{A}_{q_0}$ containing each of them. If $\widehat{\Sigma}_0(c) = \widehat{\Sigma}'_0(c')$, we say that $f_\star$ and $g_\star$ are \textbf{geometrically coherent at $(c,c')$}. Otherwise, we say that they are \textbf{geometrically noncoherent}.

%%%
                \subsubsection{Triple overlaps} \label{sec:triple}
%%%
We finally consider a third marked forward envelope $h_\star = (h : W \rightarrow M, c) \in \mathscr{U}_p$ with corresponding poset $\mathcal{I}_{h_\star} \subset \mathcal{L}_h$ and backward domain $\mathcal{D}_{h_\star} \subset W$. We define the `triple intersections' by:
\begin{align*}
    \mathcal{I}_{f_\star \cap g_\star \cap h_\star} &\coloneqq \mathcal{I}_{f_\star \cap g_\star} \cap \mathcal{I}_{f_\star \cap h_\star} \subset \mathcal{I}_{f_\star}, \\
    \mathcal{D}_{f_\star \cap g_\star \cap h_\star} &\coloneqq \mathcal{D}_{f_\star \cap g_\star} \cap \mathcal{D}_{f_\star \cap h_\star} \subset \mathcal{D}_{f_\star}.
\end{align*}

Lemma~\ref{lem:forlifts} easily implies

\begin{lem}[see \cite{BI08}, Lemma 6.12] \label{lem:ident} The following identities hold:
\begin{align*}
    \alpha_{f_\star \cap g_\star} \big( \mathcal{I}_{f_\star \cap g_\star \cap h_\star} \big) &= \mathcal{I}_{g_\star \cap f_\star \cap h_\star}, & \varphi_{f_\star \cap g_\star} \big( \mathcal{D}_{f_\star \cap g_\star \cap h_\star} \big) &= \mathcal{D}_{g_\star \cap f_\star \cap h_\star}, \\
    \alpha_{f_\star \cap g_\star}^{-1} \big(\mathcal{I}_{g_\star \cap h_\star}\big) &=\mathcal{I}_{f_\star \cap g_\star \cap h_\star}, & \varphi_{f_\star \cap g_\star}^{-1} \big(\mathcal{D}_{g_\star \cap h_\star}\big) &= \mathcal{D}_{f_\star \cap g_\star \cap h_\star}, \\
    \alpha_{g_\star \cap h_\star} \circ \alpha_{f_\star \cap g_\star}(\ell) &= \alpha_{f_\star \cap h_\star}(\ell), & \varphi_{g_\star \cap h_\star} \circ \varphi_{f_\star \cap g_\star}(x) &= \varphi_{f_\star \cap h_\star}(x),
\end{align*}
for every $\ell \in \mathcal{I}_{f_\star \cap g_\star \cap h_\star}$ and $x \in \mathcal{D}_{f_\star \cap g_\star \cap h_\star}$.
\end{lem}

Notice that by Lemma~\ref{lem:embed}, either $I_{f_\star \cap g_\star \cap h_\star} = I_{f_\star \cap g_\star}$ or $I_{f_\star \cap g_\star \cap h_\star} = I_{f_\star \cap h_\star}$ or equivalently, either $I_{f_\star \cap g_\star } \subset I_{f_\star \cap h_\star}$ or $I_{f_\star \cap h_\star } \subset I_{f_\star \cap g_\star}$. In view of the cases (I0), (I1), and (I2) above, one easily deduces:

\begin{lem} \label{lem:tripleinc}
    Up to permuting $g_\star$ and $h_\star$, the following inclusions hold:
    \begin{align}
        \mathcal{I}_{f_\star \cap g_\star } \subset \mathcal{I}_{f_\star \cap h_\star}, \qquad \mathcal{D}_{f_\star \cap g_\star } \subset \mathcal{D}_{f_\star \cap h_\star}
    \end{align}
\end{lem}

Note that this lemma depends on the hypothesis that none of the forward envelopes we consider have an initial edge corresponding to the pure stable branch of a quadratic-like singularity.

%%%
            \subsection{The preorder between envelopes} \label{sec:envorder}
%%%

We are now ready to define a weak prefoliation structure on the collection $\mathscr{U}$ of forward envelopes. As in~\cite[Section 6.2]{BI08}, we construct a preorder $\precsim$ on $\mathscr{U}_\star$ in several steps. We first define a total order $\precsim_p$ on $\mathscr{U}_p$ at every point $p \in M^\mathrm{good}$ in five steps, before extending it to every $p \in M$.

%%%
    \subsubsection{Definition of the preorder}
%%%

Let $f_\star, g_\star \in \mathscr{U}_p$ be two marked forward envelopes passing through $p \in M^\mathrm{good}$, with respective domains $U$ and $V$, and shadows $\Sigma$ and $\Sigma'$.

%%%
                \paragraph{The geometric relation $\prec^{\mathrm{geo}}$.}
%%%
This relation corresponds to the relation $>_1$ from~\cite[Section 6.2.2]{BI08}, and is defined in a similar way. We set $f_\star \prec^\mathrm{geo}_p g_\star$ if there exists $a_0 \in \mathcal{D}_{f_\star \cap g_\star}$ with $b_0 = \varphi_{f_\star \cap g_\star}(a_0)$ such that we either have
        \begin{enumerate}
            \item[(a)] $\Sigma(a_0) \subsetneq \Sigma'(b_0)$, or
            \item[(b)] $(g,b_0)$ is locally strictly above $(f,a_0)$.
        \end{enumerate}

%%%
                \paragraph{The inclusion relation $\prec^{\mathrm{inc}}$.}
%%%

This relation corresponds to the relation $>_2$ from~\cite[Section 6.2.2]{BI08} and is defined in a similar way. We say that $f_\star$ is a \textbf{strict subenvelope} of $g_\star$ if $f_\star$ and $g_\star$ are not comparable by $\prec_p^\mathrm{geo}$, and  $\mathcal{D}_{f_\star} = \mathcal{D}_{f_\star \cap g_\star}$ but $\mathcal{D}_{g_\star \cap f_\star} \subsetneq \mathcal{D}_{g_\star}$. In this situation, $f_\star$ and $g_\star$ are comparable by $\prec^\mathrm{inc}$ and we set $f_\star \prec^\mathrm{inc}_p g_\star$.

%%%
                \paragraph{The relation at corners $\prec^{\mathrm{cor}}$.}
%%%
This relation has no equivalent in~\cite{BI08} and is a feature of the presence of saddle and quadratic singularities. We start with some definitions. 

\begin{defn} \label{def:admcorner}
    A collection of forward corners $\mathfrak{c} \subset \mathfrak{C}^+$ is \textbf{admissible} if it satisfies the following conditions.
    \begin{enumerate}
        \item $\mathfrak{c}$ contains exactly one forward corner at each pure saddle singularity, and no corner at quadratic-like singularities,
        \item If $(c,c')$ is a forward double-corner, then $\{c,c'\} \cap \mathfrak{c} = \varnothing$.
    \end{enumerate}
\end{defn}

\begin{choice}
    The relation $\prec^\mathrm{cor}$ depends on the choice of an admissible collection of forward corners $\mathfrak{c} \subset \mathfrak{C}^+$.
\end{choice}
This choice will be crucial later when we extend the weak prefoliation structure to completed envelopes, see Section~\ref{sec:comp} below.

We first assume that $f_\star$ and $g_\star$ are not comparable by $\prec^\mathrm{geo} \cup \prec^\mathrm{inc}$ and have an initial corner of overlap $c$. After possibly renaming them, we further assume that they meet along $c$ in the way described in Section~\ref{sec:initcor}. Then, we set $f_\star \prec^\mathrm{cor}_p g_\star$ if one of the following holds:
\begin{enumerate}
    \item $f_\star$ and $g_\star$ are geometrically noncoherent at $c$, and $\widehat{\Sigma}(c) \subsetneq \widehat{\Sigma}'(c)$,
    \item $f_\star$ and $g_\star$ are geometrically coherent at $c$, and $c \in \mathfrak{c}$.
\end{enumerate}
Otherwise, we set $g_\star \prec^\mathrm{cor}_p f_\star$.

We now assume that $f_\star$ and $g_\star$ are not comparable by $\prec^\mathrm{geo} \cup \prec^\mathrm{inc}$ and have an initial \emph{double} corner of overlap $(c,c')$ where they meet in the way described in Section~\ref{sec:initcor}. Then, we set $f_\star \prec^\mathrm{cor}_p g_\star$ if $f_\star$ and $g_\star$ are geometrically noncoherent at $(c,c')$ and $\widehat{\Sigma}_0(c) \subsetneq \widehat{\Sigma}'_0(c')$. Otherwise, we set $g_\star \prec^\mathrm{cor}_p f_\star$.

%%%
                \paragraph{The relation along unstable components $\prec^{\mathrm{un}}$.}
%%%
This relation has no equivalent in~\cite{BI08} and is a feature of the presence of source and quadratic singularities. We need to introduce some terminology and make some auxiliary choices before defining it.

Let $L \hookrightarrow M$ be a leaf of $\eta$ in $M^\mathrm{un}$. A \textbf{forward boundary component} of $L$ is a connected component of the complement of the backward edges of $L$ in $\partial L$. We denote by $\mathcal{B}_L^+$ the set of forward boundary components of $L$. A forward boundary component might be closed or might contain one or two corners. Moreover, a boundary component might contain two consecutive flow lines of $X$ with opposite orientations, corresponding to pure unstable branches $\gamma^\pm_u$ at a quadratic-like singularity $q$. We say that the pure stable branch $\gamma_s$ at $q$ is a \textbf{spike} for $L$ if it is connected to a (pure) saddle singularity in negative times. In that case, $q$ is a saddle singularity as well. Since $X$ has no broken triple saddle connections (see Definition~\ref{def:brokensad}), $L$ has \emph{at most one} spike. The forward boundary component $\gamma^-_u$ of $L$ will be called a \textbf{distinguished boundary component} and will be denoted by $b^\dagger_L \in \mathcal{B}_L^+$. A total order on $\mathcal{B}_L^+$ for which $b^\dagger_L$ (if it exists) is maximal will be called \textbf{admissible}. We define $\mathfrak{O}^+_L$ as the collection of all admissible total orders on $\mathcal{B}_L^+$, and we set
$$\mathfrak{O}^+ \coloneqq \prod_{L} \mathfrak{O}^+_L,$$
where the product runs over all the leaves of $\eta$ in $M^\mathrm{un}$. The admissibility condition will become relevant in the proof of Lemma~\ref{lem:trans}, to ensure transitivity of the full preorder.

\begin{choice} \label{choice:or}
The relation $\prec^\mathrm{un}$ depends on the choice of a collection of admissible orders $\mathfrak{o} \in \mathfrak{O}^+$.
\end{choice}

Let $U_c \subset \mathcal{D}_{f_\star}$ be an unstable component of $f$ which is not initial if $f$ is of unstable type, and let $h \coloneqq f_{\vert U_c}$ be the corresponding unstable tile, which covers a leaf $L \hookrightarrow M$ of $\eta$ in $M^\mathrm{un}$. An \textbf{entrance edge} of $f_\star$ to $U_c$ is a forward edge of $U_c$ which belongs to $\mathcal{D}_{f_\star}$. It exists since $U_c$ is not initial, but it is not necessarily unique. However, there are at most three of those. If there are two, then they correspond to consecutive stable and unstable branches of a singularity of $X$, and if there are three, then they correspond to a succession of a stable branch of a saddle singularity $q$, a connection between $q$ and another saddle singularity $q'$, and an unstable branch of $q'$. Such an entrance edge projects to $\overline{L}$ via the universal covering map, and  to a forward boundary component $b_{f_\star} \in \mathcal{B}_L^+$. Obviously, the latter is independent on the choice of entrance edge. 

We set $f_\star \prec_p^{\mathrm{un}} g_\star$ if $f_\star$ and $g_\star$ are not comparable by $\prec^\mathrm{geo}_p \cup \prec^\mathrm{inc}_p \cup \prec^\mathrm{cor}_p$, and the following holds. There exists an initial unstable component of overlap $U_c \subset \mathcal{D}_{f_\star \cap g_\star}$ between $f_\star$ and $g_\star$ corresponding to a leaf $L \hookrightarrow M$ as above, together with forward boundary components $b_{f_\star}, b_{g_\star} \in \mathcal{B}^+_L$ induced by an entrance edge of $f_\star$ to $U_c$ and of $g_\star$ to $\varphi_{f_\star \cap g_\star}(U_c)$, respectively, which satisfy $b_{f_\star} <^\mathfrak{o} b_{g_\star}$. Here, $<^\mathfrak{o}$ denotes the total order on $\mathcal{B}^+_L$ coming from $\mathfrak{o} \in \mathfrak{O}^+$. Notice that $b_{f_\star}$ and $b_{g_\star}$ are necessarily distinct, since $f_\star$ and $g_\star$ don't have an initial (double) corner of overlap.

\begin{figure}[h]
    \centering
                \includegraphics{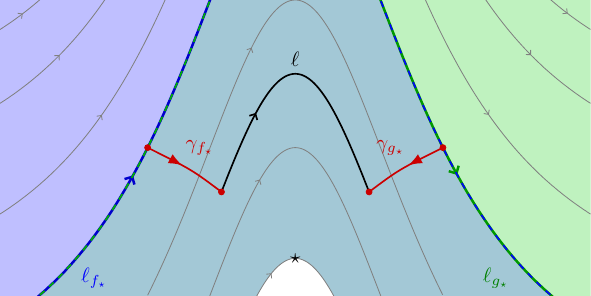}
    \caption{Overlap of $f_\star$ and $g_\star$ for $\prec^\mathrm{or}$. The blue region corresponds to $\mathcal{D}_{f_\star}$, the green region corresponds to $\mathcal{D}_{g_\star}$, and their region of overlap is in teal.}
    \label{fig:or}
\end{figure}

%%%
                \paragraph{The $X$-lines orientation relation $\prec^{\mathrm{or}}$.}
%%%
This relation corresponds to the relation $>_3$ from~\cite[Section 6.2.4]{BI08}. We define $f_\star \prec^\mathrm{or}_p g_\star$ if $f_\star$ and $g_\star$ are not comparable by $\prec^\mathrm{geo}_p \cup \prec^\mathrm{inc}_p \cup \prec^\mathrm{cor}_p \cup \prec^\mathrm{un}_p$, and the following holds. There exist maximal lower bounds $\ell_{f_\star} \in \mathcal{I}_{f_\star}$ and $\ell_{g_\star} \in \mathcal{I}_{g_\star}$ of $\mathcal{I}_{f_\star \cap g_\star} \subset \mathcal{I}_{f_\star}$ and $\mathcal{I}_{g_\star \cap f_\star} \subset \mathcal{I}_{g_\star}$, respectively, and two forward curves $\gamma_{f_\star} : [0,1] \rightarrow U$ and $\gamma_{g_\star} : [0,1] \rightarrow V$ such that
\begin{itemize}
    \item $\gamma_{f_\star}(0) \in \ell_{f_\star}$ and $\gamma_{g_\star}(0) \in \ell_{g_\star}$,
    \item $\gamma_{f_\star}$ and $\varphi_{g_\star \cap f_\star}\big({\gamma_{g_\star}}_{\vert (0,1]} \big)$ are disjoint,
    \item There exist a pure $X$-line $\ell \in \mathcal{I}_{f_\star \cap g_\star}$ and real numbers $s < t$ such that $\ell(s) \in \gamma_{f_\star}$ and $\varphi_{f_\star \cap g_\star} \big( \ell(t)\big) \in \gamma_{g_\star}$.
\end{itemize}
For short, we use the orientation on the pure $X$-lines near $\ell_{f_\star}$ and $\ell_{g_\star}$ to define the relation in this situation, see Figure~\ref{fig:or}.

Strictly speaking, this definition depends on an overall choice of orientation on $\R$: we could have replaced the condition $s < t$ with $t<s$ in the above third bullet.

\paragraph{Combining the five relations.} We are now in the position to define a relation $\precsim_p$ on $\mathscr{U}_p$ as follows. First, we set

\begin{center}
$f_\star \prec_p g_\star$ if and only if $f_\star \prec^{\heartsuit}_p g_\star$ for some $\heartsuit \in \{\mathrm{geo}, \mathrm{inc}, \mathrm{cor}, \mathrm{un}, \mathrm{or}\}$.
\end{center}

We say that $f_\star$ and $g_\star$ are \textbf{forward-equivalent} if there exists an equivalence $\varphi : f_\star \overset{\sim}{\rightarrow} g_\star$ satisfying $\Sigma' \circ \varphi = \Sigma$. In particular, $f_\star \approx_p g_\star$. Then, we define:

\begin{center}
$f_\star \precsim_p g_\star$ if and only if $f_\star \prec_p g_\star$, or $f_\star$ and $g_\star$ are forward-equivalent.
\end{center}

So far, the relation $\precsim_p$ is only defined for $p \in M^\mathrm{good}$. We extend it to every $p \in M$ as follows: $f_\star \precsim_p g_\star$ if and only if there exists $q \in M^\mathrm{good}$ and a continuous path $\gamma : [0,1] \rightarrow M$ with $\gamma(0) = p$, $\gamma(1) = q$, which lifts to both $f_\star$ and $g_\star$, and $f_{\gamma \star} \precsim_q g_{\gamma \star}$. Here, $f_{\gamma \star}$ and $g_{\gamma \star}$ are obtained from $f_\star$ and $g_\star$ by sliding the marked points along $\gamma$ in the obvious way.

The collection of homogeneous relations $\{ \precsim_p\}_{p \in M}$ naturally induces a homogeneous relation $\precsim$ on $\mathscr{U}_\star$. 

%%%
        \subsubsection{Proof of the weak prefoliation axioms}
%%%

This section is devoted to the proof of
\begin{prop} \label{prop:forord1}
The set $\mathscr{U}$ of forward envelopes equipped with the relation $\precsim$ is a weak prefoliation.
\end{prop}

We split the proof into the following four lemmas.

\begin{lem}[Definition~\ref{def:weakpref}, Axiom 2] \label{lem:coh}
The relation $\precsim$ is coherent along intersections.
\end{lem}

\begin{proof}
The relation $\precsim$ is trivially coherent along unstable components and along pure $X$-lines by definition. It is also coherent along forward curves for the same reasons as in the proof of~\cite[Lemma 6.15]{BI08}. We are left to show that it is coherent at singularities in $M^\mathrm{bad}$. If $p \notin M^\mathrm{good}$ and $f_\star \in \mathscr{U}_p$, then the marked point of $f_\star$ lies in the interior of the domain of $f$. Therefore, coherence of $\precsim$ along intersections contained in $M^\mathrm{good}$ immediately implies coherence at $p$.
\end{proof}

\begin{lem}[Definition~\ref{def:weakpref}, Axiom 0] \label{lem:tot}
For every $p \in M$, $\precsim_p$ is total. Moreover, for all $f_\star, g_\star \in \mathscr{U}_\star$,  $f_\star \simeq g_\star$ if and only if $f_\star$ and $g_\star$ are forward-equivalent.
\end{lem}

\begin{proof}
First assume $p \in M^\mathrm{good}$, and consider $f_\star, g_\star \in \mathscr{U}_p$ with our usual notations for their domains, shadows, and marked points. 

If $f_\star$ and $g_\star$ are comparable by $\prec^\mathrm{geo}_p$, say $f_\star \prec^\mathrm{geo}_p g_\star$, then $f_\star$ and $g_\star$ are not equivalent and $g_\star \not\prec^\mathrm{geo}_p f_\star$. Indeed, arguing as in~\cite[Section 6.2.2]{BI08} by considering forward curves from/to $a_0$, it is easy to see that for every $x \in \mathcal{D}_{f_\star \cap g_\star}$, $\Sigma(x) \subset \Sigma'\big( \varphi_{f_\star \cap g_\star}(x)\big)$ and $(f,x)$ is not locally strictly above $\big(g,\varphi_{f_\star \cap g_\star}(x)\big)$.

If $f_\star$ and $g_\star$ are not comparable by $\prec^\mathrm{geo}_p$, then for every $x \in \mathcal{D}_{f_\star \cap g_\star}$, $\Sigma(x) = \Sigma'\big(\varphi_{f_\star \cap g_\star}(x)\big)$ and $(f,x)$ and $\big(g,\varphi_{f_\star \cap g_\star}(x)\big)$ locally coincide.\footnote{Notice that they automatically locally coincide in the forward direction as soon as their shadows coincide, because of basic properties of forward envelopes.} We now analyze how $f_\star$ and $g_\star$ `first meet' along $\mathcal{D}_{f_\star \cap g_\star}$ depending on which case from Lemma~\ref{lem:casesI} holds:
\begin{itemize}
    \item[(I0)] $\mathcal{I}_{f_\star \cap g_\star} = \varnothing$. In that case, $f_\star$ and $g_\star$ have an initial unstable component of overlap $U_c \subset \mathcal{D}_{f_\star \cap g_\star}$. 
    
    If both $U_c$ and $\varphi\big(U_c\big)$ correspond to the initial unstable components of $f_\star$ and $g_\star$, respectively, then $f_\star$ and $g_\star$ are forward-equivalent by Proposition~\ref{prop:forenv}. 
    
    If $U_c$ corresponds the initial unstable component of $f_\star$, but $\varphi\big(U_c\big)$ does not correspond to the initial unstable component of $g_\star$, then $f_\star$ and $g_\star$ are comparable by $\prec^\mathrm{inc}$ and $f_\star \prec^\mathrm{inc} g_\star$. The case where $U_c$ and $\varphi\big(U_c\big)$ are swapped is similar.
    
    Otherwise, $f_\star$ and $g_\star$ are comparable by $\prec^\mathrm{un}$, since $U_c$ does not correspond to their initial unstable components and they don't have an initial (double) corner of overlap. Obviously, either $f_\star \prec^\mathrm{un}_p g_\star$ of $g_\star \prec^\mathrm{un}_p f_\star$ holds, but not both.
    
    \item[(I1)] \textit{$\mathcal{I}_{f_\star \cap g_\star}$ has minimal elements.} If $f_\star$ and $g_\star$ have an initial unstable component of overlap, we proceed as in the last step of the previous item. Otherwise, $f_\star$ and $g_\star$ are comparable by $\prec^\mathrm{inc}_p$. Indeed, if $\mathcal{I}_{f_\star \cap g_\star}$ has a unique minimal element, which is a minimum, then this minimum corresponds to the initial edge of either $f_\star$ or $g_\star$, depending on whether it is the minimum of $\mathcal{I}_{f_\star}$ or not. If $\mathcal{I}_{f_\star \cap g_\star}$ has strictly more than one minimal element, then it is easy to see that $g_\star \prec^\mathrm{inc}_p f_\star$.
    
    \item[(I2)] \textit{$\mathcal{I}_{f_\star \cap g_\star}$ has no minimal elements.} If $f_\star$ and $g_\star$ have an initial (double) corner of overlap, then they are comparable by $\prec^\mathrm{cor}_p$, and only one of $f_\star \prec^\mathrm{cor}_p g_\star$ or $g_\star \prec^\mathrm{cor}_p f_\star$ holds. We now assume that they do not have an initial (double) corner of overlap.
    
    Observe that if $f_\star$ and $g_\star$ are comparable by $\prec^\mathrm{cor} \cup \prec^\mathrm{inc}_p \cup \prec^\mathrm{un}_p$, then $\mathcal{I}_{f_\star \cap g_\star}$ is either empty or has a minimal element. Therefore, we need to show that $f_\star$ and $g_\star$ are comparable by $\prec^\mathrm{or}_p$. For that, we consider a forward curve $\gamma_0 : [0,1] \rightarrow U$ from the initial component of $f$ to $a$, and a forward curve $\gamma_1 : [0,1] \rightarrow V$ from the initial component of $g$ to $b$. There exist $t_0, t_1 \in [0,1]$ such that $\gamma_0(t_0)$ and $\gamma_1(t_1)$ belong $X$-lines $\ell_0 \subset U$ and $\ell_1 \subset V$ which are maximal lower bounds for $\mathcal{I}_{f_\star \cap g_\star} \subset \mathcal{I}_{f_\star}$ and $\mathcal{I}_{g_\star \cap f_\star} \subset \mathcal{I}_{g_\star}$, respectively. For $t > t_0$, $\gamma_0(t) \in \mathcal{D}_{f_\star \cap g_\star}$ and for $t > t_1$, $\gamma_1(t) \in \mathcal{D}_{g_\star \cap f_\star}$. Moreover, there exists $\epsilon > 0$ such that $\gamma_{f_\star} \coloneqq {\gamma_0}_{[t_0, t_0+ \epsilon)}$ and  $\gamma_{g_\star} \coloneqq {\gamma_1}_{[t_1, t_1+ \epsilon)}$ are strictly forward, and $\gamma_{f_\star}$ and $\varphi_{g_\star \cap f_\star}\big({\gamma_{g_\star}}_{\vert(t_1, t_1+ \epsilon)}\big)$ are disjoint. Since $f_\star$ and $g_\star$ do not have an initial unstable domain of overlap, there exists either
    \begin{enumerate} \setcounter{enumi}{-1}
        \item A pure $X$-line $\ell$ intersecting both $\gamma_{f_\star}$ and $\gamma_{g_\star}$,
        \item A special triple $(\ell_s, \ell^-_u, \ell^+_u)$ such that $\ell_s$ intersects one of $\gamma_{f_\star}$ and $\gamma_{g_\star}$ and $\ell^-_u$ or $\ell^+_u$ intersects the other,
        \item Two special triples $(\ell_s, \ell^-_u, \ell^+_u)$ and $(\ell'_s, \ell'^-_u, \ell'^+_u)$ with $\ell'_s \in \{\ell^-_u, \ell^+_u\}$ such that $\ell_s$ intersects one of $\gamma_{f_\star}$ and $\gamma_{g_\star}$ and $\ell'^-_u$ or $\ell'^+_u$ intersects the other.
    \end{enumerate}
In cases 1 and 2, we can modify $\gamma_{f_\star}$ or $\gamma_{g_\star}$ by performing suitable $T$-moves to reduce to case 0 while keeping the two curves disjoint. This implies that $f_\star$ and $g_\star$ are comparable by $\prec^\mathrm{or}_p$. It remains to show that $f_\star \prec^\mathrm{or}_p g_\star$ and $g_\star \prec^\mathrm{or}_p f_\star$ cannot hold simultaneously. If both inequalities hold, then after performing suitable $T$-moves and forward isotopies, we would obtain the following situation: there exists two forward curves $\gamma_{f_\star}$ and $\gamma_{g_\star}$ as above, together with two pure $X$-lines or pure `broken $X$-lines' (as in cases 1 and 2 above), one from $\gamma_{f_\star}$ to $\gamma_{g_\star}$ and the other one from $\gamma_{g_\star}$ to $\gamma_{f_\star}$. This would imply the existence of a pure $X$-line which is `trapped' (in positive or negative times) in a compact region of $U$, which is impossible. The details are left to the reader.
\end{itemize}

Finally, we assume that $p \notin M^\mathrm{good}$. In particular, $p \in \Delta_\mathrm{sa} \cup Q$, and the marked points of $f_\star$ and $g_\star$ belong to the interior of their domains of definition. We can slide these marked points along a stable or unstable branch at $p$ to obtain marked tiles $(f, a'), (g, b') \in \mathscr{U}_{p'}$ where $p' \in M^\mathrm{good}$. Therefore, $(f, a')$ and $(g, b')$ are comparable by $\precsim_{p'}$ by the above, so $f_\star$ and $g_\star$ are comparable by $\precsim_p$ by definition. Finally, if $f_\star \precsim_p g_\star$ and $g_\star \precsim_p f_\star$, then the same holds at some nearby $p' \in M^\mathrm{good}$ by sliding the marked point $a$ and $b$ to nearby $a'$ and $b'$, hence $(f, a') \approx_{p'} (g,b')$. Sliding the marked points back yields $f_\star \approx_p g_\star$.
\end{proof}

\begin{lem}[Definition~\ref{def:weakpref}, Axiom 1] \label{lem:refine}
For all $f_\star, g_\star \in \mathscr{U}_\star$, if $f_\star \precsim g_\star$, then $g_\star$ is locally above $f_\star$.
\end{lem}

\begin{proof}
The proof is essentially the same as the one of~\cite[Lemma 6.13]{BI08}. Let $p \in M$ and $f_\star, g_\star \in \mathscr{U}_p$. Assume first that $p \in M^\mathrm{good}$. If $f_\star$ and $g_\star$ are comparable by $\prec^\mathrm{geo}_p$, say $f_\star \prec^\mathrm{geo}_p g_\star$, then we connect the point $a_0$ given by the definition to the marked point $a$ of $f$ by a forward curve. Properties of upper enveloping tiles (see Definition~\ref{def:upenv}) and one-dimensional envelopes (see~\cite[Section 5]{BI08}) immediately imply that $g_\star$ is locally above $f_\star$. If $f_\star$ and $g_\star$ are not comparable by $\prec^\mathrm{geo}_p$, then $f_\star$ is not locally strictly above $g_\star$, and since one of these two forward tiles is locally above the other (for the same reasons as in~\cite[Lemma 6.5]{BI08}), then $g_\star$ is locally above $f_\star$. If $p \notin M^\mathrm{good}$ and $f_\star$ is locally strictly above $g_\star$, then there exists $p' \in M^\mathrm{good}$ on one of the stable or unstable branches at $p$ such that if $(f, a'), (g, b') \in \mathscr{U}_{p'}$ denote the marked tiles obtained by sliding the marked points of $f_\star$ and $g_\star$, respectively, then $(f,a')$ is locally strictly above $(g,b')$. Therefore, $(g, b') \prec^\mathrm{geo}_{p'} (f, a')$, hence $g_\star \prec f_\star$.
\end{proof}

We finish the proof of Proposition~\ref{prop:forord1} with the trickiest property: 

\begin{lem}\label{lem:trans}
The relation $\precsim$ is transitive.
\end{lem}

\begin{proof}
Let $p \in M$. Notice that if $f_\star, g_\star, h_\star \in \mathscr{U}_p$ and $f_\star$ and $g_\star$ are forward-equivalent, then $f_\star \precsim_p h_\star$ if and only if $g_\star \precsim_p h_\star$, by definition. Therefore, in view of Lemma~\ref{lem:tot}, transitivity of $\precsim_p$ is equivalent to the statement: \emph{every triple of elements in $\mathscr{U}_p$ admits an extremum (maximum or minimum)}. Moreover, by Lemma~\ref{lem:coh}, it is enough to consider the case $p \in M^\mathrm{good}$. By Lemma~\ref{lem:tripleinc}, we can assume that $\mathcal{I}_{f_\star \cap g_\star } \subset \mathcal{I}_{f_\star \cap h_\star}$ and $\mathcal{D}_{f_\star \cap g_\star } \subset \mathcal{D}_{f_\star \cap h_\star}$. We enumerate all the possible cases for the order between $f_\star$ and $g_\star$.
    \begin{enumerate}[leftmargin=*]
        \item \underline{$f_\star$ and $g_\star$ are comparable by $\prec^\mathrm{geo}$.} Assume that there exist $a_0 \in \mathcal{D}_{f_\star \cap g_\star}$ and $b_0 = \varphi_{f_\star \cap g_\star}(a_0)$ such that $\Sigma_f(a_0) \neq \Sigma_g(b_0)$, and consider $c_0 \coloneqq \varphi_{f_\star \cap h_\star}(a_0) \in \mathcal{D}_{h_\star}$. This is well-defined since $a_0 \in \mathcal{D}_{f_\star \cap g_\star} \subset \mathcal{D}_{f_\star \cap h_\star}$. Then, among the three sections $\Sigma_f(a_0)$, $\Sigma_g(b_0)$ and $\Sigma_h(c_0)$ of $\mathscr{A}_p$, there must exist a \emph{strict} maximum or minimum for the inclusion, implying that the corresponding marked tiles is an extremum in $\{f_\star, g_\star, h_\star\}$. Otherwise, there exists $a_0$ and $b_0$ as above such that $(g, b_0)$ is locally strictly above $(f, a_0)$. In that case, $p \in M^\mathrm{pure}$. Intersecting $f_\star$, $g_\star$ and $h_\star$ with a small embedded disk in $M$ passing through $p$ and transverse to $X$, we obtain three $1$-dimensional envelopes for the restriction of $\mathscr{A}$ to that disk, two of which are strictly comparable by the geometric order. We further consider two cases.
            \begin{enumerate}[leftmargin=1em]
                \item If $c_0 = \varphi_{f_\star \cap h_\star}(a_0)$ does not belong to the initial edge of $h$ (if it exists), then one of the three $1$-dimensional envelopes is strictly below or above the two others, and the corresponding marked forward envelope is an extremum in $\{f_\star, g_\star, h_\star\}$.
                \item Otherwise, $f_\star$ and $g_\star$ `begin to overlap' along the initial edge of $h_\star$. If $\Sigma_h(c_0) \neq \Sigma_f(a_0)$, then since $\Sigma_f(a_0)=\Sigma_g(b_0)$ by assumption, $h_\star$ is comparable by $\prec^\mathrm{geo}$ with both $f_\star$ and $g_\star$. Moreover, it is either below or above both of them for $\prec^\mathrm{geo}$, so $h_\star$ is an extremum of $\{f_\star, g_\star, h_\star\}$. If $\Sigma_h(c_0) = \Sigma_f(a_0) = \Sigma_g(b_0)$, then $h_\star$ is comparable by $\prec^\mathrm{inc}$ with both $f_\star$ and $g_\star$, and $h_\star$ is similarly an extremum of $\{f_\star, g_\star, h_\star\}$.
            \end{enumerate}
            
        \item \underline{$f_\star$ and $g_\star$ are comparable by $\prec^\mathrm{inc}$.} We first assume that $g_\star$ is a strict subenvelope of $f_\star$; in particular, $\mathcal{D}_{g_\star} = \mathcal{D}_{g_\star \cap f_\star}$. Notice that $\mathcal{D}_{h_\star \cap g_\star} \subset \mathcal{D}_{h_\star \cap f_\star}$ by Lemma~\ref{lem:ident}. By the previous case, we can assume that $g_\star$ and $h_\star$ are not comparable by $\prec^\mathrm{geo}$. We also have $\mathcal{D}_{g_\star} = \mathcal{D}_{g_\star \cap h_\star}$, so either $g_\star$ and $h_\star$ are order-equivalent, or $g_\star$ is a strict subenvelope of $h_\star$. Hence, either $f_\star$ or $g_\star$ is an extremum in $\{f_\star, g_\star, h_\star\}$. We now assume that $f_\star$ is a strict subenvelope of $g_\star$. Hence, $\mathcal{D}_{f_\star} = \mathcal{D}_{f_\star \cap g_\star} \subset \mathcal{D}_{f_\star \cap h_\star} = \mathcal{D}_{f_\star}$. In particular, $\mathcal{D}_{g_\star \cap f_\star} = \mathcal{D}_{g_\star \cap h_\star}$ and $\mathcal{D}_{h_\star \cap f_\star} = \mathcal{D}_{h_\star \cap g_\star}$, and by case 1 above, we can assume that $f_\star$ and $g_\star$ are not comparable by $\prec^\mathrm{geo}$, and $f_\star$ and $h_\star$ are not comparable by $\prec^\mathrm{geo}$. In this situation, $f_\star$ is a (not necessarily strict) subenvelope of $g_\star$ and $h_\star$, so $\{f_\star, g_\star, h_\star\}$ has an extremum by the previous argument.

        \item \underline{$f_\star$ and $g_\star$ are comparable by $\prec^\mathrm{cor}$.} We consider three cases.
            
            \begin{enumerate}[leftmargin=1em]
                \item $f_\star$ and $g_\star$ have a simple initial corner of overlap $c$, which is not connected to another forward corner. Since we do not consider forward envelopes starting at the pure stable branch of a quadratic-like singularity, $c$ is a pure corner. We consider two subcases.
            
                \begin{enumerate}[leftmargin=0.5em]
                    \item $\mathcal{D}_{f_\star \cap g_\star} \subsetneq \mathcal{D}_{f_\star \cap h_\star}$. Then $c$ is also the initial corner of overlap of $g_\star$ and $h_\star$. In that case, $\mathcal{D}_{g_\star \cap h_\star} = \mathcal{D}_{g_\star \cap f_\star}$ so by case 1 above, we can assume that $g_\star$ and $h_\star$ are not comparable by $\prec^\mathrm{geo}$, hence they are comparable by $\prec^\mathrm{cor}$. Let $a_0 \in U$ be a point on the $X$-line corresponding to the branch at $c$ bounding $\mathcal{D}_{f_\star \cap g_\star}$ in $\mathcal{D}_{f_\star}$, and $c_0 \coloneqq \varphi_{f_\star \cap h_\star}(a_0)$. Note that $f_\star$ and $h_\star$ play a symmetric role here. If $\Sigma_f(a_0) \neq \Sigma_h(c_0)$, then $\widehat{\Sigma}_f(c) \neq \widehat{\Sigma}_h(c)$ (see Section~\ref{sec:initcor} for the notation), and either $f_\star$ or $h_\star$ is an extremum of $\{f_\star, g_\star, h_\star\}$. Otherwise, $\widehat{\Sigma}_f(c) = \widehat{\Sigma}_h(c)$, and the order of $f_\star$ with $g_\star$ for $\prec^\mathrm{cor}$ is the same as the one of $h_\star$ with $g_\star$, implying that $g_\star$ is an extremum of $\{f_\star, g_\star, h_\star\}$.

                    \item $\mathcal{D}_{f_\star \cap g_\star} = \mathcal{D}_{f_\star \cap h_\star}$. By case 1 above, we can assume that $f_\star$ and $h_\star$ are not comparable by $\prec^\mathrm{geo}$, and moreover they are not comparable by $\prec^\mathrm{inc}$. We consider three further subcases.
                        \begin{itemize}[leftmargin=*]
                            \item $f_\star$ and $h_\star$ are comparable by $\prec^\mathrm{cor}$. Then $c$ is also the initial corner of overlap between $f_\star$ and $g_\star$, and $\mathcal{D}_{g_\star \cap f_\star} \subsetneq \mathcal{D}_{g_\star \cap h_\star}$ so we can appeal to case 3(a)i above.
                        
                            \item $f_\star$ and $h_\star$ are comparable by $\prec^\mathrm{un}$. Then one easily sees that $h_\star$ and $g_\star$ are also comparable by $\prec^\mathrm{un}$ and share the same initial unstable domain of overlap as $h_\star$ and $f_\star$. Since the entrance edges of $f_\star$ and $g_\star$ lie on the same forward boundary component of the leaf corresponding to that unstable component, but the entrance edge(s) of $h_\star$ necessarily lies on a different forward boundary component, the order of $f_\star$ with $h_\star$ for $\prec^\mathrm{un}$ is the same as the one of $g_\star$ with $h_\star$, implying that $h_\star$ is an extremum of $\{f_\star, g_\star, h_\star\}$.
                        
                            \item $f_\star$ and $h_\star$ are comparable by $\prec^\mathrm{or}$. Then $g_\star$ and $h_\star$ are comparable by $\prec^\mathrm{or}$ as well, and it follows from the definition of $\prec^\mathrm{or}$ and the behavior of $X$-lines near a forward corner that the order of $f_\star$ with $h_\star$ for $\prec^\mathrm{or}$ is the same as the one of $g_\star$ with $h_\star$, implying that $h_\star$ is an extremum of $\{f_\star, g_\star, h_\star\}$.
                        \end{itemize}
                \end{enumerate}
                
                \item $f_\star$ and $g_\star$ have a simple initial corner of overlap $c$, which is connected to another forward corner $c'$. We again consider two subcases.

                    \begin{enumerate}[leftmargin=0.5em]
                        \item \textit{$c'$ is a corner at a pure saddle singularity.} We further assume that $f_\star$ and $h_\star$ are comparable by $\prec^\mathrm{cor}$ and $c'$ is their initial corner of overlap, as all the other cases can be treated as in case 3(a).

                        Up to permuting $f_\star$, $g_\star$ and $h_\star$, we can assume that $(c,c')$ is a double corner and $f_\star$ and $h_\star$ meet along $(c,c')$. We will use the notations of Section~\ref{sec:initcor}. We consider the sections $\widehat{\Sigma}_f(c), \widehat{\Sigma}_g(c) \subset \mathscr{A}_{q}$, $\widehat{\Sigma}_g(c'), \widehat{\Sigma}_h(c') \subset \mathscr{A}_{q'}$, and $\Sigma_g(b_0), \widehat{\Sigma}_{f,0}(c), \widehat{\Sigma}_{h,0}(c') \subset \mathscr{A}_{q_0}$. Here, $b_0$ lies on the boundary of $\mathcal{D}_{g_\star \cap f_\star} \subset \mathcal{D}_{g_\star}$ and $g(b_0) = q_0$. Tracing through the definitions of these sets, one easily shows:
                        \begin{align}
                             \widehat{\Sigma}_f(c) \subset \widehat{\Sigma}_g(c) &\implies \widehat{\Sigma}_{f,0}(c) \subset \Sigma_g(b_0), \label{eq:sigma1} \\
                            \widehat{\Sigma}_g(c) \subset \widehat{\Sigma}_f(c) &\iff \Sigma_g(b_0) \subset \widehat{\Sigma}_{f,0}(c), \label{eq:sigma2}
                        \end{align}
                        and similarly,
                        \begin{align*}
                            \widehat{\Sigma}_h(c') \subset \widehat{\Sigma}_g(c') &\implies \widehat{\Sigma}_{h,0}(c') \subset \Sigma_g(b_0), \\
                            \widehat{\Sigma}_g(c') \subset \widehat{\Sigma}_h(c') &\iff \Sigma_g(b_0) \subset \widehat{\Sigma}_{h,0}(c').
                        \end{align*}
                        This immediately implies the following:
                            \begin{itemize}[leftmargin=*]
                                \item If $\widehat{\Sigma}_f(c) \subset \widehat{\Sigma}_g(c)$ and  $\widehat{\Sigma}_g(c') \subsetneq \widehat{\Sigma}_h(c')$, then $\widehat{\Sigma}_{f,0}(c) \subsetneq \widehat{\Sigma}_{h,0}(c')$.
                                \item If $\widehat{\Sigma}_g(c) \subset \widehat{\Sigma}_f(c)$ and $\widehat{\Sigma}_h(c') \subset \widehat{\Sigma}_g(c')$, then $\widehat{\Sigma}_{h,0}(c') \subset \widehat{\Sigma}_{f,0}(c)$.
                            \end{itemize}
                        By the definition of $\prec^\mathrm{cor}$, we have:
                            \begin{align}
                                \widehat{\Sigma}_h(c') \subset \widehat{\Sigma}_g(c') \implies h_\star \prec^\mathrm{cor} g_\star, \notag \\
                                \widehat{\Sigma}_g(c) \subset \widehat{\Sigma}_f(c) \implies g_\star \prec^\mathrm{cor} f_\star, \label{eq:impl1}\\
                                \widehat{\Sigma}_{h,0}(c') \subset \widehat{\Sigma}_{f,0}(c) \implies h_\star \prec^\mathrm{cor} f_\star. \label{eq:impl2}
                            \end{align}
                        It is easy to conclude that in all the possible cases for the inclusions between $\widehat{\Sigma}_f(c)$ and $\widehat{\Sigma}_g(c)$, and between $\widehat{\Sigma}_g(c')$ and $\widehat{\Sigma}_h(c')$, $\{f_\star, g_\star, h_\star\}$ has an extremum.
                        
                        \item \textit{$c'$ is a corner at a quadratic-like singularity}. If $(c',c)$ is a double forward corner, all the possible cases can be treated as 3(a), so we assume that $(c,c')$ is a double forward corner. For the same reasons, we assume that $f_\star$ and $h_\star$ are comparable by $\prec^\mathrm{cor}$ and $(c,c')$ is their initial corner of overlap, so that $\mathcal{D}_{g_\star \cap f_\star} \subseteq \mathcal{D}_{g_\star \cap h_\star}$. Here, $f_\star$ and $g_\star$ intersect as in case (C2)(i) of Section~\ref{sec:initcor}, see Figure~\ref{fig:corC2i}, while $f_\star$ and $h_\star$ intersect as in case (C1) (ii), see Figure~\ref{fig:corC1iiun} (in particular the bottom one).

                        We use the same notations as before. We also consider points $y_0 \in \mathcal{D}_{g_\star}$ and $z_0 \in \mathcal{D}_{h_\star}$ such that $\varphi_{h_\star \cap g_\star}(z_0) = y_0$, and $z_0$ lies on the $X$-line corresponding to the minimum of $\mathcal{I}_{h_\star \cap g_\star}$ in $\mathcal{I}_{h_\star}$, so that $h(z_0)$ belongs to the right unstable branch of the singularity at $c'$.

                        The implications~\eqref{eq:sigma1} and~\eqref{eq:sigma2} still hold, as well as
                        \begin{align*}
                            \Sigma_h(z_0) \subset \Sigma_g(y_0) \iff \widehat{\Sigma}_{h,0}(c') \subset \Sigma_g(b_0).
                        \end{align*}
                        Both implications follow from the definitions and crucially use the coherence of the sections at singularities (axiom 5 of Definition~\ref{def:upenv}), applied at the quadratic-like singularity at $c'$.

                        The implications~\eqref{eq:impl1} and~\eqref{eq:impl2} hold for similar reasons, as well as
                        \begin{align*}
                            \Sigma_h(z_0) \subset \Sigma_g(y_0) &\implies h_\star \prec g_\star, \\
                            \Sigma_g(y_0) \subsetneq \Sigma_h(z_0) &\implies g_\star \prec h_\star
                        \end{align*}
                        Indeed, if $\Sigma_h(z_0) \subsetneq \Sigma_g(y_0)$ (resp.~$\Sigma_g(y_0) \subsetneq \Sigma_h(z_0)$), then $g_\star$ and $h_\star$ are comparable by $\prec^\mathrm{geo}$ and $h_\star \prec^\mathrm{geo} g_\star$ (resp.~$g_\star \prec^\mathrm{geo} h_\star$). If $\Sigma_h(z_0) = \Sigma_g(y_0)$, then $g_\star$ and $h_\star$ are not comparable by $\prec^\mathrm{geo}$, and  two possibilities arise. Either $z_0$ belong to the initial edge of $h_\star$, and $h_\star \prec^\mathrm{inc} g_\star$, or $g_\star$ and $h_\star$ are comparable by $\prec^\mathrm{un}$ and their initial unstable domain of overlap corresponds to the unstable tile $L$ at $q'$ (recall that $f_\star$ and  $h_\star$ do not meet along $c$). Moreover, the connection in $c \cap c'$ is a spike for $L$, and the entrance edge for $g_\star$ is the distinguished boundary component of $L$. By our conventions for admissible orders on unstable tiles, we have $h_\star \prec^\mathrm{un} g_\star$.

                       Examining all the possible inclusions between $\widehat{\Sigma}_f(c)$ and $\widehat{\Sigma}_g(c)$, and between $\Sigma_g(y_0)$ and $\Sigma_h(z_0)$, it is easy to show that $\{f_\star, g_\star, h_\star\}$ always has an extremum.
                    \end{enumerate}

                \item $f_\star$ and $g_\star$ have an initial double corner of overlap $(c,c')$. The cases where $h_\star$ is comparable by $\prec^\mathrm{geo}$, $\prec^\mathrm{inc}$, or $\prec^\mathrm{cor}$ with either $f_\star$ or $g_\star$ can be treated as in 3(a) or 3(b), up to permutation. In the remaining cases, $h_\star$ is comparable with both $f_\star$ and $g_\star$ by one of the two relations $\prec^\mathrm{un}$ or $\prec^\mathrm{or}$, and the strategy of case 3(a)ii applies. 
            \end{enumerate}
        
        \item \underline{$f_\star$ and $g_\star$ are comparable by $\prec^\mathrm{un}$.} The only case that has not been treated yet, up to permutation, is when $f_\star$, $g_\star$, and $h_\star$ are all comparable by $\prec^\mathrm{un}$ and share the same initial unstable domain of overlap. In this situation, transitivity immediately follows from the definition of $\prec^\mathrm{un}$.
        
        \item \underline{$f_\star$ and $g_\star$ are comparable by $\prec^\mathrm{or}$.} The only case that has not been treated yet, up to permutation, is when $f_\star$, $g_\star$, and $h_\star$ are all comparable by $\prec^\mathrm{or}$, and $\mathcal{D}_{f_\star \cap g_\star} = \mathcal{D}_{f_\star \cap h_\star}$ and $\mathcal{D}_{g_\star \cap f_\star} = \mathcal{D}_{g_\star \cap h_\star}$. Without loss of generality, let us assume that $f_\star \prec^\mathrm{or} g_\star$ and $g_\star \prec^\mathrm{or} h_\star$. We consider forward curves $\gamma_{f_\star}$, $\gamma_{g_\star}$, and $\gamma_{h_\star}$ as in the definition of $\prec^\mathrm{or}$; by item (I2) in the proof of Lemma~\ref{lem:tot}, we can assume that there exist a pure $X$-line $\ell_0 \in \mathcal{I}_{f_\star \cap g_\star}$ from $\gamma_{f_\star}$ to $\gamma_{g_\star}$ and a pure $X$-line $\ell_1 \in \mathcal{I}_{g_\star \cap h_\star}$ from $\gamma_{g_\star}$ to $\gamma_{h_\star}$. Let $\ell'_1 \in \mathcal{I}_{f_\star \cap g_\star}$ denote $\alpha_{g_\star \cap f_\star}(\ell_1)$. If $\ell'_1 = \ell_0$, we readily obtain a pure $X$-line from $\gamma_{f_\star}$ to $\gamma_{h_\star}$ and $f_\star \prec^\mathrm{or} h_\star$. Otherwise, $\varphi_{g_\star \cap f_\star}({\gamma_{g_\star}}_{\vert (0,1]})$ induces a strictly forward curve intersecting $\ell_0$ and $\ell'_1$, so either $\ell_0 \lhd \ell'_1$ or $\ell'_1 \lhd \ell_0$. Let us consider the case $\ell_0 \lhd \ell'_1$ first. Using Lemma~\ref{lem:forcurv2}, it is easy to show that either $\ell'_1$ intersects $\gamma_{f_\star}$, or is part of a broken $X$-line intersecting $\gamma_{f_\star}$ as in cases 1 and 2 of item (I2) in the proof of Lemma~\ref{lem:tot}. It follows that $f_\star \prec^\mathrm{or} h_\star$. The case $\ell'_1 \lhd \ell_0$ can be treated similarly by considering $\ell'_0 \coloneqq \alpha_{f_\star \cap g_\star}(\ell_0)$ and $\ell_1$ instead, and showing that $\ell'_0$ is part of a (broken) pure $X$-line intersecting $\gamma_{h_\star}$. \qedhere
\end{enumerate}
\end{proof}

%%%
            \subsubsection{Completed forward envelopes} \label{sec:comp}
%%%

In this section, we extend the previously constructed weak prefoliation structure on $\mathscr{U}$ to the collection of \emph{completed} forward envelopes $\overline{\mathscr{U}}$. Notice that since $\mathscr{U}$ covers $M^\mathrm{good}$ (minus the stable branches at quadratic-like singularities), $\overline{\mathscr{U}}$ covers $M$.

Let $f, g \in \mathscr{U}$ be forward envelopes (with our usual notations), with completions $\bar{f}, \bar{g} \in \overline{\mathscr{U}}$, and $a \in \overline{U}$ and $b \in \overline{V}$ such that $\bar{f}(a) = \bar{g} (b) \eqqcolon p$. We write $\bar{f}_\star = (\bar{f}, a)$ and $\bar{g}_\star = (\bar{g}, b)$.

If $p \in M^\mathrm{good}$, then the extension of $\precsim_p$ to $\overline{\mathscr{U}}_p$ is entirely determined by axioms 1 and 2 of Definition~\ref{def:weakpref}. Indeed, $a$ either belongs to $U$ or to a forward edge of $\overline{U}$, and similarly for $b$. Near $p$, $\bar{f}_\star$ and $\bar{g}_\star$ locally look like a disk or a half-disk bounded by the flow line of $X$ passing through $p$. Therefore, we can distinguish two cases:
\begin{enumerate}[leftmargin=*]
    \item $\bar{f}_\star$ and $\bar{g}_\star$ locally coincide. Then we can slide the marked points along a (possibly trivial) small curve $\gamma$ from $p$ to $p' \in M$ such that the new marked points satisfy $a' \in U$ and $b' \in V$. Then, we set $\bar{f}_\star \dot{\precsim}_p \bar{g}_\star$ if $(f, a') \precsim_{p'} (g,b')$, and $\bar{g}_\star \dot{\precsim}_p \bar{f}_\star$ otherwise. Notice that this does not depend on the choice of $\gamma$, since $f$ and $g$ intersect along a path-connected set near $p$. 
    \item One of $\bar{f}_\star$ and $\bar{g}_\star$ is locally strictly above the other, since $f$ and $g$ cannot topologically cross. We set $\bar{f}_\star \dot{\precsim}_p \bar{g}_\star$ if $\bar{g}_\star$ is locally strictly above $\bar{f}_\star$, and $\bar{g}_\star \dot{\precsim}_p \bar{f}_\star$ otherwise.
\end{enumerate}

The previous arguments can easily be adapted to extend $\precsim_p$ to $\overline{\mathscr{U}}_p$ for $p \in Q^*$, because of the shape of forward tiles near a quadratic-like singularities, see Figure~\ref{fig:backcor}. Crucially, the interiors of two such tiles intersect along a \emph{path-connected} set near $p$ when they locally coincide. It is then straightforward to check that $\dot{\precsim}_p$ is an extension of $\precsim_p$ to $\overline{\mathscr{U}}_p$ and that the axioms of Definition~\ref{def:weakpref} are satisfied \underline{on $M\setminus \Delta^\mathrm{pure}_\mathrm{sa}$}.

The extension of $\precsim_p$ to $\overline{\mathscr{U}}_p$ for $p \in \Delta^\mathrm{pure}_\mathrm{sa}$ is more subtle. We consider two cases:
\begin{enumerate}
    \item $\bar{f}_\star$ and $\bar{g}_\star$ intersect along a stable or unstable branch at $p$. We then slide the marked points along this branch to obtain $(\bar{f}, a'), (\bar{g}, b') \in \overline{\mathscr{U}}_{p'}$ and we set $\bar{f}_\star \dot{\precsim}_p \bar{g}_\star$ if $(\bar{f}, a') \dot{\precsim}_{p'} (\bar{g}, b')$, and $\bar{g}_\star \dot{\precsim}_p \bar{f}_\star$ otherwise. Notice that $\dot{\precsim}_{p'}$ is already defined.
    
    \item $\bar{f}_\star$ and $\bar{g}_\star$ do not intersect along the stable and unstable branches at $p$. In that case, they form two opposite (forward) quadrants at $p$, and they only intersect at $p$ near $p$. Let $\Sigma^\llcorner_f \subset \mathscr{A}_p$ denotes the set of marked tiles $h_\star \in \mathscr{A}_p$ such that either $h_\star$ is locally geometrically below $f$, or they locally coincide (relative to the domain of $f$) and after sliding the marked points inside the interiors of the domain of $f$ and $h$, $(h, z') \in \Sigma_f(x')$. Here, $x'$ and $z'$ denote the slid marked points of $f$ and $h$, respectively, satisfying $f(x') = h(z')$. We further denote by $\widehat{\Sigma}^\llcorner_f \subset \mathscr{A}_p$ the section generated by $\Sigma^\llcorner_f$. We similarly define $\widehat{\Sigma}^\llcorner_g \subset \mathscr{A}_p$. If $\widehat{\Sigma}^\llcorner_f \subsetneq \widehat{\Sigma}^\llcorner_g$ (resp.~$\widehat{\Sigma}^\llcorner_g \subsetneq \widehat{\Sigma}^\llcorner_f$), we set $\bar{f}_\star \dot{\precsim}_p \bar{g}_\star$ (resp.~$\bar{g}_\star \dot{\precsim}_p \bar{f}_\star$). This corresponds to the `geometrically noncoherent' case in the definition of $\prec^\mathrm{cor}$. 
    
    Furthermore, if there exists $\bar{h}_\star \in \overline{\mathscr{U}}_p$ such that 
        \begin{itemize}
            \item $\bar{f}_\star$ and $\bar{h}_\star$ intersect as in case 1 and $\bar{f}_\star \dot{\precsim}_p \bar{h}_\star$ (resp.~$\bar{h}_\star \dot{\precsim}_p \bar{f}_\star$),
            \item $\bar{g}_\star$ and $\bar{h}_\star$ intersect as in case 1 and $\bar{h}_\star \dot{\precsim}_p \bar{g}_\star$ (resp.~$\bar{g}_\star \dot{\precsim}_p \bar{h}_\star$),
        \end{itemize}
        then necessarily $\widehat{\Sigma}^\llcorner_f \subset \widehat{\Sigma}^\llcorner_g$ (resp.~$\widehat{\Sigma}^\llcorner_g \subset \widehat{\Sigma}^\llcorner_f$). If moreover $\widehat{\Sigma}^\llcorner_f = \widehat{\Sigma}^\llcorner_g$, then we set $\bar{f}_\star \dot{\precsim}_p \bar{g}_\star$ (resp.~$\bar{g}_\star \dot{\precsim}_p \bar{f}_\star$).

        If $\widehat{\Sigma}^\llcorner_f = \widehat{\Sigma}^\llcorner_g$ and no such $\bar{h}_\star \in \overline{\mathscr{U}}_p$ exists, we set $\bar{f}_\star \dot{\precsim}_p \bar{g}_\star$ if $\bar{f}_\star$ is bounded by the forward corner in $\mathfrak{c}$, and $\bar{g}_\star \dot{\precsim}_p \bar{f}_\star$ otherwise.\footnote{Alternatively, we could use the forward corner at $p$ not in $\mathfrak{c}$. This choice has to be the same for all completed forward envelopes at $p$, but it does not need to be the same for all pure saddle singularities.}
\end{enumerate}

It is not a priori clear that $\dot{\precsim}$ is coherent at $p \in \Delta^\mathrm{pure}_\mathrm{sa}$ in case 1, since $f$ and $g$ can locally intersect along a \emph{disconnected} set. This corresponds to the situation that we now describe. Let us pick local coordinates $(x,y,z)$ about $p$ adapted to $\eta$ such that the stable branches at $p$ are contained in $\{x=0\}$ and the unstable branches at $p$ are contained in $\{y=0\}$. In these coordinates, we identify $f$ and $g$ with the graphs of two functions, still denoted by $f$ and $g$, over domains $D_f, D_g \subset \R^2_{x,y}$. After possibly shrinking the coordinate neighborhood of $p$, we further assume 
\begin{align*}
    D_f &= \{- \epsilon < x < 0\} \cup \{- \epsilon < y < 0\}, \\
    D_g &= \{0 < x < \epsilon\} \cup \{0 < y < \epsilon\},
\end{align*}
for some $\epsilon > 0$. Moreover, $f$ and $g$ continuously extend to $0$ along $\{x = 0\} \cup \{y=0\}$. Writing 
\begin{align*}
    Q_- &\coloneqq \{ 0 < x < \epsilon, - \epsilon < y< 0\}, \\
    Q_+ &\coloneqq \{ - \epsilon < x < 0, 0 < y < \epsilon\}, 
\end{align*}
we have
$$D_f \cap D_g = Q_-  \sqcup Q_+,$$
and this set is \emph{disconnected}. We choose arbitrary points $x_- \in Q_-$ and $x_+ \in Q_+$, and we set $f \prec_\pm g$ if one of the following holds:
\begin{enumerate}
    \item[(a)] $f_{\vert Q_\pm} \leq  g_{\vert Q_\pm}$ and $f_{\vert Q_\pm} \neq g_{\vert Q_\pm}$,
    \item[(b)] $f_{\vert Q_\pm} = g_{\vert Q_\pm}$, and $(f,x_\pm) \prec (g, x_\pm)$.
\end{enumerate}
Otherwise, we set $g \prec_\pm f$. In the latter case, $g_{\vert Q_\pm} \leq  f_{\vert Q_\pm}$ and $(g,x_\pm) \prec (f, x_\pm)$ if $f_{\vert Q_\pm} = g_{\vert Q_\pm}$, since $f$ and $g$ are forward envelopes and $(f, x_\pm) \not\approx (g, x_\pm)$. Moreover, when $f_{\vert Q_\pm} = g_{\vert Q_\pm}$, $(f, x_\pm)$ and $(g, x_\pm)$ are comparable by $\prec^\mathrm{geo}$ or $\prec^\mathrm{cor}$, and in the latter case their initial corner of overlap $c_\pm$ bounds $Q_\pm$. See Figure~\ref{fig:overlap+-}. The coherence at $p$ of the relation $\dot{\precsim}$ between $\bar{f}$ and $\bar{g}$ can be stated as follows:

\begin{figure}[t]
    \centering
            \includegraphics{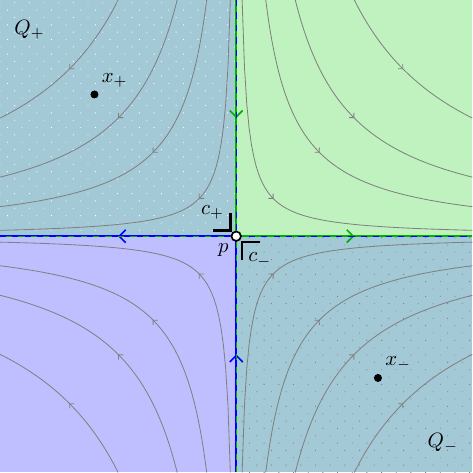}
    \caption{Overlap of $f$ and $g$ near $p$. The blue region corresponds to $D_f$, the green region corresponds to $D_g$, and their overlap is in teal. The quadrant with gray dots represents $Q_-$ while the quadrant with white dots represents $Q_+$.}
    \label{fig:overlap+-}
\end{figure}

\begin{lem} \label{lem:compatible} $f \prec_- g$ if and only if $f \prec_+ g$.
\end{lem}

\begin{proof}
    Without loss of generality, we assume $c_+ \in \mathfrak{c}$. By sliding marked points to the corners and completing, we obtain sections $\widehat{\Sigma}_f(c_\pm), \widehat{\Sigma}_g(c_\pm) \subset \mathscr{A}_p$ as in Section~\ref{sec:initcor}. Because of the shape of backward tiles passing though $p$ and of the flow of $X$ near $p$, one easily shows
    \begin{align} \label{eq:sigmafg}
        \widehat{\Sigma}_f(c_-) = \widehat{\Sigma}_f(c_+), \qquad
        \widehat{\Sigma}_g(c_-) = \widehat{\Sigma}_g(c_+).
    \end{align}
    Moreover, if $\widehat{\Sigma}_f(c_-) \subset \widehat{\Sigma}_g(c_-)$ then $f_{\vert Q_-} \leq g_{\vert Q_-}$ (and similarly for $f$ and $g$ swapped and for $c_+$). Indeed, we consider two cases.
    \begin{enumerate}
        \item There exists $y_- \in Q_-$ such that $\Sigma_f(y_-)\neq \varnothing$. Then every backward tile in $\Sigma_f(y_-)$ intersects $\{x=y=0\}$ and is geometrically below $f$. If it passes through $0$, then it induces an element in $\widehat{\Sigma}_f(c_-) \subset \widehat{\Sigma}_g(c_-)$, and this tile is also geometrically below $g$ over $Q_-$, implying $f(y'_-) \leq g(y'_-)$ for every $y'_- \in Q_-$ in the forward region of $Q_-$ containing $y_-$ by the definition of one-dimensional envelopes (see~\cite[Section 5]{BI08}). Moreover, for every $y'_- \in Q_-$ in the backward region of $Q_-$ containing $y_-$, $\Sigma_f(y'_-) \neq \varnothing$ and the same argument shows that $f(y_-) \leq g(y_-)$. Therefore, $f_{\vert Q_-} \leq g_{\vert Q_-}$.
        \item Otherwise, $\Sigma_f(y_-)= \varnothing$ for every $y_- \in Q_-$. Then no tile in $\widehat{\Sigma}_f(c_-)$ intersects $Q_- \times (-\epsilon, \epsilon)$, and $f$ is necessarily geometrically below $g$ over $Q_-$. Indeed, $f_{\vert Q_-}$ coincides with the restriction to $Q_-$ of (the underlying tile of) the forward envelope starting at $\big((\epsilon/2, 0), \varnothing\big)$ which is geometrically the lowest such tile.
    \end{enumerate}

We now show: 
\begin{align}
    \widehat{\Sigma}_g(c_\pm) \subset \widehat{\Sigma}_f(c_\pm) &\implies g\prec_\pm f, \label{eq:precpm1} \\
    \widehat{\Sigma}_f(c_\pm) \subsetneq \widehat{\Sigma}_g(c_\pm) &\implies f\prec_\pm g \label{eq:precpm2},
\end{align}
which will immediately imply the desired result. To show~\eqref{eq:precpm1}, we assume $\widehat{\Sigma}_g(c_\pm) \subset \widehat{\Sigma}_f(c_\pm)$ and we consider the order for $f$ and $g$ with respect to $\prec_-$.
\begin{enumerate}
    \item $f_{\vert Q_-} \neq g_{\vert Q_-}$. Then since $g_{\vert Q_-} \leq f_{\vert Q_-}$ by the above argument, $g \prec_- f$.
    \item $f_{\vert Q_-} = g_{\vert Q_-}$. Then $f_\star = (f, x_-)$ and $g_\star = (g, x_-)$ are comparable by either $\prec^\mathrm{geo}$ or $\prec^\mathrm{cor}$. First, for every $y_- \in Q_-$, $\Sigma_g(y_-) \subset \Sigma_f(y_-)$. Otherwise, there would exist a marked backward tile $a_\star \in \Sigma_g(y_-) \setminus \Sigma_f(y_-)$ which would geometrically coincide with $f$ and $g$ over $Q_-$, since $f \leq a \leq g$ over $Q_-$ by elementary properties of one-dimensional envelopes and their shadows (see~\cite[Section 5]{BI08}). However, sliding the marked point of $a_\star$ to the singularity at $c_-$ would produce a marked tile in $\widehat{\Sigma}_g(c_-) \setminus \widehat{\Sigma}_f(c_-)$, a contradiction. This implies that if $f_\star$ and $g_\star$ are comparable by $\prec^\mathrm{geo}$, then $g_\star \prec^\mathrm{geo} f_\star$. Otherwise, $f_\star$ and $g_\star$ are comparable by $\prec^\mathrm{cor}$ and since $c_- \notin \mathfrak{c}$, we have $g_\star \prec^\mathrm{cor} f_\star$ independently of the geometric (in)coherence of $f_\star$ and $g_\star$ at $c_-$.
\end{enumerate}
Notice that the same argument applies to $\prec_+$, since $c_+ \in \mathfrak{c}$ so the order $\prec^\mathrm{cor}$ for geometrically coherent forward envelopes meeting at $c_+$ compared to $c_-$ is reversed. Here, the choice of admissible corners `twisting' $\prec^\mathrm{cor}$ in the geometrically coherent case is essential.

The proof of~\eqref{eq:precpm2} is the same up to swapping $f$ and $g$, and the geometrically coherent case for $\prec^\mathrm{cor}$ does not arise.
\end{proof}

It is now rather straightforward to check that the collection of relations $\{ \dot{\precsim}_p\}_{p \in M}$ induces a relation $\dot{\precsim}$ on $\overline{\mathscr{U}}_\star$ which is coherent along intersections. This implies that $\dot{\precsim}$ satisfies the axioms of Definition~\ref{def:weakpref} since $\precsim$ does; the details are left to the reader. To summarize, we obtain:

\begin{prop}
The set $\overline{\mathscr{U}}$ of completed forward envelopes equipped with $\dot{\precsim}$ is a weak prefoliation. It extends $(\mathscr{U}, \precsim)$ and covers $M$.
\end{prop}

%%%%%%%%%%%%%%%%%%%%%%%%%%%%%%%%%%%%%%%%%%%%%%%%%%%%%%%%%%%%%%%%%%%%%
%%%%%%%%%%%%%%%%%%%%%%%%%%%%%%%%%%%%%%%%%%%%%%%%%%%%%%%%%%%%%%%%%%%%%
        \section{Proof of Proposition~\ref{prop:prefoil}} \label{sec:proof}

Let $\big(\mathscr{A}, {\precsim^\mathscr{A}}\big)$ be a (possibly empty) complete and backward weak prefoliation. In the previous sections, we constructed a complete and forward weak prefoliation $\big(\overline{\mathscr{U}}, \dot{\precsim}\big)$ which covers $M$. Our construction depends on auxiliary choices that are omitted from the notations. To construct an extension $\mathscr{A} \hookrightarrow \mathscr{B}$, we proceed in two steps. 
\begin{enumerate}
    \item First, we construct a weak prefoliation structure on $\mathscr{A} \cup \overline{\mathscr{U}}$ such that the inclusion $\mathscr{A} \hookrightarrow \mathscr{A} \cup \overline{\mathscr{U}}$ is an extension. 
    \item Then, we construct an extension $\mathscr{A} \cup \overline{\mathscr{U}} \hookrightarrow \mathscr{B}$ by gluing the edges and corners of tiles in $\mathscr{A}$ to suitable completed forward envelopes (the \emph{upper neighbors}), while keeping copies of completed forward envelopes. 
\end{enumerate}
%%%
        \subsection{First extension}
%%%

We define a relation $\precsim^\cup$ on $\mathscr{A} \cup \overline{\mathscr{U}}$ by adapting~\cite[Section 6.13]{BI08} and the proof of~\cite[Proposition 4.12]{BI08}. It restricts to $\precsim^\mathscr{A}$ on $\mathscr{A}_\star$ and to 
$\dot{\precsim}$ on $\overline{\mathscr{U}}_\star$. We define  $\precsim^\cup_p$ in two steps. The first one is of a geometric nature, while the second one is purely formal.

%%%
            \subsubsection{Step 1: geometric extension}
%%%

Let $p \in M$, $f_\star \in \mathscr{A}_p$, and $\bar{u}_\star \in \overline{\mathscr{U}}_p$. Here, $\bar{u}$ is the completion of a forward envelope $u \in \mathscr{U}$, $u : U \rightarrow M$ with shadow $\Sigma$.

We first use the shadow $\Sigma$ and the geometric order to define a relation $\precsim^{\cup_\circ}_p$:
    \begin{itemize}
        \item If the marked point $a \in \overline{U}$ of $\bar{u}_\star$ belongs to $U$, we set $f_\star \precsim^{\cup_\circ}_p \bar{u}_\star$ if $f_\star \in \Sigma(a)$, and $\bar{u}_\star \precsim^{\cup_\circ}_p f_\star$ otherwise.
        \item If $a \in \overline{U} \setminus U$ and one of $f_\star$ or $\bar{u}_\star$ is locally strictly above the other, their relation with respect to $\precsim^{\cup_\circ}_p$ is defined as to match the geometric order.
        \end{itemize}

We then extend $\precsim^{\cup_\circ}$ along intersections at certain boundary components of $\overline{U}$ in a coherent way. More precisely, suppose that there exists a continuous path $\gamma : [0,1] \rightarrow \overline{U}$ starting at $a \in \overline{U} \setminus U$ which lifts to $f_\star$ so that $f$ and $\bar{u}$ are comparable by $\precsim^{\cup_\circ}$ at $p' = u\circ \gamma(1)$ by Step 1, after sliding the marked points along $\gamma$ and its lift. We define the relation between $f_\star$ and $\bar{u}_\star$ with respect to $\precsim^{\cup_\circ}_p$ as to match the one with respect to $\precsim^{\cup_\circ}_{p'}$ after sliding the marked points. Notice that this step applies if the marked point of $\bar{u}_\star$ belongs to an edge connected to the initial edge of $u$ via a (double) corner, and this (double) corner lifts to $f$.

This defines a relation $\precsim^{\cup_\circ}$ on $\mathscr{A}_\star \cup \overline{\mathscr{U}}_\star$ that restricts to $\precsim^\mathscr{A}$ on $\mathscr{A}_\star$ and to 
$\dot{\precsim}$ on $\overline{\mathscr{U}}_\star$. By definition, it is coherent along intersections and is compatible with the geometric order. Moreover, it satisfies the following `pseudo-transitivity' property:

\begin{lem} \label{lem:pretrans}
    Let $a_\star, b_\star, c_\star \in \mathscr{A}_\star \cup \overline{\mathscr{U}}_\star$. If $a_\star \precsim^{\cup_\circ} b_\star$, $b_\star \precsim^{\cup_\circ} c_\star$, and $a_\star$ and $c_\star$ are comparable by $\precsim^{\cup_\circ}$, then $a_\star \precsim^{\cup_\circ} c_\star$.
\end{lem}

\begin{proof}
Let $p \in M$ and consider $a_\star, b_\star, c_\star \in \mathscr{A}_p \cup \overline{\mathscr{U}}_p$ such that $a_\star \precsim^{\cup_\circ}_p b_\star$, $b_\star \precsim^{\cup_\circ}_p c_\star$, and $a_\star$ and $c_\star$ are comparable by $\precsim^{\cup_\circ}_p$. The case $p \in M^\mathrm{un}$ is immediate. We then assume that $p \in M^\mathrm{pure}$ and we consider three cases, depending on the nature of the flow line $\gamma$ of $X$ passing through $p$.
    \begin{enumerate}[leftmargin=*]
    \item \textit{$\gamma$ is not connected to a singularity of $X$.} Notice that if $f_\star \in \mathscr{A}_p$ and $\bar{u}_\star \in \overline{\mathscr{U}}_p$ are comparable by $\precsim^{\cup_\circ}_p$, and if they only overlap along $\gamma$, then the marked point of $\bar{u}_\star$ necessarily belongs to the initial edge of $u$. Then, by considering all the possible cases depending on whether the three marked tiles belong to $\mathscr{A}_p$ or $\overline{\mathscr{U}}_p$, it is easy to check that $a_\star \precsim^{\cup_\circ}_p c_\star$.
    
    \item \textit{$\gamma$ is connected to a single singularity of $X$, which is not itself connected to a saddle singularity.} This case is slightly more complicated because of the following: if $f_\star \in \mathscr{A}_p$ and $\bar{u}_\star \in \overline{\mathscr{U}}_p$ are comparable by $\precsim^{\cup_\circ}_p$, and if they only overlap along a (concave or convex) corner $c$ containing $\gamma$, then the marked point of $\bar{u}_\star$ belongs to one of the edges at $c$, but not necessarily to the initial edge of $u$. 
    
    Let us consider the following particular case: $a_\star = \bar{u}_\star \in \overline{\mathscr{U}}_p$, $b_\star = f_\star \in \mathscr{A}_p$, and $c_\star = \bar{v}_\star \in \overline{\mathscr{U}}_p$. Moreover, the marked point of $f_\star$ belongs to a concave corner on the boundary of $f$, the marked point of $\bar{u}_\star$ belongs to the domain of $u$, but the marked point of $\bar{v}_\star$ does not belong to the domain of $v$. Arguing as in the proof of Lemma~\ref{lem:compatible}, we have that $\bar{u}_\star$ is locally geometrically below $\bar{v}_\star$. If they locally coincide, we can slide the marked points inside the forward domain of overlap of $u$ and $v$, and compare them there. They would be comparable by either $\prec^\mathrm{geo}$ or by $\prec^\mathrm{cor}$, and would be geometrically noncoherent at their initial corner of overlap. Both situations yield $\bar{u}_\star \dot{\precsim} \bar{v}_\star$ as desired.

    The other cases can be treated in a similar way.
    
    \item \textit{$\gamma$ is connected to a saddle singularity of $X$ which is itself connected to another saddle singularity.} This case is more tedious than the previous one, but can be treated in a similar way. The details are left to the reader.
    \end{enumerate}

    We are now left with the case where $p$ is a pure saddle or a quadratic-like singularity. This includes the case of a saddle connected to another saddle singularity. 
    
    Let us consider the following particular case: $p \in \Delta_\mathrm{sa}^\mathrm{pure}$ is not connected to another saddle singularity, $a_\star = \bar{u}_\star \in \overline{\mathscr{U}}_p$, $b_\star = f_\star \in \mathscr{A}_p$, and $c_\star = \bar{v}_\star \in \overline{\mathscr{U}}_p$. Moreover, the marked point of $f_\star$ belongs to a concave corner on the boundary of $f$, and $\bar{u}_\star$ and $\bar{v}_\star$ overlap with $f_\star$ near $p$ along two different edges. Then using the notations of Section~\ref{sec:comp}, $f_\star \in \widehat{\Sigma}^\llcorner_v \setminus \widehat{\Sigma}^\llcorner_u$, hence $\widehat{\Sigma}^\llcorner_u \subsetneq \widehat{\Sigma}^\llcorner_v$ and $\bar{u}_\star \precsim^\cup_p \bar{v}_\star$.
    
    The remaining cases can be treated in a similar way, and the details are left to the reader.
\end{proof}

%%%
        \subsubsection{Step 2: abstract extension}
%%%

Even after the previous step, the relation $\precsim^{\cup_\circ}_p$ between $f_\star$ and $\bar{u}_\star$ might still remain undefined if the marked points of $f_\star$ and $\bar{u}_\star$ both belong to certain edges, or if they are singularities at certain corners intersecting along edges only. We proceed in three steps:

\begin{enumerate}[leftmargin=*]
    \item Let us first assume that $p$ belongs to a pure flow line of $X$ which is not connected to any singularity of $X$. We extend $\precsim^{\cup_\circ}$ at $p$ by the following abstract procedure:
        \begin{itemize} 
            \item If there exists $a_\star \in \mathscr{A}_p \cup \overline{\mathscr{U}}_p$ such that $f_\star \precsim^{\cup_\circ}_p a_\star$ and  $a_\star \precsim^{\cup_\circ}_p \bar{u}_\star$, we set $f_\star \precsim^\cup_p \bar{u}_\star$.
            \item Otherwise, we set $\bar{u}_\star \precsim^\cup_p f_\star$.
        \end{itemize}
    By Lemma~\ref{lem:pretrans}, this definition ensures that $\precsim^\cup_p$ is a total preorder on $\mathscr{A}_p \cup \overline{\mathscr{U}}_p$.
    
    \item We now assume that $p$ belongs to a connection between two saddle singularities $q_\pm \in \Delta_\mathrm{sa}$. By sliding the marked points towards $q_\pm$ along this connection, we obtain injections
    \begin{align*}
        \varsigma_\pm : \mathscr{A}_p \cup \overline{\mathscr{U}}_p \longrightarrow \mathscr{A}_{q_\pm} \cup \overline{\mathscr{U}}_{q_\pm}
    \end{align*}
    which preserve $\precsim^{\cup_\circ}$.
    We modify the previous abstract procedure to extend $\precsim^{\cup_\circ}$ on $\mathscr{A}_p \cup \overline{\mathscr{U}}_p$. With the same notations as before, we set $f_\star \precsim^\cup_p \bar{u}_\star$ if one of the following holds:
        \begin{itemize}
            \item There exists $a_\star \in \mathscr{A}_p \cup \overline{\mathscr{U}}_p$ such that $f_\star \precsim^{\cup_\circ}_p a_\star$ and  $a_\star \precsim^{\cup_\circ}_p \bar{u}_\star$.
            \item There exists $a^\pm_\star \in \mathscr{A}_{q_\pm} \cup \overline{\mathscr{U}}_{q_\pm}$ such that $\varsigma_\pm\big(f_\star\big) \precsim^{\cup_\circ}_{q_\pm} a^\pm_\star$ and  $a^\pm_\star \precsim^{\cup_\circ}_{q_\pm} \varsigma_\pm\big(\bar{u}_\star\big)$.
        \end{itemize}
    Otherwise, we set $\bar{u}_\star \precsim^\cup_p f_\star$. We further coherently extend $\precsim^\cup$ along the flow line of $X$ passing through $p$ and at $q_\pm$ by sliding the marked points along this flow line. As in the previous item, Lemma~\ref{lem:pretrans} ensures that $\precsim^\cup_p$ is a total preorder on $\mathscr{A}_p \cup \overline{\mathscr{U}}_p$.

    \item Finally, we assume that $p \in \Delta_\mathrm{sa} \cup Q$ is a saddle or quadratic singularity. We use the same abstract procedure as in the first item to extend $\precsim^\cup$ at $p$, taking into account the extension in step 2 if $p$ is connected to another saddle singularity. As before, this yields a total preorder $\precsim^\cup_p$ on $\mathscr{A}_p \cup \overline{\mathscr{U}}_p$. We then extend $\precsim^\cup$ by coherence along the pure stable and unstable branches at $p$ (no extension is needed along a connection to another saddle singularity).
\end{enumerate}

It is now rather straightforward to prove:

\begin{lem} \label{lem:ext1}
The set $\mathscr{A} \cup \overline{\mathscr{U}}$ equipped with the relation $\precsim^\cup$ is a weak prefoliation which extends $\big( \mathscr{A}, {\precsim^\mathscr{A}}\big)$.
\end{lem}

\begin{proof}
    Routine verification left to the reader.
\end{proof}

%%%
        \subsection{Second extension}
%%%

We now define an extension $\mathscr{A} \cup \overline{\mathscr{U}} \hookrightarrow \mathscr{B}$ by adapting the strategy of the proof of~\cite[Proposition 4.11]{BI08}. To that end, we define a similar notion of \emph{upper neighbors} along the edges of tiles in $\mathscr{A}$ (see~\cite[Definition 4.10]{BI08}), and we extend every such tile by attaching its upper neighbors to its boundary components. We also keep copies of tiles in $\overline{\mathscr{U}}$ to ensure that $\mathscr{B}$ covers $M$. We then define a weak prefoliation structure on $\mathscr{B}$ in a purely formal way. As before, we will need to treat the corners with particular care.

%%%
            \subsubsection{Typology of backward corners}
%%%

\begin{figure}[t]
    \centering
        \begin{subfigure}{0.5\textwidth}
                \includegraphics{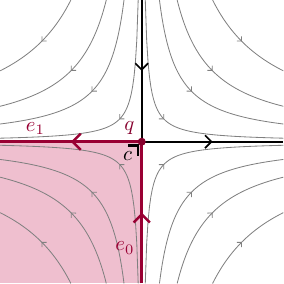}
        \centering
        \caption{Pure convex corner.}
        \end{subfigure}%
        \begin{subfigure}{0.5\textwidth}
                \includegraphics{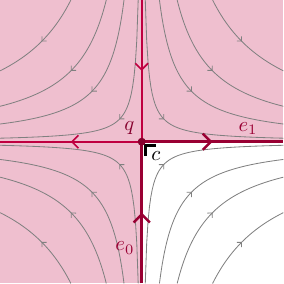}
        \centering
        \caption{Pure concave corner.}
        \end{subfigure}
        \par\bigskip
        \begin{subfigure}{0.5\textwidth}
                \includegraphics{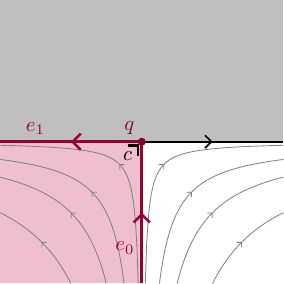}
        \centering
        \caption{Unstable convex corner.}
        \end{subfigure}%
        \begin{subfigure}{0.5\textwidth}
                \includegraphics{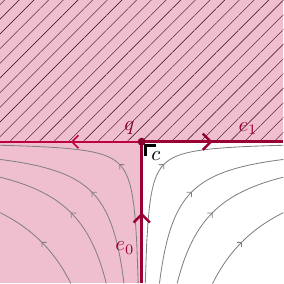}
        \centering
        \caption{Unstable concave corner.}
        \end{subfigure}
    \caption{The four types of corners of backward tiles.}
    \label{fig:fourcorners}
\end{figure}

Let $f : U \rightarrow M$ be a complete backward tile. We denote by $\mathcal{B}$ the set of boundary components of $U$, and by $\mathcal{E}$ the set of edges of $U$. A boundary component of $U$ contains between one and three edges. Two edges $e_0, e_1 \in \mathcal{E}$ are \textbf{adjacent} if $c = (f(e_0), f(e_1))$ is a (backward or forward) corner in $M$. Note the importance of the order. We further distinguish two cases:
\begin{itemize}
    \item The singularity at $c$ is a pure saddle singularity. If $c$ is a backward (resp.~forward) corner, we say that $(e_0,e_1)$ is a \textbf{pure convex (resp.~concave) corner} of $U$ (or $f$). 
    \item The singularity at $c$ is a quadratic-like singularity. If $c$ is a backward (resp.~forward) corner, we say that $(e_0,e_1)$ is an \textbf{unstable convex (resp.~concave) corner} of $U$ (or $f$). 
\end{itemize}
These four types of backward corners are illustrated in Figure~\ref{fig:fourcorners}.

We say that $(e_0,e_1)$ is an \textbf{adjoining pair of edges} if it is either a pure concave corner or an unstable (convex or concave) corner of $f$. If $e_2 \in \mathcal{E}$ and both $(e_0, e_1)$ and $(e_1, e_2)$ are adjoining pairs of edges, we call $(e_0,e_1,e_2)$ an \textbf{adjoining triple of edges}. We define a \textbf{boundary row} of $U$ (or $f$) as either
\begin{enumerate}
    \item A single edge which is not part of an adjoining pair,
    \item An adjoining pair of edges which is not part of an adjoining triple, or
    \item An adjoining triple of edges.
\end{enumerate}
We denote by $\mathcal{R}$ the set of boundary rows of $f$. Since every edge is included in a row, and every row is included in a boundary component, we have natural surjections $\mathcal{E} \twoheadrightarrow \mathcal{R} \twoheadrightarrow \mathcal{B}$. 

We also define a collection of \textbf{distinguished edges} as follows. Each boundary row $r$ contains exactly one distinguished edge, and moreover,
\begin{itemize}
    \item If $r$ is an adjoining pair of edges forming an unstable convex corner $c$, then the last edge of $r$ is distinguished.
    \item If $r$ is an adjoining pair of edges forming a forward corner $c$, then the last edge of $r$ is distinguished if $c \notin \mathfrak{c}$, and the first one is if $c \in \mathfrak{c}$. Recall that $\mathfrak{c}$ denotes the chosen admissible collection of forward corners in the definition of the weak prefoliation structure on $\mathscr{U}$.
    \item If $r$ is an adjoining triple of edges, the last edge of $r$ is distinguished.
\end{itemize}
We denote by $\mathcal{E}^\bullet \subset \mathcal{E}$ the set of distinguished edges of $f$. A row made of a single edge will be called a \textbf{trivial row}, in which case this edge is automatically distinguished.

\medskip

We now explain the relevance of these definitions. Any complete forward tile intersecting a distinguished edge of $f$ will also intersect all the other edges of $f$ in the same boundary row. Moreover, if it starts at such an edge, it will not intersect edges of $f$ in different boundary rows. We leave it to the reader to draw the corresponding pictures. Therefore, we have to be careful to define upper neighbors at edges of $f$ that are `coherent' along corners at adjoining edges.

%%%
            \subsubsection{Upper neighbors} \label{sec:upneigh}
%%%

Let $p \in M^\mathrm{pure}$, $f_\star \in \mathscr{A}_p$, $\bar{u}_\star \in \overline{\mathscr{U}}_p$, and assume that the marked points of $f_\star$ and $\bar{u}_\star$ belong to edges of the respective domains for $f$ and $\bar{u}$. Following~\cite[Definition 4.10]{BI08}, we say that $\bar{u}_\star$ is an \textbf{upper neighbor} of $f_\star$ in $\big(\mathscr{A} \cup \overline{\mathscr{U}}, {\precsim^\cup}\big)$ if $f_\star \prec^\cup_p \bar{u}_\star$, and $\bar{u}_\star$ is the lowest such element in $\mathscr{A}_p \cup \overline{\mathscr{U}}_p$, i.e.,
\begin{align}
    \forall g_\star \in \mathscr{A}_p \cup \overline{\mathscr{U}}_p, \quad f_\star \prec^\cup_p g_\star \implies \bar{u}_\star \precsim^\cup_p g_\star.
\end{align}

\begin{lem}
    With the previous notations, $f_\star$ admits an upper neighbor which is unique up to order-equivalence.
\end{lem}

\begin{proof}
    The uniqueness immediately follows from the definition. To show existence, we consider the section $\Sigma_0 \subset \mathscr{A}_p$ defined by
    \begin{align} \label{eq:sigma0}
    \Sigma_0 \coloneqq \left\{ g_\star \in \mathscr{A}_p \ \vert \ g_\star \precsim^\mathscr{A}_p f_\star \right\}.
    \end{align}
    If $p$ does not belong to the pure stable branch of a quadratic-like singularity, we define $\bar{u}_\star$ as the forward envelope starting at $\big(p, \Sigma_0\big)$ and with marked point corresponding to the initial basepoint of $\bar{u}$. We emphasize that this forward envelope might \emph{not} be the upper neighbor of $f_\star$. However, we claim that if the marked point of $f_\star$ belong to the \emph{distinguished edge} of a boundary row $r$, then $\bar{u}_\star$ is the upper neighbor of $f_\star$. To prove this, we consider three cases. 
    \begin{enumerate}[leftmargin=*]
        \item $r$ is a trivial row. By definition of $\precsim^\cup$ in Step 1 of the previous section, $f_\star \precsim^\cup \bar{u}_\star$. Then, let $g_\star \in \mathscr{A}_p \cup \overline{\mathscr{U}}_p$ such that $f_\star \prec^\cup g_\star$. We consider two cases:
            \begin{itemize}
                \item $g_\star \in \mathscr{A}_p$. Then $g_\star \notin \Sigma_0$, so $\bar{u}_\star \precsim^\cup_p g_\star$ by definition.
                \item $g_\star = \bar{v}_\star \in \overline{\mathscr{U}}_p$. By definition, there exists $\bar{w}_\star \in \overline{\mathscr{U}}_p$ such that $f_\star \precsim^\cup_p \bar{w}_\star \precsim^\cup_p \bar{v}_\star$, and the marked point $z$ of $\bar{w}$ belongs to the domain of $w$. We must have $f_\star \in \Sigma_w(z)$, otherwise $\bar{w}_\star$ would be locally geometrically below $f_\star$, so $\Sigma_0 \subset \Sigma_w(z)$. If this inclusion is strict, then $u_\star \prec^{\mathrm{geo}}_p w_\star$ and $\bar{u}_\star \dot{\precsim}_p \bar{w}_\star$, hence $\bar{u}_\star \precsim^\cup_p \bar{w}_\star \precsim^\cup_p \bar{v}_\star$. If $\Sigma_0 = \Sigma_w(z)$, then $u_\star \prec^\mathrm{inc}_p w_\star$, hence $\bar{u}_\star \precsim^\cup_p \bar{v}_\star$ as well.
            \end{itemize}
        
        \item $r = (e_0, e_1)$ is an adjoining pair of edges. We first assume that the corresponding corner $c$ is a pure forward corner. Let $g_\star \in \mathscr{A}_p \cup \overline{\mathscr{U}}_p$ such that $f_\star \prec^\cup g_\star$. As in the previous case, $g_\star \in \mathscr{A}_p$ implies $\bar{u}_\star \prec^\cup_p g_\star$. The case $g_\star \in \overline{\mathscr{U}}_p$ is slightly more complicated. As before, there exists $\bar{w}_\star \in \overline{\mathscr{U}}_p$ such that $f_\star \precsim^\cup_p \bar{w}_\star \precsim^\cup_p \bar{v}_\star$. However, the marked point of $\bar{w}_\star$ might not belong to the domain of $w$. If it does, then the previous argument still implies $\bar{u}_\star \precsim^\cup_p \bar{w}_\star$. Otherwise, we slide the marked point of $\bar{w}_\star$ along $c$, and the new marked point $z'$ now belongs to the domain of $w$. We also slide the marked point of $f_\star$ along $c$ to get a new tile $f'_\star$ passing through $p' \coloneqq w(z')$. We also define the corresponding section
        $$\Sigma'_0 \coloneqq \left\{ g'_\star \in \mathscr{A}_{p'} \ \vert \ g'_\star \precsim^\mathscr{A}_{p'} f'_\star \right\}$$
        so that $\Sigma'_0 \subset \Sigma_w(z')$. If this inclusion is strict, then $\widehat{\Sigma}_u(c) \subsetneq \widehat{\Sigma}_w(c)$, implying $\bar{u}_\star \precsim^\cup_p \bar{w}_\star$. If $\Sigma'_0 = \Sigma_w(z')$, then the definitions of $\prec^\mathrm{cor}$ and distinguished edges imply $\bar{u}_\star \precsim^\cup_p \bar{w}_\star$ as well.

        The cases where $c$ is a (forward or backward) quadratic-like corner are simpler since no forward envelope starts at the pure stable branch of the corresponding singularity. The details are left to the reader.

        \item $r = (e_0, e_1, e_2)$ is an adjoining triple of edges. There are multiple cases to consider depending on the nature of the corners corresponding to $(e_0, e_1)$ and to $(e_1, e_2)$, but all of these cases can be treated as in case 2. The details are left to the reader.
    \end{enumerate}

    If the marked point of $f_\star$ does not belong to a distinguished edge, we consider the previously constructed upper neighbor at (a point on) the distinguished edge in $r$, and we slide the marked point along the corner(s) to obtain a marked tile at $p$. By the previous arguments, it is easy to check that the latter is the upper neighbor of $f_\star$.
\end{proof}

We emphasize that by construction, the upper neighbors of $f \in \mathscr{A}$ satisfy the following coherence condition along boundary rows. Let $(e_0, e_1)$ be an adjoining pair of edges for $f$ forming a corner $c$ in $M$, and $x_0 \in e_0$, $x_1 \in e_1$ be two boundary marked points. We denote by $\bar{u}^0_\star$ and $\bar{u}^1_\star$ the upper neighbors of $(f,x_0)$ and $(f,x_1)$, respectively. Then:

\begin{lem} \label{lem:corcoh}
    Up to order-equivalence, $\bar{u}^0_\star$ and $\bar{u}^1_\star$ differ by sliding the marked point along $c$.
\end{lem}

Therefore, each boundary row $r \in \mathcal{R}$ of $f$ determines a completed forward envelope $\bar{u}^r \in \overline{\mathscr{U}}$, unique up to forward-equivalence, such that the upper neighbor of $f$ for a marked point on an edge of $r$ is $\bar{u}^r$, with a suitable marked point on the boundary component of the initial edge of $u^r$. Importantly, gluing $\bar{u}^r$ to $f$ along $r$ yields a tile.

\begin{rem}
Assume that $p \in M^\mathrm{pure}$ belongs to a closed orbit $\gamma$ of $X$ with period $T > 0$, and $\bar{u}_\star \in \overline{\mathscr{U}}_p$ is the upper neighbor of $f_\star \in \mathscr{A}_p$. In that case, $p$ belongs to the initial edge of $u$. Assume that $\bar{u}^T_\star \simeq_p \bar{u}_\star$, where we use the notation from Section~\ref{sec:prelim} for the $T$-translation. Then $u^T_\star \simeq_p u_\star$ and $\Phi^T\big(\Sigma_0\big) = \Sigma_0$, where $\Sigma_0$ is defined by~\eqref{eq:sigma0}. It immediately follows that $f^T_\star \simeq_p f_\star$. Conversely, if $f^T_\star \simeq_p f_\star$ then $\Phi^T\big(\Sigma_0\big) = \Sigma_0$ and $u^T_\star \simeq_p u_\star$ by the uniqueness of forward envelopes (see Proposition~\ref{prop:forenv}), hence $\bar{u}^T_\star \simeq_p \bar{u}_\star$.

To summarize, the upper neighbor of a tile $f_\star$ in $\mathscr{A}_p$ is order-equivalent to its $T$-translation if and only if $f_\star$ is. However, it is possible that $f^T_\star \approx_p f_\star$ without $f^T_\star \simeq_p f_\star$ ($f$ is geometrically covering itself, but not from the viewpoint of the preorder), in which case we might not even have $\bar{u}^T_\star \approx_p \bar{u}_\star$. This is because our preorders contain strictly more information than the geometric order.
\end{rem}

%%%%
            \subsubsection{Gluing upper neighbors}
%%%

We now define a complete and backward weak prefoliation $\mathscr{B}$ extending $\mathscr{A} \cup \overline{\mathscr{U}}$ by gluing the upper neighbors to the backward edges of the tiles in $\mathscr{A}$, while keeping copies of the tiles in $\overline{\mathscr{U}}$. More precisely, to each tile $f \in \mathscr{A}$ we associate a new tile $\widetilde{f}$ defined as follows. If $f$ has no edge, we set $\widetilde{f}=f$. Otherwise, for each $e \in \mathcal{E}^\bullet$, we choose a point $x_e \in e$ and write $p_e \coloneqq f(x_e)$. Then, we denote by $\bar{u}_e \in \overline{\mathscr{U}}_{p_e}$ the upper neighbor of $(f, x_e)$. It can be glued to $f$ along the row of edges containing $e$. Recall that by Lemma~\ref{lem:corcoh} above, the upper neighbor of $f$ at adjoining edges to $e$ coincides with $\bar{u}_e$ with a slid marked point. We then define $\widetilde{f}$ as the complete and forward tile obtained by gluing \emph{all} such upper neighbors along the (rows of) edges of $f$, without repetitions. By construction, there is a natural embedding $\varphi_f$ of the domain of $f$ into the domain of $\widetilde{f}$ as in Definition~\ref{def:ext}.

We now define 
$$\widetilde{\mathscr{A}} \coloneqq \big\{ \widetilde{f} \ \vert \ f \in \mathscr{A}\big\}, \qquad \mathscr{B} \coloneqq \widetilde{\mathscr{A}} \cup \overline{\mathscr{U}}.$$ 
It comes with a natural map $i : \mathscr{A} \cup \overline{\mathscr{U}} \rightarrow \mathscr{B}$ which restricts to the identity on $\overline{\mathscr{U}}$. By construction, the tiles in $\mathscr{B}$ are complete and forward, and they cover $M$. It remains to define a weak prefoliation structure $\precsim^\mathscr{B}$ on $\mathscr{B}$ making $i$ an extension. To this end, we proceed as in the proof of~\cite[Proposition 4.11]{BI08} in three steps.

\begin{enumerate}[leftmargin=*]
    \item We first define the restriction of $\precsim^\mathscr{B}$ to $\widetilde{\mathscr{A}}$. Let $\widetilde{f}_\star \in \widetilde{\mathscr{A}}_\star$, where $f \in \mathscr{A}$. We consider two cases:
        \begin{enumerate}
            \item The marked point belongs to the domain of $f$. This induces a marked tile $f_\star \in \mathscr{A}_\star$.
            \item The marked point does not belong to the domain of $f$, but to the domain of a completed forward envelope $\bar{u} \in \overline{\mathscr{U}}$. This induces a marked tile $\bar{u}_\star \in \overline{\mathscr{U}}_\star$.
        \end{enumerate}
This allows us to define a `restriction' map $\rho : \widetilde{\mathscr{A}}_\star \rightarrow \mathscr{A}_\star \cup \overline{\mathscr{U}}_\star$ as 
\begin{itemize}
\item $\rho \big(\widetilde{f}_\star\big) = f_\star$ in case (a), 
\item $\rho \big(\widetilde{f}_\star\big) = \bar{u}_\star$ in case (b).
\end{itemize}
Then, we define $\precsim^\mathscr{B}$ on $\widetilde{\mathscr{A}}$ as the pullback of $\precsim^\cup$ along $\rho$: 
$$\forall \widetilde{f}_\star, \widetilde{g}_\star \in \widetilde{\mathscr{A}}_\star, \quad \widetilde{f}_\star \precsim^\mathscr{B} \widetilde{g}_\star \iff \rho\big(\widetilde{f}_\star\big) \precsim^\cup \rho\big(\widetilde{g}_\star\big).$$

    \item We now define $\precsim^\mathscr{B}$ between elements $\widetilde{f}_\star \in \widetilde{\mathscr{A}}_\star$ and $\bar{v}_\star \in \overline{\mathscr{U}}_\star$. If $\rho\big(\widetilde{f}_\star\big) = f_\star \in \mathscr{A}_\star$, we define the order between $\widetilde{f}_\star$ and $\bar{v}_\star$ for $\precsim^\mathscr{B}$ to match the one between $f_\star$ and $\bar{v}_\star$ for $\precsim^\cup$. If $\rho\big(\widetilde{f}_\star\big) = \bar{u}_\star \in \overline{\mathscr{U}}_\star$, we set
\begin{align*}
    \widetilde{f}_\star \precsim^\mathscr{B} \bar{v}_\star \iff \bar{u}_\star \ \dot{\precsim} \ \bar{v}_\star, \\
    \bar{v}_\star \precsim^\mathscr{B} \widetilde{f}_\star \iff \bar{v}_\star \ \dot{\prec} \ \bar{u}_\star.
\end{align*}
That way, if $\bar{u}_\star \simeq \bar{v}_\star$, then $\widetilde{f}_\star \prec^\mathscr{B} \bar{v}_\star$.

    \item Finally, the restriction of $\precsim^\mathscr{B}$ on $\overline{\mathscr{U}}_\star$ is defined to agree with $\dot{\precsim}$.
\end{enumerate}

By virtue of the properties of upper neighbors, it is a rather mechanical exercise to check:

\begin{lem} \label{lem:ext2}
The set $\mathscr{B}$ equipped with the relation $\precsim^\mathscr{B}$ is a weak prefoliation, and the map $i : \mathscr{A} \cup \overline{\mathscr{U}} \rightarrow \mathscr{B}$ is an extension.
\end{lem}

\begin{proof}
Routine verification left to the reader.
\end{proof}

%%%
        \subsection{Finishing the proof}
%%%

We conclude the proof of Proposition~\ref{prop:prefoil}. Starting from a complete and backward weak prefoliation $\mathscr{A}$, the desired forward extension $\mathscr{A} \hookrightarrow \mathscr{B}$ is obtained as the composition of the extensions $\mathscr{A} \hookrightarrow \mathscr{A} \cup \overline{\mathscr{U}}$ from Lemma~\ref{lem:ext1} and $\mathscr{A} \cup \overline{\mathscr{U}} \hookrightarrow \mathscr{B}$ from Lemma~\ref{lem:ext2}. Finally, the case where $\mathscr{A}$ is \emph{forward} instead of backward is obtained by reversing the (co)orientation of $\eta$ and applying the previous results. \qed

%%%%%%%%%%%%%%%%%%%%%%%%%%%%%%%%%%%%%%%%%%%%%%%%%%%%%%%%%%%%%%%%%%%%%%%%%%%%%%%
\newpage
\appendix
%%%%%%%%%%%%%%%%%%%%%%%%%%%%%%%%%%%%%%%%%%%%%%%%%%%%%%%%%%%%%%%%%%%%%%%%%%%%%%%

%%%%%%%%%%%%%%%%%%%%%%%%%%%%%%%%%%%%%%%%%%%%%%%%%%%%%%%%%%%%%%%%%%%%%%%%%%%%%%%%%%%%%%%%%%%%%%%%%%%%%%%%%%%%%%%%%%%%%%%%%%%%%%%%%%%%%%%%%%%%%%%%%%%%%%%%%%%%%%
    \section{Key ODE lemma} \label{appsec:ODE}
%%%%%%%%%%%%%%%%%%%%%%%%%%%%%%%%%%%%%%%%%%%%%%%%%%%%%%%%%%%%%%%%%%%%%%%%%%%%%%%%%%%%%%%%%%%%%%%%%%%%%%%%%%%%%%%%%%%%%%%%%%%%%%%%%%%%%%%%%%%%%%%%%%%%%%%%%%%%%%

Let $D \subset \R^N$ be a compact subset of $\R^N$, for some integer $N$. We consider a smooth function 
$$ \begin{array}{rrcl}
F : & D \times \R \times \R & \longrightarrow & \R  \\
    & (x,y,t) & \longmapsto & F(x,y,t) = F_x(y,t) 
\end{array}$$
satisfying
\begin{enumerate}
    \item[(C1)] There exists $\epsilon > 0$ such that $$\partial_y F_x \geq \epsilon,$$
    \item[(C2)] There exists $C>0$ such that for every $(x,y,t) \in D \times \R \times \R$, $$\big\vert F_x(y,t) - F_x(y,0) \big\vert \leq C.$$
\end{enumerate}

We consider the following family of ODEs parametrized by $x \in D$:
\begin{equation*}
  (E_x) : \left\{
    \begin{aligned}
    \dot{y}(t) &= F_x\big(y(t),t\big), \\
    y(0) &= y_0.
    \end{aligned}
  \right.
\end{equation*}

In Part~\ref{part:I}, we repeatedly use the following elementary lemma.
\begin{lem} \label{lem:ODE}
    For every $x \in D$, there exists a unique $y_0(x) \in \R$ such that the maximal solution to $(E_x)$ with initial value $y_0(x)$ is defined on $\R$ and is bounded. Moreover, the map 
    $$\begin{array}{rrcl}
    y_0 : & D & \longrightarrow & \R\\
    & x & \longmapsto &  y_0(x) 
    \end{array}$$ is continuous.
\end{lem}

\begin{proof}
We first show the existence and uniqueness of $y_0 = y_0(x)$ for a fixed $x \in D$. We will write $F(y,t) = F_x(y,t)$ for simplicity. The solutions to $(E_x)$ satisfy the two key properties:
\begin{enumerate}
\item[(P1)] If $y_0 < \widetilde{y}_0$, and $y, \widetilde{y} : [0,T] \rightarrow \R$ are solutions to $(E_x)$ with initial values $y_0$ and $\widetilde{y}_0$, respectively, then 
\begin{align} \label{ineq:exp}
\forall t \in [0,T], \quad \widetilde{y}(t) - y(t) &\geq (\widetilde{y}_0 - y_0) \, e^{\epsilon t}.
\end{align}
\item[(P2)] There exists $A > 0$ such that if $y : [-T, T] \longrightarrow \R$ is a solution of $(E_x)$, then
\begin{align*}
    y(t) > A &\Longrightarrow \dot{y}(t) > 1, \\
    y(t) < -A &\Longrightarrow \dot{y}(t) < -1.
\end{align*}
\end{enumerate}

Property (P1) easily follows from condition (C1). Under the hypothesis of (P1), we readily have $$\forall t \in [0, T], \quad \widetilde{y}(t) > y(t).$$
Then, writing $z(t) \coloneqq \widetilde{y}(t) - y(t)$, we have
\begin{align*}
\dot{z}(t) &= F(\widetilde{y}(t),t) - F(y(t), t) \\
            &= \int_{y(t)}^{\widetilde{y}(t)} \partial_y F(y, t) \,dy \\
            &\geq \epsilon (\widetilde{y}(t) - y(t)) \\
            &= \epsilon z(t),
\end{align*}
and~\eqref{ineq:exp} follows.

On the other hand, property (P2) easily follows from conditions (C1) and (C2). By condition (C1), $F(\, \cdot \, , 0) : \R \rightarrow \R$ is a strictly increasing diffeomorphism and there exists $A > 0$ such that 
\begin{align*}
\forall y > A, \quad F(y,0) &> C + 1, \\
\forall y < -A, \quad F(y,0) &< -C - 1.
\end{align*}
Then, if $y(t) > A$,  
\begin{align*}
\dot{y}(t) &= F(y,t) - F(y,0) + F(y,0) \\
&\geq -C + F(y,0) \\
& > -C+ C +1 \\
&=1.
\end{align*}
The other implication holds similarly. Here, we can further assume that $A$ is independent of $x$ by simply considering
$$A \coloneqq \epsilon^{-1}\left(1 + C + \sup_{x \in D} \vert F_x(0,0)\vert \right).$$

Notice that a solution $y$ of  $(E_x)$ is defined for all negative times, and $\vert y \vert$ is bounded by $\max(A, \vert y_0\vert )$ on $\R_{\leq 0}$. However, such a solution might be unbounded on $\R_{\geq 0}$ and might blow up in finite time. We define
\begin{itemize}
\item $I_+ = I_+(x) \subset \R$ the set of $y_0 \in \R$ such that the solution $y$ to $(E_x)$ with initial value $y_0$ is \emph{not} bounded from above on $\R_{\geq 0}$,
\item $I_- = I_-(x) \subset \R$ the set of $y_0 \in \R$ such that the solution $y$ to $(E_x)$ with initial value $y_0$ is \emph{not} bounded from below on $\R_{\geq 0}$.
\item $I_0 = I_0(x) \subset \R$ the set of $y_0 \in \R$ such that the solution $y$ to $(E_x)$ with initial value $y_0$ is bounded on $\R_{\geq 0}$.
\end{itemize}
By definition, we have $$\R = I_- \sqcup I_0 \sqcup I_+.$$
By property (P2), if $y_0 \in I_\pm$ then the corresponding solution $y$ diverges to $\pm \infty$, possibly in finite time. Moreover, $I_+$ is of the form $[a_+, +\infty)$ or $(a_+, +\infty)$, $I_-$ is of the form $(-\infty, a_-]$ or $(-\infty, a_-)$, and property (C1) implies that $I_0$ contains at most one point. It is therefore enough to show that $I_\pm$ are open. Let $y_0 \in I_+$, and let $y$ be the corresponding maximal solution, defined on $(-\infty, T)$ for some $0 < T \leq +\infty$. By definition of $I_+$, there exists $0 <  t_0 < T$ such that $y(t_0) > A$. By the continuous dependence of solutions to an ODE with respect to the initial value, there exists $\delta > 0$ such that for every $\widetilde{y}_0 \in (y_0 - \delta, y_0+ \delta)$, the solution $\widetilde{y}$ to $(E_x)$ with initial value $\widetilde{y}_0$ also exists up to time at least $t_0$, and $\widetilde{y}(t_0) > A$. Therefore, by property (P2), $\widetilde{y}$ diverges to $+\infty$ (possibly in finite time), and $\widetilde{y}_0 \in I_+$. It follows that $I_+$ is open, and $I_-$ is open for similar reasons.

We now have a well-defined map $x \in D \mapsto y_0(x)$, and we are left to show that it is continuous. Let $x \in D$. If $\widetilde{y}_0 > y_0(x)$, then $\widetilde{y}_0 \in I_+(x)$ and the maximal solution $\widetilde{y}$ to $(E_x)$ with initial value $\widetilde{y}_0$ diverges to $+\infty$, possibly in finite time. Hence, there exists $t_0 >0$ such that $\widetilde{y}(t_0)$ is defined and $\widetilde{y}(t_0) > A + 1$. By the continuous dependence of solutions to a family of ODEs with respect to a parameter (the family of ODEs depend smoothly on $x$), there exists an open neighborhood $U$ of $x$ in $D$ such that for every $\widetilde{x} \in U$, every maximal solution $\widetilde{y}$ to the ODE $(E_{\widetilde{x}})$ with initial value $\widetilde{y}_0$ is defined up to time at least $t_0$ and satisfies $\widetilde{y}(t_0) > A$. Therefore, by property (P2), $\widetilde{y}$ diverges to $+\infty$ and $\widetilde{y}_0 \in I_+(\widetilde{x})$, so $\widetilde{y}_0 > y_0(\widetilde{x})$. To summarize, we have shown that for every $x \in D$ and every $\widetilde{y}_0 > y_0(x)$, there exists an open neighborhood $U$ of $x$ in $D$ such that for every $\widetilde{x} \in U$, $\widetilde{y}_0 > y_0(\widetilde{x})$. This means that $y_0$ is upper semicontinuous, and it is lower semicontinuous for similar reasons. 
\end{proof}

The same methods also show that the map $y_0 : D \rightarrow \R$ depends continuously on $F$ in the appropriate topology. More precisely, let $F$ be as before, and let $(F_n)_{n \geq 1}$ be a sequence of smooth functions $F_n : D \times \R \times \R \rightarrow \R$ satisfying the following:
\begin{itemize}
    \item For every $n \geq 1$, $F_n$ satisfies conditions (C1) and (C2) for the \underline{same} constants $\epsilon, C$ as $F$,
    \item $(F_n)_{n \geq 1}$ converges to $F$ in $C^\infty_\mathrm{loc}$.
\end{itemize}
We denote by $(E^n_x)$ the family of ODEs corresponding to $F_n$.

\begin{lem} \label{lem:ODE2}
With the above notations, let $y^n_0 : D \rightarrow \R$ denote the continuous map given by Lemma~\ref{lem:ODE} for $F_n$, $n \geq 1$. Then the sequence $(y^n_0)_{n \geq 1}$ converges uniformly to $y_0$ in the $C^0$ topology.
\end{lem}

\begin{proof}
For $n \geq 1$, we denote by $y_n(x,t)$ the unique bounded maximal solution of $(E^n_x)$, so that $y^n_0(x) = y_n(x, 0)$. By the proof of Lemma~\ref{lem:ODE}, there exists a constant $A > 0$, depending on $\epsilon$ and $C$ but \emph{not on $x$ nor $n$}, such that
$$\forall n \geq 1, \ \forall t \in \R, \qquad \vert y_n(x,t) \vert \leq A.$$

Let us assume by contradiction that $(y^n_0)_{n \geq 1}$ does not converge to $y_0$. Up to passing to a subsequence, we can assume that there exists a sequence $(x_n)_{n \geq 1}$, $x_n \in D$, such that 
\begin{itemize}
    \item $(x_n)_{n \geq 1}$ converges to some $x_\infty \in D$,
    \item $\big(y^n_0(x_n)\big)_{n \geq 1}$ converges to some $y_\infty \in \R$,
    \item $y_\infty \neq y_0(x_\infty)$.
\end{itemize}
Without loss of generality, we assume that $y_\infty > y_0(x_\infty)$. By the proof of Lemma~\ref{lem:ODE}, the solution $\widetilde{y}=\widetilde{y}(t)$ of $(E_{x_\infty})$ starting at $y_\infty$ diverges to $+ \infty$, possibly in finite time, so there exists $t_0 > 0$ such that $\widetilde{y}(t_0) = A+1$. Since $F_n$ converges to $F$ in $C^\infty_\mathrm{loc}$, the continuous dependence of solutions of a family of ODEs with respect to parameters implies that for $n \geq 1$ sufficiently large, $y_n(x_n, t_0) > A$, a contradiction.
\end{proof}

We end this appendix with a classical result about ODEs with solutions blowing up in finite time.

\begin{lem} \label{lem:ODE3}
    Let $a > 0 $ and $y : [0,a) \rightarrow \R$ be a differentiable function satisfying 
    \begin{equation*}
    \left\{
    \begin{aligned}
    \dot{y}(t) &\geq y^2(t), \\
    y(0) &\geq 1.
    \end{aligned}
  \right.
\end{equation*}
Then $a \leq 1$.
\end{lem}

\begin{proof}
    The conditions readily imply that $y$ is positive on $[0,a)$. Let us assume by contradiction that $a > 1$, so that $y$ is defined on $[0,1]$. Let $z(t) \coloneqq \frac{1}{y(t)} + t$. Then $z$ is differentiable and satisfies $\dot{z}(t) \leq 0$, hence $\frac{1}{y(t)} + t = z(t) \leq z(0) = \frac{1}{y(0)} \leq 1$. In particular, for $t=1$, we get $0 < \frac{1}{y(1)} \leq 0$, a contradiction.
\end{proof}

%%%%%%%%%%%%%%%%%%%%%%%%%%%%%%%%%%%%%%%%%%%%%%%%%%%%%%%%%%%%%%%%%%%%%%%%%%%%%%%%%%%%%%%%%%%%%%%%%%%%%%%%%%%%%%%%%%%%%%%%%%%%%%%%%%%%%%%%%%%%%%%%%%%%%%%%%%%%%%
        \section{Taut plane fields} \label{appsec:taut}
%%%%%%%%%%%%%%%%%%%%%%%%%%%%%%%%%%%%%%%%%%%%%%%%%%%%%%%%%%%%%%%%%%%%%%%%%%%%%%%%%%%%%%%%%%%%%%%%%%%%%%%%%%%%%%%%%%%%%%%%%%%%%%%%%%%%%%%%%%%%%%%%%%%%%%%%%%%%%%

In this appendix, $M$ still denotes a closed, oriented $3$-manifold. We extend the definition of tautness to arbitrary continuous cooriented plane fields on $M$ and show some equivalent characterizations. This generalization is motivated by the following two observations:
\begin{enumerate}
    \item There exist several \emph{nonequivalent} notions of tautness for $C^0$-foliations; see~\cite{CKR19}. The most robust one, \emph{everywhere tautness}, can be phrased uniquely in terms of the tangent plane field of the foliation, i.e., it does not refer to the leaves.
    \item For the purpose of translating properties of foliations into properties of contact structures and vice versa, it is convenient to work directly at the level of plane fields.
\end{enumerate}

\begin{defn} \label{def:tautdist}
A cooriented continuous plane field $\eta$ on $M$ is \textbf{(everywhere) taut} if for every $p \in M$, there exists a smooth closed curve transverse to $\eta$ and passing through $p$.
\end{defn}

\begin{defn}
A \textbf{dead-end component} for $\eta$ is a $C^1$, codimension-$0$ submanifold with boundary $N \subsetneq M$ such that $\partial N$ is tangent to $\eta$ and the co-orientation of $\eta$ points inwards along $\partial N$. 
\end{defn}

The next proposition is well-known if $\eta$ is (uniquely) integrable, and most of the proof immediately extends to the general case.

\begin{prop} \label{prop:tautdistrib}
Let $\eta$ be a cooriented continuous plane field on $M$. The following are equivalent:
\begin{enumerate}
\item $\eta$ is taut,
\item For all points $p, q \in M$, there exists a smooth curve from $p$ to $q$ positively transverse to $\eta$,
\item There exists a (smooth, or $C^1$) closed $2$-form $\omega$ such that $\omega_{\vert \eta} > 0$,
\item There exists a (smooth, or $C^1$) volume preserving vector field transverse to $\eta$,
\item There exists a continuous Riemannian metric on $M$ and a (smooth, or $C^1$) closed $2$-form $\omega$ which calibrates $\eta$,
\item There are no dead-end components for $\eta$.
\end{enumerate}
In particular, a surface tangent to $\eta$ is a stable minimal surface for a Riemannian metric as in 5.
\end{prop}

\begin{proof}
The implications 
$$\mathit{2 \Longrightarrow 1 \Longrightarrow 3 \Longleftrightarrow 4 \Longleftrightarrow 5 \Longrightarrow 6}$$ are standard and follow the same proofs as in the (uniquely) integrable case. We show $\mathit{6 \Longrightarrow 2}$

Let $p_0 \in M$ and denote by $A_+ \coloneqq A_+(p_0)$ the set of points $q \in M$ such that there exists a smooth curve from $p_0$ to $q$ positively transverse to $\eta$. Let us assume that $A_+ \neq M$. We show that $N \coloneqq \overline{A_+}$ is a dead-end component for $\eta$.

Since $\eta$ is continuous, every point $q \in M$ has a neighborhood $U \cong D^2 \times (-1,1)$ with smooth local coordinates $(x,y,z)$ in which  $\eta$ is given by 
$$ \eta = \ker \left(dz + f dx +  g dy\right),$$
for continuous functions $f, g : D^2 \times (-1,1) \rightarrow \R$ with $f(0) = g(0) = 0$ and $\vert f \vert, \vert g \vert \leq 1/10$. We call such a neighborhood $U$ a \emph{standard neighborhood of $q$}. Moreover, if $q \in A_+$ and $\gamma : [0,1] \rightarrow M$ is a smooth curve from $p_0$ to $q$ transverse to $\eta$, we can find a standard neighborhood of $q$ in which $\gamma$ is contained in $\{x= y = 0, \ z \leq 0\}$. It follows that the upper open cone in $U$ centered at $(0, 0, -1/2)$ defined by $$z > \sqrt{x^2 + y^2} -1/2$$
is contained in $A_+$ and contains $q \cong (0,0,0)$. This shows that $A_+$ is open.

Since $A_+ \neq M$ by assumption and $M$ is connected, $\partial N = \partial A_+ = \overline{A_+} \setminus A_+ \neq \varnothing.$ Let $q \in \partial N$ and $U$ be a standard neighborhood of $q$ as above. Since $\partial_z$ is positively transverse to $\eta$ in this neighborhood, there exists a map $h : D^2 \rightarrow (-1,1)$ such that
\begin{align*}
N \cap U  &= \left\{ h(x,y) \leq z \right\}.
\end{align*}
This map can be explicitly defined as
$$h(x,y) \coloneqq \inf \big\{z \in (-1,1) \, : \, (x,y,z) \in A_+ \cap U \big\}$$
and is lower semi-continuous. For $p \in U$, we consider the upper and lower half-cones $$C_\pm(p) \coloneqq p+ \left\{ \sqrt{x^2+y^2} \leq \pm z \right\}.$$
We have the following elementary properties:
\begin{itemize}
    \item If $p \in A_+ \cap U$ then $C_+(p) \subset A_+$,
    \item If $p \in U \setminus A_+$ then  $C_-(p) \cap A_+ = \varnothing$.
\end{itemize}
These imply that if $p \in \partial N$ then $C_-(p) \cap N = \{p\}$ and $\mathrm{int} \, C_+(p) \subset N$. Therefore, $h$ is $1$-Lipschitz, hence differentiable almost everywhere. We now show that $h$ is $C^1$ and tangent to $\eta$. Without loss of generality, we can assume that $h$ is differentiable at $(0,0)$. Note that $\eta(0,0,0) = \mathrm{span} \{ \partial_x, \partial_y\}$,  $h(0,0)=0$, and $\vert \partial_x h(0,0) \vert, \vert \partial_y h(0,0) \vert \leq 1$. These properties are satisfied \emph{in any standard neighborhood of $q$}, which allows us to change the $z$-direction (after possibly shrinking $U$). Therefore, $\partial_x h(0,0) = \partial_yh(0,0) = 0$ and (the graph of) $h$ is tangent to $\eta$ at $(0,0)$. Similarly, $h$ is tangent to $\eta$ at every point where it is differentiable. Since $\eta$ is continuous and $h$ is Lipschitz, $h$ is $C^1$ and tangent to $\eta$ everywhere. 

As a result, $N$ is a dead-end component for $\eta$.
\end{proof}

The fourth characterization of tautness immediately implies that it is an open condition in the $C^0$ topology. Moreover, the sixth characterization generalizes the well-known fact that a $C^1$-foliation without spherical and toroidal leaves is taut:

\begin{cor}
If $\eta$ has no closed integral surfaces which are spheres or tori, then it is taut.
\end{cor}

\begin{proof}
    Let us assume that $\eta$ is not taut and let $N$ be a dead-end component. We write $$\partial N = \bigsqcup_{1 \leq i \leq k} \Sigma_i,$$
    where the ${\Sigma_i}'s$ are closed connected orientable surfaces tangent to $\eta$. There exists a smooth vector field tansverse to $\eta$, which restricts to a nowhere vanishing vector field on $N$ transverse to $\partial N$ and inward pointing. By~\cite{P68} (see also~\cite[Proposition 1.2]{G75}), $$\sum_{i=1}^k \chi(\Sigma_i) = 0,$$ which is impossible if the genera of all of the ${\Sigma_i}$'s at least $2$.
\end{proof}

%%%%%%%%%%%%%%%%%%%%%%%%%%%%%%%%%%%%%%%%%%%%%%%%%%%%%%%%%%%%%%%%%%%%%%%%%%%%%%%%%%%%%%%%%%%%%%%%%%%%%%%%%%%%%%%%%%%%%%%%%%%%%%%%%%%%%%%%%%%%%%%%%%%%%%%%%%%%%%
        \section{Approximation lemma}
%%%%%%%%%%%%%%%%%%%%%%%%%%%%%%%%%%%%%%%%%%%%%%%%%%%%%%%%%%%%%%%%%%%%%%%%%%%%%%%%%%%%%%%%%%%%%%%%%%%%%%%%%%%%%%%%%%%%%%%%%%%%%%%%%%%%%%%%%%%%%%%%%%%%%%%%%%%%%%

In this appendix, we prove a technical lemma that is used in the proof of Proposition~\ref{prop:Reeb}. This is a generalization of~\cite[Lemma 4.3]{H22a} (see also~\cite[Appendix A]{Mas22}) which might be of independent interest. We do not need to assume that $M$ is $3$-dimensional.

\begin{lem} \label{lem:smooth}
Let $X$ be a smooth vector field on a closed manifold $M$. Let $f : M \rightarrow \R$ be a continuous function which is continuously differentiable along $X$, i.e., $X \cdot f$ is defined on $M$ and continuous. For every $\epsilon > 0$, there exists a smooth function $\widetilde{f} : M \rightarrow \R$ such that $$\big\vert f - \widetilde{f}\big\vert_{C^0} \leq \epsilon, \qquad \big\vert X\cdot (f - \widetilde{f})\big\vert_{C^0} \leq \epsilon.$$ 
\end{lem}

\begin{proof} Let $$\Delta \coloneqq \{ x \in M : X(x) =0\}$$ be the singular set of $X$, which is a compact subset of $M$. We proceed in 3 steps.

\begin{itemize}[leftmargin=*]
\item \textit{Step 1 : smoothing near $\Delta$.} We fix some arbitrary metric $g$ on $M$. For $\delta > 0$, let $U_\delta(\Delta)$ be the $\delta$-neighborhood of $\Delta$ for this metric. Recall the following standard fact: there exists some constant $C_0 > 0$ such that for every sufficiently small $\delta > 0$, there exists a smooth cutoff function $\varphi_\delta : U_\delta(\Delta) \rightarrow [0,1]$ such that 
\begin{itemize}
\item $\varphi_\delta \equiv 1$ on $U_{\delta/10}(\Delta)$,
\item $\varphi_\delta \equiv 0$ on $U_{\delta}(\Delta) \setminus U_{9\delta/10}(\Delta)$,
\item $\vert d\varphi_\delta \vert \leq \frac{C_0}{\delta}$.\footnote{To see this, first embed $M$ in some $\R^N$ and mollify a suitable Lipschitz cutoff function constructed from the distance function to the image of $\Delta$.}
\end{itemize}
Moreover, there exists a constant $C_1>0$ such that for every $x \in M$,
$$\vert X(x) \vert \leq C_1 \, \mathrm{dist}(x, \Delta).$$
It follows that there exists a constant $C_2 > 0$ such that for every $\delta > 0$, $$ \vert X \cdot \varphi_\delta \vert \leq C_2.$$ 

Let $f_0 : M \rightarrow \R$ be a smooth function so that $$\vert f-f_0\vert_{C^0} \leq \epsilon.$$
Since $X \cdot f = X \cdot f_0 = 0$ on $\Delta$, there exists $\delta > 0$ such that $$\vert X \cdot (f-f_0) \vert_{C^0(U_\delta(\Delta))} \leq \epsilon.$$
We now write $U_0 \coloneqq U_\delta(\Delta)$, $V_0 \coloneqq U_{\delta/10}(\Delta) \subset U_0$ and $\psi_0 = \varphi_\delta$ for such a $\delta > 0$. It follows from the previous paragraph that
$$\vert X \cdot (\psi_0 f - \psi_0 f_0)\vert_{C^0(U_0)} \leq \vert X \cdot \psi_0\vert_{C^0} \vert f - f_0 \vert_{C^0} + \vert X \cdot (f-f_0) \vert_{C^0(U_0)} \lesssim \epsilon,$$
where the symbol $\lesssim$ means ``less than
or equal to, up to a constant factor that only depends on $M$, $X$ and $g$''. We define $\widetilde{f}_0 \coloneqq \psi_0 f_0$.
\item \textit{Step 2: smoothing away from $\Delta$.} As a consequence of the proof of~\cite[Lemma 4.3]{H22a} (see also~\cite[Appendix A]{Mas22}), there exists
\begin{itemize}
\item A finite collection $\mathcal{U}= (U_i)_{1 \leq i \leq N}$ of open subsets of $M$ covering $M \setminus V_0$ (made of flow boxes for $X$),
\item A family of cutoff functions $(\varphi_i)_{1 \leq i \leq N}$ subordinate to $\mathcal{U}$ such that on $M \setminus V_0$, $$\sum_{i=1}^N \varphi_i = 1,$$ where $\varphi_i : M \rightarrow [0,1]$,
\item Smooth functions $f_i : M \rightarrow \R$, $1 \leq i \leq N$, such that $$\big\vert \varphi_i f - f_i\vert_{C^0} \leq  \frac{\epsilon}{N}, \qquad \big\vert X \cdot \big(\varphi_i f - f_i\big)\big\vert_{C^0} \leq \frac{\epsilon}{N}.$$
\end{itemize}
For $1 \leq i \leq N$, we define $\psi_i \coloneqq (1-\psi_0) \varphi_i$, so that $(\psi_i)_{0 \leq i \leq N}$ is a partition of unity subordinate to the open cover $\mathcal{U}' \coloneqq (U_i)_{0 \leq i \leq N}$ of $M$. We also write $\widetilde{f}_i \coloneqq (1-\psi_0) f_i$ for $1 \leq i \leq N$, so that 
$$\big \vert \psi_i f - \widetilde{f}_i \big\vert_{C^0} \leq \frac{\epsilon}{N}, \qquad \big \vert X \cdot \big(\psi_i f - \widetilde{f}_i \big) \big\vert_{C^0} \lesssim \frac{\epsilon}{N}.$$
\item \textit{Step 3 : putting things together.} Finally, we define
$$\widetilde{f} \coloneqq \sum_{i = 0}^N \widetilde{f}_i,$$ and we readily obtain \begin{align*}
\big\vert f - \widetilde{f} \big\vert_{C^0} &\leq \sum_{i=0}^N \big\vert \psi_i f - \widetilde{f}_i \big\vert_{C^0} \lesssim \big\vert \psi_0 f - \widetilde{f}_0 \big\vert_{C^0} + \epsilon \lesssim \epsilon, \\
\big\vert X \cdot \big(f - \widetilde{f} \big) \big\vert_{C^0} &\leq \sum_{i=0}^N \big\vert X \cdot \big( \psi_i f - \widetilde{f}_i \big) \big\vert_{C^0} \lesssim \big\vert X \cdot \big( \psi_0 f - \widetilde{f}_0 \big) \big\vert_{C^0} + \epsilon \lesssim \epsilon.
\end{align*} \qedhere
\end{itemize}
\end{proof}

%%%%%%%%%%%%%%%%%%%%%%%%%%%%%%%%%%%%%%%%%%%%%%%%%%%%%%%%%%%%%%%%%%%%%%%%%%%%%%%%%%%%%%%%%%%%%%%%%%%%%%%%%%%%%%%%%%%%%%%%%%%%%%%%%%%%%%%%%%%%%%%%%%%%%%%%%%%%%%
        \section{Periodic solutions to ODEs on the cylinder} \label{sec:ODEper}
%%%%%%%%%%%%%%%%%%%%%%%%%%%%%%%%%%%%%%%%%%%%%%%%%%%%%%%%%%%%%%%%%%%%%%%%%%%%%%%%%%%%%%%%%%%%%%%%%%%%%%%%%%%%%%%%%%%%%%%%%%%%%%%%%%%%%%%%%%%%%%%%%%%%%%%%%%%%%%

This appendix concerns the topology of periodic orbits of \emph{continuous} vector fields on the cylinder. In particular, we give a necessary and sufficient condition for the existence of a foliation by circles tangent to (an approximation of) such a vector field. This serves as evidence for Conjecture~\ref{conj:tight}, which would require a version of this result for continuous plane fields on a thickened $2$-torus $\mathbb{T}^2 \times (-1,1)$. Our result only applies to the intersection of such a plane field with an annulus transverse to the core torus $\mathbb{T}^2 \times \{0\}$.

\medskip

Let
$$ \begin{array}{crcl}
F : &\R \times \R & \longrightarrow & \R  \\
    & (x,t) & \longmapsto & F(x,t)
\end{array}$$
be a \emph{continuous}, bounded function satisfying
$$
\forall x, t\in \R, \quad F(x,t+1) = F(x,t),
$$
so that $F$ naturally induces a continuous function $\overline{F} : C \rightarrow \R$ on the cylinder $C \coloneqq \R \times \R \slash \Z$. The corresponding vector field $X$ on $C$ is defined by
$$X(x,t) \coloneqq \overline{F}(x,t) \, \partial_x + \partial_t$$

We are interested in the solutions to the ODE
\begin{equation*}
  (E) : \left\{
    \begin{aligned}
    \dot{y}(t) &= F\big(y(t),t\big), \\
    y(t_0) &= y_0,
    \end{aligned}
  \right.
\end{equation*}
where $(t_0, y_0) \in \R \times \R$, or equivalently, in the curves on $C$ tangent to $X$.

Clearly, if there exists a smooth closed curve transverse to $X$, then $X$ cannot be tangent to a foliation by circles on $C$, and the same is true for any vector field sufficiently $C^0$-close to $X$. Conversely,

\begin{prop} \label{prop:ODEper}
Assume that there exist no smooth closed curves transverse to $X$. Then $X$ can be $C^0$-approximated by continuous vector fields tangent to foliations by circles on $C$.

\end{prop}

\begin{proof} A simple closed curve $\gamma$ on $C$ tangent to $X$ can be written as the graph of a ($C^1$) function $f : S^1 \rightarrow \R$. Hence, we can naturally identify $\Gamma$ with a subset of $C^0(S^1, \R)$. We endow $\Gamma$ with the $C^0$ topology and the natural partial order inherited from $C^0(S^1, \R)$. Notice that $\Gamma \subset C^0(S^1, \R)$ is closed, and the bounded subsets of $\Gamma$ are relatively compact. This easily follows from the integral formulation of $(E)$, the fact that the curves in $\Gamma$ are $K$-Lipschitz, where $K \coloneqq \sup \vert F\vert$, and Arzel\`{a}--Ascoli theorem.

We construct a subset $\Gamma_\mathrm{tot} \subset \Gamma$ such that $\leq$ restricts to a total order on $\Gamma_\mathrm{tot}$, and 
\begin{align}\bigcup_{\gamma \in \Gamma_\mathrm{tot}} \gamma = C. \label{eq:union}
\end{align}
For that purpose, we will use the following claim which will be proved later:

\begin{center}
\textbf{Claim.} For every point $p \in C$, there exists $\gamma \in \Gamma$ passing through $p$.
\end{center}
One can easily deduce the following properties by suitably cutting and pasting curves in $\Gamma$:
\begin{itemize}
\item If $\gamma \in \Gamma$ and $p \in C \setminus \Gamma$ lies above (resp.~below) $\gamma$, there exists $\gamma' \in \Gamma$ passing through $p$ such that $\gamma \leq \gamma'$ (resp.~$\gamma' \leq \gamma$).
\item If $\gamma_0, \gamma_1 \in \Gamma$, $\gamma_0 \leq \gamma_1$, and $p \in C$ is in the interior of the domain in $C$ bounded by $\gamma_0$ and $\gamma_1$, then there exists $\gamma \in \Gamma$ passing through $p$ satisfying $\gamma_0 \leq \gamma \leq \gamma_1$.
\end{itemize}

Let $\mathcal{D} \subset \R \times S^1$ be a dense countable set. Using the above two bullet points, we can inductively construct a subset $\Gamma_\mathcal{D} \subset \Gamma$ on which $\leq$ restricts to a total order, and such that
$$\mathcal{D} \subset \bigcup_{\gamma \in \Gamma_\mathcal{D}} \gamma.$$
We now define 
$$\Gamma_\mathrm{tot} \coloneqq \overline{\Gamma_\mathcal{D}} \subset \Gamma,$$ 
the closure of $\Gamma_\mathcal{D}$ in $\Gamma$ for the $C^0$ topology. It is easy to check that $\leq$ restricts to a total order on $\Gamma_\mathrm{tot}$, and that~\eqref{eq:union} is satisfied. In the terminology of~\cite{BI08}, $\Gamma_\mathrm{tot}$ is a \emph{branching foliation} tangent to $X$. The methods of~\cite[Section 7]{BI08} can easily be adapted to show the following modified statement of~\cite[Theorem 7.2]{BI08}: for every $\epsilon > 0$, there exists a continuous vector field $X_\epsilon$ on $C$ such that $\vert X_\epsilon - X\vert_{C^0} \leq \epsilon$, and $X_\epsilon$ is tangent to a foliation by circles.

We are left to prove the Claim stated above. Let $p \in C = \R \times \R \slash \Z$. Without loss of generality, we can assume $p=(0,0)$. We will slightly abuse notation and write $X$ for the vector field $F(x,t) \partial_x + \partial_t$ on the strip $Z \coloneqq \R \times [0,1]$. Let $I \subset \R$ be the set of real numbers $x$ such that there exists a curve in $Z$ tangent to $X$ from $(0,0)$ to $(x, 1)$. By Kneser's theorem,\footnote{See~\cite{H02}, Chapter II, Theorem 4.1.} $I$ is a closed interval. It contains $0$ if and only if there exists $\gamma \in \Gamma$ passing through $(0,0)$. Otherwise, we can assume without loss of generality that $m \coloneqq \min I > 0$, and we construct a closed transversal to $X$, contradicting our original hypothesis. Let $\gamma_m$ be the curve in $Z$ tangent to $X$ from $(0,0)$ to $(m,1)$ which is minimal with respect to $\leq$.\footnote{See~\cite{H02}, Chapter III, Theorem 2.1.} For every $n \geq 0$, we let $$X_n \coloneqq \big(F - 2^{-n} \big) \partial_x + \partial_t = X - 2^{-n} \partial_x,$$
and we denote by $\delta_n$ the minimal curve tangent to $X_n$ starting at $(0,0)$. Note that for every $n \geq 0$, $\delta_n$ is (negatively) transverse to $X$, and $\delta_n \leq \gamma_m$. Let $x_n \in \R$ be such that $\delta_n$ intersects $\R \times \{1\}$ at $\{x_n\} \times \{1\}$. It is easy to see that $\lim_n x_n = m$. Indeed, $x_n \leq m$ and $(\delta_n)_n$ converges uniformly to a curve starting at $(0,0)$ and tangent to $X$. In particular, there exists $N \geq 0$ such that $x_N > 0$. It is then easy to modify $\delta_N$ near $\R \times \{1\}$ to obtain a smooth curve $\widetilde{\delta}_N$ from $(0,0)$ to $(0,1)$ and which is (negatively) transverse to $X$. Therefore, we have constructed a closed transversal to $X$, a contradiction.
\end{proof}

%%%%%%%%%%%%%%%%%%%%%%%%%%%%%%%%%%%%%%%%%%%%%%%%%%%%%%%%%%%%%%%%%% BIBLIOGRAPHY
%%%%%%%%%%%%%%%%%%%%%%%%%%%%%%%%%%%%%%%%%%%%%%%%%%%%%%%%%%%%%%%%%%
\newpage

\printbibliography[heading=bibintoc, title={References}]

\end{document}